\documentclass{article}
\IfFormatAtLeastTF{2026-06-01}{%
  \RemoveFromHook{package/amsthm/after}[firstaid/aliascounter]
  \AddToHook{begindocument/before}{%
  }%
}{}

\usepackage{graphicx}
\usepackage[T1]{fontenc}
\usepackage{amsfonts,amsmath,amssymb,amsthm}
\usepackage[
backend=biber,
style=alphabetic,
backref=true
]{biblatex}
\usepackage{quiver}
\usepackage{tikz,tikz-cd}
\usepackage[shortlabels]{enumitem}
\usepackage{comment}
\usepackage[dvipsnames]{xcolor}
\usepackage{needspace}
\usepackage{etoolbox}
\usepackage{multicol}
\usepackage{mathtools}
\usepackage[parfill]{parskip}

\BeforeBeginEnvironment{theorem}{\Needspace{3\baselineskip}}
\BeforeBeginEnvironment{lemma}{\Needspace{3\baselineskip}}
\BeforeBeginEnvironment{proposition}{\Needspace{3\baselineskip}}
\BeforeBeginEnvironment{corollary}{\Needspace{3\baselineskip}}
\BeforeBeginEnvironment{definition}{\Needspace{3\baselineskip}}
\BeforeBeginEnvironment{example}{\Needspace{3\baselineskip}}
\BeforeBeginEnvironment{convention}{\Needspace{3\baselineskip}}
\BeforeBeginEnvironment{remark}{\Needspace{3\baselineskip}}
\usepackage{hyperref}
\hypersetup{
    colorlinks=true,
    citecolor=blue,
    urlcolor=cyan,
    linkcolor=blue,
    }
\usepackage{geometry}
\usepackage{stmaryrd}
\renewcommand{\llparenthesis}{\mathopen{(\!(}}
\renewcommand{\rrparenthesis}{\mathclose{)\!)}}
\usepackage{xfrac}
\usepackage{float}
\usepackage{adjustbox}

\usepackage{thmtools}
\usepackage{thm-restate}
\usepackage[capitalize]{cleveref}
\crefname{appendix}{Appendix}{Appendices}

\declaretheoremstyle[
  headfont=\bfseries,
  bodyfont=\itshape,
  spaceabove=0pt,
  spacebelow=5pt,
  headpunct={.},
  postheadspace=0.5em
]{plainspaced}

\declaretheoremstyle[
  headfont=\bfseries,
  bodyfont=\normalfont,
  spaceabove=0pt,
  spacebelow=6pt,
  headpunct={.},
  postheadspace=0.5em
]{definitionspaced}

\declaretheoremstyle[
  headfont=\itshape,
  bodyfont=\normalfont,
  spaceabove=0pt,
  spacebelow=5pt,
  headpunct={.},
  postheadspace=0.5em
]{remarkspaced}

\declaretheorem[
  name=Theorem,
  style=plainspaced,
  numberwithin=section
]{theorem}

\declaretheorem[
  name=Lemma,
  style=plainspaced,
  sibling=theorem
]{lemma}

\declaretheorem[
  name=Proposition,
  style=plainspaced,
  sibling=theorem
]{proposition}

\declaretheorem[
  name=Corollary,
  style=plainspaced,
  sibling=theorem
]{corollary}

\declaretheorem[
  name=Definition,
  style=definitionspaced,
  sibling=theorem
]{definition}

\declaretheorem[
  name=Example,
  style=definitionspaced,
  sibling=theorem
]{example}

\declaretheorem[
  name=Convention,
  style=definitionspaced,
  sibling=theorem
]{convention}

\declaretheorem[
  name=Notation,
  style=definitionspaced,
  sibling=theorem
]{notation}

\declaretheorem[
  name=Remark,
  style=remarkspaced,
  sibling=theorem
]{remark}

\renewcommand{\i}{\mathrm{i}}
\renewcommand{\SS}{\mathbb{S}}
\newcommand{\Sp}{\mathbf{Sp}}
\newcommand{\N}{\mathbb{N}}
\newcommand{\Z}{\mathbb{Z}}
\newcommand{\Q}{\mathbb{Q}}
\newcommand{\R}{\mathbb{R}}
\newcommand{\C}{\mathbb{C}}
\newcommand{\E}{\mathcal{E}}
\newcommand{\F}{\mathbb{F}}
\newcommand{\TT}{\mathbb{T}}
\newcommand{\Fcal}{\mathcal{F}}
\newcommand{\Ffrak}{\mathfrak{F}}
\newcommand{\Ccal}{\mathcal{C}}
\newcommand{\Cfrak}{\mathfrak{C}}
\newcommand{\Zcal}{\mathcal{Z}}
\newcommand{\G}{\mathcal{G}}
\newcommand{\id}{\mathrm{id}}
\newcommand{\K}{\mathbb{K}}
\newcommand{\Kur}{\mathcal{K}}
\newcommand{\X}{\mathcal{X}}
\newcommand{\T}{\mathcal{T}}

\newcommand{\CP}{\mathbb{CP}}
\DeclareMathOperator{\Stab}{Stab}
\newcommand{\fix}{\mathrm{fix}}
\newcommand{\mov}{\mathrm{mov}}
\newcommand{\eq}{\mathrm{eq}}
\newcommand{\vir}{\mathrm{vir}}
\newcommand{\std}{\mathrm{std}}
\newcommand{\pre}{\mathrm{pre}}
\newcommand{\elem}{\mathrm{elem}}
\newcommand{\Ob}{\mathrm{Ob}}

\renewcommand{\O}{\mathcal{O}}
\renewcommand{\P}{\mathbb{P}}
\newcommand{\M}{\mathcal{M}}
\newcommand{\Mb}{\overline{\mathcal{M}}}
\newcommand{\A}{\mathbb{A}}

\newcommand{\m}{\mathfrak{m}}
\newcommand{\pb}{\overline{\partial}}
\renewcommand{\H}{\mathcal{H}}
\newcommand{\J}{\mathcal{J}}
\newcommand{\simple}{\mathrm{simple}}
\newcommand{\W}{\mathcal{W}}

\newcommand{\pt}{\mathrm{pt}}
\renewcommand{\X}{\mathcal{X}}
\newcommand{\V}{\mathcal{V}}

\newcommand{\Y}{\mathcal{Y}}
\newcommand{\I}{\mathcal{I}}
\newcommand{\afrak}{\mathfrak{a}}
\newcommand{\loc}{{\operatorname{loc}}}
\newcommand{\source}{{\operatorname{source}}}
\newcommand{\target}{{\operatorname{target}}}
\def\acts{\curvearrowright}

\DeclareMathOperator{\Aut}{Aut}
\DeclareMathOperator{\Def}{Def}
\DeclareMathOperator{\colim}{colim}
\DeclareMathOperator{\coker}{coker}
\DeclareMathOperator{\Hom}{Hom}
\DeclareMathOperator{\rk}{rk}
\DeclareMathOperator{\res}{Res}
\DeclareMathOperator{\Fix}{Fix}
\DeclareMathOperator{\im}{im}
\DeclareMathOperator{\Tot}{Tot}
\DeclareMathOperator{\GW}{GW}

\DeclareMathOperator{\GL}{GL}
\DeclareMathOperator{\Th}{Th}
\DeclareMathOperator{\inthom}{\underline{Hom}}
\DeclareMathOperator{\MU}{MU}
\DeclareMathOperator{\Ho}{Ho}

\DeclareMathOperator{\KO}{KO}
\DeclareMathOperator{\ind}{ind}
\DeclareMathOperator{\codim}{codim}
\DeclareMathOperator{\vdim}{vdim}
\DeclareMathOperator{\Fr}{Fr}
\DeclareMathOperator{\comb}{comb}
\DeclareMathOperator*{\limone}{lim^{\rlap{\textnormal{\scriptsize 1}}}}
\title{Equivariant Localization for $K(n)$-Local Gromov--Witten Invariants with Applications to Gopakumar--Vafa Invariants}
\author{Guangxi Yan}

\begin{document}
\maketitle
\begin{abstract}
For certain complex-oriented $K(n)$-local cohomology theories,
we establish equivariant localization formulas for global Kuranishi
charts and a Riemann--Roch formula.
We apply these formulas to compute the genus-zero Gromov--Witten
invariants of $\mathcal O_{\mathbb P^1}(-1)\oplus
\mathcal O_{\mathbb P^1}(-1)$ at odd primes, obtaining a residue
formula in terms of the formal group law of the chosen complex
orientation.
As an application, we derive a multiple-cover formula expressing
the genus-zero Morava $K$-theoretic Gromov--Witten invariants of
closed symplectic Calabi--Yau threefolds in terms of their
Gopakumar--Vafa invariants.
\end{abstract}
\tableofcontents

\section{Introduction}
Classical Gromov--Witten invariants are virtual
counts of stable maps to smooth projective
varieties or closed symplectic manifolds.
See
\cite{kontsevichGromovWittenClassesQuantum1994,
Li-Tian_algebraic,intrisic-normal-cone}
for the algebraic foundations, and
\cite{Ruan-Tian_1995,Li-Tian_1996,
Fukaya-Ono_1999,McDuff-Salamon_2012,Pardon_2016}
for the symplectic foundations.
Beyond ordinary rational cohomology,
$K$-theoretic Gromov--Witten theory was developed
by Givental and Lee
\cite{Givental_2000,Lee_2004}.

Going beyond ordinary cohomology and complex $K$-theory,
Abouzaid and Blumberg constructed Hamiltonian Floer homology
with coefficients in Morava $K$-theories
\cite{abouzaid2021arnoldconjecturemoravaktheory}.
For Gromov--Witten theory, global Kuranishi charts
allow the construction of invariants with coefficients in suitable
complex-oriented generalized cohomology
theories.
The genus-zero construction appears in
\cite{abouzaid2021complexcobordismhamiltonianloops},
and the arbitrary-genus construction is
developed in
\cite{hirschiGlobalKuranishiCharts2024b,
abouzaidGromovWittenInvariantsComplex2024a}.
Related developments include the integer-valued
Gromov--Witten-type invariants constructed by Bai and Xu
using normally complex perturbations
\cite{bai2025newtransversalityconditionorbifolds}.
Abouzaid and Bai construct complex cobordism-valued
Gromov--Witten invariants by splitting the natural map
from complex bordism to derived orbifold bordism
\cite{abouzaid2025bordismresolutionsingularities}.

In this paper, we study genus-zero
Gromov--Witten invariants with coefficients
in Morava $K$-theories at an odd prime $p$.
Our example\footnote{
The local $\mathbb P^1$ multiple-cover problem considered in this paper was also studied by Tim Large, who discussed equivariant localization for Morava $K$-theoretic Gromov--Witten invariants and the Aspinwall–Morrison multiple-cover formula in positive characteristic in seminar talks in 2022 and 2023. We learned of this problem from the recording of his talk \cite{Large2022}.
} is the local $\P^1$
\[
X_{\loc}:=
\Tot\left(
\O_{\P^1}(-1)\oplus\O_{\P^1}(-1)
\to\P^1
\right).
\]
In ordinary rational cohomology, the degree-$d$
genus-zero local Gromov--Witten invariant
is the Aspinwall--Morrison multiple-cover term $1/d^3$, which can be computed by
torus localization; see, for example,
\cite{Enumeration_of_Rational_Curves_Via_Torus_Actions,
hori2003mirror}.
However, when $p$ divides $d$, the rational
coefficient $1/d^3$ cannot be directly reduced modulo $p$.
This raises a natural question:
what are the corresponding multiple-cover
contributions in Morava $K$-theory?
We address this question for all positive degrees
$d$ and use the resulting local computations
to study genus-zero Gromov--Witten invariants
of Calabi--Yau threefolds.

\subsection{Main results}
We first establish two general results for
genus-zero Gromov--Witten invariants of
$X_\loc$ with coefficients in an HKR theory
(such as the Morava $E$-theory; see \cref{def:definition of HKR theory})
$E$ of height $n$ at an odd prime $p$,
assuming that $E^*$ has no $p$-torsion.

\noindent\textbf{Separation of the prime-to-\(p\) factor.}
Let $I_d^E$ denote the degree-$d$
invariant of $X_\loc$. We prove in \cref{prop:separating p and p-coprime part nonequivariantly}
\[
    I_{md}^E=m^{-3}I_d^E,
    \qquad p\nmid m.
\]
Thus, the computation reduces to
$p$-power degrees. 

\noindent\textbf{Residue formula for $p$-power degrees.}
Suppose the complex orientation of $E$
satisfies $[-1]_E x=-x$. Then
\[
    I_{p^a,\eq}^E
    =
    \frac{
        \tau_{p^a,u}(x^2)
        -\tau_{p^a,0}(x^2)
    }{u^2},
    \qquad
    I_{p^a}^E
    =
    [u^2]\tau_{p^a,u}(x^2),
\]
where
\[
    \tau_{d,u}(f)
    :=
    \res_x
    \frac{f(x)\omega_E(x)}
    {[d]_E x-_E u},
\]
and $\omega_E$ is the invariant differential
of the formal group law of $E$. This is \cref{thm:GW computation if [-1]x=-x}.

For example, at height one, Morava $E$-theory $E_1$ is the
$p$-completion $KU_p^\wedge$ of complex $K$-theory,
so our framework includes $p$-adic complex $K$-theory
as its height-one case.
For a specific complex orientation\footnote{
Here the complex orientation gives the formal group law $(x+y)/(1+xy)$,
which satisfy $[-1]x=-x$.
This differs from the standard complex $K$-theory
orientation, for which the degree-$p^a$ invariant
is $p^{2a}$ \cite{ChouLee2025}.
Indeed, the Thom isomorphisms, and consequently the
Euler classes and pushforwards used to define these
invariants, depend on the choice of complex orientation.} on $E_1$,
the residue formula above gives
$I_{p^a}^{E_1}=p^a$ for every $a\geq 0$.

\noindent\textbf{Evaluation in Morava $K$-theory.}
We then specialize these formulas to
Morava $K$-theory $K(n)$, using the
complex orientation induced by the
integral lift $\widetilde K(n)$.
This gives an explicit evaluation
of the local Gromov--Witten invariants
$I_d^{K(n)}\in\F_p$ for every positive
degree $d$.
\begin{restatable}{theorem}{Kn}
\label{thm:Kn-local-invariants}
Let $p\ge3$, $n\ge1$, and use the standard $p$-typical complex
orientation of $K(n)$.
For every integer $a\ge1$, the degree-$p^a$
local invariants satisfy
\begin{equation}
\label{eq:Kn-pa-invariants}
I_{p^a}^{K(n)}
=
\begin{cases}
-1\in\F_3,
& (p,n,a)=(3,1,1),\\
0,
& \text{otherwise}.
\end{cases}
\end{equation}
Consequently, for any $d\ge1$,
\begin{equation}
\label{eq:Kn-general-invariants}
I_d^{K(n)}
=
\begin{cases}
d^{-3}, & p\nmid d,\\[1mm]
-m^{-3}, &
p=3,\ n=1,\ d=3m,\ 3\nmid m,\\[1mm]
0, & \text{otherwise}.
\end{cases}
\end{equation}
\end{restatable}
As an application, we combine the
local multiple-cover computation with
the cluster formalism
of Ionel--Parker
\cite{Ionel-Parker_Gopakumar-Vafa}.
This yields a multiple-cover formula
for Morava $K$-theoretic Gromov--Witten
invariants of Calabi--Yau threefolds.

For a closed symplectic Calabi--Yau threefold $X$,
let $n_{\beta,0}$ denote its genus-zero
Gopakumar--Vafa invariants.
These invariants originate in the work of
Gopakumar and Vafa \cite{Gopakumar-Vafa_I,Gopakumar-Vafa_II}
and are determined by the ordinary rational
Gromov--Witten invariants through the genus-zero
multiple-cover formula.
Their integrality in the symplectic setting was
proved by Ionel and Parker \cite{Ionel-Parker_Gopakumar-Vafa}.

\begin{restatable}{theorem}{CY}
    \label{CY}
    Let $X$ be a closed symplectic CY threefold and let $A\in H_2(X)$ be nonzero.
    The genus zero Gromov--Witten invariant $\GW_{0,0}^{K(n)}(X,A)\in \F_p$ equals the following:
    \begin{enumerate}[(1)]
        \item If $n>1$ or $p>3$, then
        \[
            \GW_{0,0}^{K(n)}(X,A)=
            \sum_{A=d \beta,\ p\nmid d}
            n_{\beta,0}\frac{1}{d^3}\in \F_p
        \]
        In particular, if $A$ is not divisible by $p$,
        then $\GW_{0,0}^{K(n)}(X,A)$ equals the image in $\F_p$ of the corresponding
        ordinary Gromov--Witten invariant.
        \item If $(n,p)=(1,3)$, then
        \[
            \GW_{0,0}^{K(n)}(X,A)=
            \sum_{A=d \beta,\ p\nmid d}
            n_{\beta,0}\frac{1}{d^3}
            -
            \sum_{A=3m \beta,\ p\nmid m }n_{\beta,0}\frac{1}{m^3}
            \in \F_p.
        \]
    \end{enumerate}
\end{restatable}

\subsection{Ingredients and Proof strategy}
\noindent\textbf{Main ingredients.}
Our approach uses two general tools 
not specific to the local $\P^1$ Gromov--Witten moudli space:

\begin{enumerate}[(1)]
    \item \emph{Equivariant localization.}
    Virtual localization in Gromov--Witten theory
    was established in the algebraic setting by
    Graber and Pandharipande
    \cite{Graber-Pandharipande99}.
    More recently, Hirschi established a virtual
    localization formula for global Kuranishi charts
    in rational equivariant cohomology
    \cite[Theorem~5.13]{Hirschi_Properties_of_Gromov_Witten}.

    We establish Atiyah--Bott-type localization
    formulas with coefficients in a Morava
    $K$-theory or an HKR theory, for
    torus-equivariant orbifolds
    (\cref{prop:Atiyah-Bott equivariant localization formula})
    and global Kuranishi charts
    (\cref{prop:Atiyah-Bott localization for equivariant global Kuranishi charts}).

    For quotient orbifolds and global Kuranishi
    charts, the relevant fixed loci are determined
    by the induced torus action on the quotient,
    rather than by naive fixed-point conditions
    on the presenting space
    \cite[Lemmas~5.8 and 5.11]{Hirschi_Properties_of_Gromov_Witten}.
    See also
    \cref{def:fixed locus of equivariant global Kuranishi chart}.

    \item \emph{HKR character theory and Riemann--Roch.}
    The fixed-point contributions in the
    localization formula are orbifold pushforwards
    in $E$-cohomology. To compute them, we use
    $K(n)$-local ambidexterity and
    Hopkins--Kuhn--Ravenel character theory.
    We establish a Riemann--Roch formula
    for pushforwards between stably almost
    complex orbifolds and its torus-equivariant
    extension
    (\cref{thm:a Riemann--Roch,prop:Riemann Roch for equivariant HKR}).
    These formulas show that HKR character maps,
    corrected by the Euler classes of the
    HKR-moving normal bundles, commute with
    pushforwards, allowing computations
    in the rational theory $L(E)$ over
    HKR sectors.

    The homotopy-theoretic motivation comes from
    the compatibility of generalized characters
    with transfers, established by
    Hopkins--Kuhn--Ravenel
    \cite[Theorem~D]{Generalized-Group-Characters-and-Complex-Oriented-Cohomology-Theories}
    and extended to transchromatic characters
    by Stapleton
    \cite[Theorem~2.18]{Stapleton2015}.
    More recently, Ben-Moshe established
    compatibility with integration over
    $\pi$-finite spaces
    \cite[Theorem~A]{benmoshe2024highersemiadditivitytranschromatichomotopy}.
    In our setting, the oriented Gysin
    pushforwards require additional
    corrections arising from Thom classes.

    At height one, the HKR construction has the familiar
    description in terms of pairs $(x,g)$, where $g$
    belongs to the stabilizer of $x$, with the additional
    requirement that $g$ have $p$-power order.
    Its sectors are therefore precisely the usual
    twisted sectors associated to $p$-power-order
    automorphisms. Our Riemann--Roch formula is
    correspondingly analogous to Kawasaki's
    Riemann--Roch formula \cite{ojm/1200771835}:
    orbifold pushforwards are computed from characters
    on these sectors, corrected by the Euler classes
    of the moving normal bundles.
\end{enumerate}

\medskip
\noindent\textbf{Proof strategy.}
The computation for local $\P^1$ proceeds
in three stages.

\begin{enumerate}[(1)]
    \item \emph{Geometric reduction.}
    We first compare the AMS global Kuranishi
    chart with the obstruction-bundle model from algebraic geometry
    over $\Mb_{0,0}(\P^1,d)$
    (\cref{thm:equivalence of AMS charts and AG charts}).
    Equivariant localization then expresses
    the invariant as a sum over torus-fixed
    components. An Euler-class vanishing
    argument eliminates all but two types:
    single-edge components and comb components
    with equal-degree tails
    (\cref{Euler class vanishing}).
    We describe their orbifold structures
    and compute the corresponding virtual
    Euler classes
    (\cref{lem:virtual-euler-class-case1,lem:virtual-euler-class-case2}).

    \item \emph{Separation of the prime-to-$p$ factor.}
    We compare the equivariant localization
    contributions in degrees $d$ and $md$,
    where $p\nmid m$.
    This yields an equivariant identity
    (\cref{prop:separating p and p-coprime part equivariantly})
    whose nonequivariant consequence is
    \[
        I_{md}^E=m^{-3}I_d^E.
    \]
    Thus, the prime-to-$p$ part contributes
    exactly the classical multiple-cover
    factor, and the remaining calculation
    reduces to $p$-power degrees.

    \item \emph{$p$-primary computation.}
    For $d=p^a$, a further Euler-class
    vanishing eliminates the nontransitive
    HKR sectors of comb components
    (\cref{lem:Euler class vanishing after HKR}).
    Consequently, only comb components
    with a $p$-power number of tails
    can contribute.

    When $[-1]_E x=-x$, the Euler-class ratio on the single-edge
    component simplifies to $\frac{x^2}{([d]_E x)^2}$.
    The projection formula and \cref{prop:An equivariant pushforward formula computed by residue}
    then give the single-edge contribution $\tau_{d,u}(x^2)/u^2$.
    For the comb components, HKR
    Riemann--Roch reduces the orbifold
    pushforwards to sums over transitive
    sectors. We evaluate these sums using
    the combinatorics of finite abelian
    groups and trace--residue identities.
    The resulting comb contributions form
    a telescoping sum (\cref{prop:contribution from components}(2)), giving precisely
    the correction term
    $-\tau_{p^a,0}(x^2)/u^2$.
    Combining this with the single-edge
    contribution yields
    \cref{thm:GW computation if [-1]x=-x}.
\end{enumerate}

Finally, we specialize the residue formula
to the integral lift $\widetilde K(n)$
and reduce modulo $p$ to obtain
\cref{thm:Kn-local-invariants}.
The application to Calabi--Yau threefolds
then follows from the genus-zero cluster
decomposition
(\cref{decomposition into elementary clusters}).

\medskip
\noindent\textbf{Further directions.}
A natural direction is to extend the explicit calculations
in this paper to other targets with torus actions, including
toric varieties and resolutions of ADE singularities.
This would involve describing their torus-fixed moduli spaces
and evaluating the corresponding orbifold pushforwards using
the localization and Riemann--Roch formulas developed here.
For Springer resolutions, it would be interesting to seek
Morava $K$-theoretic analogues of the equivariant
quantum-cohomological calculations of Braverman, Maulik,
and Okounkov \cite{BravermanMaulikOkounkov2011}.

\subsection{Outline}
The paper is divided into two parts.
Part~I contains the computations and applications,
while Part~II develops the general foundations
underlying them. 

In Part~I, Section~2 compares the global
Kuranishi chart with the obstruction-bundle
model for local $\P^1$, and studies the
torus-fixed loci, virtual Euler classes,
and their vanishing properties.
Section~3 computes the local Gromov--Witten
invariants by separating the prime-to-$p$
factor, deriving a residue formula for
$p$-power degrees, and evaluating it in
Morava $K$-theory.
Section~4 applies these computations to
Calabi--Yau threefolds using the genus-zero
cluster decomposition.

Part~II provides the general results used
in Part~I.
Section~5 develops orbifold pushforwards,
HKR character theory, and Riemann--Roch
and equivariant localization formulas.
Section~6 reviews global Kuranishi charts
and their construction, and establishes
equivariant localization in this setting.

\subsection{Acknowledgments}
I am grateful to my advisor, Shaoyun Bai, for suggesting this problem and for his guidance and helpful discussions. I thank Ivan Smith and Noah Porcelli for helpful conversations in Cambridge. I also thank Mark McLean for a helpful discussion of this work and for his valuable suggestions. I thank Jeremy Hahn and Keita Allen for answering my questions about Morava $K$-theory and chromatic homotopy theory.

\part{Computations and Applications}
\label{part:computations}

\section[Fixed loci and virtual Euler classes for torus action on local \texorpdfstring{$\P^1$}{P\textasciicircum1}]
{Fixed loci and virtual Euler classes\\for torus action on local \texorpdfstring{$\P^1$}{P\textasciicircum1}}
    \label{Localization computation}
    \begin{definition}
      By a local $\P^1$ we mean the noncompact
      Calabi--Yau threefold
    \[  
        X_\loc=\Tot\left(\O_{\P^1}(-1)\oplus\O_{\P^1}  (-1)\right)\to \P^1
    \] with the standard
    complex structure.
     Write
    \[
        N=\mathcal O_{\mathbb P^1}(-1)\oplus \mathcal O_{\mathbb P^1}(-1).
    \]
    \end{definition}

    The definitions and conventions for global Kuranishi
    charts used here are given in Part~\ref{part:foundations};
    see \cref{def:global Kuranishi chart,conv:global Kuranishi charts}.
    The AMS construction is reviewed there as well
    (\cref{non-equivariant global Kuranishi chart construction}).

    We will focus on the moduli space
    $\Mb_{0,0}(X_\loc,d)$. The AMS construction
    produces a global Kuranishi chart $\Kur^{\operatorname{AMS}}$ (up
    to equivalences and cobordisms), which
    is too abstract to perform computation.
    However, this moduli space has already
    a model from algebraic geometry (\cref{def:AG model})
\[\begin{tikzcd}
	{H^1\left(\Sigma,f^*N\right)} \\
	{\Mb_{0,0}(\P^1,d)}
	\arrow[from=1-1, to=2-1]
	\arrow["0"', curve={height=6pt}, from=2-1, to=1-1]
\end{tikzcd}\]
    This is a derived orbifold chart, and after finding a global quotient presentation of
    $\M_{0,0}(\P^1,d)$, we have a global Kuranishi
    chart of $\Mb_{0,0}(X_\loc,d)$, called $\Kur^{\operatorname{AG}}$. This chart is optimal in the sense that
    the section is zero. We devote the
    first subsection to showing that
    $\Kur^{\operatorname{AG}}$ is
    equivalent to $\Kur^{\operatorname{AMS}}$ in the sense
    of \cref{def:equivalence of global Kuranishi charts}.
    
    In the second subsection, we take into consideration a $\TT=T^2$-action
    on $N\to \P^1$ which defines a $\TT$-action on $\Mb_{0,0}(X_\loc,d)$ naturally. To prepare for equivariant
    localization (\cref{prop:Atiyah-Bott localization for equivariant global Kuranishi charts}),
    we identify the fixed locus of $\TT$-action on
    $\Kur^{\operatorname{AG}}$. For each connected component of fixed locus, we represent its virtual
    normal bundle by sheaf cohomology, and compute
    their $\TT$-equivariant Euler characteristics.
\subsection{Comparison between the AMS chart and a derived orbifold chart}
    The main result of this subsection is \cref{thm:equivalence of AMS charts and AG charts}.
    \begin{definition}
    \label{def:naive chart}
        Define
        \[
        M_d:=\overline{\mathcal M}_{0,0}(\mathbb P^1,d).
        \]
        to be the Gromov--Witten moduli space of degree $d$
        stable maps from rational nodal curves to $\P^1$.
    \end{definition}
    \begin{lemma}
        $M_d$ is a smooth proper complex Deligne–Mumford stack.
        Moreover $R^0\pi_*f^*N=0$, and $R^1\pi_*f^*N$
        is a complex vector orbi-bundle of rank $2d-2$.
    \end{lemma}
    \begin{proof}
        This is classical in \cite{hori2003mirror}.
        By Lemma 24.4.3 there, \(\mathbb P^1\) is convex, hence \(M_d\) is a nonsingular Deligne–Mumford stack. In  \S27.5, and specifically Lemma 27.5.2, they identify the obstruction bundle for the local \(\mathbb P^1\) multiple-cover problem with
\[ R^1\pi_*f^*(\mathcal O(-1)^{\oplus2}), \]
of rank \(2d-2\).

It remains only to note that \(R^0\pi_*f^*N=0\). But this is true by
noticing that $f^*N=f^*\O(-1)\oplus f^*\O(-1)$. $f^*\O(-1)$ has either
negative or zero degree on each irreducible component of $\Sigma$. Any holomorphic section
of this line bundle vanishes on negative-degree components and is
constant on zero-degree components. Since there is at least one
negative-degree component and $\Sigma$ is connected, we have
$R^0\pi_*f^*N=0$.
    \end{proof}
    \begin{definition}
        Define
        \[
            \Ob_d=R^1\pi_*f^*N.
        \]
        as a complex vector orbi-bundle over $M_d$.
    \end{definition}
Let \(i:\mathbb P^1\hookrightarrow X_{\loc}\) denote the zero section. 
There is a canonical splitting along the zero section
\[
    i^*TX_{\loc} \cong T\mathbb P^1\oplus N.
\] 
Consequently, for a stable map \(f:\Sigma\to\mathbb P^1\), the deformation–obstruction theory of the map \(i\circ f:\Sigma\to X_{\loc}\) splits into a tangent part, corresponding to deformations as a stable map to \(\mathbb P^1\), and a normal part given by \(R\Gamma(\Sigma,f^*N)\).
Since \(\mathbb P^1\) is convex, the tangent part is unobstructed and its deformation space is the tangent space of \(M_d\). On the other hand, by the preceding lemma,
$$ H^0(\Sigma,f^*N)=0, \qquad H^1(\Sigma,f^*N)=(\operatorname{Ob}_d)_{[f]}. $$
Thus the deformation and obstruction spaces of the local-\(\mathbb P^1\) problem at \([f]\) are respectively
$$ T_{[f]}M_d \qquad\text{and}\qquad (\operatorname{Ob}_d)_{[f]}. $$
In particular, the virtual tangent bundle of this obstruction-bundle model is
$$ T^{\mathrm{vir}} = TM_d-\operatorname{Ob}_d. $$
This is the reason that the derived orbifold chart
$$ (M_d,\operatorname{Ob}_d,0) $$
is the natural obstruction-bundle model for the local-\(\mathbb P^1\) moduli problem.
    \begin{remark}
         In fact, $M_d$ admits
    a global quotient presentation. Fix an integer $k\ge 3$,
    and the universal curve is denoted $\Ccal_d \xrightarrow{\pi} M_d$, with universal stable map $f$.
    \[\begin{tikzcd}
	\Ccal_d & {\P^1} \\
	{M_d}
	\arrow["f", from=1-1, to=1-2]
	\arrow["\pi", from=1-1, to=2-1]
    \end{tikzcd}\]
    Then we can define a line bundle over $\Ccal$,
    \begin{equation}
    \label{eq:polarizing line bundle in AG}
        L_k=\left(
            \omega_{\Ccal_d/M_d}\otimes f^*\O_{\P^1}(3)
        \right)^{\otimes k},
    \end{equation}
    which turns out to be relatively very ample for sufficiently
    large $k$. Write $D:=k(3d-2)$ for the degree of $L_k$ on each domain
    curve, using $\deg\omega_\Sigma=-2$, and set $r=D+1$.
    Here $d$ remains the stable-map degree. Let $E_k=\pi_*L_k$;
    since $H^1(L_k)=0$ on each fiber, Riemann--Roch gives $\operatorname{rank}E_k=r$.
    Let $\Fr(E_k)$ be the framings of $E_k$.
    For \(k\) sufficiently large, the stabilizer of every point of \(M_d\) acts faithfully on the corresponding fiber of \(E_k\). Consequently the frame bundle \(\Fr(E_k)\) is a smooth complex manifold and
        \[
            M_d\simeq [\Fr(E_k)/\GL(r,\C)]
        \]
    \end{remark}
   \begin{definition}
   \label{def:AG model}
   \begin{enumerate}[(1)]
       \item  We call $(M_d,\Ob_d,0)$ the \emph{derived orbifold chart}
       of $\Mb_{0,0}(X_\loc,d)$.
       \item Choose a Hermitian metric on $E_k$,
        and let $\Fr_U(E_k)$
        be the unitary framing bundle.
        Then we can define a global Kuranishi chart
        \begin{equation}
            \label{eq:AG chart}
            \Kur^{\mathrm{AG}}=\left(U(r),\Fr_U(E_k), \widetilde{\Ob}_d,0\right)
        \end{equation}
        Here, $\widetilde{\Ob}_d$ is the lift of $\Ob_d\to M_d$ to $\Fr_U(E_k)$.
   \end{enumerate}
   \end{definition}
    On the other hand, we also have an
    abstract global Kuranishi chart
    $\Kur=(G,\T,\E,s)$
    from the AMS construction. We prove that this chart is equivalent to
    $\Kur^{\mathrm{AG}}$ in the sense of \cref{def:equivalence of global Kuranishi charts}(2)
    (see \cref{lem:4 moves} for the four elementary moves).
    Since changing the AMS auxiliary data gives an equivalent chart,
    it suffices to carry out the comparison for one choice.
    
    The polarizing line bundles in \eqref{eq:polarizing line bundle in AMS}
    and \eqref{eq:polarizing line bundle in AG} suggest the following choice of auxiliary data.
    \begin{itemize}
        \item 
        $L=p^*\O_{\P^1}(3)$, where $p\colon X_\loc\to \P^1$ is the projection.
        Equip $L$ with the Hermitian metric and Chern connection pulled back
        from the standard Hermitian metric on $\O_{\P^1}(3)$.\footnote{The curvature form of this connection vanishes in vertical
    directions, so it does not satisfy the ambient taming condition of
    \cite[Definition~4.23]{abouzaidGromovWittenInvariantsComplex2024a}.
    What matters here is that, on every non-contracted component of a stable
    holomorphic map $u$, the pullback $u^*\Omega$ is nonnegative and is positive
    wherever the differential of $u$ is injective. Together with stability on
    contracted components, this gives the required positivity of the polarizing
    line bundle, which persists near the zero locus.}
        \item The auxiliary degree is $D=\deg_{\Sigma}\bigl((\omega_\Sigma\otimes u^*L)^{\otimes k}\bigr)=k(3d-2)$, and $r=D+1$.
        \item Although
    in general construction we need the finite dimensional
    approximation scheme
    to be dense in (compactly supported sections of) $\Omega^{0,1}_{\Ccal^\circ/\Fcal}\otimes_\C T M$, now we only need them to be dense in vertical direction. In fact, we can choose
    a sufficiently large complex $U(r)$-representation $W$ with a complex-linear perturbation map
    \[
        \lambda^v\colon W\to \Gamma\left(
        \Ccal\times \P^1, \Omega^{0,1}_{\Ccal^\circ/\Fcal}\otimes_\C
        N
        \right).
    \]
     This is because the base $\P^1$ is convex, so there's no obstruction
    bundle on $T\P^1$ direction. Also, at $s^{-1}(0)$ all ($J$-holomorphic) curves fall into $\P^1\subset X_\loc$. 
    By pulling back along $p\colon X_\loc\to \P^1$ and
    then inclusion $p^*N\hookrightarrow TX_\loc$, this $\lambda^v$ is lifted 
    to 
    \[
        \lambda\colon W\to \Gamma\left(
        \Ccal\times X_\loc,\Omega^{0,1}_{\Ccal^\circ/\Fcal}\otimes_\C
        T X_\loc
        \right).
    \]
    $\lambda$ lies in the vertical direction and is constant
    on vertical fibers.
    \item The constant domain metric should be identical to
    the Hermitian metric on $E_k$ for holomorphic curves
    in the zero section,
    $f\colon \Sigma\to \P^1$.
    \end{itemize}
   For this choice, the AMS global
    Kuranishi chart is 
    \begin{equation}
    \label{eq:AMS chart}
        \Kur^{\operatorname{AMS}}=\left(U(r),\T,H_{0,0,D}\oplus \underline{W}\oplus \underline{\H},s\right).
    \end{equation}
    Here $s$ is the section of \cref{def: AMS construction}, specialized to
    $g=h=0$, auxiliary degree $D$, and the perturbation space $W$ chosen above.

    For the purpose of showing $\Kur^{\mathrm{AG}}$ and
    $\Kur^{\operatorname{AMS}}$ are equivalent, we introduce
    an intermediate global Kuranishi chart. Over $\Fr_U(E_k)$, there 
    is a map of vector bundles
    \[
        W\xrightarrow{\lambda} \Gamma(\Sigma, \Omega^{0,1}_{\Sigma}\otimes f^* N)
        \twoheadrightarrow H^1(\Sigma,f^* N)=\widetilde{\Ob}_d|_{[f]}
    \]
        Let $K_1=\ker \left(\underline{W}\to\widetilde{\Ob}_d\right)$.
    \begin{definition}
        The intermediate global Kuranishi chart is
        \begin{equation}
            \label{eq:intermediate chart}
            \Kur_1=\left(
            U(r), \Tot(K_1),\underline{W},s_1
            \right)
        \end{equation}
        Here, the section $s_1$ is induced by the vector
        bundle inclusion $K_1\hookrightarrow \underline{W}$.
        More strictly speaking, we may choose $\T_1$ 
        to be a sufficiently small neighborhood of
        the zero section
        $\Fr_U(E_k)\subset \Tot(K_1)$
    \end{definition}
    \begin{remark}
        A point in $\Tot(K_1)$ consists of the following data
        \begin{itemize}
            \item $f\colon \Sigma\to \P^1$ is a
            holomorphic map from a rational nodal curve 
            $\Sigma$ to $\P^1$.
            \item $F$ is a unitary framing of
            $\left(\omega_\Sigma\otimes f^*\O_{\P^1}(3)\right)^{\otimes k}$.
            \item $e\in K_1|_{(f,F)}$, i.e. a vector $e\in W$ which lies in $\ker(W\to\widetilde{\Ob}_d)$.
        \end{itemize}
    \end{remark}
    By $e\in K_1$, its image under $\lambda$ $\lambda(e)$ lies
    in the image of 
    \[
        \Omega^0(\Sigma,f^*N)\xrightarrow{\pb_{f^*N}}
        \Omega^{0,1}(\Sigma, f^* N)
    \]
    Noticing that $\ker(\pb_{f^*N})=H^0\left(
    \Sigma, f^*N
    \right)=0$, there exists a unique $\xi\in \Omega^0(\Sigma, f^* N)$ satisfying the equation
    \[
        \pb_{f^* N}\xi +\lambda(e)|_{f}=0
    \]
    Thus, we get a map $u=(f,\xi)\colon \Sigma\to X_\loc$,
    a unitary framing $F$, and $e\in W$ satisfying the perturbed
    Cauchy--Riemann equation. 
    \begin{definition}
        We define a morphism of global Kuranishi charts 
        $i\colon \Kur_1\to \Kur^{\operatorname{AMS}}$
        as follows:
        \begin{enumerate}[(1)]
            \item The groups are the same: both $U(r)$.
            \item The aforementioned discussion gives a map\footnote{
        For the moment it is only a map between sets.}
        \[
            i\colon \Tot(K_1)\to \T
        \]
        \[
            (f,F,e)\mapsto \left(\phi=\phi_F, (f,\xi),F,e\right)
        \]
        Here $\phi_F\colon \Sigma\to \P^D$ is the map from framing $F$.
        \item The map of vector bundles is just inclusion of a summand
        \[
            \underline{W}\hookrightarrow H_{0,0,D}\oplus \underline{W}
            \oplus\underline{\H}
        \]
        \[
            e\mapsto (0,e,0).
        \]
        \end{enumerate}
    \end{definition}
    \begin{proposition}
     \label{prop:an equivalence prop}
         $i$ is an equivalence of global
        Kuranishi charts.
    \end{proposition}
    We leave the proof to \cref{app:Proof of an equivalence}. Once we have the equivalence, it is easy
    to observe a further equivalence:
\[\begin{tikzcd}
	{\underline{W}} & {\coker(K_1\hookrightarrow W)} \\
	{\Tot(K_1)} & {\Fr_U(E_k)}
	\arrow[from=1-1, to=1-2]
	\arrow[from=1-1, to=2-1]
	\arrow[from=1-2, to=2-2]
	\arrow["{s_1}"', curve={height=6pt}, from=2-1, to=1-1]
	\arrow[from=2-1, to=2-2]
	\arrow["0"', curve={height=6pt}, from=2-2, to=1-2]
\end{tikzcd}\]
In this diagram, notice that $\coker(K_1\hookrightarrow W)\simeq H^1(\Sigma,f^* N)=\widetilde{\Ob}_d$. Therefore, we get
    \begin{theorem}
    \label{thm:equivalence of AMS charts and AG charts}
     $\Kur=(G,\T,\E,s)$ is equivalent to the chart
            $\left(U(r),\Fr_U(E_k), \widetilde{\Ob}_d,0\right)$.
\end{theorem}
  
\subsection{Fixed-locus and Euler-class computations}

\begin{notation}
Let $\mathcal X$ be a closed stably almost complex
orbifold of dimension $d$, and let $E$ be a
Morava $K$-theory or an HKR theory.
For $a\in H^*(\mathcal X,E)$, write
\[
  \int_{\mathcal X}a
  :=(\pi_{\mathcal X})_!(a)
  \in E^{*-d},
\]
where $\pi_{\mathcal X}:\mathcal X\to\mathrm{pt}$
and the pushforward is defined in
\cref{def:orbifold proper pushforward}.

For a $\mathbb T$-equivariant orbifold,
we similarly write
\[
  \int_{\mathcal X}^{\mathrm{eq}}a
  :=\pi_{!,\mathrm{eq}}(a),
\]
with the equivariant pushforward defined
in \cref{def:equivariant pushforward}.
\end{notation}
    Consider the following map
    \begin{equation}
    \label{eq:torus action on C^2}
        (t_0,t_1)\in T^2, (z_0,z_1)\in \C^2
        \mapsto (t_0z_0,t_1z_1)\in \C^2.
    \end{equation}
    We construct some $T^2$-actions:
    \begin{definition}
        \begin{enumerate}[(1)]
        \item $T^2$-action on $\P^1$. For
        $(t_0,t_1)\in T^2$
        $[z_0:z_1]\in \P^1$, it is sent to
        $[t_0z_0:t_1z_1]\in \P^1$.
        \item The standard $T^2$-action on $\O_{\P^1}(-1)$.
        Notice that $\Tot(\O_{\P^1}(-1))=\operatorname{Bl}_0\C^2$. Then blowing up
        $\C^2$ at $0$ in \eqref{eq:torus action on C^2},
        we get a $T^2$-action on $\Tot(\O_{\P^1}(-1))$ which
        is linear on fibers. We denote this equivariant
        lift of $\O(-1)$ by $\O(-H)$.
        \item Let $\C_{\alpha_0},\C_{\alpha_1}$ be the one-dimensional $T^2$ complex representations given by
        \[
            (t_0,t_1)\mapsto t_0^{-1}, \quad
            (t_0,t_1)\mapsto t_1^{-1}.
        \]
        For $k,l\in \Z$, we write $\C_{k\alpha_0+l\alpha_1}
        =(\C_{\alpha_0})^{\otimes k}\otimes (\C_{\alpha_1})^{\otimes l}$. These are
        all one-dimensional $T^2$ representations.
        \item For $k,l\in \Z$, we have a $T^2$-equivariant
        vector bundle $\O_{\P^1}(-H+k\alpha_0+l\alpha_1)$, 
        defined as the tensor product
        \[
            \O(-H)\otimes (\C_{\alpha_0})^{\otimes k}\otimes (\C_{\alpha_1})^{\otimes l}
        \]
        Since the underlying vector bundles of
        $\C_{\alpha_0},\C_{\alpha_1}$ are trivial line bundles, $\O_{\P^1}(-H+k\alpha_0+l\alpha_1)$ is also an equivariant lift of $\O(-1)$.
    \end{enumerate}
    \end{definition}
    \begin{definition}
    \label{def:convention for 0 and infty}
          We write
$p_0=[1:0]$, $p_1=[0:1]$,
and use the affine coordinate
$z=\frac{z_1}{z_0}$
on the chart $\{z_0\neq 0\}\subset \mathbb{P}^1$.
Then $p_0$ corresponds to $0/1=0$, and
$p_1$ corresponds to $1/0=\infty$.
    \end{definition}
    \begin{definition}
    \label{def:the right T2 action on N}
        We define the $T^2$-action on $X_{\loc}$
        using the following $T^2$ equivariant lift of $\O(-1)^{\oplus 2}$
        \[
            \O(-H+\alpha_0)\oplus \O(-H+\alpha_1).
        \]
    \end{definition}
    \begin{remark}
    \label{rmk:S1 action convention on O(-1)oplusP(-1)}
        The action factors through $q\colon T^2\to S^1$, given by
        $q(t_0,t_1)=t_1/t_0$.
        The corresponding
        $S^1$-action on $\O(-1)\oplus \O(-1)$
        has the following weights at $0$
        and $\infty$:
\[
\begin{array}{c|cc}
& 0 & \infty \\ \hline
\mathcal{O}(-1)_1 & 1 & \chi^{-1} \\
\mathcal{O}(-1)_2 & \chi & 1
\end{array}
\]
Here $\chi$ stands for the $S^1$-representation
inducing $\O(1)\to BS^1=\CP^\infty$\footnote{
Contrary to what the notation suggests, this means the representation $S^1\to \C^*$ given by $t\mapsto
t^{-1}$.
}.
Writing $u=e_E(\chi)$ and using $\alpha_i$ also for the corresponding
Euler classes, we have $q^*u=\alpha_1-_E\alpha_0$.
By naturality of equivariant pushforward, the $T^2$-equivariant invariant
is the pullback along $Bq$ of the $S^1$-equivariant invariant.
    \end{remark}
    By \cref{sec:Equivariant lift of AMS construction}, the $T^2$-action defines
    an equivariant global Kuranishi chart
    for $\Mb_{0,0}(X_\loc,d)$, which allows us to
    compute the Gromov--Witten invariant using equivariant localization (\cref{prop:Atiyah-Bott localization for equivariant global Kuranishi charts}). 
    Our approach follows the same idea as \cite[section 2.3]{Enumeration_of_Rational_Curves_Via_Torus_Actions}
    or 
    \cite[section 27.3--27.7]{hori2003mirror}. 

    Let $e_{\eq}(\Ob_d)$ be the $T^2$-equivariant
    Euler class of $\Ob_d$. 
    By \cref{prop:Atiyah-Bott equivariant localization formula}
    \begin{equation}
    \label{eq:concrete equivariant localization formula}
        \int_{\Mb_{0,0}(\P^1,d)}e_\eq\left(
            \Ob_d
        \right)=\sum_{F}\int_F \frac{e_{\eq}(\Ob_d)}{e_{\eq}(\nu_F)}
    \end{equation}
    where $F$ ranges over the components of fixed loci (defined
    as $T^2$-equivariant orbifolds by \cref{def:fixed points of orbifolds}) of the $T^2$-action on $\Mb_{0,0}(\P^1,d)$. 

   By \cite[Section 3.2]{Enumeration_of_Rational_Curves_Via_Torus_Actions}, components of fixed locus can be
   described by graphs $\Gamma$.
   \begin{definition}
       We consider decorated trees
       $\Gamma$ with the following
       data:
       \begin{itemize}
           \item Each vertex $v$ is
           assigned either $p_0$ or $p_1$
           in $\P^1$. Let $\deg(v)\in \N$
           be its degree. 
           \item Each edge $e$ is assigned
           a degree $d_e\in \N$.
           \item If two vertices $v_1,v_2$ are adjacent, then their assigned
           points are different.
       \end{itemize}
   \end{definition}
   Given such a tree $\Gamma$, we can construct
   stable maps in the fixed locus in the following
   way:
   \begin{itemize}
       \item For each vertex $v$ such that $\deg(v)\ge 3$, associate
       a $\Sigma_v\in \Mb_{0,\deg(v)}$.
       The marked points of $\Sigma_v$ correspond
       to the edges incident to $v$.
       The map $u_v\colon \Sigma_v\to \P^1$ is entirely
       to $p_0$ or $p_1$, depending
       on which point is associated to $v$.
       \item For each edge $e$, associate
       a copy of $\P^1$ to it, named $\P^1_e$.
       The map on $\P^1_e$ is modeled by
       \[
            u_e\colon \P^1_e\to \P^1,\quad 
            z\mapsto z^{d_e}.
       \]
       We still use $p_0,p_1\in \P^1_e$ to denote the
       point mapped to $p_0,p_1\in \P^1$
       under $u_e$.
       \item For a vertex $v$ with $\deg(v)\ge 3$ that is an endpoint of an edge $e$, assume
       $v$ is assigned $p_\epsilon$, and the marked
       point corresponding to $e$ is $q$.
       Then
       we identify the $p_\epsilon\in \P^1_e$ and $q\in \Sigma_v$ to glue the two curves
       together. 

       For a vertex with $\deg(v)=2$, with two
       edges $e_1,e_2$ we glue
       $\P^1_{e_1}$ and $\P^1_{e_2}$ at
       the corresponding $p_\epsilon$ assigned to
       $v$. 

       For a vertex with $\deg(v)=1$, we don't
       glue.
   \end{itemize}
    For all such incident pairs $(v,e)$, we glue
    curves by identifying points, getting a
    rational nodal curve $\Sigma$.
    Since the
    maps $u_v,u_e$ coincide on identified points,
    there is a holomorphic map
    \[
        u\colon \Sigma\to \P^1.
    \]
    It is easy to see this $u$ lies in the
    $T^2$-fixed locus. Actually, the inverse
    is true.
    Notice that a tree $\Gamma$ actually gives
    a family of stable maps, because the choice
    of $\Sigma_v$ allows moving
    in $\Mb_{0,\deg(v)}$. 
    \begin{definition}
        \begin{enumerate}[(1)]
            \item Write 
                \[
                    \Mb_\Gamma=\prod_{v}\Mb_{0,\deg(v)}.
                \] 
            \item Let $\Aut\Gamma$ be
            the automorphism group
            of $\Gamma$. We define
            a map
           \[
\begin{aligned}
\Aut(\Gamma)
&\rightarrow
\Aut\left(\prod_e \Z/d_e\Z\right),\\
\gamma
&\mapsto
\left(
(i_e)\mapsto (i_{\gamma^{-1}(e)})
\right).
\end{aligned}
\]
            \item Define $\A_\Gamma=
            \Aut(\Gamma)\ltimes \prod_e \Z/d_e\Z$,
            so there is an exact
            sequence
            \begin{equation}
            \label{eq:A_Gamma extension sequence}
                 0\to \prod_e \Z/d_e\Z
                \to 
                \A_\Gamma
                \to 
                \Aut(\Gamma)\to 1.
            \end{equation}
            \item Define
            an $\A_\Gamma$-action on $\Mb_\Gamma$ by
            using the $\Aut\Gamma$-action on $\Mb_\Gamma$
            (and the projection $\A_\Gamma\to \Aut(\Gamma)$).
        \end{enumerate}
    \end{definition}
    \begin{lemma}[{\cite[Section 3.2]{Enumeration_of_Rational_Curves_Via_Torus_Actions}\cite[Section 27.3]{hori2003mirror}}]
       The components of $T^2$-fixed locus
       of $\Mb_{0,0}(\P^1,d)$ correspond
       bijectively
       to decorated trees $\Gamma$ with
       $\sum_{e\in\Gamma}d_e=d$.
       Let $F_\Gamma$ be the component corresponding to $\Gamma$, then
       there is an isomorphism of orbifolds
       \[
            F_\Gamma\simeq [\Mb_\Gamma/\A_\Gamma]
       \]
   \end{lemma} 
   Now that we have identified components 
   of the
   fixed locus, we need the deformation
   and obstruction information on each component,
   i.e. the $e_{\eq}(\Ob_d)$ and $e_\eq(\nu_F)$
   term in \eqref{eq:concrete equivariant localization formula}.
   
   Let us recall the computation in 
   the proof of \cite[Lemma 27.5.3]{hori2003mirror}. 
     For a stable map $u\colon \Sigma\to \P^1$ in $F_\Gamma$,
       let $b_1,\ldots,b_m$ be the nodes joining the pieces associated to $\Gamma$.
       Let $\pi\colon \widetilde{\Sigma}\to \Sigma$ be the partial normalization
       at these nodes, leaving the internal nodes of contracted subcurves unresolved,
       and write $\widetilde{u}=u\circ \pi$.
       There is an exact sequence of sheaves
    \[
        0\to\O_\Sigma\to \pi_*\O_{\widetilde{\Sigma}}\to \bigoplus_{i=1}^m \O|_{b_i}\to 0
    \]
 Taking the tensor product with $u^*N$\footnote{
 Recall that $N$ admits the $T^2$-action given in
 \cref{def:the right T2 action on N}
 } and then taking sheaf cohomology,
    we have:
    \begin{equation}
    \label{eq:sheaf cohomology long exact sequence}
        \begin{tikzcd}[column sep=small]
0 \arrow[r]
& H^0(\Sigma,u^*N) \arrow[r]
& H^0(\widetilde{\Sigma},\widetilde u^*N) \arrow[r,"\phi"]
& \displaystyle\bigoplus_{i=1}^m N|_{u(b_i)}
  \arrow[d]
\\
&&&
H^1(\Sigma,u^*N) \arrow[r]
& H^1(\widetilde{\Sigma},\widetilde u^*N) \arrow[r]
& 0
\end{tikzcd}
    \end{equation}

    We will identify some terms in this long exact sequence:
      \begin{itemize}
        \item $H^0(\Sigma,u^*N)=0$. 
        
        Recall $N=\O(-1)\oplus \O(-1)$.
        $u^*(\O(-1))$
        has negative degree on non-contracted components, and is the trivial line bundle on
        contracted components. Thus, a section must be zero on non-contracted components and constant on contracted components. Since $\Sigma$
        is connected, sections can only be zero.
        \item Similarly, $H^0(\widetilde{\Sigma},\widetilde u^*N)$ is identified with
        \[\bigoplus_{\widetilde{\Sigma}_j}N|_{u\left(\widetilde{\Sigma}_j\right)}\]
        where the sum runs over the contracted connected components
        $\widetilde{\Sigma}_j$, mapped to $p_0$ or $p_1$.
        These are genus-zero nodal curves, so
        \[
            H^0(\widetilde{\Sigma}_j,\O_{\widetilde{\Sigma}_j})=\C,
            \qquad
            H^1(\widetilde{\Sigma}_j,\O_{\widetilde{\Sigma}_j})=0.
        \]
    \end{itemize}
   Under this identification, the map 
   \[
    \phi \colon H^0\left(\widetilde{\Sigma}, \widetilde{u}^*\left(N\right)\right) 
    \to \bigoplus_{i=1}^m N|_{u(b_i)}
   \]
   becomes the restriction of sections from each contracted component
   to its attachment nodes.
Notice that $H^1\left(\Sigma, u^*\left(N\right)\right)$ is the obstruction bundle $\Ob_d$.
    By the exact sequence above, it has a subbundle $\coker(\phi)$.
   \begin{lemma}[{cf. the proof of \cite[Lemma~27.5.3]{hori2003mirror}}]
        \label{Euler class vanishing}
     Let $F$ be a component of the fixed locus. 
     Then $e_{\eq}(\Ob_d)|_F=0$ unless the domain has one of the following two forms: 
     \begin{itemize} 
         \item $\Sigma\cong\CP^1$, and $u\colon\Sigma\to\CP^1$ is a degree-$d$ cover totally ramified over $0$ and $\infty$. 
         \item $\Sigma$ consists of one contracted rational nodal curve $\Sigma_0$ with $h$ marked points, where $h$ is a power of $p$, mapped to either $0$ or $\infty$, together with $h$ non-contracted components $\Sigma_1,\ldots,\Sigma_h\cong\CP^1$ attached at these marked points. Each restriction $u|_{\Sigma_i}$ is a degree-$l$ cover totally ramified over $0$ and $\infty$. In particular, $d=h\cdot l$.
     \end{itemize}
   \end{lemma}
\begin{proof}
We refine the vanishing argument in the proof of
\cite[Lemma~27.5.3]{hori2003mirror} by keeping track of the finite
group presenting the fixed component.

Let
$N=\mathcal O_{\mathbb P^1}(-1)\oplus\mathcal O_{\mathbb P^1}(-1)$
and $G=\A_\Gamma$. For $\epsilon\in\{0,\infty\}$, let
$L_\epsilon\subset N|_{p_\epsilon}$ be the one-dimensional
$T^2$-fixed subspace determined by the \cref{def:the right T2 action on N}, which is used in
\cite[Lemma~27.5.3]{hori2003mirror}.

Let $B=\{b_1,\ldots,b_m\}$, and let $C$ range over the connected components
of $\widetilde{\Sigma}$. Throughout this proof, contracted components mean
the maximal connected contracted subcurves associated to vertices of $\Gamma$.
\eqref{eq:sheaf cohomology long exact sequence} gives
\begin{equation}
    \label{eq:Def-Ob exact sequence computation of Ob}
    \bigoplus_C H^0(C,u|_C^*N)
\xrightarrow{\phi}
\bigoplus_{b\in B}N|_{u(b)}
\longrightarrow
\operatorname{Ob}_d|_u\to
\bigoplus_C H^1(C,u|_C^*N)\to 0,
\end{equation}
and hence $\coker(\phi)\hookrightarrow (\Ob_d)_u$.
If $C$ is non-contracted, then $H^0(C,u|_C^*N)=0$. If $C_v$ is
contracted to $p_\epsilon$, with attachment nodes $B_v\subset B$, then on the
$T^2$-fixed line $L_\epsilon$ its contribution to $\phi$ is the
diagonal map
\[
L_\epsilon\longrightarrow L_\epsilon^{B_v},
\qquad
w\longmapsto (w)_{b\in B_v}.
\]

Thus it is enough to show, outside the two stated cases, that
$\operatorname{coker}(\phi)$ contains an equivariant subbundle
whose Euler class vanishes. Multiplicativity then gives
$e_{\mathrm{eq}}(\operatorname{Ob})|_F=0$. For the rest of this proof,
Steps 1 and 2 rule out two cases, and Step 3 shows that what remains
are those indicated by the lemma.

\smallskip
\noindent
\underline{Step 1: nodes joining two non-contracted components.}
Suppose that one node $b$ joins two non-contracted components. Neither
branch contributes to $H^0(C,u|_C^*N)$. Hence, writing $u(b)=p_\epsilon$,
the $L_\epsilon$-contributions at the nodes in $G\cdot b$ form a summand
$V_b$ of $\operatorname{coker}\phi$.

The action of the torus extension $\widetilde{\A}_\Gamma$ on $V_b$
factors through the graph action
$\rho:\widetilde{\A}_\Gamma\to H=\operatorname{Aut}(\Gamma)$, since
$\ker\rho$ fixes both branches and acts trivially on the weight-zero line
$L_\epsilon$. Since $\phi$ takes the difference of the branch values,
the stabilizer $H_b=\operatorname{Stab}_H(b)$ acts on the $b$-term by
its branch-permutation sign character $\bar\sigma_b$, and transitivity
on the orbit gives
\[
V_b\cong\rho^*W_b,
\qquad
W_b=\operatorname{Ind}_{H_b}^{H}\C_{\bar\sigma_b}.
\]
Over $\Mb_\Gamma$, this summand defines an equivariant subbundle of
$\operatorname{Ob}_d|_{F_\Gamma}$, pulled back from $W_b$ along
$[\Mb_\Gamma/\widetilde{\A}_\Gamma]\to BH$.
By \cref{lem:euler-induced-sign}, $e_E(W_b)=0$.
Consequently, multiplicativity yields
\[
e_{\mathrm{eq}}(\operatorname{Ob}_d)|_{F_\Gamma}=0.
\]

Hence we may assume that every node in $B$ is incident to a contracted
component.

\smallskip
\noindent
\underline{Step 2: more than one contracted component.}
Assume that there are at least two contracted components. Since
$\Gamma$ is a tree, choose a contracted
vertex $v$ which is a leaf of the minimal subtree containing all
contracted vertices. Among the branches at $v$, exactly one points
toward another contracted vertex, while, by stability, at least one
other branch does not.

Therefore the stabilizer $G_v=\operatorname{Stab}_G(v)$ does not act
transitively on the set $B_v$ of nodes incident to $v$. 
At each node in $B_v$, the contracted and non-contracted branches cannot
be exchanged, so the action on the node terms is the ordinary permutation
action.
It follows that
the reduced permutation representation
$\C[B_v]/\mathbb C\mathbf 1$\footnote{
 $B_v$ is a $G_v$-set, and
 $\C[B_v]$ means the $G_v$-representation spanned by the basis
 $B_v$. $\C\mathbf{1}$ is the vector
 $(1,1,\ldots,1)$, whose coordinate
 is $1$ for each basis vector.
}
contains a nonzero
$G_v$-invariant vector. 
Notice that we have the summand
inclusion
\[
    \iota_v\colon \C[B_v]\otimes N|_{u(v)}
    \hookrightarrow
    \bigoplus_{b\in B}N|_{u(b)}
\]
such that $\C\mathbf{1}\otimes N_{u(v)}=\iota_v^{-1}(\im \phi)$.
Tensoring with $L_{u(v)}$ and summing over $G/G_v$ gives a nonzero
$G$-invariant class in $\operatorname{coker}\phi$: different translates
are supported on disjoint attachment-node sets, and the relations
contributed by $\phi$ are diagonal within each such set. Since the torus
extension acts through the graph action on these weight-zero terms,
this class defines a nowhere-vanishing equivariant section over $F_\Gamma$.
Hence
\[
e_{\mathrm{eq}}(\operatorname{Ob}_d)|_{F_\Gamma}=0.
\]

Thus a possibly nonvanishing fixed locus has at most one contracted
component.

\smallskip
\noindent
\underline{Step 3: classification of the remaining graphs.}
If there is no contracted component, Step~1 implies that there are no
nodes. Hence $\Sigma\simeq\mathbb P^1$, and $u$ is a $T^2$-fixed
degree-$d$ cover, hence is totally ramified over $0$ and $\infty$.
This is the first exceptional case.

Suppose instead that there is exactly one contracted component $C_v$,
mapped to $p_\epsilon$. By Step~1 every node in $B$ is an attachment to $C_v$.
Let $B_v$ be these attachment nodes; then, as in Step~2,
there is a $T^2$-fixed subspace of $\coker \phi$
\[
W_v
=
L_\epsilon\otimes
\left(
\mathbb C[B_v]/\mathbb C\mathbf 1
\right).
\]
Since $G$ is finite, $W_v^G=0$ if and only if $G$ acts transitively on
$B_v$. Thus, unless the action on the branches is transitive, $W_v$
contains a trivial $G$-representation and the Euler class vanishes.

In the transitive case, graph automorphisms preserve the degree labels
of the non-contracted components. Hence all branches have the same
degree.

Set $h=|B_v|$ and
\[
V_{\mathrm{std}}=\mathbb C[B_v]/\mathbb C\mathbf1.
\]
Since $L_\epsilon$ has torus weight zero and deck transformations act
trivially on it, the equivariant subbundle defined by $W_v$ is pulled
back from $V_{\mathrm{std}}$ along
\[
[\Mb_\Gamma/\widetilde{\A}_\Gamma]\longrightarrow BS_h.
\]
Suppose that $h$ is not a power of $p$. A Sylow $p$-subgroup
$S\subset S_h$ then has at least two orbits on $B_v$, so
$V_{\mathrm{std}}|_S$ contains a trivial complex line. Thus
\[
\operatorname{res}^{S_h}_S e_E(V_{\mathrm{std}})=0.
\]
By the $p$-local stable Sylow splitting, restriction
$E^*(BS_h)\to E^*(BS)$ is injective, and hence
$e_E(V_{\mathrm{std}})=0$; see
\cite[Section~6.2]{Ragnarsson2008RetractiveTransfers}.
Pullback and multiplicativity give
\[
e_{\eq}(\Ob_d)|_{F_\Gamma}=0.
\]
Therefore, a possibly nonvanishing component must have $h$ a power
of $p$, giving the second exceptional case.
\end{proof}

\begin{lemma}[Euler class of an induced sign representation]
\label{lem:euler-induced-sign}
Let $E$ be a complex-oriented $p$-local cohomology theory, where $p$ is odd.
Let $H$ be a finite group, $K\subset H$ a subgroup, and
$\sigma:K\to\{\pm1\}$ a character. Then
\[
e_E\!\left(\operatorname{Ind}_K^H\C_\sigma\right)=0
\quad\text{in }E^*(BH).
\]
\end{lemma}
\begin{proof}
Choose a Sylow $p$-subgroup $S\subset H$. For each double-coset
representative $g\in S\backslash H/K$, the group $S\cap gKg^{-1}$
has odd order, so the conjugate character $\sigma^g$ restricts trivially
to it. Mackey decomposition formula gives
\[
\operatorname{Res}_S^H\operatorname{Ind}_K^H\C_\sigma
\cong
\bigoplus_{[g]\in S\backslash H/K}
\C\!\left[S/(S\cap gKg^{-1})\right].
\]
Each summand contains a trivial complex line, so the Euler class of
this restricted representation vanishes. But
\[
E^*(BH)\longrightarrow E^*(BS)
\]
is injective, and hence the original Euler class vanishes;
see \cite[Section~6.2]{Ragnarsson2008RetractiveTransfers}.
\end{proof}

By Lemma~\ref{Euler class vanishing}, it remains to analyze two types of fixed loci. We refer to the two geometric types as the \emph{single-edge component}
and the \emph{comb component}, respectively.

The computations for these two cases are taken from \cite{Graber-Pandharipande99} and \cite[Section 27]{hori2003mirror}.
\subsubsection{Case I: the single-edge component.}
\label{sec:single-edge component examination}
In this case, there is a unique
degree $d$ branched cover $u\colon \Sigma\to \P^1$
modelled by $z\mapsto z^d$ (using the coordinate
$z=z_1/z_0$ in \cref{def:convention for 0 and infty}) up to isomorphism, with $\A_\Gamma=C_d$. Not considering
the $T^2$-action, this component is isomorphic to $[\pt/C_d]$.

 However, if taking the $T^2$-action into
 consideration, then \cref{def:working definition of group action on orbifold} comes
 into play. The map
 $u\colon \Sigma\to \P^1$
    is modeled by $z\mapsto z^d$. 
    For $(t_0,t_1)\in T^2$, we can formally
    lift the action of $(t_0,t_1)$
    on $\P^1$ to $\Sigma$ by
    \[
        z\mapsto \left(\frac{t_1}{t_0}\right)^{1/d}
        z
    \]
    which is multivalued when $d>1$\footnote{
    In other words, the lift depends
    on the choice of a path (up to homotopy) from the unit element
    to $(t_0,t_1)$ in $T^2$.
    }.
    For $d>1$, this means the $T^2$-action on
    $[\pt/C_d]$ cannot be defined by
    a $C_d\times T^2$-action on $\pt$. 
    
    As $t=t_1/t_0\in S^1$ circulates counterclockwise along $S^1$,
    the lifted map is an isotopy from
    identity to $e^{{2\pi \i}/{d}}\in C_d$.
    This motivates us to define the following
    diagram of groups in which each row is exact (cf. sequence \eqref{eq:working definition exact sequence of groups} in \cref{def:working definition of group action on orbifold})
\begin{equation}
\label{eq:single-edge component working model}
    \begin{tikzcd}
	1 & {C_d} & {\widetilde{G}} & {T^2} & 1 \\
	1 & {C_d} & {S^1} & {S^1} & 1
	\arrow[from=1-1, to=1-2]
	\arrow[from=1-2, to=1-3]
	\arrow["\id", from=1-2, to=2-2]
	\arrow[from=1-3, to=1-4]
	\arrow[from=1-3, to=2-3]
	\arrow["\lrcorner"{anchor=center, pos=0.125}, draw=none, from=1-3, to=2-4]
	\arrow[from=1-4, to=1-5]
	\arrow["{(t_0,t_1)\mapsto t_1/t_0}", from=1-4, to=2-4]
	\arrow[from=2-1, to=2-2]
	\arrow[from=2-2, to=2-3]
	\arrow["{t\mapsto t^d}"', from=2-3, to=2-4]
	\arrow[from=2-4, to=2-5]
\end{tikzcd}
\end{equation}
    The RHS square is a pullback diagram.
    In the sense of \cref{def:working definition of group action on orbifold}, this defines
    a $T^2$ or $S^1$-action on $[\pt/C_d]$.
    The latter one comes from reducing
    the $T^2$-action on $\P^1$
    to the multiplication action 
    $S^1\acts \C\cup\{\infty\}=\P^1$,
    related via 
    \[
        T^2\to S^1,\quad (t_0,t_1)\mapsto
        t_1/t_0
    \]
    This gives the model for the $T^2$ or $S^1$-action on the single-edge component $[\pt/C_d]$, depending on which action we choose on $\P^1$. 
    The quotient stack $[\pt/C_d]/T^2 $ (resp., $[\pt/C_d]/S^1$)
    is actually $[\pt/\widetilde{G}]$ (resp., $[\pt/S^1]$\footnote{
    This stands for the middle $S^1$ in the
    second row of \eqref{eq:single-edge component working model}.
    }). The restrictions of $\Ob_d$ and $\nu_F$ to the single-edge
    component are thus equivalently $\widetilde{G}$ (resp., $S^1$)
    representations.
    
    Let $\alpha_0,\alpha_1\in H^2(BT^2,E)$ and $x\in 
    H^2(BS^1,E)$ (the middle $S^1$) be the hyperplane classes.
    We also denote their pullbacks to
    $H^2(B\widetilde{G},E)$ by $\alpha_0,\alpha_1,x$, and the same
    for restriction of $x$ to $H^2(BC_d,E)$
    along $C_d\hookrightarrow S^1$. Recall that 
    $[n]_E x:=\underbrace{x +_E x +_E \dots +_E x}_{n \text{ times}}$
    for $n>0$, and $[-n]_E(x)$ is its inverse with regard to $F_E(-,-)$.
\begin{lemma}[Essentially {\cite[Exercise 27.2.6\&27.4.1]{hori2003mirror}}]
\label{lem:virtual-euler-class-case1}
        The Euler classes\footnote{
        Since the $T^2$-quotient of $[\pt/C_d]$ is $[\pt/\widetilde{G}]$,
        the $T^2$-equivariant Euler classes lie in $H^*(B\widetilde{G},E)$.
        } of $\Ob_d$ and $\nu_F$ are
        \[
            e_{\eq}(\Ob_d)=\prod_{i=1}^{d-1}[i]_E(x)\cdot [-i]_E(x)
        \]
        and
        \[
            e_{\eq}(\nu_F)=\prod_{i=2}^{d}[i]_E(x)\cdot [-i]_E(x)
        \]
        where we view $x$ as a class in $H^2(B\widetilde{G},E)$ or $H^2(BS^1,E)$.
        Thus, after Euler localization with respect to $T^2$ (or the target
        $S^1$), we have
        \[
        \frac{e_{\eq}(\Ob_d)}{e_\eq(\nu_F)}
        =\frac{x\cdot [-1]x}{[d]x\cdot [-d]x }.
        \]
\end{lemma}
    Since there is only one point in the
    (coarse moduli of) single-edge component, it suffices
    to identify $\Ob_d$ and $\nu_F=T\Mb_{0,0}\left(\P^1,d\right)$
    as $\widetilde{G}$ or $S^1$-representations. 
    These weight computations (i.e. the proof
    of \cref{lem:virtual-euler-class-case1}) are taken from
    \cite[Section 27.3--27.4]{hori2003mirror}.
    We only present them here for readers'
    convenience.
\begin{proof}
    W.L.O.G. we consider the $S^1$-equivariant
    version, as the $\widetilde{G}$-equivariant
    case is just a pullback. We distinguish two $S^1$s 
    as $S^1_\source$ which acts on $\Sigma$ (the middle $S^1$
    in \eqref{eq:single-edge component working model})
    and $S^1_\target$ which acts on $\P^1$ (the right
    $S^1$ in \eqref{eq:single-edge component working model}).
    
    We use 
    \cite[Exercise 27.2.3]{hori2003mirror}(\cref{lem:from pole weights to cohomology weights}), which enables getting
    the 
    weights of cohomology of an equivariant line bundle 
    from its weights of fibers at $0,\infty$. Note that the $S^1_\target$-weights of
    $\O(-1)_1,\O(-1)_2$ at $0,\infty$ are discussed in
    \cref{rmk:S1 action convention on O(-1)oplusP(-1)}.
    \begin{itemize}
        \item For $\Ob_d\simeq H^1\left(\Sigma,u^* N\right)$, by
        Serre duality, we have 
        \[
            \Ob_d\simeq 
            H^0\left(
            \Sigma,  u^*(\O(-1)_1)^\vee\otimes\omega_\Sigma       
            \right)^\vee
            \oplus
            H^0\left(
            \Sigma,  u^*(\O(-1)_2)^\vee\otimes\omega_\Sigma       
            \right)^\vee
        \]
        The $S^1_\source$-weights of relevant bundles
        at $0,\infty$ are as follows:
   \[
\begin{array}{c|cc}
& 0 & \infty \\ \hline
u^*\mathcal{O}(-1)_1^\vee & 1 & \chi^{\otimes d} \\
u^*\mathcal{O}(-1)_2^\vee & \chi^{\otimes (-d)} & 1 \\
T_\Sigma & \chi^{-1} & \chi \\
\omega_\Sigma & \chi & \chi^{-1} \\
u^*\mathcal{O}(-1)_1^\vee \otimes \omega_\Sigma
& \chi & \chi^{\otimes(d-1)} \\
u^*\mathcal{O}(-1)_2^\vee \otimes \omega_\Sigma
& \chi^{\otimes(1-d)} & \chi^{-1}
\end{array}
\]
Note that $\deg_\Sigma u^*(\O(-1)_1)^\vee\otimes\omega_\Sigma
=\deg_\Sigma u^*(\O(-1)_2)^\vee\otimes\omega_\Sigma
=d-2$.
By \cref{lem:from pole weights to cohomology weights},
$H^0\left(\Sigma,  u^*(\O(-1)_1)^\vee\otimes\omega_\Sigma       
\right)$ is a direct sum of weights
            $\chi, \chi^{\otimes2},\ldots, \chi^{\otimes{d-1}}$.

            The space $H^0\left(\Sigma,  u^*(\O(-1)_2)^\vee\otimes\omega_\Sigma       \right)$ is a direct sum of
            $\chi^{-1},\chi^{\otimes(-2)},\ldots,\chi ^{\otimes(1-d)}$.
            Since $x=e(\chi)$ in $H^2(BS^1,E)$,
            the Euler class $e_{\eq}(\Ob_d)$ is
            $\prod_{i=1}^{d-1}[i]_E(x)\cdot [-i]_E(x)$.
        \item For $\nu_F$,
        there is an exact sequence
    \[
        0\to H^0(\Sigma,T_\Sigma) \to H^0\left(\Sigma,u^*T_{\P^1}\right)
        \to \nu_F\to \Def\Sigma=0 
    \]
    The $S^1_\source$-weights over $0,\infty$ are
   \[
\begin{array}{c|cc}
& 0 & \infty \\ \hline
T_\Sigma
& \chi^{-1} & \chi \\[2pt]
u^*T_{\mathbb P^1}
& \chi^{\otimes(-d)} & \chi^{\otimes d}
\end{array}
\]
Thus (by \cref{lem:from pole weights to cohomology weights}) $H^0(\Sigma,T_\Sigma)$ is a direct sum of $\chi,1,\chi^{-1}$,
and $H^0(\Sigma, u^*T_{\P^1})$ is a direct sum of
$\chi^{\otimes(-d)},\chi^{\otimes(-d+1)},\ldots,\chi^{\otimes (d-1)},
\chi^{\otimes d}$. Thus $\nu_F$, as the quotient, is
a direct sum of $\chi^{\otimes(-d)},\chi^{\otimes(-d+1)},\ldots,\chi^{\otimes (-2)}, \chi^{\otimes 2},\ldots,\chi^{\otimes (d-1)},
\chi^{\otimes d}$.
The formula of $e_{\eq}(\nu_F)$ follows.
\end{itemize}   
\end{proof}
The cited lemma in this proof is the following.
\begin{lemma}[{\cite[Exercise 27.2.3]{hori2003mirror}}]
\label{lem:from pole weights to cohomology weights}
    Let $\Sigma=\P^1$ with weight $(1,1)$ $T^2$-action. If a line bundle $\O_{\Sigma}(l)$
    with a $T^2$-action has weights $W_0,W_1$ at the two fixed points $p_0,p_1$, then, for $l>0$, the weights of the $T^2$-action on $H^0(\Sigma,\O_{\Sigma}(l))$ are \[
        \frac{i W_0+(l-i)W_1}{l}
    \]
    for $i=0,\ldots,l$.
    For $l=0$, $H^0(\Sigma,\O_{\Sigma}(l))$ is one-dimensional of weight $W_0=W_1$;
    for $l<0$, it is zero.
\end{lemma}
Note that we can safely replace $T^2$-action
by $S^1$-action in \cref{def:the right T2 action on N}.
\subsubsection{Case II: the comb components.}
\label{sec:comb components examination}
For the comb component, recall $h$ is the number of attachment points
of the central contracted subcurve, and $l$ is the degree
of each edge. We call the central subcurve $\Sigma_0$ with
marked points $b_1,b_2,\ldots,b_h$, and call the curve corresponding
to each edge $\Sigma_i,i=1,2\ldots,h$. $\Sigma_i$ is attached to
$\Sigma_0$ at $b_i\in \Sigma_0$. We use $S_h$ to denote the $h$-th permutation group.

The domain of a stable map in the comb component is determined
by the central curve $[\Sigma_0]\in \Mb_{0,h}$. 
The sequence \eqref{eq:A_Gamma extension sequence} takes the form
\[
    1\to (C_{l})^h\to G\to S_h\to 1
\]
It is interpreted as:
\begin{itemize}
    \item $(C_{l})^h$ is the group of automorphisms
of the domain curve that \emph{fixes the nodes $b_i$},
corresponding to rotating $\Sigma_i$ by $2\pi/l$.
    \item Also, there is a global action of $S_h$ on $\Mb_{0,h}$ by permuting
marked points.
\end{itemize}
We also write $G=(C_{l})^h\rtimes S_h$. 
Let $M=\Mb_{0,h}$,
then
the comb component is isomorphic to 
\begin{equation}
\label{eq:model for comb components}
    [M/G]=\left[\Mb_{0,h}/\left((C_{l})^h\rtimes S_h\right)\right].
\end{equation}
where $G$ acts on $\Mb_{0,h}$ through $G\to S_h$. 
Taking the $\TT=T^2$ or $\TT=S^1$-action into consideration, we have (\cref{def:working definition of group action on orbifold})
an exact sequence
\begin{equation}
    \label{eq:define widetilde G}
    1\to G\to \widetilde{G}\to \TT\to 1
\end{equation}
\begin{remark}
\label{rmk:semi-direct product and Sn action}
    For $\TT=S^1$, set
    \[
        A_l:=\left\{
            (z_1,\ldots,z_h;t)\in (S^1)^h\times S^1
            \mid z_1^l=\cdots=z_h^l=t
        \right\},
    \]
    with coordinatewise multiplication. Then $\widetilde{G}=A_l\rtimes S_h$,
    where $S_h$ permutes the coordinates $z_i$, and we have a split exact sequence
    \begin{equation}
        \label{eq:another way to define widetildeG}
        1\to A_l\to \widetilde{G}\to S_h\to 1.
    \end{equation}
    The projection $\widetilde{G}\to S^1$ sends $(z_1,\ldots,z_h;t;\sigma)$
    to $t$; the corresponding extension by $G$ does not split when $l>1$.
    The action on $\Mb_{0,h}$ factors through $\widetilde{G}\to S_h$.

    Similarly, for $\TT=T^2$, set
    \[
        A_l^{T^2}:=\left\{
            (z_1,\ldots,z_h;t_0,t_1)\in (S^1)^h\times T^2
            \mid z_1^l=\cdots=z_h^l=t_1/t_0
        \right\}.
    \]
    Then $\widetilde{G}=A_l^{T^2}\rtimes S_h$, with projection to $T^2$
    given by $(t_0,t_1)$.
\end{remark}
\begin{lemma}
\label{lem:virtual-euler-class-case2}
     Euler classes of $\Ob_d$ and $\nu_F$ are as follows
    \[
        e_{\eq}(\Ob_d)=e_\eq\left(N|_{p_\epsilon}\otimes V_{\std}\right)
        \cdot 
        e_\eq\left(\bigoplus_{i=1}^h H^1\left(\Sigma_i, u|_{\Sigma_i}^*N\right)\right)
    \]
    \[
        e_{\eq}(\nu_F)=e_\eq\left(\bigoplus_{j=1}^h 
        (T_{b_j}\Sigma_j\otimes T_{b_j} \Sigma_0)\right)
        \cdot e_\eq\left(\bigoplus_{i=1}^h E_i\right)
        \cdot e_\eq (T_{p_\epsilon}\P^1).
    \]
    Here $V_{\std}=\left(\bigoplus_{i=1}^h \C b_i\right)/\C\mathbf{1}$ is the standard representation of $S_h$,
    and 
    \begin{equation}
        \label{eq:definition of Ei}
         E_i=\coker\left(
        H^0(\Sigma_i, T\Sigma_i(-b_i))
        \to H^0\left(\Sigma_i, u_i^* T_{\P^1}(-b_i)\right)
        \right)
    \end{equation}
\end{lemma}
\begin{remark}
    These classes lie in $H^*\left((\Mb_{0,h})_{h\widetilde{G}},E\right)$. 
\end{remark}
\begin{proof}
\begin{itemize}
    \item For $\Ob_d$, from the deformation-obstruction exact
    sequence \eqref{eq:Def-Ob exact sequence computation of Ob}
    we have
    \[
        H^0\left(\Sigma_0,u|_{\Sigma_0}^* N\right)\to 
        \bigoplus_{i=1}^h N|_{u(b_i)}\to
        \Ob_d
        \to \bigoplus_{i=1}^h H^1\left(
        \Sigma_i, u|_{\Sigma_i}^* N
        \right)
        \to 0
    \]
    The first term is just $N|_{p_\epsilon}$ because $\Sigma_0$ is contracted.
    The second term can be identified with the tensor product of $N|_{p_\epsilon}$ and
    the permutation representation of $S_h$ on the $h$ attachments. Thus,
    the quotient of the second term by the image of the first term is
    $N|_{p_\epsilon}\otimes V_{\std}$.
    \item For $\nu_F$,
    use the 
    deformation exact sequence
   \[
\begin{tikzcd}[column sep=1.5em]
0
\arrow[r] &
\Aut\Sigma
\arrow[r] &
H^0(\Sigma,u^*T_{\P^1})
\arrow[r] &
T_{[\Sigma]}\Mb_{0,0}(\P^1,d)
\arrow[r] &
\Def\Sigma
\arrow[r] &
H^1(\Sigma,u^*T_{\P^1})
\arrow[r] &
0 .
\end{tikzcd}
\]
The normal bundle $\nu_F$ is the $T^2$-moving part of 
$T_{[\Sigma]}\Mb_{0,0}(\P^1,d)$.

We have the following computations:
\begin{enumerate}[(a)]
    \item $\coker\left(\Aut\Sigma\to H^0(\Sigma,u^* T_{\P^1})\right)
    =T_{p_\epsilon}\P^1+\sum_{i=1}^h E_i$ in equivariant orbifold $K$-theory. In fact, consider the normalization sequence for $H^0(\Sigma,u^* T_{\P^1})$
    \[
        0\to \bigoplus_{i=1}^h H^0(\Sigma_i, u_i T_{\P^1}(-b_i))
        \to H^0(\Sigma,u^*T_{\P^1})
        \to T_{p_\epsilon}\P^1\to 0
    \]
    and also 
    \[
        \Aut(\Sigma)=\bigoplus_{i=1}^h H^0(\Sigma_i, T_{\Sigma_i}(-b_i)).
    \]
    we get the exact sequence:
    \[
        0\to \bigoplus_{i=1}^h E_i\to \coker\left(\Aut\Sigma\to H^0(\Sigma,u^* T_{\P^1})\right)\to T_{p_{\epsilon}}\P^1.
    \]
    \item $H^1(\Sigma,u^*T_{\P^1})=0$. We omit the proof.
    \item The moving part of $\Def(\Sigma)$ is exactly
    \[
        \bigoplus_{j=1}^h \left(
        T^*_{b_j}\Sigma_j\otimes T_{b_j}^* \Sigma_0
        \right)^\vee
    \]
    Here $b_j$ are the nodes between the central subcurve and
    the $h$ tails. It corresponds to smoothing of nodes\footnote{
    The $T^2$-fixed part of $\Def(\Sigma)$ is the deformation space of
    the central stable $h$-pointed curve.
    This is irrelevant to $T^2$-action because
    $T^2$ only acts on non-contracted spheres, i.e. the $h$
    tail components.
    }.
\end{enumerate}
\end{itemize}
Therefore, in $K$-theory we have
\[
    [\nu_F]=\left[\left(\oplus_{j=1}^h 
        T_{b_j}\Sigma_j\right)\otimes T_{b_j} \Sigma_0
        \right]+[T_{p_\epsilon}\P^1]
        +\sum_{i=1}^h E_i.
\]
and the Euler class is the product.
\end{proof}
We next study vanishing for nontransitive HKR sectors;
see \cref{rmk:definition of p-inertia stack}(1) for their definition.
\begin{lemma}
\label{lem:Euler class vanishing after HKR}
    For an HKR sector $\alpha\colon \Lambda_r\to (C_l)^h\rtimes S_h$, the class $e_{\operatorname{HKR},\eq}^L(\Ob_d)$ is zero
    unless $\Lambda_r$ acts transitively on the marked
    points set via $\alpha\colon \Lambda_r\to (C_l)^h\rtimes S_h\
    \twoheadrightarrow S_h$.
\end{lemma}
\begin{proof}
    Using the exact sequence \eqref{eq:Def-Ob exact sequence computation of Ob}, we get a subbundle of $\Ob_d$ 
    \[
       \coker\phi=\left(\left(\bigoplus_{i=1}^h\C b_i\right)/ \C\mathbf{1}\right)\otimes N|_{p_\epsilon}
    \] 
    where $p_\epsilon\in \P^1$ is
    the point that the central subcurve lies over.
    Since $\coker\phi$ is trivial over $\Mb_{0,h}$, we
    can view $\coker\phi$ as a $\widetilde{G}$-representation,
    or equivalently a bundle over $[\pt/\widetilde{G}]$, and
    $\coker\phi\to [M/\widetilde{G}]$ is given by the pullback
    along $[M/\widetilde{G}]\to [\pt/\widetilde{G}]$. From now
    on, we will forget $M$ and only consider $\coker \phi$ as
    a representation.
    
    Let's recall definition of $e_{\operatorname{HKR},\eq}^L$. 
    By \cref{equivariant HKR character}, we first view $\coker \phi$
    as a $Z_{\widetilde{G}}(\alpha)\times \Lambda_r$ representation.
    We get its Euler class in $H^*\left(BZ_{\widetilde{G}}(\alpha)\times B\Lambda_r,E\right)$, and then send this class to $H^*(BZ_{\widetilde{G}}(\alpha),L(E))$ to get $e_{\operatorname{HKR},\eq}^L$.

    If $\Lambda_r$ acts nontransitively on the set of marked points,
    define $W=\left(\bigoplus_{i=1}^h\C b_i/ \C\mathbf{1}\right)^{\Lambda_r}$ to be the vectors fixed by $\Lambda_r$. $W$ is a $Z_{\widetilde{G}}(\alpha)$ representation
    since $Z_{\widetilde{G}}(\alpha)$ commutes with $\im\alpha$.
    Let $c$ be the number of $\Lambda_r$-orbits on the $h$ marked points.
    The invariant vectors in $\bigoplus_{i=1}^h\C b_i$ are constant on each orbit.
    Since taking invariants is exact over $\C$, quotienting by $\C\mathbf{1}$ gives
    $\dim_\C W=c-1$. Nontransitivity implies $c\ge2$, so $W$ is nonzero.
    Denote by $\O(-1)_\epsilon$ the suitable
    $\O(-1)$ summand in $N$ whose $p_\epsilon$ fiber has trivial $\TT$ weight, we get a 
    $Z_{\widetilde{G}}(\alpha)\times \Lambda_r$ sub-representation
    \[
        W\otimes \O(-1)_\epsilon|_{p_\epsilon}\subset \coker \phi.
    \]
    This sub-representation is also nonzero, since tensoring with a line preserves dimension.
    We examine the $Z_{\widetilde{G}}(\alpha)\times \Lambda_r$ action
    on it:
    \begin{itemize}
        \item $\Lambda_r$ acts trivially on $W\otimes \O(-1)_\epsilon|_{p_\epsilon}$
        by definition.
        \item The identity connected component $Z_{\widetilde{G}}(\alpha)^\circ $ of $Z_{\widetilde{G}}(\alpha)$ acts
        trivially, too. This is because the map 
        $Z_{\widetilde{G}}(\alpha)\subset \widetilde{G}\to \TT$
        is an isomorphism on Lie algebras, and we choose $\O(-1)_\epsilon$
        with trivial $\TT$ weight at $p_\epsilon$.
    \end{itemize}
    Therefore, let $K=Z_{\widetilde{G}}(\alpha)/Z_{\widetilde{G}}(\alpha)^\circ$ which is a finite group, then the $Z_{\widetilde{G}}(\alpha)\times \Lambda_r$-representation on $W\otimes \O(-1)_\epsilon|_{p_\epsilon}$ is a pullback of 
    a $K$-representation on $W\otimes \O(-1)_\epsilon|_{p_\epsilon}$
    along the projection $Z_{\widetilde{G}}(\alpha)\times \Lambda_r\to Z_{\widetilde{G}}(\Lambda_r)\to K$. Consider the following
    commutative diagram
\[\begin{tikzcd}
	{H^*\left(BZ_{\widetilde{G}}(\alpha)\times B\Lambda_r,E\right)} & {H^*\left(BZ_{\widetilde{G}}(\alpha),L(E)\right)} \\
	{H^*\left(BK\times B\Lambda_0,E\right)} & {H^*(BK,L(E))}
	\arrow[from=1-1, to=1-2]
	\arrow[from=2-1, to=1-1]
	\arrow[from=2-1, to=2-2]
	\arrow[from=2-2, to=1-2]
\end{tikzcd}\]
Starting with $e(W\otimes \O(-1)_\epsilon|_{p_\epsilon})\in H^*(BK,E)$ in the
bottom left, there are two directions:
\begin{itemize}
    \item Going up and then right: the image in $H^*(BZ_{\widetilde{G}}(\alpha)\times B\Lambda_r,E)$ is again an Euler class, and since
    $W\otimes \O(-1)_\epsilon|_{p_\epsilon}\subset \coker \phi$ is a $Z_{\widetilde{G}}(\alpha)\times \Lambda_r$ sub-representation, $e(W\otimes \O(-1)_\epsilon|_{p_\epsilon})$ is a factor of $e(\coker\phi)$, both
    in $H^*(BZ_{\widetilde{G}}(\alpha)\times B\Lambda_r,E)$ and
    $H^*(BZ_{\widetilde{G}}(\alpha),L(E))$.
    \item Going right and then up: 
    by \cref{lem:pullback of finite group representations vanish in rational theory} we get zero in $H^*(BZ_{\widetilde{G}}(\alpha),L(E))$.
\end{itemize}
Combining these two, we get $e_{\operatorname{HKR},\eq}^L(\coker\phi)$ is zero. Since $\coker\phi\subset \Ob_d$ is a $\widetilde{G}$ sub-representation, we get $e_{\operatorname{HKR},\eq}^L(\Ob_d)=0$
for a similar reason.
\end{proof}
The lemma used during the proof will be useful later.
\begin{lemma}
\label{lem:pullback of finite group representations vanish in rational theory}
    Let $\phi\colon G\to K$ be a homomorphism from a compact Lie group $G$ to
    a \emph{finite} group K. Let $V$ be a nonzero complex $K$-representation.
    Then the Euler class $e(\phi^* V)\in H^*(BG,E)$ is mapped
    to $0$ in $H^*(BG,L(E))$.
\end{lemma}
\begin{proof}
    Since $K$ is finite and $L(E)$ is rational, restriction to a point induces
    an isomorphism $H^*(BK,L(E))\to L(E)^*$.
    The restriction of $V$ is a nonzero trivial complex vector bundle over a point,
    whose Euler class is zero. Hence $e(V)=0$ in $H^*(BK,L(E))$.
    Pulling back along $B\phi$ proves the claim.
\end{proof}

The bundles on the RHS of \cref{lem:virtual-euler-class-case2}
carry a complicated $\widetilde{G}=A_l\rtimes S_h$-action:
$A_l$ acts on the individual tail contributions, while $S_h$ permutes them.
However, after applying the HKR character map, we can use the simpler
structure of $Z_{\widetilde{G}}(\alpha)$ described below to identify
the relevant representations.

For $\alpha\colon \Lambda_r\to (C_l)^h\rtimes S_h$, let $P\subset S_h$
be its image in $S_h$. 
\begin{lemma}
\label{lem:identifying centralizer}
    Assume $\alpha$ acts transitively on $h$ marked points. Then:
\begin{enumerate}[(1)]
    \item we have a diagram where each row is exact
        \[\begin{tikzcd}
	1 & {S^1} & {Z_{\widetilde{G}}(\alpha)} & P & 1 \\
	1 & {A_l} & {\widetilde{G}} & {S_h} & 1
	\arrow[from=1-1, to=1-2]
	\arrow[from=1-2, to=1-3]
	\arrow[hook, from=1-2, to=2-2]
	\arrow[from=1-3, to=1-4]
	\arrow[hook, from=1-3, to=2-3]
	\arrow[from=1-4, to=1-5]
	\arrow[hook, from=1-4, to=2-4]
	\arrow[from=2-1, to=2-2]
	\arrow[from=2-2, to=2-3]
	\arrow[from=2-3, to=2-4]
	\arrow[from=2-4, to=2-5]
\end{tikzcd}\]
The inclusion $S^1\hookrightarrow A_l$ is
$z\mapsto (z,\ldots,z;z^l)$.
\item The action $P\acts \{b_1,b_2,\ldots,b_h\}$
is free and transitive. In particular,
$|P|=h$.
\end{enumerate}
\end{lemma}
\begin{proof}
\begin{enumerate}[(1)]
    \item Since $\Lambda_r$ is commutative, $P$ is contained
        in the image of $Z_{\widetilde{G}}(\alpha)$. Since $\alpha$
        acts transitively on $h$ points, any centralizer of $P$
        inside $S_h$ must lie in $P$, which means image of $Z_{\widetilde{G}}(\alpha)$ is exactly $P$.

        By the product rule of the semi-direct product $A_l\rtimes S_h$,
        the elements of $Z_{\widetilde{G}}(\alpha)\cap A_l$ are precisely
        those fixed by $P$ under the permutation action
        (\cref{rmk:semi-direct product and Sn action}).
        Since $P$ acts transitively, we have
        \[
            Z_{\widetilde{G}}(\alpha)\cap A_l
            =\{(z,\ldots,z;z^l):z\in S^1\}
            \simeq S^1.
        \]
        \item That $P$-action is transitive is our assumption.
        Let $\sigma\in P$ and let $b_i$ be one of the marked points,
        and assume $\sigma(b_i)=b_i$. Then for any $\tau\in P$,
        since $P$ is abelian, we have $\sigma(\tau(b_i))=\tau(\sigma(b_i))
        =\tau(b_i)$, i.e. $\tau(b_i)$ is a fixed point of $\sigma$.
        Since $P$ acts transitively, we have $\sigma=\id$. This means
        the $P$-action is free.
\end{enumerate}    
\end{proof}
On an HKR sector $\alpha$, we can represent classes in \cref{lem:virtual-euler-class-case2} explicitly as pullbacks of representations.
Let's fix some notations of representations and vector bundles
\begin{notation}
    \begin{enumerate}[(1)]
        \item As in case I, we view the $S^1$-representation
        $\chi$\footnote{
        Recall $\chi$ is multiplying by $t^{-1}$; see \cref{rmk:S1 action convention on O(-1)oplusP(-1)}.
        }as a $Z_{\widetilde{G}}(\alpha)$-representation
        which is the pullback along $Z_{\widetilde{G}}(\alpha)\subset\widetilde{G}\to S^1$.
        \item For the $S^1\subset Z_{\widetilde{G}}(\alpha)$,
        we use $V_i,i\in \Z$ to represent the $S^1$-representation
        of weight $i$.
        \item We view $V_\std$ as a 
        $Z_{\widetilde{G}}(\alpha)$-representation
        via $Z_{\widetilde{G}}(\alpha)\to P\subset S_h$
    \end{enumerate}
\end{notation}
\begin{remark}
    There are two $S^1$ involved: one appears in $\widetilde{G}\to S^1$
    which is acting on the target $\P^1$; another appears in $S^1\to Z_{\widetilde{G}}(\alpha)$ which acts on the tail spheres $\Sigma_i,1\le i\le h$. We denote representations of the first $S^1$ by $\chi^j$,
    and denote representations of  the second $S^1$ by
    $V_j$. The composition $S^1\to Z_{\widetilde{G}}(\alpha)\subset
    \widetilde{G}\to S^1$ is of degree $l$, as the stable map
    has degree $l$ on each tail $\Sigma_i\to \P^1$ 
    (so that two $S^1$-actions are compatible).
\end{remark}
As in case I, we view $\chi$ as a
$Z_{\widetilde{G}}(\alpha)$-representation which
is the pullback $Z_{\widetilde{G}}(\alpha)\subset \widetilde{G}
\to S^1$, 

\begin{lemma}
\label{lem:case2 virtual bundle reduction after HKR}
    Restricted to a transitive HKR sector $[\Mb_{0,h}^{\im\alpha}/Z_{\widetilde{G}}(\alpha)]$ of comb components whose central subcurve is over
    $p_\epsilon$, we have in $K$-theory:
    \begin{enumerate}[(1)]
        \item $\Ob_d=V_\std\otimes (\chi^{-\epsilon}+\chi^{1-\epsilon})
        +\operatorname{Ind}_{S^1}^{Z_{\widetilde{G}}(\alpha)}\left(\bigoplus_{i=1}^{l-1}V_i\oplus V_{-i}\right)$.
        \item $\nu_F=\bigoplus_{j=1}^h 
        \left(T_{b_j}\Sigma_j\otimes T_{b_j} \Sigma_0\right)+ \operatorname{Ind}_{S^1}^{Z_{\widetilde{G}}(\alpha)}\left(\sum_{i=1}^{l-1}V_{\eta i}+\sum_{i=2}^l V_{-\eta i}\right)
        +\chi^{\eta}$,
         where $\eta=2\epsilon-1$.
         \item The difference satisfies
         \[
             (\Ob_d-\nu_F)|_{[\Mb_{0,h}^{\im\alpha}/Z_{\widetilde{G}}(\alpha)]}
             =V_\std+
             \operatorname{Ind}_{S^1}^{Z_{\widetilde{G}}(\alpha)}
        \left(
            V_{-\eta}
        \right)
        -\bigoplus_{j=1}^h 
        \left(T_{b_j}\Sigma_j\otimes T_{b_j} \Sigma_0\right)
        -\chi-\chi^{-1}.
         \]
    \end{enumerate}
\end{lemma}
\begin{proof}
It suffices to identify the $Z_{\widetilde{G}}(\alpha)$-representation
on every bundle in \cref{lem:virtual-euler-class-case2}.
    \begin{enumerate}[(1)]
        \item The first term follows from the fact that
        the weight of $N|_{p_\epsilon}$ is $\chi^{-\epsilon}\oplus\chi^{1-\epsilon}$.

        For the second term,
        in view of \cref{lem:identifying centralizer} and $|P|=h$
        acting regularly on $h$-marked points, one can write
        \[
            \bigoplus_{i=1}^h H^1\left(\Sigma_i, u_i^*N\right)
            \simeq
            \operatorname{Ind}_{S^1}^{Z_{\widetilde{G}}(\alpha)}
            \left(
            H^1\left(\Sigma_1, u_1^*N\right)
            \right),
        \]
        viewing $H^1\left(\Sigma_1, u_1^*N\right)$ as a $S^1$-representation. Its weight computation is the same
        as that of the obstruction bundle in the proof of \cref{lem:virtual-euler-class-case1}, which is
        $\sum_{i=1}^{l-1}\left(V_i+V_{-i}\right)$.
        \item The first term is intact.

        For the second term, we also have
            $\bigoplus_{i=1}^h E_i=\operatorname{Ind}_{S^1}^{Z_{\widetilde{G}}(\alpha)}
            \left(E_1\right)$, where
        $E_1$ is defined in \eqref{eq:definition of Ei}.
        Using the coordinate $z$ on $\Sigma_1$ and $w$ on $\P^1$
        with $w=z^l$, we claim
        \begin{itemize}
            \item $H^0(\Sigma_1, T\Sigma_1(-b_1))$
            is of rank $2$ with generators $z\partial_z,z^2\partial_z$.
            The $S^1$-weights are $V_0,V_{-\eta}$ respectively.
            \item $H^0\left(\Sigma_1, u_1^* T_{\P^1}(-b_1)\right)$
            has rank $2l$ with generators
            $z\partial_w,z^2\partial_w,\ldots, z^{2l}\partial_w$.
            The $S^1$-weight of $z^k\partial_w$ is $V_{\eta(l-k)}$.
            \item The map $du_1$ sends 
            \[
                z\partial_z\mapsto l\cdot z^l\partial_w
            \]
            \[
                z^2\partial_z\mapsto l\cdot z^{l+1}\partial_w
            \]
        \end{itemize}
        Therefore, in $E_1$ there are only weights
        $V_{\eta j}$ with $j\in \{-l,-l+1\,\ldots,-2\}\cup 
        \{1,2,\ldots,l-1\}$. The second term follows.

        The last term $T_{p_\epsilon}\P^1$ has weight $\chi^\eta$.
        \item Subtracting (2) from (1), we get directly
        \[
        \begin{aligned}
        (\Ob_d-\nu_F)|_{[\Mb_{0,h}^{\im\alpha}/Z_{\widetilde{G}}(\alpha)]}
        &=V_{\std}\otimes(\chi^{-\epsilon}+\chi^{1-\epsilon})
        -\bigoplus_{j=1}^h
        \left(T_{b_j}\Sigma_j\otimes T_{b_j} \Sigma_0\right)\\
        &\quad+\operatorname{Ind}_{S^1}^{Z_{\widetilde{G}}(\alpha)}
        \left(
            V_{-\eta}-V_{-l\cdot \eta}
        \right)
        -\chi^\eta.
        \end{aligned}
        \]
        We claim that $\operatorname{Ind}_{S^1}^{Z_{\widetilde{G}}(\alpha)}(V_{-l\cdot\eta})=\chi^{-\eta}\left(1+V_{\std}\right)$.
        Using projection formula and $\operatorname{Res}_{S^1}^{Z_{\widetilde{G}}(\alpha)}(\chi)=V_l$:
        \begin{align*}
            \operatorname{Ind}_{S^1}^{Z_{\widetilde{G}}(\alpha)}(V_{-l\cdot\eta})\otimes\chi^\eta
            &\simeq 
            \operatorname{Ind}_{S^1}^{Z_{\widetilde{G}}(\alpha)}
            \left(V_{-l\cdot\eta}\otimes \operatorname{Res}_{S^1}^{Z_{\widetilde{G}}(\alpha)}(\chi^\eta)\right)
            \\
            &\simeq \operatorname{Ind}_{S^1}^{Z_{\widetilde{G}}(\alpha)}
            \left(V_{-l\cdot\eta}\otimes V_{l\cdot\eta}\right)
            \\
            &\simeq\operatorname{Ind}_{S^1}^{Z_{\widetilde{G}}(\alpha)}(1)
            \\
            &\simeq 1+V_\std .    
        \end{align*}
        Plug this into the original difference ($\eta=2\epsilon-1$): 
        \[
        V_{\std}(\chi^{-\epsilon}+\chi^{1-\epsilon}-\chi^{1-2\epsilon})
        -\bigoplus_{j=1}^h 
        \left(T_{b_j}\Sigma_j\otimes T_{b_j} \Sigma_0\right)
        +\operatorname{Ind}_{S^1}^{Z_{\widetilde{G}}(\alpha)}
        \left(
            V_{-\eta}
        \right)
        -\chi^{\eta}-\chi^{-\eta}
        \]
        For $\epsilon=0,1$, $\chi^{-\epsilon}+\chi^{1-\epsilon}-\chi^{1-2\epsilon}=1$.
    \end{enumerate}
\end{proof}
\begin{remark}
    Consider the bundle
        $T_{b_1}\Sigma_1\otimes T_{b_1}\Sigma_0\in K([\Mb_{0,h}/S^1])$,
        then 
        \[
        \bigoplus_{j=1}^h 
        \left(T_{b_j}\Sigma_j\otimes T_{b_j} \Sigma_0\right)
        \simeq \left(
            [\Mb_{0,h}^{\im\alpha} /S^1]\to [\Mb_{0,h}^{\im\alpha}/Z_{\widetilde{G}}(\alpha)]
        \right)_{!}\left(T_{b_1}\Sigma_1\otimes T_{b_1}\Sigma_0\right)
        \]
        is a $K$-theoretic pushforward.
        This can be viewed as a ``twisted Ind''.
\end{remark}

\section{Computation}
Throughout this section, let $E$ be an HKR theory with
$p$-torsion-free coefficient ring $E^*$. We will compute the
Gromov--Witten invariants of $X_\loc$ with coefficients in $E$.
This covers the Morava $K$-theory
coefficient
because one can choose a suitable $E$ with
a morphism from $E\to K(n)$
so the result in $K(n)$ is just that in $E$
mod $p$. We assume $p\ge 3$ (\cref{rmk: homotopy commutative requirement}).

We denote the $S^1$-equivariant
Gromov--Witten invariants of $X_\loc$ by
$I_{d,\eq}\in E^*\llbracket u\rrbracket$, and
denote the nonequivariant one by $I_d\in E^0$, which
is also the constant term of $I_{d,\eq}$.
Notions like $x+_Ey$ and $[n]x$ are defined using 
the formal group law  of a fixed complex orientation
of $E$.

\subsection{Reducing from general \texorpdfstring{$d$}{d} to \texorpdfstring{$p$}{p}-powers}
With rational cohomology coefficients,
the Gromov--Witten invariant of $X_\loc$
is $1/d^3$. We show that such a factor
also appears with coefficients in $E$, but only
for the $p$-coprime part. Precisely,
we prove:
\begin{proposition}
\label{prop:separating p and p-coprime part nonequivariantly}
    Assume that $E^*$ has no $p$-torsion.
    Let $d=p^a m$ with $p\nmid m$. Then
    \[
        I_{p^a m}= \frac{1}{m^3} I_{p^a}.
    \]
    Note that $m$ is invertible in $E^0$, so $\frac{1}{m^3}$ is defined.
\end{proposition}
We actually prove the $S^1$-equivariant version of it. 
\begin{proposition}
    \label{prop:separating p and p-coprime part equivariantly}
    Assume that $E^*$ has no $p$-torsion.
    Let $d,m$ be positive integers and assume $p\nmid m$. Then
    the degree-$d$ and degree-$m\cdot d$ $S^1$-equivariant Gromov--Witten invariants of $X_\loc$ satisfy the following relation
    \[
        Q(u)\cdot I_{md,\eq}(u)=
        \frac{1}{m}\cdot Q\left(\left[\frac{1}{m}\right] u\right)\cdot 
        I_{d,\eq}\left(\left[\frac{1}{m}\right] u\right)
    \]
    where $Q(u):=u\cdot [-1]u$.
\end{proposition}
Taking the $u^2$ coefficients of both sides implies \cref{prop:separating p and p-coprime part nonequivariantly}.

    Consider the following on the target $\P^1$,
        $\varphi_m\colon \P^1\to \P^1,\quad w\mapsto w^m$
    which induces a morphism 
    \begin{equation}
        \label{eq:deg m covering map}
        \varphi_m\colon \Mb_{0,0}(\P^1,d)\to \Mb_{0,0}(\P^1,md)
    \end{equation}
    by postcomposition. $\varphi_m$ is equivariant with regard to
    the group homomorphism $S^1\xrightarrow{t\mapsto t^m}S^1$.
    $\varphi_m$ maps points in $S^1$-fixed components
    of degree $d$ to those of degree $md$. In fact,
    \begin{itemize}
        \item It sends the single edge component to the single edge component.
        \item For comb components whose number of marked points
         is a power of $p$ (\cref{Euler class vanishing}),
         since $d$ and $md$ are
    divisible by the same power of $p$, there are
    the same number of such comb components, with $\varphi_m$
    comparing them.
    \end{itemize}
    
    We examine $\varphi_m$ on comb components.
    Nonequivariantly, these components are $[M/G]$ and $[M'/G']$,
    as in \eqref{eq:model for comb components}, where
    \[
        M=M'=\Mb_{0,h},\qquad
        G=(C_l)^h\rtimes S_h,\qquad
        G'=(C_{lm})^h\rtimes S_h.
    \]
    Their $S^1$-quotient stacks are $[M/\widetilde{G}]$ and
    $[M'/\widetilde{G'}]$, where
    \[
        \widetilde{G}=A_l\rtimes S_h,\qquad
        \widetilde{G'}=A_{lm}\rtimes S_h,
    \]
    as in \eqref{eq:another way to define widetildeG}.
    The map $\varphi_m$ on these comb components is described as follows.
    \begin{itemize}
        \item The map $M\to M'$ is the identity on $\Mb_{0,h}$,
        since postcomposition does not change the contracted central curve.
        \item The homomorphism $G\to G'$ is induced by the inclusion $C_l\hookrightarrow C_{lm}$.
        \item The homomorphism $\widetilde{G}\to\widetilde{G'}$
        fits into the following commutative diagram:
\[\begin{tikzcd}
	1 & {A_l} & {\widetilde{G}} & {S_h} & 1 \\
	1 & {A_{lm}} & {\widetilde{G'}} & {S_h} & 1
	\arrow[from=1-1, to=1-2]
	\arrow[from=1-2, to=1-3]
	\arrow[from=1-2, to=2-2]
	\arrow[from=1-3, to=1-4]
	\arrow[from=1-3, to=2-3]
	\arrow[from=1-4, to=1-5]
	\arrow["\id", from=1-4, to=2-4]
	\arrow[from=2-1, to=2-2]
	\arrow[from=2-2, to=2-3]
	\arrow[from=2-3, to=2-4]
	\arrow[from=2-4, to=2-5]
\end{tikzcd}\]
The left vertical arrow is
\[
    (z_1,\ldots,z_h;t)\longmapsto(z_1,\ldots,z_h;t^m).
\]
See \cref{rmk:semi-direct product and Sn action}.
    \end{itemize}
    Thus $\varphi_m$ induces a homeomorphism between the coarse moduli
    spaces of these comb components. The following lemma gives a stronger comparison.
    \begin{lemma}
    \label{lem:comparison of m and md}
        Assume $p\nmid m$, then under \eqref{eq:deg m covering map}
        \begin{enumerate}[(1)]
            \item $\varphi_m$ induces a bijection on conjugacy classes
            of HKR sectors 
            \[
                \Hom(\Lambda_r,G)/G\to \Hom(\Lambda_r,G')/G'
            \]
            \item For an HKR sector $\alpha\colon \Lambda_r\to G$, we
            use $\alpha$ itself to denote the HKR sector $(G\to G')\circ \alpha$. \emph{If $\alpha$ acts transitively on $h$ marked points}, then the natural homomorphism
            \[
                Z_{\widetilde{G}}(\alpha)\to Z_{\widetilde{G'}}(\alpha)
            \]
            is an isomorphism.
        \end{enumerate}
    \end{lemma}
    \begin{proof}
        \begin{enumerate}[(1)]
            \item We need a fact from group theory, \cref{lem:a group theory lemma}, which
     is cited later.
 Write $l=p^b l_0$ with $p\nmid l_0$, and define
 $H=(C_{p^b})^h\rtimes S_h$. There are two exact sequences
 \[
    1\to (C_{l_0})^h\to G\to H\to 1
 \]
 \[
    1\to (C_{l_0 m})^h\to G'\to H\to 1.
 \]
 By \cref{lem:a group theory lemma}, the projection
 $G\to H$ and $G'\to H$ induce bijections of HKR sectors up
 to conjugacy. Thus $G\hookrightarrow G'$ also induces a bijection
 of HKR sectors up to conjugacy.
 \item By \cref{lem:identifying centralizer} we have
\[\begin{tikzcd}
	1 & {S^1} & {Z_{\widetilde{G}}(\alpha)} & P & 1 \\
	1 & {S^1} & {Z_{\widetilde{G'}}(\alpha')} & P & 1
	\arrow[from=1-1, to=1-2]
	\arrow[from=1-2, to=1-3]
	\arrow["\id", from=1-2, to=2-2]
	\arrow[from=1-3, to=1-4]
	\arrow[from=1-3, to=2-3]
	\arrow[from=1-4, to=1-5]
	\arrow["\id", from=1-4, to=2-4]
	\arrow[from=2-1, to=2-2]
	\arrow[from=2-2, to=2-3]
	\arrow[from=2-3, to=2-4]
	\arrow[from=2-4, to=2-5]
\end{tikzcd}\]
Notice that the left vertical arrow is the identity because in
\cref{lem:identifying centralizer} the inclusion $S^1\to A_l$ is
\[
    z\longmapsto(z,\ldots,z;z^l),
\]
and the map $A_l\to A_{lm}$ is
\[
    (z_1,\ldots,z_h;t)\longmapsto(z_1,\ldots,z_h;t^m).
\]
Therefore, the middle map is also an isomorphism.
\end{enumerate}
\end{proof}
\begin{lemma}
    \label{lem:a group theory lemma}
    Let $K\lhd G$ be a normal subgroup, and let $H=G/K$. If
    $p\nmid |K|$, then the projection induces a bijection of HKR sectors:
    \[
        \Hom(\Lambda_r,G)/G\xrightarrow{\simeq }\Hom(\Lambda_r,H)/H.
    \]
\end{lemma}

     This is a consequence of the Schur–Zassenhaus theorem (for example, see \cite[Theorems 7.41 and 7.42]{Rotman95}):
\begin{quote}
    A normal Hall(meaning $|K|$ and $|G/K|$ are coprime) subgroup $K$ of a finite group $G$ has a complement;
if either $K$ or $G/K$ is solvable, then any two complements
are conjugate.
\end{quote}
Complements up to conjugacy
are in bijection with the sections of projection
$G\twoheadrightarrow G/K$ up to conjugacy, so the theorem means that there's a unique section up to conjugacy.
\begin{proof}[Proof of \cref{lem:a group theory lemma}]
The lemma is implied by Schur–Zassenhaus theorem in the
following way. 
\begin{itemize}
    \item Surjectivity: given a homomorphism $\overline{\alpha}\colon \Lambda_r\to H$,
we pull the exact sequence $1\to K\to G\to H\to 1$ to
\begin{equation}
\label{eq:the pullback exact sequence}
    1\to K\to \overline{G}\to \Lambda_r\to 1.
\end{equation}
Then $|K|$ is coprime to $|\Lambda_r|$ because $|\Lambda_r|$
is a power of $p$, and $\Lambda_r$ is abelian, thus solvable.
Applying Schur–Zassenhaus theorem to \eqref{eq:the pullback exact sequence}, we get a section $\alpha\colon \Lambda_r
\to G$. This proves surjectivity of HKR sectors.
\item Injectivity: assume there are two HKR sectors $\alpha_1,\alpha_2
\colon \Lambda_r\to G$ such that their compositions
with the projection $G\to G/K$, denoted by $\overline{\alpha_1},
\overline{\alpha_2}\colon \Lambda_r\to H$, are conjugate. W.L.O.G., we
can assume that $\overline{\alpha_1}=\overline{\alpha_2}$.
Then they give two sections in \eqref{eq:the pullback exact sequence}, and are conjugate by the second sentence of 
Schur–Zassenhaus theorem.
\end{itemize}
\end{proof}
\begin{proof}[Proof of \cref{prop:separating p and p-coprime part equivariantly}]
    After equivariant localization and HKR character, we reduce
    the Gromov--Witten invariants to a sum of equivariant pushforward maps
    (\cref{prop:Atiyah-Bott equivariant localization formula,prop:Riemann Roch for equivariant HKR})
    \[
       I_{d,\eq,\alpha}^{\comb}:= \left(\left[(\Mb_{0,h})^{\im\alpha}/Z_{\widetilde{G}}(\alpha)\right]
        \to [\pt/S^1]\right)_{!,\eq}^{L(E)}
        \left(
            \chi_{n,p;\eq}\left(\frac{e_\eq(\Ob_d)}{e_\eq(\nu_F)}\right)
            \frac{1}{e_{\operatorname{HKR},\eq}^L(\nu_{\Mb_{0,h}^{\im\alpha}/\Mb_{0,h}})}
        \right)
    \]
    Here the character value is restricted to the $\alpha$-sector before
    applying the pushforward.
    By \cref{lem:comparison of m and md}, $\varphi_m$ induces
    an isomorphism of $\left[(\Mb_{0,h})^{\im\alpha}/Z_{\widetilde{G}}(\alpha)\right]$ for degree $d$ and $md$. 
    
    We compare
    the two virtual bundles $\Ob_d-\nu_F-\nu_{\Mb_{0,h}^{\im\alpha}/\Mb_{0,h}}$
    for degree $d$ and $md$ over the same $\left[\Mb_{0,h}^{\im\alpha}/Z_{\widetilde{G}}(\alpha)\right]$.
    The last terms $\nu_{\Mb_{0,h}^{\im\alpha}/\Mb_{0,h}}$ are the same.
    By \cref{lem:case2 virtual bundle reduction after HKR} we have
    \[
             (\Ob_d-\nu_F)|_{[\Mb_{0,h}^{\im\alpha}/Z_{\widetilde{G}}(\alpha)]}
             =V_\std+
             \operatorname{Ind}_{S^1}^{Z_{\widetilde{G}}(\alpha)}
        \left(
            V_{-\eta}
        \right)
        -\bigoplus_{j=1}^h 
        \left(T_{b_j}\Sigma_j\otimes T_{b_j} \Sigma_0\right)
        -\chi-\chi^{-1}.
         \]
     \[
             (\Ob_{md}-\nu_F)|_{[\Mb_{0,h}^{\im\alpha}/Z_{\widetilde{G'}}(\alpha)]}
             =V_\std+
             \operatorname{Ind}_{S^1}^{Z_{\widetilde{G'}}(\alpha)}
        \left(
            V_{-\eta}
        \right)
        -\bigoplus_{j=1}^h 
        \left(T_{b_j}\Sigma_j\otimes T_{b_j} \Sigma_0\right)
        -\chi'-(\chi')^{-1}.
         \]
    where $\chi,\chi'$ are the $S^1$-characters satisfying $\varphi^*_m(\chi')=\chi^{\otimes m}$. 
    We discover that
    \[
    \begin{aligned}
        &\chi+\chi^{-1}+(\Ob_d-\nu_F-\nu_{\Mb_{0,h}^{\im\alpha}/\Mb_{0,h}})|_{[\Mb_{0,h}^{\im\alpha}/Z_{\widetilde{G}}(\alpha)]}\\
        &\quad=\chi'+(\chi')^{-1}+(\Ob_{md}-\nu_F-\nu_{\Mb_{0,h}^{\im\alpha}/\Mb_{0,h}})|_{[\Mb_{0,h}^{\im\alpha}/Z_{\widetilde{G'}}(\alpha)]}
    \end{aligned}
    \]
    The comparison diagram for equivariant pushforward is
\[\begin{tikzcd}
	{[\Mb_{0,h}^{\im\alpha}/Z_{\widetilde{G}}(\alpha)]} & {[\Mb_{0,h}^{\im\alpha}/Z_{\widetilde{G'}}(\alpha)]} \\
	{[\pt/S^1]} & {[\pt/S^1]}
	\arrow["\id"', from=1-1, to=1-2]
	\arrow["{\pi_{d}}", from=1-1, to=2-1]
	\arrow["{\pi_{md}}"', from=1-2, to=2-2]
	\arrow["{\varphi_m\colon t\mapsto t^m}", from=2-1, to=2-2]
\end{tikzcd}\]
Let $u=e(\chi),u'=e(\chi')$. Then
\begin{align*}
    Q(u')\cdot I_{md,\eq,\alpha}^{\comb}(u')
    &=
    (\pi_{md})_{!,\eq}\left(
    e_{\operatorname{HKR},\eq}^L\left(
    \chi'+(\chi')^{-1}+(\Ob_{md}-\nu_F-\nu_{\Mb_{0,h}^{\im\alpha}/\Mb_{0,h}})
    \right)
    \right)
    \\
    &=(\varphi_m)_!\circ (\pi_d)_{!,\eq}
    \left(e_{\operatorname{HKR},\eq}^L
    \left(
        \chi+\chi^{-1}+(\Ob_d-\nu_F-\nu_{\Mb_{0,h}^{\im\alpha}/\Mb_{0,h}})
    \right)
    \right)
    \\
    &=(\varphi_m)_!\left(
    Q(u)\cdot I_{d,\eq,\alpha}^{\comb}(u)
    \right)
\end{align*}
For $\varphi_m\colon BS^1\to BS^1$, since $L(E)$ is rational, we have
\[
    (\varphi_m)_!\circ \varphi_m^*(f(u'))=\frac{1}{m}\cdot f(u')
\]
But $\varphi_m^*(f(u'))=f([m]_E u)$ with $p\nmid m$ (so we can define$\left[\frac{1}{m}\right]u\in E^*\llbracket u\rrbracket$), we get the formula
\[
    (\varphi_m)_!f(u)=\frac{1}{m}\cdot f\left(\left[\frac{1}{m}\right] u'\right)
\]
Therefore
\[
    Q(u')\cdot I_{md,\eq,\alpha}^{\comb}(u')
    =\frac{1}{m}\cdot Q\left(\left[\frac{1}{m}\right] u'\right)
    I_{d,\eq,\alpha}^{\comb}\left(
    \left[\frac{1}{m}\right] u'\right)
\]
Summing over all $\alpha$ and all components gives the stated identity
in $L(E)^*((u))$. Since $E^*$ is $p$-torsion-free and $L(E)^*$ is
faithfully flat over $E^*[1/p]$
\cite[Corollary~6.8(ii)]{Generalized-Group-Characters-and-Complex-Oriented-Cohomology-Theories},
the map $E^*\llbracket u\rrbracket\to L(E)^*((u))$ is injective.
Thus the identity holds in $E^*\llbracket u\rrbracket$.
\end{proof}
\subsection{\texorpdfstring{$p$}{p}-primary computation assuming \texorpdfstring{$[-1]x=-x$}{[-1]x=-x}}
After \cref{prop:separating p and p-coprime part nonequivariantly},
we assume in this subsection that $d=p^a$ for $a\in \N$. $E$
is an HKR theory of height $n$.
\begin{definition}
    For $d=p^a$
    \begin{enumerate}[(1)]
        \item Define \[
    A_{d,u}:=E^*\llbracket x, u\rrbracket /
    ([d]_Ex-_Eu)
    =E^*\llbracket x, u\rrbracket /
    ([d]_Ex-u)
\]
$A_{d,u}$ is an $E^*\llbracket u\rrbracket$-algebra and is free
of rank $d^n=p^{an}$ as an $E^*\llbracket u\rrbracket$-module.
\item For a power series $f(x)\in E^*\llbracket x\rrbracket$
which represents the element $[f]\in A_{d,u}$, we define
\[
    \tau_{d,u}([f]):= \res_x\left(
    \frac{f(x)\omega(x)}{[d]_Ex-_E u}
    \right)
    \in E^*\llbracket u\rrbracket.
\]
Here $\omega\in E^*\llbracket x\rrbracket$ is the $1$-form characterized by
\cref{lem:invariant form}.
    \end{enumerate}
\end{definition}
Then our main result is
\begin{theorem}
    \label{thm:GW computation if [-1]x=-x}
    Assume that $E^*$ has no $p$-torsion and that the complex orientation on $E$
    satisfies $[-1]_Ex=-x$. Then the equivariant
    Gromov--Witten invariants of $X_\loc$
    in degree $p^a$ are
    \begin{equation}\label{eq:equivariant computation in total}
        I_{p^a,\eq}=\frac{
        \tau_{p^a,u}(x^2)-\tau_{p^a,0}(x^2)
        }{u^2}
    \end{equation}
    In particular,
    \begin{equation}\label{eq:nonequivariant computation in total}
        I_{p^a}=[u^2]\left(\tau_{p^a,u}(x^2)\right).
    \end{equation}
\end{theorem}

The underlying $S^1$-fixed component contribution
is as follows.
\begin{proposition}
\label{prop:contribution from components}
Under the hypotheses of \cref{thm:GW computation if [-1]x=-x}, the following formulas hold.
    \begin{enumerate}[(1)]
        \item The single-edge component
        contribution is
        \[
            I_{p^a,\eq}^{\operatorname{single}}=\frac{\tau_{p^a,u}(x^2)}{u^2}
            \in E^*\llparenthesis u
            \rrparenthesis.
        \]
        \item For $1\le b\le a$, let $I_{p^a,\eq}^{\mathrm{comb},b}$ be the corresponding
        contribution from both\footnote{
        For one comb graph, there are two components of this shape, one whose central curve lies over $p_0$
        and another whose central curve lies over $p_1$.
        } comb components with $h=p^b$ tails
        of degree $l=p^{a-b}$. Put
        \[
            \Theta_b=\tau_{p^b,0}(x^2),\quad
            \Theta_0=0.
        \]
        Then
        \[
            I_{p^a,\eq}^{\mathrm{comb},b}
            =\frac{p^{(n-1)(a-b)}}{u^2}\left(
            -\Theta_b+p^{n-1}\Theta_{b-1}\right).
        \]
    \end{enumerate}
\end{proposition}
\begin{proof}[Proving \cref{thm:GW computation if [-1]x=-x} using \cref{prop:contribution from components}]
    The sum of comb contributions is a telescope
\[
    \sum_{b=1}^a\frac{p^{(n-1)(a-b)}}{u^2}\left(
            -\Theta_b+p^{n-1}\Theta_{b-1}\right)
            =
            -\frac{\Theta_a}{u^2}+\frac{p^{a(n-1)}}{u^2}\Theta_0
            =-\frac{\tau_{p^a,0}(x^2)}{u^2}.
\]
Combining this with the single-edge contribution gives the stated formula
after Euler localization. Since $I_{p^a,\eq}\in E^*\llbracket u\rrbracket$,
multiplication by $u^2$ and injectivity of Euler localization give
\[
    u^2I_{p^a,\eq}
    =\tau_{p^a,u}(x^2)-\tau_{p^a,0}(x^2)
    \quad\text{in }E^*\llbracket u\rrbracket.
\]
Thus the numerator is divisible by $u^2$, and comparison of the coefficients
of $u^2$ gives $I_{p^a}=[u^2]\tau_{p^a,u}(x^2)$.
\end{proof}
We will devote the rest of this subsection to
proving \cref{prop:contribution from components}.
\subsubsection{Single-edge component computation}
We have the following
knowledge of this component.
\begin{itemize}
    \item The orbifold structure is $[\pt/C_{p^a}]$ with
    a $S^1$-action described (in the sense of \cref{def:working definition of group action on orbifold}) by 
    \[
    1\to C_{p^a}\to S^1\xrightarrow[t\mapsto t^{p^a}]{\pi} S^1\to 1.
    \]
    which is the lower row of \eqref{eq:single-edge component working model}. We use $S^1_\source$ and $S^1_\target$ to represent
    the middle and right $S^1$.
    \item The equivariant pushforward $\pi_{!,\eq}\colon H^*(BS^1_\source,E)\to H^*(BS^1_\target,E)$ is computed in
    \cref{prop:An equivariant pushforward formula computed by residue}. We postpone
    the proof to later.
    \item The equivariant Euler class of $\Ob_d-\nu_F$ is given
    in \cref{lem:virtual-euler-class-case1}
    \[
        e_{\eq}(\Ob_d-\nu_F)=\frac{x\cdot [-1]x}{[d]x\cdot [-d]x}
    \]
    By assumption $[-1]x=-x$, this is also $\frac{x^2}{\left([d]x\right)^2}$.
\end{itemize}
\begin{proof}[Proof of \cref{prop:contribution from components}(1)]
    Since $\pi^*(u)=[d]_E x$, the projection formula gives
\begin{align*}
    I_{d,\eq}^{\operatorname{single}}
    &=\pi_{!,\eq}\left(\frac{x^2}{(\pi^*u)^2}\right)\\
    &=\frac{1}{u^2}\pi_{!,\eq}(x^2)\\
    &=\frac{1}{u^2}\res_x\left(
    \frac{x^2\omega(x)}{[d]_E x-_E u}\right)\\
    &=\frac{\tau_{d,u}(x^2)}{u^2},
\end{align*}
where the third equality applies
\cref{prop:An equivariant pushforward formula computed by residue} to $x^2$.
\end{proof}
Now we supply the technical computations. The following formula
is an equivariant generalization of \cite[Proposition 9.2]{strickland2000knlocaldualityfinitegroups}.
\begin{proposition}
\label{prop:An equivariant pushforward formula computed by residue}
    Assume that $E^*$ has no $p$-torsion.
    Let $f(x)\in E^*\llbracket x\rrbracket$ represent
    an element of $H^*\left(M_{h\widetilde{G}},E\right)$.
    Then the equivariant pushforward map is computed by
    \[
        \pi_{!,\eq}\left(f(x)\right)
        =\res_x\left(
        \frac{f(x)\omega(x)}{[p^a]_Ex-_E u}
        \right)\in E^*\llbracket u\rrbracket \simeq H^*(BS^1_\target,E).
    \]
    (Recall $\res$ in \cref{def:residue}.)
\end{proposition}

\begin{proof}
    We can verify the formula after applying the equivariant HKR character map.
    Choose $r$ sufficiently large. The classes in HKR
    theory are
    $
        \alpha\colon \Lambda_r\to C_{p^a}.
    $
    Since $C_{p^a}$ is abelian, conjugation
    fixes every $\alpha$, so there are
    $p^{an}$ classes in total.
    \[
        \Fix_{n,p}M=\bigsqcup_{\alpha
        \in \Hom(\Lambda_r,C_{p^a})}
        \pt
    \]
    with the trivial $\widetilde{G}=S^1$-action. 
    The equivariant $n$-fold
    $p$-inertia stack is 
    \[
        \bigsqcup_{\alpha
        \in \Hom(\Lambda_r,C_{p^a})} [\pt/S^1].
    \]

    The universal map (\cref{def:equivariant universal map and equivariant normal bundle}) on each component $\alpha$ is
    $BS^1\times B\Lambda_r\to BS^1$, induced by the group homomorphism
    $S^1\times \Lambda_r\xrightarrow{\id_{S^1}\times \iota}
    S^1\times S^1\xrightarrow{\text{multiplication}} S^1$.
    Write $H^*(BS^1,E)=E^*\llbracket x\rrbracket$ and let 
    $a_\alpha\in H^*(B\Lambda_r,E)$ be the pullback of $x$
    under $ B\Lambda_r\xrightarrow{\alpha} BC_{p^a}\to BS^1$.
    We can also view $a_\alpha$ as an element of $L(E)^*$ by
    the map $H^*(B\Lambda_r,E)\to L(E)^*$, and still denote
    this image by $a_\alpha$. They are exactly the roots of
    $[p^a]x=0$ in $L(E)^*$ by \cref{lem:e(chi) becomes exactly the roots of [p^a]x=0}.
    
    We trace the diagram \eqref{eq:equivariant HKR Rieman Roch}
    to compare the images of the two sides:
    \begin{enumerate}[(1)]
        \item First we trace along this direction:
        \[\begin{tikzcd}
    {H^*(M_{h\widetilde{G}},E)} && {H^*(\Fix_{n,p}M_{h\widetilde{G}},L(E))} \\
	&& {H^*(B\TT,L(E))}
    \arrow["{\chi_{n,p;\eq}}", from=1-1, to=1-3]
	\arrow["{\pi_{!,\eq}^{L(E)}}", from=1-3, to=2-3]
\end{tikzcd}\]
            The HKR
    character is
    \[
        f(x)\in E^*\llbracket x\rrbracket=H^*(M_{h\widetilde{G}},E)
        \mapsto \oplus_{\alpha\in \Hom(\Lambda_r,C_{p^a})}
        f(x+_Ea_\alpha)\in H^*(\Fix_{n,p}M_{h\widetilde{G}},L(E)).
    \]
    Now we compute the $L(E)$-pushforward 
    \[
        H^*\left(\Fix_{n,p}M_{h\widetilde{G}},L(E)\right)
        \to H^*(B\TT,L(E))
    \]
    On each component $\alpha$, it is the pushforward along 
    $\pi\colon B\widetilde{G}\to B\TT, z\mapsto z^{p^a}$.
    Since $L(E)$ is of rational coefficient, there is an
    $y\in H^*(B\TT,L(E))$ such that $[p^a]y=u$. So
    $f(x+_Ea_\alpha)=\pi^* f(y+_E a_\alpha)$.
    By \cref{cor:projection formula}, we have
    \[
        \pi_{!,\eq}^{L(E)}(f(x+_Ea_\alpha))
        =f(y+_E a_\alpha) \pi_{!,\eq}^{L(E)}(1)
        =\frac{1}{p^a}f(y+_E a_\alpha).
    \]
    Here $\pi_{!,\eq}^{L(E)}(1)=\frac{1}{p^a}$ is because $L(E)$
    is rational. To conclude, we get
    \begin{equation}
        \label{eq:first diagram tracing}
        \sum_{\alpha\in \Hom(\Lambda_r,C_{p^a})}
        \frac{1}{p^a}f(y+_Ea_\alpha)\in H^*(B\TT,L(E)).
    \end{equation}
        \item On the other hand, we consider the image of
        \[
            \res_x\left(
        \frac{f(x)\omega(x)}{[p^a]_Ex-_Eu}
        \right)\in E^*\llbracket u\rrbracket
        \]
        in $H^*(B\TT,L(E))=L(E)^*\llbracket u\rrbracket$.
        By $u=[p^a]y$, we have
        \begin{align*}
            \res_x\left(
        \frac{f(x)\omega(x)}{[p^a]x-_Eu}
        \right)&=\res_x\left(
        \frac{f(x)\omega(x)}{[p^a]x-_E[p^a]y}
        \right)\\
        &=_{z=x-_Ey}\res_z\left(\frac{f(z+_E y)\omega(z)}{
        [p^a]z
        }\right)
        \end{align*}
        The last equality holds because $\omega$ is invariant under a change of variable
        $z=x-_E y$.
        By construction of $L(E)$, $a_\alpha$ are roots of $[p^a]z=0$
        in $L(E)$,
        and they are all roots. Also $a_\alpha-a_\beta$ are invertible,
        by construction. 
        We have
        \begin{equation}
            \label{eq:second diagram tracing}
            \begin{aligned}
            \res_z\left(\frac{f(z+_E y)\omega(z)}{
        [p^a]z
        }\right)
        &=\operatorname{Tr}_{L(E)^*\otimes_{E^*} \left( E^*\llbracket y,z\rrbracket /([p^a]z) \right)}
        \left(
        \frac{f(z+_E y)w(z)}{([p^a]z)'}
        \right)\\
        &=\sum_{\alpha}\frac{f(a_\alpha+_E y)w(a_\alpha)}{([p^a]z)'|_{z=a_\alpha}}.
            \end{aligned}
        \end{equation}
        The first equality is by \cref{lem:residue=trace}
        and the second holds because $[p^a]z=0$ has distinct roots.
    \end{enumerate}
    Comparing \eqref{eq:first diagram tracing} and \eqref{eq:second diagram tracing}, we are left to prove that
    \[
        \frac{w(a_\alpha)}{([p^a]x)'|_{x=a_\alpha}}=\frac{1}{p^a}
    \]
    By \cref{lem:invariant form}
    \[
        w([p^a]x)([p^a]x)'dx=w([p^a]x)d([p^a]x)=p^a w(x)dx
    \]
    Let $x=a_\alpha$ we get
    $p^a w(a_\alpha)=([p^a]x)'|_{x=a_\alpha}$ because
    $[p^a]a_\alpha=0$ and $w(0)=1$.
    Thus the two sides have the same image in $H^*(BS^1_\target,L(E))$.
    Since $E^*$ has no $p$-torsion and $L(E)^*$ is faithfully flat over
    $E^*[1/p]$
    \cite[Corollary~6.8(ii)]{Generalized-Group-Characters-and-Complex-Oriented-Cohomology-Theories},
    the map
    \[
        E^*\llbracket u\rrbracket\longrightarrow L(E)^*\llbracket u\rrbracket
    \]
    is injective. Hence the two sides are equal in $E^*\llbracket u\rrbracket$.
\end{proof}

\begin{remark}
Put $d=p^a$. The argument
of \cite[Proposition~9.2]{strickland2000knlocaldualityfinitegroups} gives
the nonequivariant pushforward $\epsilon$ of $[\pt/C_d]\to\pt$ as
    \begin{equation}
        \label{eq:nonequivariant pushforward formula by residue}
        f(x)\mapsto
        \res_x\left(
        \frac{f(x)\omega(x)}{[p^a]_Ex}
        \right)\in E^*.
    \end{equation}

Let
\[
\pi\colon[\pt/S^1_\source]\longrightarrow[\pt/S^1_\target]
\]
be induced by $t\mapsto t^d$. This is the equivariant projection
associated to the $S^1_\target$-action on $[\pt/C_d]$. Pulling back this
action along $S^1_\source\to S^1_\target$ splits the defining extension,
giving the $2$-cartesian square
\[
\begin{tikzcd}
    {[\pt/(S^1_\source\times C_d)]}
    & {[\pt/S^1_\source]} \\
    {[\pt/S^1_\source]}
    & {[\pt/S^1_\target]}
    \arrow["m", from=1-1, to=1-2]
    \arrow["{\pi_1}"', from=1-1, to=2-1]
    \arrow["\pi", from=1-2, to=2-2]
    \arrow["\pi"', from=2-1, to=2-2]
\end{tikzcd}
\]
where $m$ is induced by $(t,\zeta)\mapsto t\zeta$, so that
$m^*(x)=x+_Ey$.

By change-of-group naturality of the equivariant pushforward in
\cref{def:equivariant pushforward}, inherited from transverse base change
(\cref{prop:properties of f_!}(2)), we have
\[
\pi^*\circ\pi_{!,\eq}=(\pi_1)_{!,\eq}\circ m^*.
\]
The left vertical map corresponds to the trivial $S^1_\source$-action on
$[\pt/C_d]$, so its pushforward applies $\epsilon$ in the $C_d$-factor.
Consequently,
\[
\begin{aligned}
\pi^*\!\left(\pi_{!,\eq}(f)\right)
&=\res_y\frac{f(x+_Ey)\omega(y)}{[d]_Ey}\\
&=\res_z\frac{f(z)\omega(z)}{[d]_Ez-_E[d]_Ex},
\end{aligned}
\]
using $z=x+_Ey$ and translation invariance of $\omega$.

Since $\pi^*(u)=[d]_Ex$, this is the pullback of the asserted formula.
The previously established finite freeness of
$A_{d,u}\cong E^*\llbracket x\rrbracket$ over $E^*\llbracket u\rrbracket$
makes $\pi^*$ injective, completing the proof.
\end{remark}
\subsubsection{Comb components computation}
Recall $\Lambda_r=\left(\Z/p^r\right)^n$ for a sufficiently
large $r$ (which only depends on $a$). Write $h=p^b$, $l=p^{a-b}$. $M=\Mb_{0,h}$,
$G=(C_l)^h\rtimes S_h$ and $\widetilde{G}=A_l\rtimes S_h$ as
in \eqref{eq:model for comb components} and \eqref{eq:define widetilde G}.

By \cref{lem:Euler class vanishing after HKR}, only transitive
HKR sectors $\alpha\colon \Lambda_r\to G$ could possibly 
have a nonzero contribution. Such sectors are classified as follows.
\begin{lemma}
    There is a bijective correspondence between
    \begin{itemize}
        \item Conjugacy classes of transitive sectors $\alpha\colon\Lambda_r\to (C_l)^h\rtimes S_h$.
        \item A subgroup $H\subset \Lambda_r$ with $[\Lambda_r\colon H]
        =h$
        and a group homomorphism $\phi\colon H\to C_l$.
    \end{itemize}
    Also, for sufficiently large $r$ (depending only on $h,l$), any fixed such $H$ has exactly $l^n$ homomorphisms $\phi$.
\end{lemma}
 We use the following fact. Let $E$ be the union of $h$ free right
    $C_l$-sets. Then there is a (non-canonical) isomorphism
    $\Aut_{C_l}(E)\simeq (C_l)^h\rtimes S_h$. A homomorphism
    $\alpha\colon \Lambda_r\to (C_l)^h\rtimes S_h$ is a left $\Lambda_r$-action on $E$ that commutes with the $C_l$-action. Correspondingly,
    for two $\alpha_1,\alpha_2$, they are conjugate iff $E_1,E_2$
    as $\Lambda_r\times C_l$-sets are isomorphic. Thus we can replace
    conjugacy classes of sectors by $\Lambda_r$ actions on $E$ up
    to isomorphism.
\begin{proof}

    Choose one free $C_l$-orbit $O_1$ from the $h$ orbits, and define $H\subset \Lambda_r$
    to be the subgroup that maps $O_1$ to itself; $H$ acts on $O_1$
    and commutes with the $C_l$-action on $O_1$.
    For $h_1\in H$, the action
    on $O_1$ is multiplication by $\phi(h_1)\in C_l$. This
    gives the homomorphism $\phi$.
    By an argument similar to \cref{lem:identifying centralizer}, $H$ actually
   preserves every orbit, and $\alpha|_H$
    factors through $(C_l)^h$; in fact it factors through
    the diagonal $C_l\to (C_l)^h$. This means for $h_1\in H$ and any
    orbit $O_i$, the rotations $\phi_i(h_1)\in C_l$ are the same.
    So $(H,\phi)$ does not depend on the choice of $O_1$.

    Conversely, for any $H$ and $\phi$ satisfying
    conditions, there is a transitive $\alpha\colon \Lambda_r
    \to (C_l)^h\rtimes S_h$. We define
    \[
        E_{H,\phi}=(\Lambda_r\times C_l)/\sim.
    \]
    where 
    \[
        (\lambda+h_1, c)\sim (\lambda, \phi(h_1)c),\quad h_1\in H.
    \]
    $E_{H,\phi}$ admits commuting left $\Lambda_r$-actions and
    right $C_l$-actions, and the $C_l$-action is free. $E_{H,\phi}$
    is the disjoint union of $h$ free $C_l$-orbits because $E_{H,\phi}/C_l\simeq \Lambda_r/H$. 

    Now we prove the constructions above are inverses to each other
    \begin{itemize}
        \item  $E\to (H,\phi)\to E_{H,\phi}$. In the first construction
        we have $E\simeq \operatorname{Ind}_{H}^{\Lambda_r}(C_l)$\footnote{This map is given by the adjoint inclusion
        $C_l=O_1\hookrightarrow \operatorname{Res}_{H}^{\Lambda_r} E$
        as $H\times C_l$ sets.}
        as $\Lambda_r\times C_l$ sets. While in the second construction
        we have $E_{H,\phi}=\operatorname{Ind}_{H}^{\Lambda_r}(C_l)$
        by definition. So $E\simeq E_{H,\phi}$.
        \item $(H,\phi)\to E_{H,\phi}\to (H',\phi')$. In the construction $(H,\phi)\to E_{H,\phi}$, since
        $E_{H,\phi}/C_l\simeq \Lambda_r/H$, we have
        $\Lambda_r$
        acts on the $C_l$-orbits like left action on $\Lambda_r/H$.
        If we take $O_1$ to be the orbit corresponding to $\id[H]$, i.e.
        \[
            \id[H]=\{(0,c)\mid c\in C_l\}\subset E_{H,\phi}.
        \]
        $H'$ that preserves $O_1$ is exactly $H$. The action of $H$
        on $\id[H]$ is
        \[
        h_1\cdot(0,c)=(h_1,c)\sim (0,\phi(h_1)c)
        \]
        which is exactly $\phi$.
    \end{itemize}
    Finally, if $r$ is sufficiently large, any $H$ with $[\Lambda_r:H]=h$
    is a direct sum of $n$ cyclic $p$-groups whose order is divisible by $l$ (because $h,l$ are all $p$-powers).
    So there are $l^n$ homomorphisms $\phi$.
\end{proof}
Notice that $H$ is the kernel of $\Lambda_r\to S_h$. As before, we
denote $P$ to be the image of $\Lambda_r\to S_h$. 
\begin{notation}
We call the following data a \emph{decorated flag}:
    a chain of normal subgroups
        \[
            1=P_0\lhd P_1\lhd \cdots \lhd P_s =P.
        \]
        such that $P_j/P_{j-1}$ are cyclic $p$-groups, along
        with faithful characters:
        \[
            \theta_j\colon P_j/P_{j-1}\to \C^\times.
        \]
\end{notation}
For transitive $\alpha$, we can relate points in $\Mb_{0,p^b}^{\im\alpha}$
with equivalence classes of decorated flags.
\begin{lemma}
\label{lem:computation of imalpha fixed points and normal bundle}
We assume $p$ is odd.
    \begin{enumerate}[(1)]
        \item The connected components of $\Mb_{0,h}^{\im\alpha}$
        correspond bijectively to
        \[
            \{\text{decorated flags}\}/(\theta_s\sim \theta_{s}^{-1})
        \]
        The equivalence relation is to replace $\theta_s\colon P_s/P_{s-1}\to \C^\times$
        by its inverse and all others are the same.
        \item At each point $x\in \Mb_{0,p^b}^{\im\alpha}$ corresponding to $1=P_0\lhd\cdots \lhd P_s=P$ and characters $\theta_j$, we have
        \[
            [T_x\Mb_{0,p^b}]=[V_\std]-[\theta_s]-[\theta_{s}^{-1}]
            +\sum_{j=1}^{s-1}\left([\operatorname{Ind}_{P_j}^P\theta_j]
            -[\operatorname{Ind}_{P_j}^P\theta_j^{-1}]\right),
        \]
        as a $P$-representation; its $Z_{\widetilde{G}}(\alpha)$-representation structure\footnote{
        $Z_{\widetilde{G}}(\alpha)$ acts trivially on $\Mb_{0,p^b}^{\im\alpha}$: it acts via $Z_{\widetilde{G}}(\alpha)
        \to P$ (see the diagram in \cref{lem:identifying centralizer}(1)) but $P$ is exactly the image of $\alpha$.
        } comes from the pullback $Z_{\widetilde{G}}(\alpha)\to P$.
    \end{enumerate}
\end{lemma}
\begin{proof}
\begin{enumerate}[(1)]
    \item Let $[\Sigma_0]\in \Mb_{0,p^b}^{\im\alpha}$ be a stable rational nodal
    curve with $p^b$ marked points in the $\im\alpha$-fixed point;
    the domain of corresponding stable map is $\Sigma=\Sigma_0\vee\Sigma_1\vee\cdots\vee \Sigma_{p^b}$, i.e.
    attaching a sphere $\Sigma_i$ at each marked point.
    Then $P$ acts
    on $\Sigma$.
    We construct the ``dual tree''
    $T$ of $\Sigma$ whose vertices correspond to spheres and
    edges correspond to nodes. The leaves of $T$ correspond
    to the non-contracted components of $\Sigma$. So $P$
    acts on $T$, and the action is transitive and free on leaves.

    We claim:
    \begin{quote}
        There is a unique fixed vertex of $T$ under the $P$-action.
    \end{quote}
    View $T$ as a $P$-space. Since it is contractible, by Smith theory,
    the fixed point $T^P$ is nonempty. Choose a point $x\in T^P$,
    then $x$ is either a vertex or the midpoint of an edge. But
    the latter is impossible when $p$ is odd: otherwise this
    defines a homomorphism $P\to C_2$ corresponding to whether
    or not flipping the edge. By $|P|=p^b$ is odd, this homomorphism
    is trivial, meaning that $P$ never exchanges the two parts separated by
    this edge, contradicting the assumption that $P$ acts on
    leaves transitively. Thus $x$ is a vertex. If there are
    two different vertices $v_0,v_1\in T^P$, then the path connecting
    $v_0,v_1$ is fixed by $T^P$, which implies again that an edge
    is fixed and thus a contradiction. Therefore, there is
    a unique vertex fixed by $P$, denoted by $v_0$.

    We choose a path $v_0v_1\ldots v_s$ from $v_0$ to a leaf $v_s$.
    By transitivity of $P$ on leaves, for any vertex $v$ there
    are some $i\in \{1,2,\ldots, s\}$ and $\sigma\in P$ such that
    $\sigma(v_i)=v$; $i$ is equal to the distance from $v_0$ to $v$.
    This means $P$ acts transitively on depth-$i$ points of $T$.
    
    Define 
    \[
        P_{s-j}:=\{\sigma\in P\mid \sigma(v_j)=v_j\}.
    \]
    Then $P_{s-j}$ acts on the sphere corresponding to $v_j$;
    we call it the ``$v_j$-sphere''. We distinguish two
    cases:
    \begin{itemize}
        \item If $0<j<s$. $P_{s-j}$ fixes the node corresponding to the edge $v_{j}v_{j-1}$,
    and $P_{s-j-1}$ is the subgroup of maps that fix two
    nodes corresponding to edge $v_jv_{j-1}$ and $v_jv_{j+1}$. By transitivity of $P_{s-j}$ on the depth-$(j+1)$ vertices adjacent to $v_j$ and commutativity of $P$, we know that $P_{s-j-1}$ fixes the entire $v_j$-sphere\footnote{
    Since there are at least $3$ nodes, it suffices to show $P_{s-j-1}$ fixes all nodes. It already fixes the unique node connecting to $v_{j-1}$; what is new is that it fixes the nodes connecting to depth $j+1$-spheres}. Thus $P_{s-j}/P_{s-j-1}$ acts faithfully on a $\P^1$ fixing two points, which gives the character $\theta_{s-j}$.
    \item If $j=0$, $v_0$ is already the root, so there are only
    $|P_s/P_{s-1}|$ nodes that connect to the lower spheres.
    The $P_s/P_{s-1}$ action on $\P^1$ is a rotation, and the nodes
    lie on the equator with identical distance.
    So the
    character is up to an $\theta_s\sim \theta_s^{-1}$ ambiguity
    because we can apply $z\mapsto \frac{1}{z}$ on the $v_0$-sphere
    that fixes all the nodes.\footnote{
    For $0<j<s$, the node joining the $v_j$-sphere
    to the $v_{j-1}$-sphere is distinguished.
    We define $\theta_{s-j}$ to be the character
    of $P_{s-j}/P_{s-j-1}$ on the tangent line
    to the $v_j$-sphere at this node.
    This distinguishes $\theta_{s-j}$ from
    its inverse. No uniqueness of the coordinate
    is required.
    }
    \end{itemize}
    \item Our convention is that on the $v_j$-sphere, choose a local
    coordinate near the $v_jv_{j-1}$-node, then the $P_{s-j}$-action
    is rotation by $\theta_{s-j}\colon P_{s-j}\to \C^\times$.
    We write $C=\Sigma_0$, and let $\widetilde{C}$ be its normalization.
    Each non-leaf vertex $v$ of $T$ corresponds to an irreducible component
    $C_v$ of $C$; we denote $S_v\subset C_v$ the union of marked
    points and nodes. For each node $q$, let $q_+,q_-$ be the preimages
    in $\widetilde{C}$.

    By algebraic geometry, $T_x\Mb_{0,p^b}$ fits into
    an exact sequence
    \[
        0\to \bigoplus_{v}H^1(C_v, T_{C_v}(-S_v))
        \to T_x\Mb_{0,p^b}
        \to \bigoplus_{\text{ Nodes}}
        T_{q_+}\widetilde{C}\otimes T_{q_-}\widetilde{C}\to 0
    \]
    We keep in mind that they are $P$-representations and the following
    computation is in the representation ring of $P$.
    \begin{enumerate}[(a)]
        \item $\bigoplus_{\text{ Nodes}}
        T_{q_+}\widetilde{C}\otimes T_{q_-}\widetilde{C}= \sum_{j=1}^{s-1}\operatorname{Ind}_{P_j}^P\theta_j.$.
        
        We group the summation by $P$-orbits of nodes, and a
        set of representatives consists of the nodes $q_j$ between the spheres
        corresponding to $v_{j-1}$ and $v_{j}$ for $1\le j\le s-1$. The stabilizer is
        $\Stab_P(q_j)=P_{s-j}\subset P$.
        $P_{s-j}$ acts on $T_{q_j}C_{v_j}$ by $\theta_{s-j}$ by convention
        and acts trivially on $T_{q_j}C_{v_{j-1}}$. Since $P/P_{s-j}$
        acts transitively and freely on the $q_j$-orbit, we have
        \[
            \bigoplus_{q\in P\cdot q_j} T_{q_+}\widetilde{C}\otimes T_{q_-}\widetilde{C}
            \simeq \operatorname{Ind}_{P_{s-j}}^P \theta_{s-j}.
        \]
        Therefore, the sum over $1\le j\le s-1$ is 
        \[
            \sum_{j=1}^{s-1}\operatorname{Ind}_{P_j}^P\theta_j.
        \]
        \item $\bigoplus_{v}H^1(C_v, T_{C_v}(-S_v))= V_{\std}-\theta_s-\theta_s^{-1}-\sum_{j=1}^{s-1}\operatorname{Ind}_{P_j}^P\theta_j^{-1}$.

        For each $C_v$ we have an exact sequence
        \[
            0\to H^0(C_v,T_{C_v})\to \bigoplus_{q\in S_v}T_q C_v
            \to H^1\left(C_v, T_{C_v}(-S_v)\right)\to 0.
        \]
        Again we group them by $P$-orbits with a representative
        set given by $v_j$,$0\le j\le s-1$.
        \begin{itemize}
            \item If $j\neq 0$, then there are $1+|P_{s-j}/P_{s-j-1}|$
            nodes. We have equalities in representation ring
            \[
                H^0(C_{v_j},T_{C_{v_j}})=1+\theta_{s-j}+\theta_{s-j}^{-1}
                \in R(P_{s-j}).
            \]
            \[
                \bigoplus_{q\in S_{v_j}}T_qC_{v_j}=
                \theta_{s-j}+\C\cdot (P_{s-j}/P_{s-j-1})
                =\theta_{s-j}+\operatorname{Ind}_{P_{s-j-1}}^{P_{s-j}}1.
            \]
            The sum over the whole $v_j$-orbit is
            \[
                \operatorname{Ind}_{P_{s-j}}^P\left(\operatorname{Ind}_{P_{s-j-1}}^{P_{s-j}}1-1-\theta_{s-j}^{-1}\right)
                =\left(\operatorname{Ind}_{P_{s-j-1}}^P 1-\operatorname{Ind}_{P_{s-j}}^P 1\right)
                -\operatorname{Ind}_{P_{s-j}}^P\theta_{s-j}^{-1}.
            \]
            \item If $j=0$, then there are only $|P_{s-j}/P_{s-j-1}|$
            nodes. We have
            \[
                H^0(C_{v_j},T_{C_{v_j}})=1+\theta_{s-j}+\theta_{s-j}^{-1}
                \in R(P_{s-j}).
            \]
            \[
                \bigoplus_{q\in S_{v_j}}T_qC_{v_j}=
                \operatorname{Ind}_{P_{s-j-1}}^{P_{s-j}}1
            \]
            The difference is
            \[
                \operatorname{Ind}^P_{P_{s-1}}1-1-\theta_s-\theta_{s}^{-1}.
            \]
        \end{itemize}
        The total is
        \begin{align*}
            &\sum_{j=1}^{s-1}\left(\left(\operatorname{Ind}_{P_{s-j-1}}^P 1-\operatorname{Ind}_{P_{s-j}}^P 1\right)
                -\operatorname{Ind}_{P_{s-j}}^P\theta_{s-j}^{-1}\right)
                +\operatorname{Ind}^P_{P_{s-1}}1-1-\theta_s-\theta_{s}^{-1}\\
                =&\operatorname{Ind}_1^P 1- \operatorname{Ind}_{P_{s-1}}^P 1
                -\sum_{j=1}^{s-1}\operatorname{Ind}_{P_j}^P\theta_j^{-1}
                +\operatorname{Ind}^P_{P_{s-1}}1-1-\theta_s-\theta_{s}^{-1}\\
                =&V_{\std}-\theta_s-\theta_s^{-1}-\sum_{j=1}^{s-1}\operatorname{Ind}_{P_j}^P\theta_j^{-1}
        \end{align*}
        In the final step we used $\operatorname{Ind}_1^P1-1=V_{\std}$.
    \end{enumerate}
    The combination of (a) and (b) proves the desired formula.
\end{enumerate}    
\end{proof}
\begin{corollary}
    $\Mb_{0,h}^{\im\alpha}$ is discrete, so $T_x\Mb_{0,h}$ coincides
        with $\nu_{M^{\im\alpha}/M}$ (\cref{def:normal bundle of alpha fixed locus}).
\end{corollary}
Under the assumption $[-1]_Ex=-x$, the Euler class of $\Ob_d-\nu_F-\nu_{\Mb_{0,h}^{\im\alpha}/\Mb_{0,h}}$ at each point of $\Mb_{0,h}^{\im\alpha}$, which is a $Z_{\widetilde{G}}(\alpha)\times \Lambda_r$-representation\footnote{
Apparently it is only a $Z_{\widetilde{G}}(\alpha)$-representation,
but we can pull back along $Z_{\widetilde{G}}(\alpha)\times \Lambda_r\to
Z_{\widetilde{G}}(\alpha)$. To define this map, we notice
the image of $\Lambda_r\to G\subset \widetilde{G}$ lies in $Z_{\widetilde{G}}(\alpha)$ as $\Lambda_r$ is commutative, and that
image of $\Lambda_r$ commutes with $Z_{\widetilde{G}}(\alpha)$
by definition. The map is defined by multiplication in $Z_{\widetilde{G}}(\alpha)$.
}
admits a simple expression.

In the following proposition and its proof, Euler classes are understood after
pullback along the map of classifying spaces induced by the multiplication homomorphism
\[
    Z_{\widetilde G}(\alpha)\times\Lambda_r
    \longrightarrow Z_{\widetilde G}(\alpha),
\]
where $\Lambda_r$ maps into $Z_{\widetilde G}(\alpha)$ via $\alpha$,
followed by extension of the $B\Lambda_r$-coefficients to $L(E)^*$,
as in the equivariant HKR character construction.

\begin{proposition}
    Assume $[-1]x=-x$ and $d=p^a$. At $x\in \Mb_{0,h}^{\im\alpha}$ corresponding to $P_j,\theta_j$, we
    have
    \[
        \frac{e_{\eq}(\Ob_d)}{e_\eq (\nu_F)e(\nu_{\Mb_{0,h}^{\im\alpha}/\Mb_{0,h}})}=
        (-1)^s\frac{e(\theta_s)^2}{u^2}
        \in S_{S^1}^{-1}H^*\left(BZ_{\widetilde{G}}(\alpha),L(E)\right).
    \]
\end{proposition}
\begin{proof}
    We identified (as $Z_{\widetilde{G}}(\alpha)$-representations) $\Ob_d-\nu_F$ and $\nu_{\Mb_{0,p^b}^{\im\alpha}/\Mb_{0,p^b}}$
    in \cref{lem:case2 virtual bundle reduction after HKR}(3)
    and \cref{lem:computation of imalpha fixed points and normal bundle}(2).
    Since $[-1]x=-x$, for any vector bundle $V$ we have $e(V^\vee)=(-1)^{\rk V}e(V)$.
    We interpret $e_\eq(\Ob_d-\nu_F)$ as $e_\eq(\Ob_d)/e_\eq(\nu_F)$ in the displayed localized ring.
    The smoothing Euler factors are invertible there by torus localization.
   
    \begin{itemize}
        \item $\Ob_d-\nu_F$.  At each point inside $\Mb_{0,p^b}^{\im\alpha}$ we have
        \[
            \bigoplus_{j=1}^{p^b} 
            \left(T_{b_j}\Sigma_j\otimes T_{b_j} \Sigma_0\right)
            \simeq \operatorname{Ind}_{S^1}^{Z_{\widetilde{G}}(\alpha)}V_\eta
        \]
        So (at each point)
        \[
            e_\eq(\Ob_d-\nu_F)=e_\eq(V_\std)(-1)^{p^a}e(\chi)^{-1}e(\chi^{-1})^{-1}
            =\frac{e(V_\std)}{u^2}\in S_{S^1}^{-1}H^*\left(BZ_{\widetilde{G}}(\alpha),L(E)\right).
        \]
        \item $\nu_{\Mb_{0,p^b}^{\im\alpha}/\Mb_{0,p^b}}$.
        The cancellation is obvious
        \[
            e_\eq(\nu_{\Mb_{0,p^b}^{\im\alpha}/\Mb_{0,p^b}})
            =e(V_\std)e(\theta_s)^{-1}e(\theta_s^{-1})^{-1}(-1)^{s-1}.
        \]
    \end{itemize}
    Therefore, 
    \[
        \frac{e_\eq(\Ob_d-\nu_F)}{e_\eq(\nu_{\Mb_{0,p^b}^{\im\alpha}/\Mb_{0,p^b}})}
        =(-1)^s\frac{e_\eq(\theta_s)^2}{u^2}.
    \]
\end{proof}
After localizing at $S_{S^1}$, the equivariant HKR Riemann--Roch formula
\eqref{eq:equivariant HKR Rieman Roch} gives
\[
    \left(e_{\operatorname{HKR,\eq}}^L\right)^{-1} \cdot\chi_{n,p;\eq}(-)
    \colon S_{S^1}^{-1}H^*\left((\Mb_{0,h})_{h\widetilde{G}},E\right)
    \to \bigoplus_{[\alpha],x\in \Mb_{0,h}^{\im\alpha}}
    S_{S^1}^{-1}H^*\left(BZ_{\widetilde{G}}(\alpha),L(E)\right).
\]
\begin{corollary}
    The $(\alpha,x)$-component of the image of $e_\eq(\Ob_d-\nu_F)$
    under $\left(e_{\operatorname{HKR,\eq}}^L\right)^{-1} \cdot\chi_{n,p;\eq}(-)$
    is 
    \[
        (-1)^s\frac{a_{\theta_s}^2}{u^2}
    \]
    as an element of $S_{S^1}^{-1}H^*\left(BZ_{\widetilde{G}}(\alpha),L(E)\right)$.
    Consequently, its $L(E)$-equivariant pushforward is
    \[
        \pi_{!,\eq}^{L(E)}\left((-1)^s\frac{a_{\theta_s}^2}{u^2}\right)
        =(-1)^s\frac{a_{\theta_s}^2}{p^a u^2}
        \in S_{S^1}^{-1}H^*(BS^1,L(E)).
    \]
    
    Here $u$ is the generator of $H^*(BS^1,L(E))$. 
    $a_{\theta_s}$ is the image of $e(\theta_s)$\footnote{
    $\theta_s$ is a $P$-representation.
    } under $H^*(BP,E)\to H^*(B\Lambda_r,E)\to L(E)$, where the first map is pullback and
    the second comes from definition of $L(E)$.
\end{corollary}
For $x\in \Mb_{0,p^b}^{\im\alpha}$, the $\alpha$-component of
$\chi_{n,p;\eq}$, restricted to $x$, is computed by pullback to
$[\pt/Z_{\widetilde{G}}(\alpha)]\times [\pt/\Lambda_r]$,
followed by the coefficient map
$H^*\left(BZ_{\widetilde{G}}(\alpha)\times B\Lambda_r,E\right)\to H^*\left(BZ_{\widetilde{G}}(\alpha),L(E)\right)$.
\begin{proof}
    The expression in the preceding proposition is the image of
    \begin{equation}
        \label{eq:the whole HKR class}
        (-1)^s\frac{\left(e_\eq(\theta_s)+_E e_\eq(\alpha^*\theta_s)\right)^2}{u^2}
    \end{equation}
    under the map $S_{S^1}^{-1}H^*(BZ_{\widetilde G}(\alpha)\times B\Lambda_r,E)\longrightarrow S_{S^1}^{-1}H^*(BZ_{\widetilde G}(\alpha),L(E))$,
    induced by $H^*(B\Lambda_r,E)\longrightarrow L(E)^*$.
    Here $e_\eq(\theta_s)\in H^*\left(BZ_{\widetilde{G}}(\alpha),E\right)$
    is a pullback,
    and $e_\eq(\alpha^*\theta_s)\in H^*\left(B\Lambda_r,E\right)$.
    However, $\theta_s$ is a pullback from a $P$-representation with $P$
    finite, so by \cref{lem:pullback of finite group representations vanish in rational theory} $e(\theta_s)$ becomes zero under the map
    $H^*\left(BZ_{\widetilde{G}}(\alpha),E\right)\to H^*\left(BZ_{\widetilde{G}}(\alpha),L(E)\right)$.
    The $B\Lambda_r$-Euler class $e_\eq(\alpha^*\theta_s)$ maps to $a_{\theta_s}$.
    Hence \eqref{eq:the whole HKR class} maps to $(-1)^s\frac{a_{\theta_s}^2}{u^2}$
    under $S_{S^1}^{-1}H^*(BZ_{\widetilde G}(\alpha)\times B\Lambda_r,E)\longrightarrow S_{S^1}^{-1}H^*(BZ_{\widetilde G}(\alpha),L(E))$.

    The relevant map for $L(E)$ equivariant pushforward is
    $Z_{\widetilde{G}}(\alpha)\to S^1$ \footnote{
    This is the composition $Z_{\widetilde{G}}(\alpha)\subset \widetilde{G}
    \to S^1$, where the RHS map is from \eqref{eq:define widetilde G}.
    }, whose kernel has order $p^a$. Notice that $e_{\theta_s}$ lives
    in $L(E)^*$ and $u\in H^*(BS^1,L(E))$, so
 \[\pi_{!,\eq}\left((-1)^s\frac{a_{\theta_s}^2}{u^2}\right)
        = (-1)^s\frac{a_{\theta_s}^2}{u^2}\pi_{!,\eq}(1)
        =(-1)^s\frac{a_{\theta_s}^2}{p^a\cdot u^2}.
    \]
\end{proof}
What remains is pure algebra. Put $\Lambda=(\Z_p)^n$. The summation is indexed over
the decorated flags. 
\begin{equation}\label{eq:indexing data}
\begin{aligned}
&H=H_0\lhd H_1\lhd \cdots \lhd H_s=\Lambda,
    \qquad |\Lambda/H|=h,\\
&\phi\colon H\to C_l,\quad
\theta_j\colon H_j/H_{j-1}\to \C^\times.
\end{aligned}
\end{equation}
\begin{remark}
    We have
two components for one $b$ as the central sphere can lie over either $p_0$
or $p_1$; changing from $p_0$ to $p_1$ amounts to replacing $\chi$ by $\chi^{-1}$ and $\eta$ by $\eta^{-1}$ in Euler classes. Under assumption $[-1]x=-x$,
$(-1)^s\frac{a_{\theta_s}^2}{p^a\cdot u^2}$ is unchanged. This
cancels with the ambiguity $\theta_s\sim \theta_s^{-1}$ in \cref{lem:computation of imalpha fixed points and normal bundle}(1),
so the summation is indexed over \eqref{eq:indexing data}
\emph{without} identifying $\theta_s\sim \theta_s^{-1}$.
\end{remark}
We write $\lambda
=\theta_s\colon H_s/H_{s-1}\to \C^\times$. Actually, $\lambda$ lands
in $C_{p^b}\subset \C^\times$. We denote $\Lambda^\vee[p^b]$ to be
all characters $\lambda\colon \Lambda\to C_{p^b}\subset \C^\times$. In the following
definition and lemma, we fix $l$ as a power of p,
but allow $h$ to change.
\begin{definition}
\begin{enumerate}[(1)]
    \item  Let $S_b$ be the sum of $(-1)^s a_{\theta_s}^2$
    over \eqref{eq:indexing data} with $h=p^b$. Let $B_k$ be the sum
    of $(-1)^s$ over the same indexing data 
    \eqref{eq:indexing data}
    with $h=p^k$.
    \item Define
    \[
        M_b=\sum_{\lambda\in \Lambda^\vee[p^b]}a_\lambda^2.
    \]
\end{enumerate}
\end{definition}
These terms relate as follows. We 
\begin{lemma}
\label{lem:computing Bk and Sb}
    \begin{enumerate}[(1)]
        \item For $k\ge 1$, $B_k=-\sum_{j=1}^k(p^{nj}-p^{n(j-1)})B_{k-j}$.
        Using this recurrence we get
       $B_0=l^{n}, B_k=l^{n}(1-p^n),k\ge 1$.
        \item $S_b=l^n(-M_b+p^n M_{b-1})$.
    \end{enumerate}
\end{lemma}
\begin{proof}
    \begin{enumerate}[(1)]
        \item For $k\ge 1$ and a decorated flag
        $H=H_0\lhd H_1\lhd \cdots \lhd H_s=\Lambda$ we let
        $|H_s/H_{s-1}|=p^j$. A faithful character
        $H_s/H_{s-1}\to \C^\times$ is exactly a surjective character
        $\Lambda\to C_{p^j}$, the total number of which
        is $p^{nj}-p^{n(j-1)}$. As the remaining data form exactly
        a decorated flag if we view $H_{s-1}=\ker(\theta_s)$ as the new
        $\Lambda$, we get
        \[
            B_k=-\sum_{j=1}^k(p^{nj}-p^{n(j-1)})B_{k-j}.
        \]
        The $-1$ comes from an extra level by adding $H_{s-1}\subset H_{s}$
        to $H_0\lhd \cdots H_{s-1}$. 

        We can use induction to prove $B_0=l^n$, $B_k=l^n(1-p^n),k\ge 1$ satisfy
        this recurrence. The starting case $B_0=l^n$ is because
        $H=H_0\simeq (\Z_p)^n$, and $B_0$ is the number of homomorphisms
        $H\to C_l$ (not necessarily faithful).
        \item If in \eqref{eq:indexing data} we let $|H_s/H_{s-1}|=p^j$,
        then we can reorder the sum by first choosing a surjective
        $\lambda\colon \Lambda\to C_{p^j}$, and then summing $(-1)^s$ over $H_0\lhd\cdots \lhd H_{s-1}=\ker(\lambda)$.
        We see that $M_{j}-M_{j-1}$ is the sum of $a_\lambda^2$ over surjective
        characters $\lambda\colon \Lambda\to C_{p^j}$, and the latter
        count is $-B_{b-j}$.
        Therefore (let $M_{0}=0$)
        \[
            S_b=-\sum_{j=1}^b B_{b-j}(M_j-M_{j-1})
        \]
        Substituting $B_0=l^n$ and $B_k=l^n(1-p^n)$ for $k\ge 1$, we obtain $S_b=l^n(-M_b+p^n M_{b-1})$.
    \end{enumerate}
\end{proof}
For the sum $M_b=\sum_{\lambda\in \Lambda^\vee [p^b]}a_\lambda^2$,
we notice that $a_\lambda\in L(E)^*$ are exactly the $p^{bn}$ distinct roots of $[p^b]x=0$ in $L(E)^*$. We let $A=E^*$ and $B=E^*\llbracket x\rrbracket/([p^b]x)$.
\begin{lemma}
\label{lem:computing Mb}
    We have 
    \[
         M_b=\operatorname{Tr}_{B/A}(x^2)=
         p^b \res_x\left(\frac{x^2 \omega(x)}{[p^b]x}\right).
    \]
    The first equality is understood in $L(E)^*$ via the coefficient map,
    and the second holds in $A=E^*$.
\end{lemma}
\begin{proof}
    Because
    $B\otimes_{E^*}L(E)^*\simeq \prod_{\lambda\in \Lambda^\vee[p^b]} L(E)^*[x]/(x-a_\lambda)$ ($a_\lambda-a_{\lambda'}$ is invertible in $L(E)^*$ for $\lambda\ne\lambda'$), we have
    \[
        M_b=\operatorname{Tr}_{B\otimes_{E^*}L(E)^*/L(E)^*}(x^2)
        =\operatorname{Tr}_{B/A}(x^2)
    \]
    Here the trace commutes with base change because $B$ is finite free over $A$.
    This proves the first half.

    For the second half, we apply \cref{lem:residue=trace}, where
    $f(x)=[p^b]x$. We write $\omega(x)=w(x)dx$.
    By \cref{lem:invariant form} we have
    \[
        w([p^b]x)d\left([p^b]x\right)=p^b w(x)dx.
    \]
    Equivalently $\left([p^b]x\right)'=p^b \frac{w(x)}{w([p^b]x)}$.
    Therefore, by \cref{lem:residue=trace}
    \begin{align*}
        \operatorname{Tr}_{B/A}(x^2)&=
        \res_x\left(\frac{x^2 \left([p^b]x\right)'}{[p^b]x}\right)\\
        &=p^b\res_x\left(\frac{x^2 w(x)}{[p^b]x}\frac{1}{w([p^b]x)}\right)
    \end{align*}
    Since the initial term of $w(x)$ is $1$, we can write $\frac{1}{w([p^b]x)}$ as a power series of $[p^b]x$. We can drop
    positive degree terms because $\frac{x^2 w(x)}{[p^b]x}\cdot [p^b]x$
    has no $x^{-1}$ terms. So finally we have
    \[
        \operatorname{Tr}_{B/A}(x^2)=p^b\res_x\left(\frac{x^2 w(x)}{[p^b]x}\right).\qedhere
    \]
\end{proof}
Recall the notation that 
$\Theta_b=\tau_{p^b,0}(x^2)=\res_x\left(\frac{x^2 w(x)}{[p^b]x}\right)$.
\begin{proof}[Proof of \cref{prop:contribution from components}(2)]
    We want to compute the sum of $(-1)^s\frac{a_{\lambda}^2}{p^a u^2}$
    over \eqref{eq:indexing data}, which is $\frac{S_b}{p^a u^2}$ in our
    notation. The following identity is understood in $L(E)^*$, with
    $\Theta_b$ and $\Theta_{b-1}$ denoting their images under the coefficient map.
    By \cref{lem:computing Bk and Sb}(2) and \cref{lem:computing Mb} (we let $l=p^{a-b}$)
    \[
        S_b=p^{(a-b)n}(-M_b+p^nM_{b-1})=p^{(a-b)n}(-p^b\Theta_b
        +p^n p^{b-1}\Theta_{b-1})
    \]
    So 
    \[
        \frac{S_b}{p^a u^2}=\frac{1}{u^2}p^{(a-b)n-a}(-p^b\Theta_b
        +p^n p^{b-1}\Theta_{b-1})
        =\frac{1}{u^2}p^{(a-b)(n-1)}\left(-\Theta_b+p^{n-1} \Theta_{b-1}\right)
    \]
    The coefficient injection above extends to
    \[
        S_{S^1}^{-1}E^*\llbracket u\rrbracket
        \hookrightarrow L(E)^*((u)),
    \]
    since every nontrivial $S^1$-Euler class becomes invertible in
    $L(E)^*((u))$. Hence the asserted contribution formula holds in
    $S_{S^1}^{-1}E^*\llbracket u\rrbracket$.
    \end{proof}
\subsection{Evaluation in Morava \texorpdfstring{$K$}{K}-theory}
In this subsection, we evaluate \eqref{eq:nonequivariant computation in total} in Morava $K$-theories. The main result is \cref{thm:Kn-local-invariants}.

We choose $E=\widetilde{K}(n)$ to be the integral
Morava $K$-theory with $\pi_* E=\Z_p [v_n,v_n^{-1}]$, getting
the Gromov--Witten invariants $I_{p^a}^E \in \Z_p$. The same
invariant in Morava $K$-theory is then $I_{p^a}^E\mod p$. 

Recall that $p$ is odd and the formal group law of
$E=\widetilde{K}(n)$ is  $[p]x=px +_Ev_n x^{p^n}$.
The following lemma is about the invariant form
of the mod $p$ Morava $K$-theory.
\begin{lemma}
\label{lem:Kn-series-and-differential}
The invariant differential of a Morava $K$-theory $K(n)$ (not $\widetilde{K}(n)$)
$\omega(x)=g(x)\,dx$ is given by
\begin{equation}
\label{eq:Kn-invariant-differential}
g(x)=
\begin{cases}
1, & n\ge 2,\\[2mm]
\displaystyle
\sum_{j\ge 0}
v_1^{(p^j-1)/(p-1)}x^{p^j-1},
& n=1.
\end{cases}
\end{equation}
Moreover, the condition $[-1]x=-x$ is true for odd $p$.
\end{lemma}
\begin{proof}
To compute the invariant differential,
work first over the rationalized coefficient
ring of $\widetilde K(n)$. We use $F$ to denote
the formal group law of $\widetilde{K}(n)$, and write $q=p^n$.
Write $\log_F(x)=\sum_{j\ge 0}\ell_jx^{q^j}$,
$\ell_0=1$.
Applying the formal logarithm to
$[p]_F(x)=px+_Fv_nx^q$ gives
\[p\log_F(x)=\log_F(px)+\log_F(v_nx^q).\]
Comparison of coefficients yields
$(p-p^{q^j})\ell_j=
v_n^{q^{j-1}}\ell_{j-1}$, $j\ge1$.
Consequently,
\[
    \ell_j=
    \frac{v_n^{(q^j-1)/(q-1)}}
    {\prod_{i=1}^j(p-p^{q^i})}.
\]

Since $\omega(x)=\log_F'(x)\,dx$,
the coefficient of $x^{q^j-1}$ in
$g(x)$ is
\[  
q^j\ell_j
=
p^{(n-1)j}
v_n^{(q^j-1)/(q-1)}
\prod_{i=1}^j
(1-p^{q^i-1})^{-1}.
    \]
These coefficients are $p$-adically
integral. Reducing modulo $p$ gives
\eqref{eq:Kn-invariant-differential}.

Finally, all powers occurring in
$\log_F(x)$ are odd, since $p$ is odd.
Thus $\log_F(-x)=-\log_F(x)$, implying
$[-1]_F(x)=-x$ over the integral
coefficient ring and hence modulo $p$.
\end{proof}

\begin{lemma}
\label{lem:u2-coefficient-extraction}
Let $F$ be a formal group law over a ring $R$ in which
$2$ is invertible. Assume
$[-1]_F(x)=-x$, and write its invariant
differential as $\omega(x)=g(x)\,dx$.

Then, in $R((x))[[u]]$,
\begin{equation}
\label{eq:u2-coefficient-extraction}
[u^2]\frac{1}{x-_Fu}
=
\frac{1}{x^3g(x)^2}
+
\frac{g'(x)}{2x^2g(x)^3}.
\end{equation}
\end{lemma}

\begin{proof}
By \cref{lem:invariant form}, we know
\begin{equation}
    \label{eq:identify of g(x)}
     g(F(x,y))\cdot \partial_y F(x,y)=g(y)
\end{equation}
whose special case by setting $y=0$ is $\frac{1}{g(x)}=\partial_2F(x,0)$.

Since $[-1]_F(x)=-x$, the invariant
differential is even, so $g'(0)=0$.
Differentiating \eqref{eq:identify of g(x)} once more with respect
to $y$ and setting $y=0$ yields
\[
\frac{\partial^2F}{\partial y^2}(x,0)
=-\frac{g'(x)}{g(x)^3}
\]
Therefore,
\[
x-_Fu=F(x,-u)
=
x-\frac{u}{g(x)}
-\frac{g'(x)}{2g(x)^3}u^2
+O(u^3).
\]

Taking the multiplicative inverse,
we get
\[
\begin{aligned}
\frac{1}{x-_Fu}
={}&\frac{1}{x}+\frac{u}{x^2g(x)}\\
&+
\left(
\frac{1}{x^3g(x)^2}
+\frac{g'(x)}{2x^2g(x)^3}
\right)u^2
+O(u^3).
\end{aligned}
\]
Taking the coefficient of $u^2$ proves the formula.
\end{proof}
\Kn*
\begin{proof}
By \eqref{eq:nonequivariant computation in total},
\[
I_{p^a}^{K(n)}
=
[u^2]\tau_{p^a,u}(x^2).
\]
We evaluate the residue after reducing
the coefficients modulo $p$.

Put
\[
Q=p^{na},\qquad
c=v_n^{(Q-1)/(p^n-1)}.
\]
By \cref{lem:Kn-series-and-differential},
$[p^a]_F(x)=cx^Q$.

Define
\[
R(x):=[u^2]\frac{1}{x-_Fu}.
\]
By \cref{lem:u2-coefficient-extraction},
\[
R(x)=
\frac{1}{x^3g(x)^2}
+\frac{g'(x)}{2x^2g(x)^3}
\in x^{-3}K(n)^*[[x]].
\]
Hence
\begin{equation}
\label{eq:Kn-residue-reduction}
I_{p^a}^{K(n)}
=
\res_x x^2g(x)R(cx^Q)\,dx.
\end{equation}

Write $R(x)=\sum_{j\ge-3}r_jx ^j$.
A monomial $x^e$ occurring in $g(x)$
can contribute to the residue only if
\[
2+e+Qj=-1
\]
for some integer $j\ge -3$. In particular,
\begin{equation}
\label{eq:Kn-divisibility-condition}
Q\mid e+3.
\end{equation}
If $n\ge2$, then $g(x)=1$. Thus
$e=0$, and the divisibility condition
would require $Q\mid3$. This is
impossible because $Q=p^{na}\ge9$.

Suppose now that $n=1$.
The exponents appearing in $g(x)$ are
\[
e=0,\qquad e=p^k-1\quad(k\ge1).
\]
For $e=p^k-1$, the condition becomes $p^a\mid p^k+2$.
This is impossible because $p$ is odd.

For $e=0$, the condition is
$p^a\mid3$, which forces $p=3$
and $a=1$.
It remains to calculate this case.
We have
\[
[3]_F(x)=v_1x^3,
\qquad
g(x)=1+v_1x^2+v_1^4x^8+\cdots.
\]
By \eqref{eq:u2-coefficient-extraction},
\[
R(x)=x^{-3}-v_1x^{-1}+O(x).
\]
Therefore,
\[
\begin{aligned}
I_3^{K(1)}
&=\res_x x^2g(x)R(v_1x^3)\,dx\\
&=v_1^{-3}[x^6]g(x)-[x^0]g(x)\\
&=-1.
\end{aligned}
\]
This proves \eqref{eq:Kn-pa-invariants}.

Finally, $I_1^{K(n)}=1$. Applying the
prime-to-$p$ degree reduction formula
\cref{prop:separating p and p-coprime
part nonequivariantly} to $\widetilde K(n)$ and reducing modulo $p$ gives
\[
I_{p^am}^{K(n)}
=m^{-3}I_{p^a}^{K(n)}
\qquad (p\nmid m).
\]
Combining these formulas proves
\eqref{eq:Kn-general-invariants}.
\end{proof}

\section{An application to Calabi--Yau threefolds}
We use the cluster formalism from \cite{Ionel-Parker_Gopakumar-Vafa} to transfer local computation results
to Gromov--Witten invariants of Calabi--Yau threefolds,
with coefficients in $K(n)$.

From now on, all Gromov--Witten invariants in this section are
understood to be \(K(n)\)-valued. For a closed symplectic Calabi--Yau threefold \(X\) and
a nonzero class \(A\), we write
\[
\GW_{0,0}^{K(n)}(X,A)=\langle 1, [\Mb_{0,0}(X,A)]_\vir\rangle\in K(n)^0(\mathrm pt)=\F_p.
\]
where $[\Mb_{0,0}(X,A)]_\vir$
denotes the virtual fundamental class.

Here $n_{\beta,0}$ denotes the genus-zero
Gopakumar--Vafa invariant of $X$, as in the Introduction.
We prove the following theorem in this section.
\CY*
The proof depends on \cref{decomposition into elementary clusters}.
\begin{proof}[Proof of \cref{CY}]
    It suffices to identify the values
    of each term in \cref{decomposition into elementary clusters}.
    The coefficients $e_{\beta,0}(X)$ are identified with $n_{\beta,0}$ below,
    and the local contributions $I_d^{K(n)}$ are computed in \cref{thm:Kn-local-invariants}.
\end{proof}

We define the genus-zero \(K(n)\)-Gromov--Witten potential by
\[
\GW^{K(n)}_0(X,q)
=
\sum_{\omega(A)>0}
\GW^{K(n)}_{0,0}(X,A)q^A
\in
\mathbb F_p[[q^A]]_{A\in H_2(X,\Z)}.
\]
Note it is the
generating function for the \(K(n)\)-valued invariants.
For the local curve
\[
\operatorname{Tot}\bigl(\mathcal O_{\mathbb P^1}(-1)
\oplus
\mathcal O_{\mathbb P^1}(-1)\bigr),
\]
we write
\[
\GW^{K(n),\mathrm{elem}}_0(q)
=
\sum_{d\geq 1}
\GW^{K(n)}_{0,0}(X_{\loc},d[\mathbb P^1])q^d
\in \mathbb F_p[[q]],
\]
where \(q=q^{[\mathbb P^1]}\).

We argue that \cite[Theorem 8.1]{Ionel-Parker_Gopakumar-Vafa} also applies in this
setting of global Kuranishi charts
and coefficients in $K(n)$.
\begin{proposition}
\label{decomposition into elementary clusters}
There are integers \(e_{\beta,0}(X)\) such that
\[
\GW^{K(n)}_0(X,q)
=
\sum_{\omega(\beta)>0}
e_{\beta,0}(X)\,
\GW^{K(n),\mathrm{elem}}_0(q^\beta)
\]
in \(\mathbb F_p[[q^A]]\), where the integer
\(e_{\beta,0}(X)\) is viewed in \(\mathbb F_p\).
\end{proposition}
These $e_{\beta,0}(X)$ are the same as those in \cite[Theorem 8.1]{Ionel-Parker_Gopakumar-Vafa}. They are equal to Gopakumar-Vafa invariants
$n_{\beta,0}$, following from \cite[Proposition 3.4]{Ionel-Parker_Gopakumar-Vafa}. 

Although Ionel--Parker prove their decomposition formula for rational Gromov--Witten invariants,
the construction itself gives a cluster decomposition of the underlying moduli space,
equivalently of the zero locus $s^{-1}(0)/G$ of a global Kuranishi chart. 
The passage from this decomposition to $K(n)$-valued invariants uses the support property of virtual fundamental classes: the $K(n)$-oriented virtual fundamental class is supported on $s^{-1}(0)$,
and hence is unchanged after restricting the chart to any $G$-invariant open neighborhood of $s^{-1}(0)$. 
Therefore the same cluster regions decompose the $K(n)$-valued virtual pushforward.
We first express this cluster decomposition for fixed $J$ in the language of global Kuranishi charts. We then use the parametrized construction to prove invariance of the virtual contributions of cluster regions whose boundaries remain disjoint from the moduli space.

\subsection{Cluster formalism and global Kuranishi charts}
Let $\J$ be the space of all smooth $\omega$-tame almost complex structures
on $X$, equipped with the $C^\infty$ topology.
Let $\Mb(X)=\{(\Sigma,J)\}$, where $\Sigma$ is a rational curve (genus zero nodal curve) and $\Sigma\to X$ is $J$-holomorphic,
be the space of all $J$-holomorphic curves for
some $J\in \J$, and let $\M_\simple(X)$ be the subspace of simple
maps. Fix an $E>0$, and let $\Mb^E(X)$ be
the subspace of $\Mb(X)$ with $\omega[\Sigma]\le E$.
Let $\pi\colon\Mb^E(X)\to\J$ be the projection $(\Sigma,J)\mapsto J$.
For $S\subset\J$, write $\Mb^{S,E}(X)=\pi^{-1}(S)$; for $S=\{J\}$, write $\Mb^{J,E}(X)$, and let $\M_\simple^{J,E}(X)$ be its subspace of simple maps.
By \cite[Lemma 1.2 and Corollary A.5]{Ionel-Parker_Gopakumar-Vafa}, there is a Baire subset $\J_E^*\subset \J$
such that for $J\in \J_E^*$, all simple $J$-holomorphic curves
with $\omega[\Sigma]\le E$
are smooth, embedded, and regular, and any two curves do not intersect. More generally, we can define $\J^E_\mathrm{isol}$
to be the set of $J\in \J$ such that 
$\M_\simple^{J,E}(X)$ is discrete in Gromov topology, all simple curves in this fiber are smooth, embedded and any two of them are disjoint; no regularity condition is required.
The deformation of almost complex structures
happens entirely in $\J^E_\mathrm{isol}$.

Let $\mathrm{Subsets}(X)$ be the space of nonempty closed subsets of $X$, endowed
with Hausdorff topology. Fix a Riemannian metric on $X$ and denote the induced Hausdorff
distance by $d_H$. We write
\[
    B(C,\epsilon):=\{A\in\mathrm{Subsets}(X):d_H(A,C)<\epsilon\}.
\]
There is a map $c\colon\Mb(X)\to \mathrm{Subsets}(X)$,
sending a $J$-holomorphic curve $\Sigma\to X$ to its image set in $X$.
\begin{definition}
    Define
    \[
        \Ccal(X)=c\left(
            \Mb(X) 
        \right)
    \]
    \[
        \Ccal^E(X)=c\left(
            \Mb^E(X)
        \right)
    \]
    For each $J\in\J$, write $\Ccal^{J,E}(X)=c\left(\Mb^{J,E}(X)\right)$.
\end{definition}
For $J\in\J^E_\mathrm{isol}$, any $J$-holomorphic curve of energy at most $E$ admits a factorization through an embedded simple
curve, so it has the same image under $c$ as that simple curve \cite[Lemma 1.5]{Ionel-Parker_Gopakumar-Vafa}. If a curve $u\colon \Sigma\to X$
represents homology class $A\in H_2(X,\Z)$, we define
the level of this curve to be the number of prime factors, counted with multiplicity, of $\operatorname{lcm}\{k\in\Z_{>0}:A\in kH_2(X,\Z)\}$.
Let $\Ccal^E(X)_i$ be the image of curves
whose level is at most $i$. The image $\Ccal^E(X)$ is filtered by levels
\[
    \Ccal^E(X)_0\subset \Ccal^E(X)_1\subset \cdots\subset \Ccal^E(X)_m\subset \cdots
\]
For fixed $J\in\J^E_\mathrm{isol}$, let $\Ccal^{J,E}(X)_i$ be the image under $c$ of maps in $\Mb^{J,E}(X)$ whose level is at most $i$.
This fixed-$J$ filtration satisfies the following \cite[Lemma 1.6]{Ionel-Parker_Gopakumar-Vafa}
\begin{itemize}
    \item For any open neighborhood $U$ of $\Ccal^{J,E}(X)_{m-1}$, $\Ccal^{J,E}(X)_m-U$ is a finite set.
    \item $\Ccal^{J,E}(X)_m=\Ccal^{J,E}(X)$ for some
    finite $m$. In other words,
    the filtration terminates after finitely many steps.
\end{itemize}
The following definition of clusters is from \cite[Definition 2.2]{Ionel-Parker_Gopakumar-Vafa}
\begin{definition}[$E$-Cluster]
    For $E>0$, an $E$-cluster $\O=\left(C,\epsilon,J\right)$
    consists of an almost complex structure $J$, an embedded
    $J$-holomorphic curve, and a radius $\epsilon>0$ satisfying the
    following conditions:
    \begin{enumerate}[(1)]
        \item all non-constant $J$-holomorphic maps $u\colon\Sigma\to X$
        with energy $\le E$ and $d_H(u(\Sigma),C)<\epsilon$
        represent the class $k[C]\, (k\ge 1)$
        in homology and have genus $g\ge g(C)=0$.
        \item The only possibility for $k=1$ is $\Sigma=C$.
        \item There are no $J$-holomorphic curves $u\colon \Sigma\to X$ with energy $\omega[\Sigma]\le E$
        such that $d_H(u(\Sigma),C)=\epsilon$.
    \end{enumerate}
\end{definition}

A cluster $\O=\left(C,\epsilon,J_0\right)$ determines the open subset
\[
    \Mb^E_{\O}(X)=c^{-1}\bigl(B(C,\epsilon)\bigr)\cap\Mb^E(X)
\]
of the universal moduli space.
For a parameter set $S\subset\J$ satisfying $d_H(A,C)\ne\epsilon$ for every $J\in S$ and $A\in\Ccal^{J,E}(X)$, its restriction
$\Mb^{S,E}_{\O}(X):=\Mb^E_{\O}(X)\cap\pi^{-1}(S)$ is open and closed in $\Mb^{S,E}(X)$.
Its fiber over $J$ is $\Mb^{J,E}_{\O}(X):=\Mb^E_{\O}(X)\cap\pi^{-1}(J)$.
By \cite[Lemma 2.1]{Ionel-Parker_Gopakumar-Vafa}, the absence of boundary curves is an open condition on $J$.
    
    For $J\in\J^E_\mathrm{isol}$ and an embedded $J$-holomorphic curve $C$, there exists $\epsilon_C>0$
    such that $(C,\epsilon,J)$ is an $E$-cluster for all but countably many $\epsilon\in(0,\epsilon_C)$, \cite[Lemma 2.3]{Ionel-Parker_Gopakumar-Vafa}. Furthermore, by \cite[Proposition 2.4]{Ionel-Parker_Gopakumar-Vafa}, there are finitely many
    clusters $\O_i=(C_i,\epsilon_i,J)$ such that
    \[
        \Mb^{J,E}(X)=\bigsqcup_i \Mb^{J,E}_{\O_i}(X).
    \]
    The choice is as follows: first choose all simple curves $C_{i_0}$ in
    $\Ccal^{J,E}(X)_0$ and take suitable $\epsilon_{i_0}$. Once we have chosen
    finitely many $(C_{i_m},\epsilon_{i_m},J)$ that cover
    $\Ccal^{J,E}(X)_m$, $\Ccal^{J,E}(X)_{m+1}-\bigcup_{i_m}c\left(\Mb^{J,E}_{\O_{i_m}}(X)\right)$ is finite. Then we can take
    curves in this finite set, and choose suitable $\epsilon$. Finally,
    we obtain finitely many clusters disjoint from each other.

More generally, suppose that over a parameter set $S\subset\J$ the regions defined by finitely many pairwise disjoint Hausdorff balls give a decomposition
\[
    \Mb^{S,E}(X)=\bigsqcup_i\Mb^{S,E}_{\O_i}(X)
\]
into open and closed regions, as in \cite[Section 2]{Ionel-Parker_Gopakumar-Vafa}.

In the presence of a cluster decomposition of images,
there is an induced decomposition of global Kuranishi charts,
and thus a decomposition of virtual fundamental classes.
Fix $J\in S$ and let $(G,\T,\E,s)$ be a global Kuranishi chart for $\Mb^{J,E}(X)$.
Taking the fiber over $J$ of the above decomposition gives
$\Mb^{J,E}(X)=\bigsqcup_i\Mb^{J,E}_{\O_i}(X)$.
Notice that any point $(\phi,u,e,F)\in\T$ already comes with a map $u\colon\Sigma_\phi\to X$.
\begin{proposition}
\label{Cluster decomposition of global Kuranishi charts}
    \begin{enumerate}[(1)]
        \item Let $\T_{\O_i}$ be the subset of $\T$ such that $d_H(u(\Sigma_\phi),C_i)<\epsilon_i$. Then $\T_{\O_i}$ is a $G$-equivariant open subset of $\T$, and each $\left(G,\T_{\O_i},\E|_{\T_{\O_i}},s|_{\T_{\O_i}}\right)$
         is a global Kuranishi chart for
         the closed and open subset $\Mb_{\O_i}^{J,E}(X)\subset \Mb^{J,E}(X)$.
        \item The virtual fundamental classes defined
        by these global Kuranishi charts satisfy
        \[
            \left[\Mb^{J,E}(X)\right]_\vir=\sum_{i}\left[
            \Mb^{J,E}_{\O_i}(X)
            \right]_\vir
        \]
        as a map $H^*_G\left(\T,\K\right)\to \K^*$. In other words, given $a\in H^*_G(\T,\K)$, there is an identity
        \begin{equation}
        \label{vfc decomposition}
            \left\langle
                a, \left[\Mb^{J,E}(X)\right]_\vir
            \right\rangle
            =
            \sum_{i}
            \left\langle
                a|_{\T_{\O_i}},
                \left[\Mb^{J,E}_{\O_i}(X)\right]_\vir
            \right\rangle
        \end{equation}
    \end{enumerate}  
\end{proposition}    
\begin{proof}
    There is a map
        \[
            \mathfrak c_{\mathcal T}:\mathcal T\longrightarrow \operatorname{Subsets}(X),
            \qquad
        (\phi,u,e,F)\longmapsto u(\Sigma_\phi).
        \]
        We claim that $\mathfrak c_\T$ is
        continuous. $\operatorname{Subsets}(X)$
        is endowed with the Hausdorff topology.
        The topology of $\T$ can be defined following \cite[Definition 4.19]{abouzaidGromovWittenInvariantsComplex2024a}:
        define the pre-thickened moduli space
        $\T^\pre=\{(\phi,u,e)\}$ by forgetting
        the framing $F$ in \cref{def: AMS construction}. Recall $\Ccal=\Ccal_{g,h,d}$ is the universal curve
        introduced in \hyperref[sec:projective parametrization of domains]{the projective parametrization of domains}.
        For $(\phi,u,e)\in \T^\pre$, we have a fiber inclusion $\Sigma_\phi\hookrightarrow \Ccal$
        and $u\colon \Sigma_\phi\to X$.
        Taking the image of $\Sigma_\phi\to \Ccal\times X$
        and $e\in V_\mu$ defines a map
        \[
            \T^\pre\to
            \operatorname{Subsets}(\Ccal\times X)\times V_\mu
        \]
        It is injective, so we can use the Hausdorff topology on $\operatorname{Subsets}(\Ccal\times X)$ and the topology on $V_\mu$ (it's a vector space) to define a topology on $\T^\pre$.
        We give $\T$ the subspace topology from the full
        $\GL(d-g+1)$-framing bundle over $\T^{\pre,\mathrm{reg}}$,
        in which it is an open subset.

        Using this topology on $\T$, $\mathfrak c_\T$
        is a composition of continuous maps
        \[
           \T\to \T^\pre\to
            \operatorname{Subsets}(\Ccal\times X)\times V_\mu
            \xrightarrow{\pi_1}
            \operatorname{Subsets}(\Ccal\times X)
            \xrightarrow{(\pi_X)_*}
            \operatorname{Subsets}(X)
        \]
     \begin{enumerate}[(1)]
        \item  The formal definition of $\T_{\O_i}$ is
        \[
        \begin{aligned}
            \T_{\O_i}
            &=\{(\phi,u,e,F)\in\T:d_H(u(\Sigma_\phi),C_i)<\epsilon_i\}\\
            &=\mathfrak{c}_{\T}^{-1}\bigl(B(C_i,\epsilon_i)\bigr).
        \end{aligned}
        \]
        $\T_{\O_i}$ is open because $\mathfrak c_\T$ is continuous, and $\T_{\O_i}$ are pairwise
        disjoint because $\O_i$ are pairwise disjoint.
        The definition implies that $\T_{\O_i}\cap s^{-1}(0)$ is the preimage of $\Mb_{\O_i}^{J,E}(X)$,
        thus the homeomorphism $s^{-1}(0)/G\cong \Mb^{J,E}(X)$
        restricts to a homeomorphism $s^{-1}(0)\cap \T_{\O_i}/G\cong \Mb_{\O_i}^{J,E}(X)$. So
        $\left(G,\T_{\O_i},\E|_{\T_{\O_i}},s|_{\T_{\O_i}}\right)$
         is a global Kuranishi chart for
         the closed and open subset $\Mb_{\O_i}^{J,E}(X)\subset \Mb^{J,E}(X)$
         \item From 
         \[
            \Mb^{J,E}(X)=\bigsqcup_i \Mb_{\O_i}^{J,E}(X)
         \]
         there are $s^{-1}(0)\subset \bigsqcup_i \T_{\O_i}$.
         So the disjoint union
         $\bigsqcup_i \left(G,\T_{\O_i},\E|_{\T_{\O_i}},s|_{\T_{\O_i}}\right)$ is an equivalent
         global Kuranishi chart of $(G,\T,\E,s)$ with the obvious
         open embedding
         \[
             \bigsqcup_i \left(G,\T_{\O_i},\E|_{\T_{\O_i}},s|_{\T_{\O_i}}\right)
             \subset (G,\T,\E,s).
         \]
         The invariance of virtual fundamental
         classes under equivalences implies \eqref{vfc decomposition}
    \end{enumerate}
\end{proof}

\subsection{Cluster isotopy and global Kuranishi charts}
The global Kuranishi chart can be decomposed into clusters $\O=
(C,\epsilon,J)$ where $C$ is an embedded $\P^1\subset X$.
The normal bundle $N_{C/X}$ has degree $-2$.
To compare the cluster with the elementary model of \cite[Section~3]{Ionel-Parker_Gopakumar-Vafa},
we choose $J_1$ so that a neighborhood of $C$ is biholomorphic to a neighborhood of the zero section
in the standard local $\P^1$.
Transporting the AMS chart data through this identification gives the local series defined above.
A generic path $J_t$ from $J_0$ to $J_1$, fixed on $C$ and supported near $C$, then gives the following theorem \cite[Theorem~7.1]{Ionel-Parker_Gopakumar-Vafa}.
For two clusters $\O_1,\O_2$ with the same core curve, by $\GW(\O_1)\approx\GW(\O_2)$ we mean
\[
    \GW(\O_1)-\GW(\O_2)=\sum_{\O_i}\pm \GW (\O_i)
\]
where $\O_i$ are finitely many clusters whose levels
are strictly greater than the level of $\O_1,\O_2$
(they have the same level).
\begin{proposition}[Cluster isotopy, {\cite[Theorem 7.1]{Ionel-Parker_Gopakumar-Vafa}}]
\label{Cluster isotopy}
    For a cluster $\O=(C,\epsilon_0,J_0)$ where $C$
    is embedded, smooth, and regular, we have
    \[
        \GW(\O)\approx \mathrm{sign}(C,J_0)\GW^\elem_0(q^{[C]}).
    \]
\end{proposition}
We use the deformation and local models of \cite[Sections~6--7]{Ionel-Parker_Gopakumar-Vafa}.
Restricting the parametrized global chart to an isolated cluster region gives, in each homology class, a cobordism between the endpoint charts.
Additivity and cobordism invariance then give the simple-isotopy and wall-crossing relations of \cite[Lemmas~7.2--7.4]{Ionel-Parker_Gopakumar-Vafa}.
The resulting wall-crossing signs multiply to $\mathrm{sign}(C,J_0)$.

The following proof of \cref{decomposition into elementary clusters} is the same
as in \cite[Lemma 8.2, Theorem 8.3]{Ionel-Parker_Gopakumar-Vafa}
\begin{proof}[Proof of \cref{decomposition into elementary clusters}]
    Choosing $J\in\J_E^*$ for the initial cluster decomposition,
    \cref{Cluster decomposition of global Kuranishi charts} reduces the claim to showing that
    \[
        \GW_0(\O,q^{[C]})=\sum_{d\ge 1}e_d\GW^\elem_0(q^{d[C]})
    \]
    for some $e_d\in \Z$.
    Applying \cref{Cluster isotopy} once, we have
    \[
        \GW_0(\O)=\mathrm{sign}(C,J_0)\GW_0^\elem(q^{[C]})
        +\sum_{\O_i}\pm\GW_0(\O_i)
    \]
    where $\O_i$ are finitely many higher-level clusters whose core classes are positive multiples of $[C]$, necessarily $k[C]$ with $k\ge2$ in genus zero.
    The isotopy construction can be arranged so that the higher-level clusters have regular cores.
    We can therefore apply \cref{Cluster isotopy} recursively to these clusters.
    Each iteration at least doubles the core area, so finite branching and the fixed energy bound $E$ force termination.
    Then $e_d$ is the signed count of all occurrences of clusters with core class $d[C]$
    in this recursion, each weighted by the product of all error-term signs along the
    path from the initial cluster and the sign of the elementary contribution at that occurrence.
    Since the recursion uses the same cluster decompositions, isotopies, and signs as in \cite[Lemma~8.2]{Ionel-Parker_Gopakumar-Vafa}, summing these integer weights over the initial clusters gives the coefficients $e_{\beta,0}(X)$ of \cite[Theorem~8.1]{Ionel-Parker_Gopakumar-Vafa}.
\end{proof}

\newpage
\part{General Foundations}
\label{part:foundations}

\section{Algebraic topology preparations}
\subsection{Basics of stable homotopy theory}
\subsubsection{Complex oriented
ring spectra and Thom isomorphism}
We recall some basic knowledge of stable homotopy theory.
This can be found in \cite[Part III]{Adamsbook} or \cite{Rugyak1998}.

Let $\Sp$ be the stable homotopy category whose objects
are spectra. On $\Sp$, there are a suspension functor $\Sigma$
and a looping functor $\Omega$ which are inverses to
each other.  $\oplus$ denotes the direct sum and $\wedge$ denotes the smash product in $\Sp$. For $E,F\in \Sp$, we can define the internal Hom object $ \inthom(E,F)\in \Sp$. 

Examples of spectra include suspension spectra $\Sigma^\infty X$ of a pointed topological space $X$, and $\Sigma^\infty_+X:=\Sigma^\infty(X\sqcup \pt)$; the sphere spectrum $\SS=\Sigma^\infty S^0$, and its suspension and looping spectra 
\[
    \SS^n=\begin{cases}
        \Sigma^n\SS=\Sigma^\infty S^n& n\ge 0\\
        \Omega^{-n}\SS& n<0
    \end{cases}
\]

$\Sp$ becomes a closed symmetric monoidal category 
with $\wedge$ and $\inthom$, with $\SS$ being
the unit object. It is also a triangulated
category. In particular, for $E,F\in \Sp$, their Hom set
is an abelian group, which we denote by $[E,F]$.

For $E\in \Sp$, its homotopy groups consist of
homotopy classes of morphisms $\SS^n\to E$:
\[
    \pi_n(E)=E_n=E^{-n}:=[\SS^n,E].
\]
Each of these is an abelian group. For a topological space $X$,
define its $E$-homology and $E$-cohomology as
\[
    H_*(X,E)=\pi_*(\Sigma^\infty_+X\wedge E)
\]
\[
    H^*(X,E)=\pi_{-*}\inthom(\Sigma^\infty_+X,E)
\]
and we can also define the reduced homology $\widetilde H_*(X,E)$ and cohomology $\widetilde H^*(X,E)$
of a pointed space $X$ by substituting $\Sigma^\infty_+X$
with $\Sigma^\infty X$ in the aforementioned expression.
\begin{definition}[{\cite[Part III section 10]{Adamsbook}}]
\label{def:ring_spectrum}
    A \emph{ring spectrum} is $E\in \Sp$ equipped with morphisms
    $\mu\colon E\wedge E\to E$ and $\eta\colon \SS\to E$ such that up
    to homotopy, $E$ becomes a ring where $\mu$
    is the ring multiplication and $\eta\colon \SS\to E$
    is the unit map. 
    Furthermore, if $\mu$ is commutative
    up to homotopy, then we call $E$ a \emph{homotopy commutative ring spectrum}.
\end{definition}

If $E$ is a homotopy commutative ring spectrum, then $\pi_* E$ is a graded-commutative ring. The identity element $1\in \pi_0E$ is
presented by the map $\SS\to E$.
The $E$-cohomology of a topological space
$X$ carries a graded commutative product: given $a\colon \Sigma^\infty_+X\to \Sigma^k E,b\colon \Sigma^\infty_+X\to \Sigma^lE$,
define $a\cup b$ to be
\[
    \Sigma^\infty_+X\xrightarrow{\Delta_X}
    \Sigma^\infty_+X\wedge\Sigma^\infty_+ X
    \xrightarrow{a\wedge b}
    \Sigma^k E\wedge \Sigma^l E
    \simeq \Sigma^{k+l}E\wedge E
    \xrightarrow{\mu} \Sigma^{k+l} E
\]

Throughout this section, we always assume $E$ is
a homotopy commutative ring spectrum. 

Let $\xi\colon V\to B$ be a real or complex vector bundle over
a topological space $B$. The Thom space of $\xi$ is
\[
    \Th(\xi)=D(V)/S(V).
\]
Here $D(V),S(V)$ are the unit disk bundle and unit sphere bundle once $V$ is equipped with an inner product. The Thom space
is also denoted $B^\xi$ or $B^V$.
\begin{definition}[{\cite[Part III section 10]{Adamsbook}}]
    Let $k$ be the real rank of $\xi$. An $E$-orientation of the vector bundle $\xi\colon V\to B$
    is a \textit{homotopy class} of morphisms
    \[
        [\mathfrak o]\in [\Th(\xi), \Sigma^k E]
    \]
    such that restricting to each fiber $V|_{b}$, the morphism
    \[
        S^k\simeq S^{V|_{b}}\hookrightarrow
        \Th(V)\xrightarrow{\mathfrak o}\Sigma^kE
    \]
    in $\pi_k(\Sigma^kE)\simeq \pi_0(E)$ represents
    $1$ or $-1$ in $\pi_*E$. The class $[\mathfrak o]\in \widetilde H^k(\Th(\xi),E)$ is called a Thom class.
\end{definition}
Let $\gamma_1\to \CP^\infty$ be the tautological complex
line bundle and denote its Thom space by $\MU(1)$. Note
that $\MU(1)\simeq \CP^\infty$.
\begin{definition}[{\cite[Part II (2.1)]{Adamsbook}}]
\label{def:complex-orientation}
    Let $E$ be a ring spectrum. A complex orientation of
    $E$ is a homotopy class of morphisms
    \[
        \CP^\infty=\MU(1)\to \Sigma^2E
    \]
    such that its restriction to $S^2=\CP^1\subset \CP^\infty$ represents $1\in \pi_0(E)=\pi_2(\Sigma^2 E)$.
    $E$ along with such a morphism is called
    a complex-oriented homotopy commutative ring spectrum.
\end{definition}
Because $\CP^\infty $ is the classifying space of complex line bundles, \cref{def:complex-orientation} means every complex line bundle admits an $E$-orientation.
In fact, 
every complex vector bundle admits a distinguished $E$-orientation:
\begin{lemma}[Consequence of {\cite[Part II Lemma 4.6]{Adamsbook}}]
\label{lem:multiplicativity of the complex orientation}
    Let $E$ be a complex-oriented homotopy commutative ring spectrum. Then for any complex vector bundle $\xi\colon V\to B$
    we can associate a Thom class $u_\xi\in \widetilde H^{2\rk_\C {\xi}}(\Th(\xi),E)$ satisfying the following:
    \begin{enumerate}[(1)]
        \item If $f\colon C\to B$ is a pullback, 
        then
        $u_{f^*\xi}=f^*u_{\xi}$.
        \item Let $\xi_i\colon V_i\to B_i,i=1,2$ be two complex vector bundles, and let $\xi_1\boxplus \xi_2$ be the exterior direct sum.
        Then $u_{\xi_1\boxplus\xi_2}=u_{\xi_1}\cup u_{\xi_2}$. 
        This expression makes sense because $\Th(\xi_1\boxplus\xi_2)
        =\Th(\xi_1)\wedge \Th(\xi_2)$.
        \item On the tautological line bundle $\gamma_1\to \CP^\infty$,
        the class $u_{\gamma_1}$ coincides with the complex orientation
        in \cref{def:complex-orientation}.
    \end{enumerate}
\end{lemma}
The following computation of cohomology groups can
be found in \cite[Part II Lemma 2.5, 4.3]{Adamsbook}.
\begin{lemma}
\label{lem:E-cohomology of BU(n)}
    \begin{enumerate}[(1)]
        \item 
        $H^*(\CP^n,E)=E^*\llbracket x\rrbracket / x^{n+1}$ and $H^*(\CP^\infty,E)=E^*\llbracket x\rrbracket$
        \item $H^*(BU(n))=E^*\llbracket c_1^E,c_2^E,\ldots,c_n\rrbracket$
        and $H^*(BU,E)=E^*\llbracket c_1^E,c_2^E,\ldots\rrbracket$
    \end{enumerate}
\end{lemma}
The classes $c_i^E\in H^*(BU(n),E)$ are generalizations
of Chern classes, and we still call them Chern classes. The Whitney sum formula is also true: 
\begin{lemma}[{\cite[Part II Lemma 4.3.(v)]{Adamsbook}}]
    Let $V\to B,W\to B$ be two complex vector bundles
    of rank $n$ and $m$, and $V\oplus W\to B$ is their direct
    sum. Then
    \[
        c_k^E(V\oplus W)=\sum_{i+j=k}c_i^E(V)\cdot c_j^E(W)
    \]
\end{lemma}
\subsubsection{Formal group laws and residues}
What's different from ordinary Chern classes is the tensor
product of line bundles. Consider the universal map classifying
the tensor product of two line bundles $\CP^\infty\times\CP^\infty\to \CP^\infty$,
We define a formal power series $F_E(x,y)\in E^*\llbracket x,y\rrbracket$:
\[
    E^*\llbracket x\rrbracket =H^*(\CP^\infty)
    \to H^*(\CP^\infty\times \CP^\infty,E)
    =E^*\llbracket x,y\rrbracket
\]
\[
    x\mapsto F_E(x,y)
\]
This power series is a formal group law \cite[Part II Lemma 2.7]{Adamsbook}.
\begin{definition}[{\cite[Part II Chapter 1]{Adamsbook}}]
\label{def:formal group law}
    Let $R_*$ be a graded commutative ring. A \emph{formal group law} with coefficients in $R_*$ is
    $F(x,y)=\sum_{i,j\ge 0}a_{ij}x^i y^j$ with $a_{ij}\in R_{2i+2j-2}$ satisfying the following:
    \begin{enumerate}[(1)]
        \item $F(x,0)=F(0,x)=x$;
        \item $F(x,y)=F(y,x)$;
        \item $F(x,F(y,z))=F(F(x,y),z)$.
    \end{enumerate}
    For the complex oriented $E$, $F_E$ is called
    the formal group law associated to $E$. 
\end{definition}
Throughout this paper, we only consider formal group laws $F_E$ associated to a complex oriented and homotopy commutative $E$. We also write $F_E(x,y)$ as $x+_Ey$, and write
\[
    x_1+_E x_2+_E\cdots +_Ex_n
    :=F_E(x_1,F_E(x_2,\cdots, F_E(x_{n-1},x_n))).
\]
In particular, if $x_i=x,\forall i$, we write 
$[n]_E x:=\underbrace{x +_E x +_E \dots +_E x}_{n \text{ times}}$.
If $E^*$ is a graded field of characteristic $p$, we say that $E$ has finite height $n\ge 1$ if
\[
    [p]_E x=a x^{p^n}+\cdots
\]
for some nonzero $a$.
For a complete local Noetherian ring $E^*$ whose residue field has characteristic $p$,
the height of $E$ means the height of $F_E$ after reduction modulo the maximal ideal.
When $p$ is invertible in the coefficient ring, we use height $0$.

\begin{lemma}
\label{lem:invariant form}
    Let $R$ be a commutative ring, and let $F\in R\llbracket x,y\rrbracket$ be a formal group law.
    There is a unique 1-form $\omega\in R\llbracket
    x\rrbracket dx$ whose initial term is
    $1\cdot dx$ such that
    \[
        g(F(x,y))d(F(x,y))
        =g(x)dx+g(y)dy.
    \]
    where $\omega=g(x)dx$.
    In particular, for constant
    $a$, there is
    \[
        g(F(x,a))dF(x,a)=g(x)dx
    \]
\end{lemma}
Informally, $\omega$ is invariant under
translation, with ``translation'' understood
as $x\mapsto F(x,a)$ for constant $a$.
We call it the \emph{invariant form}
of the formal group law $F$.
\begin{proof}
     Define a $1$-form as follows: 
    \[
        \omega(x)=g(x)dx=\frac{1}{\partial_2 F(x,0)}dx.
    \]
    Notice that coefficients of $g(x)$
    lie in $R$ because the constant term of $\partial_2F(x,0)$ is $1$.

    For formal variables $x,y,z$, we have
    \[
        F(x,F(y,z))=F(F(x,y),z).
    \]
    Differentiate w.r.t. $z$:
    \[
        \partial_2 F(x,F(y,z))\cdot
        \partial_2F(y,z)
        =\partial_2 F(F(x,y),z).
    \]
    Evaluate at $z=0$:
    \[
        \partial_2 F(x,y)\cdot\partial_2 F(y,0)
        =
        \partial_2 F (F(x,y),0),
    \]
    i.e.
    \[
        \partial_2F(x,y)g(F(x,y))=g(y).
    \]
    Symmetrically,
    \[
        \partial_1F(x,y)g(F(x,y))=g(x).
    \]
    Therefore,
    \begin{align*}
        g(F(x,y))dF(x,y)&=
        g(F(x,y))\left(\partial_1 F(x,y)
        dx+\partial_2 F(x,y)dy\right)\\
        &=g(x)dx+g(y)dy.
    \end{align*}
    So this $\omega=g(x)dx$ satisfies the requirement.

    On the other hand, if $\omega=g(x)dx$
    is invariant under translation, i.e.
    \[
        g(F(x,y))\partial_1F(x,y)=g(x)
    \]
    for formal variables $x,y$. Let
    $x=0$ and recall the normalization $g(0)=1$,
    \[
        g(y)\partial_1F(0,y)=g(0)=1.
    \]
    Thus $g(x)$ is unique and of the form $g(x)
    =\frac{1}{\partial_2F(x,0)}$.
\end{proof}
Next we define trace and residue. 
\begin{definition}
    \label{def:trace}
    Let $A\to B$ be a ring extension such that
    $B$ is finite and free as an $A$-module.
    Then we can define the trace
    \[
        \operatorname{Tr}_{B/A}\colon B\to A
    \]
    as the trace of a matrix presentation
    of $M_b$, multiplication by $b$.
\end{definition}
\begin{definition}
\label{def:residue}
For $f(x)=\sum_{i\in \Z}a_i x^i$
such that $a_i\to 0$ when $i\to -\infty$
in $\m$-adic topology, we define
the residue of $f(x)dx$
by
\[
    \res\left(f(x)dx\right):=a_{-1}
\]
\end{definition}
The next lemma relates trace and residue
for $A$ being a Noetherian complete local ring with maximal ideal $\m$,
and $B=A\llbracket x\rrbracket / f(x)$ with $f(x)$ a Weierstrass series. Such
$B$ is free over $A$. 
Also $f(x)$ admits an inverse in
the ring of ``Laurent series satisfying 
the condition in \cref{def:residue}'':
since $f(x)$ is Weierstrass,
    there is a factorization
    \begin{equation}\label{eq:weierstrass preparation}
        f(x)=(x^n+a_{n-1}x^{n-1}+\cdots a_0) 
        (b_0+b_1 x+b_2x^2+\cdots )
    \end{equation}
    such that $b_0\in A^\times$ and
    $a_0,a_1,\ldots,a_{n-1}\in\m$, the maximal ideal.
    We \emph{define}
    \[
        \frac{1}{f(x)}:=
        \frac{1}{b_0 x^n}\left(\sum_{i=0}^\infty
        \left(-a_{n-1}x^{-1}-\cdots -a_0x^{-n}\right)^i
        \right)\left(\sum_{j=0}^\infty
        \left(-\frac{b_1}{b_0}x-\frac{b_2}{b_0}x^2-\cdots\right)^j
        \right)
    \]
    This is well-defined because the coefficient of each $x^i$ converges in the $\m$-adic topology
    of $A$. It is easy to see this expression
    satisfies the condition in \cref{def:residue}.
\begin{lemma}
\label{lem:residue=trace}
    For such $A$ and $B=A\llbracket x\rrbracket /(f(x))$, let $g(x)\in A\llbracket x\rrbracket $, we have
    \[
        \operatorname{Tr}_{B/A}\left([g(x)]\right)
        =\res\left(\frac{g(x)f'(x)}{f(x)}dx\right)
    \]
\end{lemma}
\begin{proof}
    Write $f(x)=p(x)u(x)$ as in \eqref{eq:weierstrass preparation}.
    Then
    \[
        \res\left(\frac{g(x)f'(x)}{f(x)}dx\right)=\res\left(
            g(x)\frac{p'(x)}{p(x)}
        \right)
        +
        \res\left(
            g(x)\frac{u'(x)}{u(x)}
        \right).
    \]
    Since $\frac{1}{u(x)}$ has no negative
    degree terms, $\res\left(
            g(x)\frac{u'(x)}{u(x)}
        \right)=0$. So we have reduced
        to the case when $f(x)$ is a monic polynomial.

        Write $f(x)=x^n+a_{n-1}x^{n-1}+\cdots a_0$. We can assume $g(x)$ is
        a polynomial since both sides only
        depend on $g(x)\mod f(x)$.
        For a fixed $g(x)$, the trace
        and the residue are polynomials
        in $a_i$, so we can assume
        $f(x)$ splits as a product
        of linear functions with
        \[
            f(x)=(x-\alpha_1)\cdots (x-\alpha_n)
        \] 
        and $\alpha_i-\alpha_j,i\neq j$
        is invertible.
        Then, we have an isomorphism
        \[
            B\simeq A[x]/(x-\alpha_1)\oplus\cdots\oplus
            A[x]/(x-\alpha_n).
        \]
        So the polynomial $g(x)$ acts
        on the $i$-th $A$ by just multiplying $g(\alpha_i)$. This implies
        \[
            \operatorname{Tr}(g(x))=\sum_{i=1}^n g(\alpha_i).
        \]
        On the other hand
        \[
            \res\left(g(x)\frac{f'(x)}{f(x)}\right)
            =\sum_{i=1}^n\res\left(\frac{g(x)}{x-\alpha_i}
            \right)
            =\sum_{i=1}^n g(\alpha_i)
        \]
        So both sides are equal.
\end{proof}
\subsubsection{Morava
\texorpdfstring{$K$}{K}-theories, \texorpdfstring{$E$}{E}-theories and HKR theories}
We introduce the Morava $K$-theories $K(n)$ and Morava $E$-theories $E_n$ or more generally any HKR theory needed
in computation. For the machinery
used to define them, we refer to
\cite[Lecture 22]{LurieChromatic252x} or \cite{Wurgler-MoravaKtheory}. Let $p$ be a prime number and $n\ge 1$.
\begin{remark}
\label{rmk: homotopy commutative requirement}
    At odd primes, Morava $K$-theory admits a \textit{homotopy} commutative ring spectrum structure.
    At $p=2$, $K(n)$ admits a homotopy associative multiplication
    \cite[Theorem 1.5]{Wurgler-MoravaKtheory}. See also the discussion following Definition~5 in Lecture~22 of \cite{LurieChromatic252x}. $E_n$ is always homotopy commutative
    regardless of $p$. We will always assume $p\ge 3$ when
    using Morava $K$-theory.
\end{remark}
We summarize the properties of Morava $K$-theory here:
\begin{lemma}[{\cite[Lemma 5.13]{abouzaid2021complexcobordismhamiltonianloops}\cite[Theorem 1.3]{Wurgler-MoravaKtheory}}]
\label{lem:properties of K(n)}
    Let $E=K(n)$ denote a Morava $K$-theory.
    \begin{enumerate}[(1)]
        \item $\pi_*K(n)=\F_p[v_n,v_n^{-1}]$, where $|v_n|=2(p^n-1)$.
        \item Kunneth formula: if $X,Y$ are topological spaces, then
        \[
            H^*(X\times Y,E)\simeq H^*(X,E)\otimes_{E^*}H^*(Y,E)
        \]
        \item Its formal group law satisfies
        \[
            [p]_{K(n)}=v_{n}^{}x^{p^n}
        \]
        In particular, $F_{K(n)}$ is of height $n$.
    \end{enumerate}
\end{lemma}
Besides Morava $K$-theory, we will also use ``HKR theories''. The Morava $E$-theory is
an example.
\begin{definition}[{\cite{Generalized-Group-Characters-and-Complex-Oriented-Cohomology-Theories}\cite{strickland2000knlocaldualityfinitegroups}}]
\label{def:definition of HKR theory}
A complex oriented homotopy commutative ring spectrum
$E$ is called an ``HKR theory'' if it satisfies:
\begin{enumerate}[(1)]
    \item $E^*$ is a complete local Noetherian
    ring with maximal ideal $\m$.
    \item Characteristic of the field $E^*/\m$ 
    is $p$.
    \item $p^{-1}E^*$ is nonzero.
    \item The formal group law
    on $E^*/\m$ induced by the formal group law
    of $E$ has height $n$. Equivalently,
    in the original ring $E^*$ write
    $[p]x=px+\sum_{i\ge 2}a_ix^i$, then
    $a_1(=p),a_2,\ldots,a_{p^n-1}\in \m$ while
    $a_{p^n}$ is invertible.
\end{enumerate}
\end{definition}
The mod $p$ Morava $K$-theory $K(n)$ itself is not an HKR theory. There are the following examples of
HKR theories (see discussions after Theorem C in \cite{Generalized-Group-Characters-and-Complex-Oriented-Cohomology-Theories}
for more examples)
\begin{itemize}
    \item The Morava $E$-theory $E_n$ where
    \[\pi_* E_n=W(\F_{p^n})\llbracket
u_1,u_2,\ldots,u_{n-1}\rrbracket [u,u^{-1}]
\]
with $|u_i|=0$ and $|u|=2$. In the standard $p$-typical coordinate, the formal
group law satisfies
\[
[p]_{E_n}=px+_{E_n} u_1 u^{p-1}x^p+_{E_n}u_2u^{p^2-1}x^{p^2}+_{E_n}\cdots+_{E_n}u_{n-1}u^{p^{n-1}-1}x^{p^{n-1}}+_{E_n} u^{p^n-1}x^{p^n}.
\]
\item The integral lift $\widetilde{K}(n)$
of Morava $K$-theory where 
\[
    \pi_*\widetilde{K}(n)
    =
    \Z_p[v_{n}^{\pm 1}]
\]
where $|v_n|=2(p^n-1)$ with
\[
    [p]x=px+_{\widetilde{K}(n)} v_n x^{p^n}
\]
\end{itemize}
Notice the $[p]x$, and more generally $[p^m]x$, of an HKR theory
satisfy the Weierstrass  preparation theorem:
the coefficient of $x^{p^n}$ is a unit, and the coefficients
of lower terms lie in the maximal ideal $\m$ (and because
$\pi_*E$ is complete w.r.t $\m$). In particular (let $C_d$ be the cyclic
group of order $d$), we have that $H^*(BC_{p^m},E)=E^*\llbracket x\rrbracket /[p^m]x$ is free over $E^*$ of rank $p^{mn}$
\cite[Corollary 5.8]{Generalized-Group-Characters-and-Complex-Oriented-Cohomology-Theories}.

\subsection{Orbifolds and group actions}
\subsubsection{Orbifolds as groupoids}
We adopt the viewpoint that orbifolds are (presented by)
a special kind of Lie groupoid, see 
\cite{Moerdijk-Pronk1997,Moerdijk-orbifolds_as_groupoids2002,Lerman-orbifolds_as_stacks} for
discussions of this idea. 
\begin{definition}
    A Lie groupoid $\mathcal{G}$ consists of the following data:
    \begin{enumerate}[(1)]
        \item $\G_0,\G_1$ are two smooth manifolds.
        \item $s,t\colon \G_1\to \G_0$ are two smooth submersive maps. $e\colon \G_0\to \G_1$ and 
        $\mathrm{rev}\colon \G_1\to \G_1$ are smooth
        maps.
        \item Let $\G_1\times_{t,\G_0,s}\G_1$ be the
        fiber product of
        \[\begin{tikzcd}
	& {\G_1} \\
	{\G_1} & {\G_0}
	\arrow["t"', from=1-2, to=2-2]
	\arrow["s", from=2-1, to=2-2]
\end{tikzcd}\]
There is a smooth map $\mathrm{comp}\colon\G_1\times_{t,\G_0,s}\G_1\to \G_1$
For $a,b\in \G_1$ such that $t(a)=s(b)$ (we call such
$a,b$ composable), i.e. $(a,b)
\in \G_1\times_{t,\G_0,s}\G_1$, we also denote 
$a\circ b=\mathrm{comp}(a,b)$.
    \end{enumerate}
    and we require:
    \begin{enumerate}[(a)]
        \item $s\circ e=t\circ e=\id_{\G_0}$, $\mathrm{rev}\circ e=e$.
        \item $\mathrm{rev}\circ\mathrm{rev}=\id_{\G_1}$, 
        $s\circ \mathrm{rev}=t$ and $t\circ \mathrm{rev}=s$.
        \item For $a,b\in \G_1, t(a)=s(b)$, $s(a\circ b)=s(a)$, $t(a\circ b)=t(b)$.
        \item For $a,b,c\in \G_1$ such that $t(a)=s(b),t(b)=s(c)$, we have
        \[
            (a\circ b)\circ c =a\circ (b\circ c).
        \]
        \item For $a\in \G_1, s(a)=x,t(a)=y$, we have
        \[
            e(x)\circ a=a\circ e(y)=a
        \]
        \item For $a\in \G_1$, we have
        \[
            a\circ \mathrm{rev}(a)=\id_{s(a)}.
        \]
        \[
            \mathrm{rev}(a)\circ a=\id_{t(a)}.
        \]
    \end{enumerate}
\end{definition}
This definition is saying that $(\G_0,\G_1)$
is a groupoid internal to smooth manifolds. Points in $\G_0$
and $\G_1$ are objects and morphisms of the category 
and $s,t$ are the source and terminal of morphisms.
$e$ stands for the identity morphism, $\mathrm{rev}$
stands for the inverse morphism, and $\mathrm{comp}$
means composition of morphisms. 

Given two Lie groupoids $\G,\H$,
a functor $F$ from $\G$ to $\H$ consists of smooth maps $F_i\colon \G_i\to \H_i,i=0,1$ such that all structure maps ($s,t,e,\mathrm{rev},\mathrm{comp}$) commute with $F_0,F_1$.
\begin{definition}
    A functor $F\colon \G\to \H$ is a Morita equivalence
    iff
    \begin{enumerate}[(1)]
        \item Essential surjectivity: there is a map from $\G_0\times_{F_0,\H_0,s}\H_1$ to $\H_0$
        sending $x\in \G_0,a\in \H_1, F_0(x)=s(a)
        $ to $t(a)$. Then
        \[
            \G_0\times_{F_0,\H_0,s}\H_1\xrightarrow{t} \H_0
        \]
        is surjective and submersive.
        \item Full faithfulness: the following diagram is cartesian
        \[\begin{tikzcd}
	{\G_1} & {\H_1} \\
	{\G_0\times \G_0} & {\H_0\times \H_0}
	\arrow["{F_1}", from=1-1, to=1-2]
	\arrow["{s\times t}", from=1-1, to=2-1]
	\arrow["{s\times t}", from=1-2, to=2-2]
	\arrow["{F_0\times F_0}", from=2-1, to=2-2]
\end{tikzcd}\]
    \end{enumerate}
\end{definition}
\begin{remark}
\label{rmk:2-category of lie groupoid}
    In fact, ``Lie groupoids'' should be realized as a
    2-category. According to the discussion in \cite[331--332]{Lerman-orbifolds_as_stacks}, there are three ways to describe this 2-category:
    \begin{enumerate}[(1)]
        \item The Lie groupoids, functors, and natural transformations of functors (defined the same
        as in categories) form a 2-category and let $W$
        be Morita equivalences, and consider
        localization of Lie groupoids 2-category at $W$. This localized category
        is denoted $\mathbf{Gpoid}$.
        \item Lerman in \cite[Definition 3.25,3.26]{Lerman-orbifolds_as_stacks}
        constructed a 2-category $\mathbf{Bib}$ whose objects
        are Lie groupoids, 1-morphisms are called ``bibundles'' and 2-morphisms are diffeomorphisms
        of bibundles. 
        \item Geometric stacks over manifolds.
    \end{enumerate}
    These three ways give equivalent 2-categories
    \cite{Pronk1996,Lerman-orbifolds_as_stacks}.
\end{remark}
We say a Lie groupoid $\G$ is \emph{proper} if $s\times t\colon \G_1\to\G_0\times \G_0$ is proper; and \emph{\'etale} if $s$ and $t$ are both local diffeomorphisms. $\G$ is a \emph{foliation} groupoid if for every $x\in \G_0$, $G_x:=\{a\in \G_1|s(a)=t(a)=x\}$ is discrete. An \'etale Lie groupoid is a foliation groupoid, and
every foliation Lie groupoid
is Morita equivalent to an \'etale one \cite[1.5(c)]{Moerdijk-orbifolds_as_groupoids2002}.
Orbifolds are proper foliation/\'etale Lie groupoids
up to Morita equivalences. More precisely:
\begin{definition}
    We define the 2-category of orbifolds to be the full subcategory of proper foliation/\'etale Lie groupoids
    of the 2-category in \cref{rmk:2-category of lie groupoid}
\end{definition}
Let $\G$ be a Lie groupoid, we can associate to it a simplicial space $N_\bullet(\G)$ with $N_0(\G)=\G_0,N_1(\G)
=\G_1$ and $N_n(\G)=\G_1\times_{t,G_0,s}\G_1\times \cdots \times_{t,G_0,s}\G_1$ has $n$ copies of $\G_1$. $N_n(\G)$
consists of composable $a_1,a_2,\ldots,a_n\in \G_1$.
The maps $d_i\colon N_n(\G)\to N_{n-1}(\G)$ are defined as follows
\[
    d_0=(a_2,a_3,\ldots,a_n), d_n=(a_1,a_2,\ldots,a_{n-1})
\]
\[
    d_i=(a_1,\ldots,a_{i-1},a_i\circ a_{i+1},a_{i+2},\ldots,
    a_{n}), 1\le i\le n-1
\]
and the maps $s_i\colon N_{n}(\G)\to N_{n+1}(\G)$ are inserting $e(s(a_i))=e(t(a_{i-1}))$ between $a_{i-1}$ and $a_i$.
\begin{definition}\label{def:classifying space}
    Let $\X$ be an orbifold presented by a Lie groupoid $\G$.
    \begin{enumerate}[(1)]
        \item We define $B\G:=|N_\bullet(\G)|$, i.e.
        the geometric realization of $N_\bullet(\G)$.
        The homotopy type of $B\G$ doesn't change
        under Morita equivalences. We use $\Ho(\X)$
        to denote the homotopy type.
        \item For a ring spectrum $E$, the $E$-homology
        and $E$-cohomology of $\X$ are defined to
        be those of $B\G$:
        \[
            H_*(\X,E):= H_*(\Ho(\X),E).
        \]
        \[
            H^*(\X,E):=H^*(\Ho(\X),E).
        \]
    \end{enumerate}
\end{definition}
\subsubsection{Orbifolds with group actions and equivariant
localization}
Let $K$ be a compact Lie group, and let $\X$ be an orbifold.
\begin{definition}[General definition of group action]
\label{def:weak action}
    A $K$-action on $\X$ consists of a 1-morphism
    $a\colon K\times \X\to \X$ together with invertible 2-morphisms
    
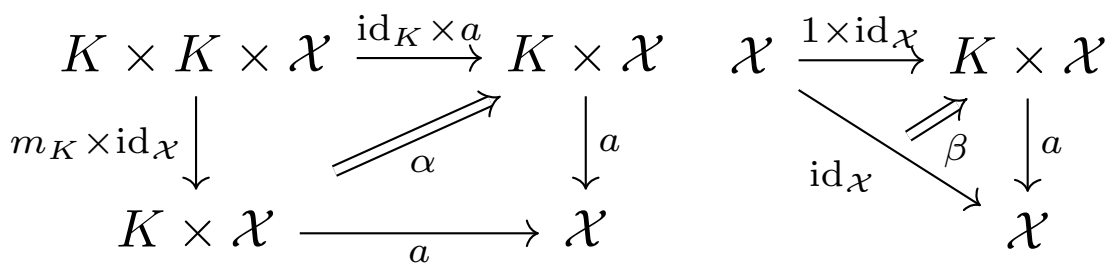
\begin{figure}[H]
        \centering
        \begin{tikzcd}
	{K\times K\times \X} & {K\times \X} \\
	{K\times \X} & \X
	\arrow["{\id_K\times a}", from=1-1, to=1-2]
	\arrow["{m_K\times \id_\X}"', from=1-1, to=2-1]
	\arrow["a", from=1-2, to=2-2]
	\arrow["\alpha"', between={0.2}{1}, Rightarrow, from=2-1, to=1-2]
	\arrow["a"', from=2-1, to=2-2]
\end{tikzcd}
\begin{tikzcd}
	\X & {K\times \X} \\
	& \X
	\arrow["{1\times \id_\X}", from=1-1, to=1-2]
	\arrow[""{name=0, anchor=center, inner sep=0}, "{\id_\X}"', from=1-1, to=2-2]
	\arrow["a", from=1-2, to=2-2]
	\arrow["\beta"', between={0.2}{1}, Rightarrow, from=0, to=1-2]
\end{tikzcd}
        \label{fig:associativity condition for group action
        on orbifolds}
    \end{figure}
    They further satisfy coherence conditions. For $k,k_1,k_2\in K$,
    $x\in \X$ (or more generally $x\colon U\to \X$
    for a test manifold $U$), write $k\cdot x:=a(k,x)$,
    $\alpha_{k_1,k_2,x}\colon (k_1\cdot k_2)\cdot x\Rightarrow k_1\cdot(k_2\cdot x)$ and $\beta_x\colon 
    x\Rightarrow 1\cdot x$. Then we further require
    \begin{enumerate}[(1)]
        \item Pentagon condition for $\alpha$: the following
        diagram commutes
        \[
\begin{tikzcd}[column sep=1em]
& ((k_1k_2)k_3)\cdot x
  \arrow[rr, Rightarrow, "\alpha_{k_1k_2,k_3,x}"]
   \arrow[dl, equal] 
&& (k_1k_2)\cdot (k_3\cdot x)
  \arrow[dr, Rightarrow, "\alpha_{k_1,k_2,k_3\cdot x}"]
\\
(k_1(k_2k_3))\cdot x
  \arrow[drr, swap, Rightarrow, "\alpha_{k_1,k_2k_3,x}"]
&&&&
k_1\cdot \bigl(k_2\cdot (k_3\cdot x)\bigr)
\\
&& k_1\cdot ((k_2k_3)\cdot x)
  \arrow[urr, swap, Rightarrow, "{\id_{k_1}\cdot \alpha_{k_2,k_3,x}}"]
\end{tikzcd}
\]
        \item The unit identities:
        \begin{figure}[H]
    \centering
    \[
    \begin{tikzcd}[column sep=large, row sep=large]
    & k\cdot x
      \arrow[dl, equal]
      \arrow[dr, Rightarrow, "\beta_{k\cdot x}"]
    \\
    (1k)\cdot x
      \arrow[rr, Rightarrow, swap, "\alpha_{1,k,x}"]
    &&
    1\cdot (k\cdot x)
    \end{tikzcd}
    \qquad
    \begin{tikzcd}[column sep=large, row sep=large]
    & k\cdot x
      \arrow[dl, equal]
      \arrow[dr, Rightarrow, "{\id_k\cdot \beta_x}"]
    \\
    (k1)\cdot x
      \arrow[rr, Rightarrow, swap, "\alpha_{k,1,x}"]
    &&
    k\cdot(1\cdot x)
    \end{tikzcd}
    \]
    \label{fig:unit-coherence}
\end{figure}
    \end{enumerate}
\end{definition}
\begin{definition}
     Let $K$ be a compact Lie group acting on $\X$ (\cref{def:weak action}). We can form the simplicial
     object\footnote{This diagram is commutative in the 2-categorical sense. The cohering 2-morphisms come
     from $\alpha,\beta$ in \cref{def:weak action}}.
\[
\begin{tikzcd}[column sep=1.2em]
\cdots
\arrow[r, shift left=1.5ex]
\arrow[r, shift left=0.5ex]
\arrow[r, shift right=0.5ex]
\arrow[r, shift right=1.5ex]
&
K^3\!\times\!\X
\arrow[r, shift left=1.1ex]
\arrow[r]
\arrow[r, shift right=1.1ex]
&
K^2\!\times\!\X
\arrow[r, shift left=0.8ex]
\arrow[r]
\arrow[r, shift right=0.8ex]
&
K\!\times\!\X
\arrow[r, shift left=0.6ex]
\arrow[r, shift right=0.6ex]
&
\X .
\end{tikzcd}
\]
The geometric realization/homotopy colimit of this diagram is called
the \emph{quotient stack} of $K\acts X$, denoted by
$[\X/K]$. 
\end{definition}
There's
a fiber sequence of stacks
\[\begin{tikzcd}
	\X & {[\X/K]} \\
	& {[*/K]}
	\arrow[from=1-1, to=1-2]
	\arrow["\pi", from=1-2, to=2-2]
\end{tikzcd}\]
Taking the homotopy type of these stacks gives a fiber
sequence of spaces
\[\begin{tikzcd}
	{\Ho(\X)} & {\Ho(\X/K)} \\
	& BK
	\arrow[from=1-1, to=1-2]
	\arrow["\pi", from=1-2, to=2-2]
\end{tikzcd}\]

Although one can define group actions on orbifolds using
stack language, it is easier to work with ``group actions
on a specific Lie groupoid presentation'', although
this notion does not transfer under Morita
equivalences.
We restrict ourselves to 
global quotient orbifolds $[M/G]$\footnotemark, where $G$ is a compact Lie group and $M$ admits an almost-free $G$-action, meaning
that the stabilizer groups of $p\in M$ are discrete (and thus finite). 
\footnotetext{The action groupoid is $\G_0=M,\G_1=G\times M$, $s(g,x)=x,t(g,x)=gx, e(x)=(\id_G,x), \mathrm{rev}(g,x)=
    (g^{-1},gx)$ and $\mathrm{comp}((g,x),(h,gx))=(hg,x)$. }

If $[M/G]$ is a global quotient orbifold, $N\subset G$ is a normal
subgroup such that $N$ acts on $M$ freely. We construct
\begin{equation}
\label{quotient construction}
    (G,M)\rightsquigarrow \left(G/N, M/N\right)
\end{equation}
It's easy to see 
$[M/G]$ and $\left[(M/N)/(G/N)\right]$ are isomorphic orbifolds. Conversely, any two isomorphic 
global quotient orbifolds $[M/G]$ and $[N/H]$ can be connected
by a zigzag of \eqref{quotient construction}. In fact, there is a manifold $P$ with commuting
$G\times H$ action, such that the actions of $G$ and $H$
on $P$ are both free, and $M\simeq P/H,N\simeq P/G$
as $G$ or $H$-manifolds; this is the special case
of a bibundle representing an equivalence, see 
\cite[Definition 3.40]{Lerman-orbifolds_as_stacks} or 
\cite[Theorem 2.26]{Behrend-Xu-differentiable_stacks}.   

In our application, group actions come in the following
form:
\begin{definition}[Working definition of group action]
\label{def:working definition of group action on orbifold}
    Let $K$ be a compact Lie group, and let $[M/G]$ be a global quotient
    orbifold. A $K$-action on \textit{the presentation} $(G,M)$\footnote{Warning: this notion is not necessarily transferable under equivalences like
    \eqref{quotient construction}. But if $P$ admits a commutative $G\times H$
    action such that the individual actions of $G$ and $H$ are free,
    we can transfer such a definition from $M=P/H$
    to $P$. Maybe one can compare two group actions on two global quotient presentations by lifting them to $P$.} consists of the following data:
    \begin{enumerate}[(1)]
        \item An extension sequence
        \begin{equation}
        \label{eq:working definition exact sequence of groups}
            1\to G\to \widetilde{G}\to K\to 1
        \end{equation}
        \item A $\widetilde G$-action on $M$ extending the
        $G$-action.
    \end{enumerate}
\end{definition}

If the action comes from \cref{def:working definition of group action on orbifold}, the quotient stack is
$[M/\widetilde{G}]$, i.e. the Lie groupoid $\widetilde{G}\times M\to M\times M$. The fiber sequence
is
\begin{equation}
\label{eq:homotopy quotient fiber sequence}
    \begin{tikzcd}
	{M_{hG}} & {M_{h\widetilde G}} \\
	& BK
	\arrow[from=1-1, to=1-2]
	\arrow[from=1-2, to=2-2,"\pi"]
\end{tikzcd}
\end{equation}
We will
define a notion of fixed points
relevant for equivariant localization.
For $p\in M$, let $\widetilde{G}_p\subset \widetilde{G}$ be
its stabilizer group. The following definition
appears in section 5.2 of \cite{Hirschi_Properties_of_Gromov_Witten}
\begin{definition}
\label{def:fixed points of orbifolds}
    Given a $K$-action on $(G,M)$, define 
    \[
        M^\fix=
        \{
            p\in M|\Stab{\widetilde{G}}(p)\to K
            \text{ is surjective.}
        \}
    \]
    $G$ preserves the set $M^\fix$.
    We define the fixed point $[M/G]^K:=[M^\fix/G]$. 
\end{definition}
\begin{remark}
    This definition is exactly the ``fixed locus''
    in \cite[Definition A.9]{virtual_localization}.
    Let $a\colon K\times \X\to \X$ be a weak action.
    Define the \emph{action inertia stack} $I_K(\X)$
    to be a fiber product
    \[\begin{tikzcd}
	{I_K(\X)} & \X \\
	{K\times \X} & {\X\times \X}
	\arrow[from=1-1, to=1-2]
	\arrow["p", from=1-1, to=2-1]
	\arrow["{\Delta\X}", from=1-2, to=2-2]
	\arrow["{(a,\pi_2)}", from=2-1, to=2-2]
\end{tikzcd}\]
Then $\X^K$ is defined to be the locus of $x\in X$
such that $I_K(\X)\times_{\X}x\xrightarrow{p\times_{\X}x}
K$ is an epimorphism of stacks\footnote{
More precisely, we define the functor $\X^K$ as follows:
for any test manifold $U$, $\X^K(U)$ is the full subcategory of $\X(U)$ 
whose objects are
$x\in \X(U)$, i.e. $x\colon U\to \X$ such that the
pullback of $\pi\colon I_K(\X)\to K\times \X$ along $x\colon U\to \X$,
$I_K(\X)\times_{\X}U\to K\times U$ is an epimorphism
of stacks.}. Therefore, the
\cref{def:fixed points of orbifolds} is presentation independent.
\end{remark}
\begin{remark}
    For connected $K$, $M^\fix$ is a submanifold
    of $M$, so we can talk about normal vector bundle
    $\nu_{M^\fix/M}$. For example, this is the case for our main application where $K=\TT$ is a torus. 
\end{remark}

We provide another description of $M^\fix$.
Let $M/G$ be the coarse quotient space,
then $\widetilde G$ action on $M$ descends
to a $K$-action on $M/G$. Let $(M/G)^K$
be the points fixed by $K$ in the ordinary
sense, then $M^\fix$ is the preimage
of $(M/G)^K$ under the quotient map 
$M\to M/G$.

Now specialize to the case $K=\TT=\TT^m$ is a torus of rank
$m$,
with an extension $1\to G\to \widetilde{G}\to \TT\to 1$. 
$H^*(B\TT,E)=E^*\llbracket \alpha_1,\ldots,\alpha_m\rrbracket$ by \cref{lem:E-cohomology of BU(n)}. 
We have the following version of Atiyah-Bott
localization for $E$-cohomology, where
$E$ is either Morava $K$-theory or $E$-theory.
See \cite[Lemma 3.6]{bai2024equivariantformalitycomplexorientedtheories} for the case of a torus acting on a compact manifold.
\begin{lemma}
    \label{isomorphism upto torsion lemma}
    Suppose that $[M/G]$ is compact.
    The restriction map $H^*_{\widetilde 
    G}(M,E)\to H^*_{\widetilde 
    G}(M^\fix,E)$ becomes an isomorphism after tensoring with
    $E^*(\alpha_1,\ldots,\alpha_m)$, the fraction field of $H^*(B\TT,E)$.
\end{lemma}

We postpone the proof of this lemma in a streamlined account
of equivariant localization in \cref{section: equivariant localization}; it is contained in
\cref{lem:combined form equivariant localzation}.
We will apply this localization lemma to global Kuranishi charts coming from the Gromov--Witten moduli problem when
the symplectic manifold admits a Hamiltonian torus
action. 
\subsubsection{Vector orbi-bundle and stably almost complex structure}
Let $\X$ be an orbifold with a Lie groupoid presentation
$\G=(\G_0,\G_1)$. A real or complex vector orbi-bundle
$\V$ is a vector bundle $V_0\to \G_0$ along with 
an isomorphism of vector bundles $\varphi\colon s^*V_0\to t^* V_0$ on $\G_1$ such that
\begin{enumerate}[(1)]
    \item at $e(x)$, $\varphi\colon V_0|_x\to V_0|_x$
    is the identity.
    \item For $a\in \G_1$, $\varphi(\mathrm{rev}(a))=\varphi(a)^{-1}$
    \item If $a,b\in \G_1$ are composable, then
    the composition $\varphi(b)\circ \varphi(a)=\varphi(\mathrm{comp}(a,b))$.
\end{enumerate}
An example is the tangent orbi-bundle of $\X$,
denoted $T\X$.
Let $\G$ be a proper \'etale Lie groupoid, 
then define $V_0= T\G_0$, and the isomorphism
is $t_*\circ (s_*)^{-1}$ because $s,t$
are local diffeomorphisms.

If $\V$ is a vector orbi-bundle over $\X$,
then on the classifying space $\Ho(\X)$ we can construct
a vector bundle with the same fiber as $\V$, which
is denoted $\Ho\V$.
The local model is, let $V\to M$ be a $G$-equivariant vector
bundle, then we can construct the vector bundle
$EG\times_G V\to EG\times_G M
=M_{hG}$.

\begin{definition}
\label{def:stably almost complex structures on vector orbi-bundles}
\begin{enumerate}[(1)]
    \item
    Let $\V\to\X$ be a real vector orbi-bundle.
    A \emph{stably almost complex structure} on $\V$ is an
    equivalence class of triples
    \[
        (\W_{\C},k,J),
    \]
    where $\W_{\C}\to\X$ is a complex vector orbi-bundle,
    $k\geq 0$, and $J$ is a complex structure on
    \[
        \V\oplus (\W_{\C})_{\R}\oplus\underline{\R}^{k}.
    \]

    \item
    Let $\xi=\V_1-\V_2$ be a virtual real vector orbi-bundle
    (with fixed presentation). A \emph{stably almost complex structure}
    on $\xi$ is an equivalence class of quadruples
    \[
        (k,\W_{1,\C},\W_{2,\C},\psi),
    \]
    where $\W_{1,\C},\W_{2,\C}\to\X$ are complex vector orbi-bundles,
    $k\geq0$, and
    \[
        \psi\colon
        \V_1\oplus(\W_{1,\C})_{\R}\oplus\underline{\R}^{k}
        \xrightarrow{\;\simeq\;}
        \V_2\oplus(\W_{2,\C})_{\R}
    \]
    is an isomorphism of real vector orbi-bundles.

    In both cases, the equivalence relation is generated by:
    \begin{enumerate}[(i)]
        \item \emph{stabilization by $\R^2\simeq\C$}:
        \[
            k \mapsto k+2,\qquad
            J \mapsto J\oplus J_{\C},\qquad
            \psi \mapsto \psi\oplus\psi_0,
        \]
        where $J_{\C}$ is the standard complex structure on $\R^2$
        and $\psi_0(x,y)=x+iy$;

        \item \emph{stabilization by complex bundles}:
        for any complex vector orbi-bundle $\F_{\C}\to\X$,
        \[
            \W_{\C} \mapsto \W_{\C}\oplus\F_{\C},\qquad
            J \mapsto J\oplus J_{\F},\qquad
            \psi \mapsto \psi\oplus\id_{\F_{\C}}.
        \]
    \end{enumerate}

    Canonical permutation isomorphisms of direct summands are understood.
\end{enumerate}
\end{definition}
Similarly, a stably almost complex structure on $\V$
induces a stably almost complex structure on 
the usual vector bundle $\Ho\V\to \Ho(\X)$.
\subsection{Dualities and proper pushforward maps}
\subsubsection{Thom isomorphism}
For a vector bundle $\xi\colon V\to B$ of real rank $k$, consider the morphism
$\Th(V)\to \Th(V)\wedge \Sigma^\infty_+B$ which
is identity on $\Th(V)$ and projection to $\Sigma^\infty_+B$.
If $V$ admits an $E$-orientation, then
we can construct
\[
    \Th(V)\wedge E\to 
    \Th(V)\wedge\Sigma^\infty_+B\wedge E\xrightarrow{u_\xi\wedge \id_{\Sigma^\infty_+B}}
    \Sigma^{k} E\wedge \Sigma^\infty_+B\wedge E
    \simeq \Sigma^{k}\Sigma^\infty_+B\wedge (E\wedge E)
    \xrightarrow{\id\wedge \mu}\Sigma^{k}_+ B\wedge E.
\]
In particular, if $\xi$ is a complex vector bundle and $E$ is
complex-oriented, then we have the preceding
map. This can be generalized to the case of a virtual
vector bundle with a stably complex structure.
Let $B$ be a space and let $V,W\to B$ be two real vector
bundles. The virtual vector bundle is $\xi=V-W$.
A stably complex structure on $\xi$ consists of a complex
vector bundle $U\to B$ and constant bundles
$\underline{\R}^k$, $\underline{\C}^l$with an isomorphism
\[
    V\oplus \underline{\C}^l \oplus \underline{\R}^k\simeq
    W\oplus U
\]
This is subject to an equivalence relation. We define
the Thom spectrum of $\xi=V-W$ to be
\[
    \Th(\xi)=B^\xi:=\Sigma^{-2l-k}\Sigma^\infty\Th(U)
\]
We can get $\Th(\xi)\wedge E\to \Sigma^{\rk V-\rk W}_+B\wedge E$ by taking desuspension.
\begin{theorem}[Thom isomorphism, {\cite[Theorem-Definition V.1.3]{Rugyak1998}}]
    Assume $\xi=V-W$ is a stably complex virtual bundle,
    $E$ is complex-oriented and homotopy commutative and $k=\rk V-\rk W$. 
    The
    map $\Th(\xi)\wedge E\to \Sigma^{k}_+ B\wedge E$
    is an equivalence. In particular, we have equivalence
    of homology and cohomology
    \[
        \widetilde{H}_{*+k}(\Th(\xi),E)\simeq H_*(B,E).
    \]
    \[
        \widetilde{H}^{*+k}(\Th(\xi),E)\simeq H^*(B,E).
    \]  
\end{theorem}
Notice this theorem is for topological spaces. If
it is a virtual vector bundle over an orbifold $\X$,
then we apply Thom construction and Thom isomorphisms
for the corresponding virtual vector bundle over $\Ho(\X)$.
\subsubsection{Pontryagin-Thom collapse for orbifolds}
Let's recall Pontryagin-Thom collapse for manifolds,
also called the Gysin map, wrong-way map, or transfer map.
Let $f\colon X\to Y$ be a smooth map of closed manifolds,
there is a vector bundle $\pi\colon E\to Y$ and an embedding
$i\colon X\to E$ such that $f=\pi\circ i$. Let $\nu=
N_{X/E}$ be the normal bundle of $X$ in $E$. Then we have
\[
    Y^E=Th(E)\to \Th(\nu)
\]
If $E$ is a trivial vector bundle (which is possible if
we add a suitable vector bundle $F\to Y$), we can desuspend
and get a morphism in $\Sp$, 
\[
    \mathrm{PT}_f\colon \Sigma^\infty_+Y\to X^{f^*TY-TX}
\] because there is an equivalence of virtual vector bundles $\nu-f^*E=f^*TY-TX$. The homotopy class of this morphism is independent of the choice of $E$ and
embedding $i$.

This construction can be generalized to a representable
morphism of closed orbifolds \cite[Theorem 3.3(1)]{Ebert-Giansiracusa2011}. A morphism $f\colon \X\to \Y$
is called representable if for every $x\in \X, y=f(x)\in\Y$,
the induced homomorphism of stabilizer groups
is injective. The proof is also by
constructing a vector orbi-bundle $\E\to \Y$
such that $\X$ embeds\footnote{
A morphism of orbifolds is an embedding iff it is representable, injective on tangent spaces and is
injective on underlying spaces.
} in $\E$; see \cite[Corollary 3.24]{Pardon2023-Orbifold_bordism}.
\begin{lemma}
\label{lem:Pontryagin-Thom for orbifolds}
    Let $f\colon \X\to \Y$ be a representable\footnote{To see why representability is essential in this lemma,
we consider  homomorphisms of finite groups $u\colon G\to H$
inducing $f\colon [\pt/G]\to [\pt/H]$. $f$ being representable
is equivalent to $u$ being injective. In this case,
the homotopy fiber of $BG\to BH$ is $H/G$, which
is a finite set. We can find a finite dimensional vector bundle $E\to BH$ with an embedding $BG\to E$.
Pontryagin-Thom collapse applies here because finite points
in Euclidean spaces have normal bundles.

If $G\to H$ is surjective with nontrivial kernel $K$, the homotopy fiber
of $BG\to BH$ is $BK$. In the injective case we could embed the finite set
$H/G$ in Euclidean space. Here we do not see how to embed $BK$ in the same way.} morphism of closed orbifolds.
    Then there is a morphism
    \[
        \mathrm{PT}_f\colon
        \Sigma_+^\infty\Ho(\Y) \to \Ho(\X)^{\Ho(f^*T\Y-T\X)}
    \]
    When $\X$ and $\Y$ are manifolds, this recovers the classical
    Pontryagin--Thom collapse.
\end{lemma}
\begin{definition}
\label{def:PT pushforwards of orbifolds}
    Let $f\colon \X\to \Y$ be a representable morphism of 
    closed \emph{stably almost complex} orbifolds,
    and let $E$ be any complex oriented theory, then
    there is a map
    \[
        f_!^{\operatorname{PT}}\colon 
        H^*(\Ho(\X),E)\to H^{*+\dim \Y-\dim\X}\left(
        \Ho(\Y),E\right).
    \]
    defined by $\operatorname{PT}_f$ and Thom isomorphism
    $\widetilde{H}^*\left(\Ho(\X)^{\Ho(f^*T\Y-T\X)},E\right)
    \simeq H^{*+\dim \X-\dim \Y}\left(\Ho(\X),E\right)$.
\end{definition}
We call this map the ``Pontryagin-Thom pushforward''. It is now
only defined for representable $f$, but later for $E$ being a
Morava $K$-theory or HKR theory, we will see this condition can be dropped.
We use the notation $f_!^{\operatorname{PT}}$ to distinguish this pushforward from
$f_!$ in \cref{def:orbifold proper pushforward}.
\subsubsection{Ambidexterity and Cheng's theorem}
\label{sec:ambidexterity}
For a closed orbifold $\X$, Cheng constructed a morphism
of spectra \cite{Cheng-Deligne_Mumford_Stacks2011}.
\begin{equation}
\label{eq:duality lambda}
    \lambda_\X\colon \Sigma^\infty\Ho(\X)_+\to \inthom \left(
    \Ho(\X)^{-\Ho T\X},\SS
    \right).
\end{equation}

In the closed case, the construction is as follows: consider the diagonal
morphism $\Delta_\X\colon \X\to \X\times \X$. $\Delta_\X$
is representable because at $x\in \X$, the homomorphism of
stabilizer groups is also a diagonal
$\Stab_x\to \Stab_x\times\Stab_x$, and this is injective.
The normal bundle of $\X\to \X\times \X$ is $T\X$.
By a slight variant of \cref{lem:Pontryagin-Thom for orbifolds}, we get
\[
    \Sigma^\infty\Ho(\X)_+\wedge \Ho(\X)^{-\Ho T\X}
    \to \Sigma^\infty\Ho(X)_+
\]
$\lambda_\X$ is the composition of this map with $\Ho(\X)_+\to S^0$. 

\setlength{\parfillskip}{0pt plus 1fil}
Cheng proved the following duality result
generalizing ambidexterity of $BG$ for finite groups 
(\cite{Greenlees1996TheTS}).
For possibly noncompact orbifolds, write $H_*^{\mathrm{lf}}(\X;E)$ for
locally finite orbifold $E$-homology.\footnote{
Henceforth, for $E=K(n)$ or an HKR theory, we use
$H_*(\X;E):=\pi_* L_{K(n)}\!\left(E\wedge\Sigma^\infty_+\Ho(\X)\right)$.
}
\begin{theorem}[{\cite[Theorem 1.0.1]{Cheng-Deligne_Mumford_Stacks2011}\cite[Corollary 1.4]{cheng2013poincaredualitymoravaktheory}\cite[Theorem 5.16]{abouzaid2021complexcobordismhamiltonianloops}}]
\label{thm:duality for orbifolds}
    Let $E$ be $K(n)$ or an HKR theory, and let $\X=[M/G]$ be a smooth
    global quotient orbifold, with $M$ closed, open, or compact with boundary.
    \begin{enumerate}[(1)]
    \item The map $\lambda_\X$ induces an isomorphism
    \[
        H_*^{\mathrm{lf}}(\X;E)
        \simeq H^{-*}\left(\Ho(\X)^{-\Ho T\X},E\right).
    \]
    \item If $\X$ is equipped with a stably almost complex structure,
    the Thom isomorphism gives
    \[
        \lambda_\X\colon H_*^{\mathrm{lf}}(\X;E)
        \xrightarrow{\;\simeq\;}
        H^{\dim \X-*}\left(\Ho(\X),E\right).
    \]
    \end{enumerate}
\end{theorem}
\begin{remark}
\label{rmk:definition by slant product with diagonal class}
    There is an equivalent way to define $\lambda_\X$ for closed stably
    almost complex $\X$. Because
    $\Delta\colon \X\to \X\times \X$ is representable,
    $\Delta_!^{\operatorname{PT}}$ is defined. Denote
    \[
        \delta_\X:=\Delta_!^{\operatorname{PT}}(1)\in H^{\dim \X}
        (\Ho(\X\times\X),E).
    \]
    Then we can define $\lambda_\X\colon H_*^{\mathrm{lf}}(\Ho(\X),E)\to H^{\dim\X-*}
    (\Ho(\X),E)$ by taking slant product with $\delta_\X$.
\end{remark}
    $\lambda_\X$ generalizes the following
    special cases:
    \begin{enumerate}[(1)]
        \item If $G=1$ and $\X=M$ is a closed manifold,
        the map $\lambda_M\colon \Sigma^\infty M_+
        \to \Hom\left(M^{-TM},\SS\right)$ coincides
        with Atiyah duality.
        \item If $M=\pt$ and $G$ is a finite group,
        in which case $\dim\X=0$. We have an isomorphism
        $H^*(BG,E)\simeq H^{-*}(BG,E)$. This is
        the same isomorphism as \cite[Corollary 1.2]{Greenlees1996TheTS}.
    \end{enumerate}

Let $E$ be $K(n)$ or an HKR theory, and let $\X,\Y$ be global quotient
stably almost complex orbifolds. Let $f\colon\X\to\Y$ be a smooth
morphism, not necessarily representable. Assume either that $\X,\Y$
have no boundary and $f$ is proper, or that they are compact with boundary
and $f(\partial\X)\subseteq\partial\Y$.

The preceding duality gives
\[\begin{tikzcd}
	{H^*(\Ho(\X),E)} & {H_{\dim \X-*}^{\mathrm{lf}}(\X;E)} & {H_{\dim\X-*}^{\mathrm{lf}}(\Y;E)} & {H^{*-\dim\X+\dim\Y}(\Ho(\Y),E)}
	\arrow["{\lambda^{-1}_\X}", from=1-1, to=1-2]
	\arrow["{f_*}", from=1-2, to=1-3]
	\arrow["{\lambda_\Y}", from=1-3, to=1-4]
\end{tikzcd}\]
Here $f_*$ is induced by properness in the boundaryless case and by the
map of pairs in the compact-with-boundary case.
\begin{definition}
\label{def:orbifold proper pushforward}
    Under these assumptions, define
    $f_!:=\lambda_\Y\circ f_*\circ (\lambda_\X)^{-1}$
    and call it the pushforward map in $E$-cohomology.
\end{definition}
This definition coincides with Pontryagin-Thom pushforward maps in
representable cases. A relevant proposition is \cite[Lemma 5.22]{abouzaid2021complexcobordismhamiltonianloops}. 
\begin{lemma}
\label{lem:coincidence with Thom isomorphisms}
   If $f\colon \X\to \Y$ is a representable morphism of closed
   stably almost complex orbifolds, then the two pushforwards 
   $f_!^{\operatorname{PT}}$(\cref{def:PT pushforwards of orbifolds}) and $f_!$(\cref{def:orbifold proper pushforward}) coincide.
\end{lemma}
\begin{proof}
   Let $\delta_\X$ and $\delta_\Y$ be as defined in \cref{rmk:definition by slant product with diagonal class}. Consider the commutative
   diagram
\[\begin{tikzcd}
	\X & {\X\times \Y} \\
	\Y & {\Y\times \Y}
	\arrow["{\Gamma_f}", from=1-1, to=1-2]
	\arrow["f"', from=1-1, to=2-1]
	\arrow["{f\times \id_\Y}", from=1-2, to=2-2]
	\arrow["{\Delta_\Y}", from=2-1, to=2-2]
\end{tikzcd}\]
We use $a \backslash\delta_\X$ to denote the slant product
for $a\in H_*(\Ho(\X),E)$. 
The transverse base-change of representable PT pushforward maps gives
\[
    (f\times \id_\Y)^*\delta_\Y=(\Gamma_f)_!^{\operatorname{PT}}(1)
    =(\id_\X\times f)_!^{\operatorname{PT}}\delta_\X.
\]
For $\alpha\in H_*(\Ho(\X),E)$ we have
\[
    \alpha\backslash (f\times \id_\Y)^*\delta_\Y
    =f_*\alpha\backslash \delta_\Y
    =\lambda_\Y(f_* \alpha)
\]
and, by compatibility of the slant product with PT pushforward in the second factor,
\[
    \alpha\backslash (\id_\X\times f)_!^{\operatorname{PT}}\delta_\X
    =f_!^{\operatorname{PT}}(\alpha\backslash \delta_\X)
    =f_!^{\operatorname{PT}}(\lambda_\X(\alpha))
\]
So we have $\lambda_\Y(f_* \alpha)=f_!^{\operatorname{PT}}(\lambda_\X(\alpha))$. Since $\lambda_\X$ is an isomorphism, this gives
\[
    f_!^{\operatorname{PT}}=\lambda_\Y\circ f_*\circ\lambda_\X^{-1}=f_!
\]
by \cref{def:orbifold proper pushforward}, for representable $f$.
\end{proof}
\begin{proposition}
Let $\X,\Y,\Zcal$ be closed stably almost complex orbifolds.
\label{prop:properties of f_!}
    \begin{enumerate}[(1)]
        \item If $f\colon \X\to\Y$ and $g\colon \Y\to \mathcal Z$
        are proper, then $(g\circ f)_!=g_!\circ f_!$
        \item (\cite[Proposition 8.6]{strickland2000knlocaldualityfinitegroups} for groupoids version) Let $f\colon \X\to \Zcal$ and $u\colon \Y\to \Zcal$
        be morphisms of orbifolds that are transverse\footnote{If $x\in \X,y\in \Y,z\in \mathcal Z$
        such that $f(x)=u(y)=z$, then the image of $df\colon T_x\X\to T_z\mathcal Z$ and $dg\colon \colon T_y\Y\to T_z\mathcal Z$ spans $T_z\mathcal Z$.}.
        Then $\mathcal{W}=\X\times_{\Zcal}\Y$
        is also an orbifold. Let $g\colon \mathcal W\to \Y$ and $v\colon \mathcal W\to \X$ be the projection
        maps.
        \[\begin{tikzcd}
	{\mathcal W} & \Y \\
	\X & {\mathcal Z}
	\arrow["g", from=1-1, to=1-2]
	\arrow["v", from=1-1, to=2-1]
	\arrow["u", from=1-2, to=2-2]
	\arrow["f", from=2-1, to=2-2]
\end{tikzcd}\]
We have
\[
    g_!\circ v^*=u^*\circ f_!
\]
    \end{enumerate}
\end{proposition}
\begin{proof}
    \begin{enumerate}[(1)]
        \item It follows from definition
        directly.
        \item Let $\Gamma_f\colon \X\to \X\times \Zcal$ be the graph
        of $f$, then $\Gamma_f$ is representable and
        $f$ is the composition $\X\xrightarrow{\Gamma_f}\X\times \Zcal
        \xrightarrow{\pi_\Zcal} \Zcal$.
\[\begin{tikzcd}
	\W & {\X\times \Y} & \Y \\
	\X & {\X\times \Zcal} & {\Zcal}
	\arrow["{j=(v,g)}", from=1-1, to=1-2]
	\arrow["v", from=1-1, to=2-1]
	\arrow["\pi_\Y", from=1-2, to=1-3]
	\arrow["{\id_\X\times u}", from=1-2, to=2-2]
	\arrow["u", from=1-3, to=2-3]
	\arrow["{\Gamma_f}", from=2-1, to=2-2]
	\arrow["{\pi_\Zcal}", from=2-2, to=2-3]
\end{tikzcd}\]
Then it suffices to prove (2) for each square. We only present a sketch
of the proof.
\begin{enumerate}[(a)]
    \item $(\id_\X\times u)^*\circ (\Gamma_f)_!=j_!\circ v^*$.\hfill

    Since $\Gamma_f$ is representable, by \cite[Corollary 3.24]{Pardon2023-Orbifold_bordism} there is a complex vector orbi-bundle
    $\pi\colon \mathcal{V}\to \X\times \Zcal$ and an embedding
    $i\colon \X\to \V$ such that $\Gamma_f=\pi\circ i$. By \cref{lem:coincidence with Thom isomorphisms}, $\pi_!,i_!$ both
    coincide with classical PT map.
    Since this is already representable transverse base-change, we get
    $(\id_\X\times u)^*\circ (\Gamma_f)_!=j_!\circ v^*$.
    
    \item $u^*\circ (\pi_Z)_!=(\pi_Y)_!\circ (\id_\X\times u)^*$.\hfill

    Let $[\X]_E=\lambda_X^{-1}(1)$ be the
    dual of $1$, then on $\X\times \mathcal{Y}$ or on $\Y\times \Zcal$, the $\pi_!$ ($\pi=\pi_\Y$ or $\pi_\Zcal$)
    is just the slant product with $[\X]$. Because the
    slant product is functorial in pulling back cohomology classes,
    we get $u^*\circ (\pi_Z)_!=(\pi_Y)_!\circ (\id_\X\times u)^*$.
\end{enumerate}
    \end{enumerate}
\end{proof}
\begin{corollary}[Projection formula]
\label{cor:projection formula}
    Let $f\colon \X\to \Y$ be proper, $a\in H^*(\X,E), b\in H^*(\Y,E)$.
    Then $f_!(a\cup f^* b)=f_!(a)\cup b$
\end{corollary}
\begin{proof}
    Consider the following transverse pullback diagram
\[\begin{tikzcd}
	\X & {\X\times \Y} & {} \\
	\Y & {\Y\times\Y}
	\arrow["{(\id_\X,f)}", from=1-1, to=1-2]
	\arrow["f"', from=1-1, to=2-1]
	\arrow["{f\times\id_\Y}", from=1-2, to=2-2]
	\arrow["{\Delta_\Y}"', from=2-1, to=2-2]
\end{tikzcd}\]
We choose class $a\times b\in H^*(\X\times\Y,E)$, then
by \cref{prop:properties of f_!} (2) we have
\[
    f_!\circ (\id_\X,f)^*(a\times b)=\Delta^*\circ (f\times \id_\Y)_!
    (a\times b).
\]
The LHS is $f_!(a\cup f^*(b))$ and the RHS
is $f_!(a)\cup b$.
\end{proof}

\subsection{HKR theory and pushforward formulae}
In \cref{sec:ambidexterity}, we introduced pushforward maps
of proper morphisms between orbifolds for an HKR theory $E$.
To actually compute pushforward maps, we employ
the character map in \cite{Generalized-Group-Characters-and-Complex-Oriented-Cohomology-Theories} which identifies
rationalization of $E$-cohomology of (the homotopy type of) an orbifold
with $L(E)$-cohomology of its $n$-fold $p$-inertia stack, where
$L(E)$ is a ring spectrum constructed from $E$ such that
$L(E)^*$ is a $\Q$-algebra. 

After that, we will establish
a Riemann--Roch type formula for orbifolds relating the $E$-pushforward
with the $L(E)$-pushforward, which is essentially\footnote{
This isn't equality on the nose. See \cref{rmk:warning on formal group law of L(E)}
} the
usual pushforward maps in ordinary cohomology.
This will be compared with Kawasaki's Riemann--Roch formula
computing holomorphic Euler characteristics of equivariant
vector bundles \cite{ojm/1200771835}. For transfer formulas for HKR characters and
transchromatic characters\footnote{
It relates
theories of height $n$ and $t$. The original HKR theory is
$t=0$ case.
}, see \cite[Theorem~D]{Generalized-Group-Characters-and-Complex-Oriented-Cohomology-Theories}
and \cite[Theorems~3.7--3.8]{BenMoshe2025Transfers}, respectively.
These formulas motivate our Riemann--Roch formula for orbifolds.

Finally, we define an ad hoc equivariant pushforward map
and develop the corresponding equivariant HKR theory.
\subsubsection{HKR characters: definition} 

\begin{definition}[{\cite{Generalized-Group-Characters-and-Complex-Oriented-Cohomology-Theories}}]
\label{def:p-inertia stack construction}
    \begin{enumerate}[(1)]
        \item Let $\X$ be a compact orbifold (possibly with boundary),
        we define the \emph{$p$-inertia stack} to be
        \[
            \I_p \X=\colim_{m\to \infty}\inthom([\pt/C_{p^m}],\X)
        \]
        Here $[*/C_{p^m}]$ denotes the quotient orbifold of
        $C_{p^m}\acts \pt$, and $\inthom$ is the internal Hom stack.
        The colimit is along the diagram
        $\cdots\to C_{p^3}\to C_{p^2}\to C_p$ mapping $1\in C_{p^m}$
        to $1\in C_{p^{m-1}}$.
        Equivalently, $\I_p\X=\{x\in \X,g\in \Stab_{x}| \text{the order of }g\text{ is a power of }p\}$. 
        \item The \emph{$n$-fold $p$-inertia stack} of $\X$ is
        defined to be $(\I_p)^n\X$,
        i.e. performing $\I_p$ $n$ times to $\X$. Equivalently,
        \[
            (\I_p)^n\X=\colim_{m\to\infty}\inthom
            \left(
            [\pt/(C_{p^m})^n],\X
            \right)
        \]
    \end{enumerate}
\end{definition}
\begin{remark}
\begin{enumerate}[(1)]
\label{rmk:definition of p-inertia stack}
    \item The original definition in \cite{Generalized-Group-Characters-and-Complex-Oriented-Cohomology-Theories} is
    for a finite group $G$ and finite $G$-CW complex $X$.
    Define a $G$-equivariant space
        \[
            \Fix_{n,p}X=\bigsqcup_{\alpha\colon \Z_p^n\to G}X^{\im\alpha}
        \]
        Here $\alpha$ is a continuous homomorphism, and this means
        $\alpha$ factors through some $(\Z/p^m)^n$ because $G$ is finite; in particular, there are only finitely many $\alpha$.
        Such $\alpha$ are called \emph{HKR sectors}.
        $X^{\im\alpha}$ is the subspace of $X$ fixed by $\im\alpha\subset G$. $\Fix_{n,p}X$ is equipped with a $G$-action
        sending $x\in X^{\im\alpha}$ to $gx\in X^{\im (g\cdot\alpha\cdot g^{-1})}$.
    If $\X$ is   ${[M/G]}$ where $G$ is finite and $M$ is a compact manifold,
    then there is an isomorphism of orbifolds
    \[
        \I_p^n[M/G]\simeq [\Fix_{n,p}M/G].
    \]
    \item Because $\X$ is compact, there are only finitely many isotopy
    types, so elements in them have bounded orders. Therefore,
    the $\inthom
            \left(
            [\pt/(C_{p^m})^n],\X
            \right)$ in definition becomes constant after a finite $m$.
    Later on, we will implicitly take such an $m$, assuming that
    $m$ is large enough to correctly define $\I_p^n\X$ as $\inthom
    \left([\pt/C_{p^m}]^n,\X\right)$.
\end{enumerate}
\end{remark}
Now let's define $L(E)$. 
Fix an HKR theory $E$ of height $n$, and denote $\Lambda_r=(C_{p^r})^n$.
There is an inverse diagram of surjective homomorphisms mapping $1$
to $1$
\[
    \cdots\to\Lambda_r\to \Lambda_{r-1}\to\cdots\to \Lambda_1
\]
The $E$-cohomologies of $B\Lambda_r$ thus form a directed diagram,
and define
\[
    L'=\colim_{r\to\infty}H^*\left(B\Lambda_r,E\right)
\]
Recall that $H^*(BS^1,E)\simeq E^*\llbracket x\rrbracket$ (\cref{lem:E-cohomology of BU(n)}). For every nontrivial group homomorphism
$f\colon \Lambda_r\to S^1$, we get $f^*x\in H^2(B\Lambda_r,E)$, and
$f^*x$ can be viewed as an element of $L'$ via natural inclusion.
\begin{definition}[{\cite{Generalized-Group-Characters-and-Complex-Oriented-Cohomology-Theories}}]
    Let $S$ be the multiplicative set generated by $f^*x$ for all 
    $r$ and nontrivial homomorphisms $f\colon\Lambda_r\to S^1$.
    Define the \textit{graded ring} $L(E)^*$ to be
    \[
        L(E)^*=S^{-1}L'
    \]
    We denote by $L(E)$ the Eilenberg-MacLane spectrum
    of the graded ring $L(E)^*$, and multiplication of $L(E)^*$
    makes $L(E)$ a commutative homotopy
    ring spectrum.
\end{definition}
It turns out $L(E)^*$ is a $\Q$-algebra, so its height is $0$. 
\begin{remark}
\label{rmk:def of L(E) coincides with localized cohomology}
    For each $r$, consider the nontrivial characters $\chi\colon\Lambda_r\to S^1$, each inducing a complex line bundle
    over $B\Lambda_r$. Let $S_r\subset H^*(B\Lambda_r,E)$ be the multiplicative set generated by the Euler characters $e_E(\chi)$
    for nontrivial $\chi$. Then we have an equivalent definition of $L(E)$
    \[
        L(E)^*=\colim_{r\to \infty}S_r^{-1}H^*(B\Lambda_r,E).
    \]
    The localized cohomology $S_r^{-1}H^*(B\Lambda_r,E)$
    appears again in \cref{section: equivariant localization}.
\end{remark}
Note that there is a natural morphism of ring spectra $E\to L(E)$
whose effect on cohomology is $H^*(X,E)\xrightarrow{a\mapsto a\otimes 1} H^*(X,E)\otimes_{E^*}L(E)^*\to H^*(X,L(E))$.
$L(E)$ carries a natural complex orientation inherited from
$E$.
\begin{remark}
\label{rmk:warning on formal group law of L(E)}
Although $L(E)$ is a rational theory, its formal group law
isn't simply given by $x+y$. As a consequence, the pushforward of $L(E)$-cohomology isn't the $\Q$-cohomology pushforward. They are related by an integration formula
which is
in the same spirit as computing $K$-theoretic pushforwards
using Grothendieck--Riemann--Roch formula.
\end{remark}
We record the following property of $L(E)$. Suppose $r$ is sufficiently
large. Note that there are exactly $p^{n a}$ homomorphisms
$\chi\colon \Lambda_r\to C_{p^a}$.
\begin{lemma}
\label{lem:e(chi) becomes exactly the roots of [p^a]x=0}
    The images of elements $e(\chi)\in H^*(B\Lambda_r,E)\to L(E)^*$
    are pairwise distinct and are exactly the roots of the equation\footnote{
    Take the Weierstrass factorization over $E^*$,
    $[p^a]_Ex=U_a(x)W_a(x)$, where
    $U_a(x)\in E^*\llbracket x\rrbracket^\times$ and
    $W_a(x)\in E^*[x]$ is the monic Weierstrass polynomial.
    By roots we mean roots of $W_a(x)$.
    Quotients by $[p^a]_Ex$ after coefficient extension are understood
    as base changes of the finite free algebra
    $E^*\llbracket x\rrbracket/([p^a]_Ex)\cong E^*[x]/(W_a(x))$.
    } $[p^a]x=0$.
\end{lemma}
The HKR character map's construction is as follows. Choose a sufficiently large
$m$ such that $\I_p^n\X=\inthom
            \left(
            [\pt/\Lambda_m],\X
            \right)$.
The universal map is $\I_p^n\X\times [\pt/\Lambda_m]
\to \X$. Taking (Borel) $E$-cohomology, we get
\[
\begin{aligned}
    H^*(\X,E)&\to H^*\left(
    \I_p^n\X\times [\pt/\Lambda_m],E
    \right)\\
    &\simeq H^*(\I_p^n\X,E)\otimes_{E^*}
    H^*(B\Lambda_m,E)\\
    &\to H^*(\I_p^n\X,E)\otimes_{E^*}L(E)^*
    \to H^*(\I_p^n\X,L(E)).
\end{aligned}
\]
The Kunneth formula is true because $H^*(B\Lambda_m,E)$
is a free $E^*$ module.
\begin{definition}[{\cite{Generalized-Group-Characters-and-Complex-Oriented-Cohomology-Theories}}]
        Let $\X$ be a compact orbifold. The \emph{HKR character}
        of $E$ is defined as the composition above:
        $\chi_{n,p}(\X)\colon H^*(\X,E)\to H^*(\I_p^n \X,L(E))$
\end{definition}
\begin{theorem}[{\cite[Theorem C]{Generalized-Group-Characters-and-Complex-Oriented-Cohomology-Theories}}]
\label{thm:HKR isomorphism}
    The character $\chi_{n,p}(\X)$ induces an isomorphism
    \[
       \chi_{n,p}(\X) \colon  L(E)^*\otimes_{E^*} H^*(\X,E)\to H^*\left(
        \I^n_p \X,L(E)
        \right).
    \]
    In fact, it is an isomorphism of graded rings.
\end{theorem}
\begin{remark}
    There is a slight difference between the statement here and
    {\cite[Theorem C]{Generalized-Group-Characters-and-Complex-Oriented-Cohomology-Theories}}, as there it is only for $[X/G]$ where $X$
    is a finite $G$-CW complex. However, we can find an orbi-CW filtration\footnote{This is guaranteed by \cite[Proposition 3.5]{Pardon2023-Orbifold_bordism}.}
    of $\X$ and prove by induction (by Mayer–Vietoris sequences).
\end{remark}
Although Kunneth isomorphism for $E$ is generally not true,
we prove a lemma of Kunneth type after taking tensor product
with $L(E)^*$, using the identification \cref{thm:HKR isomorphism}.
\begin{lemma}
\label{lem:Kunneth formula after base change}
    Let $\X,\Y$ be compact orbifolds. Then there is an isomorphism
    of rings
    \[
        L(E)^*\otimes_{E^*} H^*\left(\X\times \Y, E\right)
        \simeq 
        \left(L(E)^*\otimes_{E^*}H^*\left(\X, E\right) \right)
        \otimes_{L(E)^*}
        \left(L(E)^*\otimes_{E^*}H^*\left(\Y, E\right) \right)
    \]
\end{lemma}
\begin{proof}
    By \cref{thm:HKR isomorphism}, it is equivalent to proving
    \[
        H^*\left(\I^n_p (\X\times \Y),L(E)
        \right)
        \simeq 
        H^*\left(\I^n_p \X,L(E)\right)\otimes_{L(E)^*}
        H^*\left(\I^n_p \Y,L(E)\right).
    \]
    Since $\I ^n_p\X=\Hom([\pt/\Lambda_r], \X)$ for
    sufficiently large $r$, and
    \[
        \Hom\left([\pt/\Lambda_r],\X\times \Y\right)
        \simeq \Hom\left([\pt/\Lambda_r],\X\right)
        \times \Hom\left([\pt/\Lambda_r],\Y\right),
    \]
    we have
    $\I^n_p(X\times \Y)\simeq \I^n_p\X\times \I^n_p\Y$. Then
    the Kunneth isomorphism is true because $L(E)$ is rational.
\end{proof}
\subsubsection{A Riemann--Roch type formula}
Let $\X$ be a closed orbifold with a stably almost complex
structure. Take $m$ sufficiently large such that
$\I_p^n\X=\inthom\left([\pt/\Lambda_m],\X\right)$, and denote
the universal map by 
\[q=q_\X \colon [\pt/\Lambda_m]\times
\I_p^n \X\to \X.\]
$q^* T\X$ is a vector orbi-bundle over $\I_p^n\X\times [\pt/\Lambda_m]$, and it is equivalently
a $\Lambda_m$-equivariant vector orbi-bundle over $\I_p^n\X$,
where $\I_p^n\X$ has the trivial $\Lambda_m$ action.
\begin{definition}
\label{def:normal bundle of alpha fixed locus}
    Define $\nu_\X=(q^*T\X)^\mov$, the vector orbi-bundle of nontrivial irreducible $\Lambda_m$-representations in
    $q^*T\X$, also called the ``moving part''. 
\end{definition}
The trivial representation part/fixed part is actually $T\I_p^n \X$.
$T\I_p^n\X$ and $\nu_\X$ carry natural stably
complex structures such that $q^*T\X=T\I_p^n\X\oplus \nu_\X$ respects the stably
complex structures\footnote{
The stably complex structure on $T\X$ requires two complex
vector orbi-bundles $\V_\C,\W_\C$ and a constant real vector
bundle $\underline{\R}^k$ with an isomorphism
\[
    T\X\oplus \V_\C\oplus\underline{\R}^k \simeq
    \W_\C.
\]
Pulling these bundles by $q^*$, (the pulled back bundle of)
these bundles are also $\Lambda_m$-equivariant vector 
orbi-bundles over $\I_p^n\X$. Since the fixed and moving
parts of a complex vector orbi-bundle are still complex,
we get
\[
T\I_p^n\X\oplus (\V_\C)^\fix\oplus \underline{\R}^k\simeq (\W_\C)^\fix
\]
\[
    \nu_\X\oplus (\V_\C)^\mov\simeq (\W_\C)^\mov.
\]
}. The $E$-cohomological equivariant
Euler class of $T\I_p^n\X$ is in $H^*(I_p^n\X,E)$,
and that of $\nu_\X$ is in $H^*(\I_p^n\X\times [\pt/\Lambda_m],E)$.
\begin{definition}
    Define $e^{L}_{\mathrm{HKR}}(\nu_\X)$\footnote{
    $\nu_\X$ can have different ranks on different connected components of $\I_p^n\X$, so
    rigorously we should define $e^{L}_{\mathrm{HKR}}(\nu_\X|_C)$
    for each connected component $C$.
    }
    to be the image of the equivariant Euler
    class of $\nu_\X$ under the map
    $H^*(\I_p^n\X\times [\pt/\Lambda_m],E)
    \to H^*(\I_p^n\X,L(E))$.
\end{definition}
\begin{lemma}
    $e^{L}_{\mathrm{HKR}}(\nu_\X)$ is
    invertible in $H^*(\I_p^n\X,L(E))$.
\end{lemma}
\begin{proof}
    We appeal to \cref{prop:invertibility of Euler class&isomorphism after S_A localization}, where we take $A=\Lambda_m$.
\end{proof}

Assume $\X,\Y$ are both closed orbifolds with
stably almost complex structures, and $f\colon
\X\to \Y$ is a smooth map. Naturality of
$\I_p^n$ induces a map $\I_p^n f\colon
\I_p^n\X\to \I_p^n\Y$.
We have two
pushforward maps, one for $E$-cohomology
(height $n$)
and one for $L(E)$-cohomology (height $0$).
The goal of this section is the following
compatibility diagram, which relates the $E$ and $L(E)$ pushforward
maps.
\begin{theorem}\label{thm:a Riemann--Roch}
    There's a commutative diagram:
    \[
\begin{tikzcd}[column sep=9em, row sep=large]
	{H^*\left(\X,E\right)}
	&
    {H^*(\I_p^n\X,L(E))}
	\\
	{H^*\left(\Y,E\right)}
	&
    {H^*(\I_p^n\Y,L(E))}
	\arrow[
        "{{e_{\mathrm{HKR}}^L(\nu_\X)}^{-1}\cdot\chi_{n,p}(\X)}",
		from=1-1,
		to=1-2
	]
	\arrow["{f_!}"', from=1-1, to=2-1]
    \arrow["{(\I_p^n f)_!}", from=1-2, to=2-2]
	\arrow[
        "{{e_{\mathrm{HKR}}^L(\nu_\Y)}^{-1}\cdot\chi_{n,p}(\Y)}"',
		from=2-1,
		to=2-2
	]
\end{tikzcd}
\]
The horizontal maps are $a\in H^*(\X,E)\mapsto
{{e_{\mathrm{HKR}}^L(\nu_\X)}^{-1}\cdot\chi_{n,p}(\X)}(a)$.
\end{theorem}
We use the following strategy to prove \cref{thm:a Riemann--Roch}. We prove the representable case in \cref{lem:RR representable case} and the $\Y=\pt$ case in \cref{lem:RR Y is a point case}, and
then combine them to prove \cref{thm:a Riemann--Roch}.

The following auxiliary lemma compares the $\I_p^n$ and $\chi$
of $W$ and $\Zcal$ for a vector orbi-bundle $\pi\colon W\to \Zcal$. 
Similar to \cref{def:normal bundle of alpha fixed locus},
under the pullback 
$q_\Zcal\colon \I_p^n\Zcal\times[\pt/\Lambda_r]\to \Zcal$,
there is a splitting $q_\Zcal W= W^\fix\oplus W^\mov$\footnote{They
are all bundles over $\I_p^n \Zcal\times [\pt/\Lambda_r]$, or equivalently, $\Lambda_r$-equivariant bundles over $\I_p^n \Zcal$.}.
It is easy to see that $\I_p^n W$ is a vector orbi-bundle over $\I_p^n \Zcal$ with fiber $W^\fix$.
\begin{lemma}
    We have
        \[
            \chi_{W}(u_E(W))=u_L(W^\fix)\cdot e_{\operatorname{HKR}}^L(W^\mov).
        \]
        and
        \begin{equation}
            \label{eq:comparing Z and W}
            \chi_{W}\left(\Th_W(a)\right)
            =\Th_{W^\fix}\left(\chi_\Zcal(a)\cdot e_{\operatorname{HKR}}^L(W^\mov)\right).
        \end{equation}
        Here $u_E(W)\in \widetilde{H}^*(\Th(W),E)$ and
        $u_L(W^\fix)\in \widetilde{H}^*(\Th(W^\fix),L(E))$\footnote{
        The relevant vector bundle is $W^\fix\to \I_p^n \Zcal\times [\pt/\Lambda_r]$.
        } are the Thom classes in $E$ and $L(E)$ cohomology.
\end{lemma}
\begin{proof}
    To compute $\chi_W$ we consider
    \[  
        q_W\colon \I_p^n W\times [\pt/\Lambda_r]\to
        W
    \]
    where $\I_p^n W=W^\fix$. By splitting of $W$ into $W^\fix$, $W^\mov$
    after pullback, we have
    \[
        q^*_W  (u_E(W))=u_E(W^\fix)e_E(W^\mov)\in \widetilde{H}^*
        \left(\Th(W^\fix)\times [\pt/\Lambda_r],E\right).
    \]
    Then after $\widetilde{H}^*
        \left(\Th(W^\fix)\times [\pt/\Lambda_r],E\right)
        \to \widetilde{H}^*
        \left(\Th(W^\fix),L(E)\right)$, it becomes
        $u_L(W^\fix)e_{\operatorname{HKR}}^L(W^\mov)$. This is the first
        formula

        To prove \eqref{eq:comparing Z and W}, we set
        $\pi\colon W\to \Zcal$ and $\pi^\fix\colon W^\fix\to \I_p^n \Zcal$.
        Then,
        \begin{align*}
            \chi_W(\Th_W(a))&=\chi_W(u_E(W)\pi^* a)\\
            &=u_L(W^\fix)\cdot e_{\operatorname{HKR}}^L(W^\mov)\cdot 
            (\pi^\fix)^*\chi_\Zcal(a)\\
            &=\Th_{W^\fix}\left(\chi_\Zcal(a)\cdot e_{\operatorname{HKR}}^L(W^\mov)\right).
        \end{align*}    
\end{proof}
For simplicity, we use $\chi_\X,\chi_\Y$ to denote $\chi_{n,p}$
of $\X$ and $\Y$, and use $e_\X\in H^*\left(\I_p^n \X,L(E)\right)$ to denote $e_{\operatorname{HKR}}^L(\nu_\X)$.
\begin{lemma}
\label{lem:RR representable case}
    \cref{thm:a Riemann--Roch} is true if $f$ is representable. 
    That is, for representable $f\colon \X\to \Y$, we have
    \[
        \frac{\chi_\Y\left(f_!(a)\right)}{e_\Y}
        =(\I_p^nf)_!\left(
        \frac{\chi_\X(a)}{e_\X}
        \right).
    \]
\end{lemma}
\begin{proof}
    By \cite[Corollary 3.24]{Pardon2023-Orbifold_bordism}, we can find
    a vector orbi-bundle $V\to \Y$ and an embedding
    $\X\to V$ such that $f$ is the composition $\X\to V\to \Y$.
    Let $N$ be the normal bundle of $\X\subset V$. For
    simplicity we write $e_N=e_{\operatorname{HKR}}^L(N^\mov)$
    and $e_V=e_{\operatorname{HKR}}^L(V^\mov)$.

    Applying \eqref{eq:comparing Z and W} to the vector orbi-bundles
    $N\to \X$ and $V\to \Y$, we get
    \begin{equation}
        \label{eq:Thom N to X}
         \chi_N\left(\Th_N(a)\right)=\Th_{N^\fix}\left(
        \chi_\X(a)e_N
        \right)
    \end{equation}
    and
    \[
        \chi_V\left(\Th_V(b)\right)=\Th_{V^\fix}\left(
        \chi_\Y(b)e_V
        \right).
    \]
    In the second case, since $\Th_V$ is an isomorphism, let
    $z=\Th_V(b)$, then
    \begin{equation}
        \label{eq:inverse of Thom V to Y}
        \chi_\Y\circ (\Th_V)^{-1}(z)=e_V^{-1}\cdot (\Th_{V^\fix})^{-1}\circ\chi_V(z)
    \end{equation}
    Let $j_*\colon \widetilde{H}^*\left(\Th(N),E\right)\to \widetilde{H}^*\left(\Th(V),E\right)$ be the canonical map
    induced by inclusion. We have
    \begin{align*}
        \chi_\Y(f_!(a))&=\chi_\Y\left(
        (\Th_V)^{-1}\circ j_*\circ \Th_N (a)
        \right)\\
        &\overset{\eqref{eq:inverse of Thom V to Y}}{=}
        e_V^{-1}\cdot (\Th_{V^\fix})^{-1}\circ \chi_V\circ j_*\circ \Th_N(a)\\
        &=e_V^{-1}\cdot (\Th_{V^\fix})^{-1}\circ (\I_p^n j)_*\circ
        \chi_N\circ \Th_N(a)\\
        &\overset{\eqref{eq:Thom N to X}}{=}
        e_V^{-1}\cdot (\Th_{V^\fix})^{-1}\circ (\I_p^n j)_*\circ
        \Th_{N^\fix}\left( \chi_\X(a)e_N\right)
    \end{align*}
    We notice that $(\Th_{V^\fix})^{-1}\circ (\I_p^n j)_*\circ
        \Th_{N^\fix}$ is exactly the pushforward $(\I_p^n f)_!$ associated to
        \[
            \I_p^n f\colon \I_p^n\X\longrightarrow\I_p^n\Y.
        \]
        Thus, we get
        \[
            \chi_\Y(f_!(a))=e_V^{-1}(\I_p^nf)_!\left(
            \chi_X(a)e_N
            \right)
        \]
        By projection formula (\cref{cor:projection formula}) and $N^\mov-f^* V^\mov=f^*T\Y^\mov-T\X^\mov$,
         \[
        \frac{\chi_\Y\left(f_!(a)\right)}{e_\Y}
        =(\I_p^nf)_!\left(
        \frac{\chi_\X(a)}{e_\X}
        \right).
        \]
        Notice that we can divide by $e_\X,e_\Y,e_N,e_V$
        because they are invertible.
\end{proof}
In the next case, we deal with maps $\X\to \pt$
and $\I_p^n\X\to \pt$. We denote the corresponding
pushforward maps by $t_\X^E$ and $t_{\I_p^n\X}^L$.
\begin{lemma}
\label{lem:RR Y is a point case}
    \cref{thm:a Riemann--Roch} is true if $\Y=\pt$. Equivalently,
    \[
        (E\to L(E))\left(t_\X^E\left(a\right)\right)=
        t_{\I_p^n \X}^L\left(
        \frac{\chi_\X(a)}{e_\X}
        \right)
    \]
    Here $(E\to L(E))$ denotes the ring homomorphism
    $E^*\to L(E)^*$.
\end{lemma}
\begin{proof}
    There are two traces on $H^*(\X,E)\otimes_{E^*}L(E)^*\simeq H^*\left(\I_p^n\X,L(E)\right)$.
    \begin{itemize}
        \item $\tau(c)=(E\to L(E))t_\X^E\left(\chi_\X^{-1}(c)\right)$.
        It is defined by base-changing the trace $t_\X^E$ on $H^*(\X,E)$
        to $H^*(\X,E)\otimes_{E^*}L(E)^*$.
        \item $\sigma(c)=t_{\I_p^n\X}^L\left(\frac{c}{e_\X}\right)$.
    \end{itemize}
    The lemma is equivalent to showing $\tau=\sigma$.

    We need the following inputs. Recall $\delta_\X=(\Delta_\X)_!^E(1)$.
    \begin{enumerate}[(a)]
        \item $K:=\chi_{\X\times \X}(\delta_\X)=(\Delta_{\I_p^n\X})_!(e_\X)$.
        
        Consider the representable map $\Delta_\X\colon \X\to \X\times \X$. We have
        \[
            \I_p^n(\X\times \X)=\I_p^nX\times \I_p^nX
        \]
        \[
            e_{\X\times \X}=e_\X \boxtimes e_\X.
        \]
        Apply \cref{lem:RR representable case}:
        \[
            \chi_{\X\times \X}((\Delta_\X)_!(1))=
            e_{\X\times \X}\cdot (\Delta_{\I_p^n\X})_!\left(\frac{\chi_\X(1)}{e_{\X}}\right)
            =(\Delta_{\I_p^n\X})_! (e_\X)
        \]
        The last equality is by projection formula, because
        $(\Delta_{\I_p^n\X})^*e_{\X\times \X}=(\Delta_{\I_p^n\X})^*(e_\X\boxtimes e_\X)=e_\X^2$.
        \item For all $c\in H^*\left(\I_p^n\X,L(E)\right)$, $(c\boxtimes 1) K=(1\boxtimes c)K$.

        This is immediate by projection formula:
        \[
            (c\boxtimes 1)(\Delta_{\I_p^n\X})_!(e_\X)
            =(\Delta_{\I_p^n\X})_!\left((\Delta_{\I_p^n\X})^*(c\boxtimes 1)\cdot e_\X\right)
            =(\Delta_{\I_p^n\X})_!(c \cdot e_\X).
        \]
        The same is true for $(1\boxtimes c)K$. So they are equal.
        \item $(\tau\otimes 1)\left((c\boxtimes 1) K\right)=c$.

        This can be understood already in $H^*(\X,E)$: let $p_1,p_2\colon \X\times \X\to \X$ be the projections to the first and second entries,
        then it suffices to show for $a\in H^*(\X,E)$
        \begin{equation}
            \label{eq:inverse kernel property}
            (p_2)_!\left(p_1^*a \cup \delta_\X\right)= a.
        \end{equation}
        Again by projection formula
        \[
            (p_2)_!\left(
            p_1^* a \cup (\Delta_\X)_!(1)\right)
            =(p_2)_!\left(
            (\Delta_\X)_!(\Delta_\X)^*\circ  p_1^*(a)\right)
            =(p_2\circ \Delta_\X)_!\circ (p_1\circ \Delta_\X)^*(a)
            =a,
        \]
        because $p_1\circ \Delta_\X=p_2\circ \Delta_\X=\id_\X$.
        \item $(1\otimes \sigma)\left((1\boxtimes c) K\right)=c$.

        Let $q_1\colon \I_p^n\X\times \I_p^n\X\to \I_p^n\X$ be
        projection to the first entry. By definition of $\sigma$,
        \begin{align*}
            (1\otimes \sigma)\left((1\boxtimes c) K\right)
            &=(q_1)_!\left(
                \frac{(1\boxtimes c) K}{1\boxtimes e_\X}
            \right)\\
            &=(q_1)_!\left(
                \frac{(\Delta_{\I_p^n})_!(c e_\X)}{1\boxtimes e_\X}
            \right)\\
            &=(q_1)_!\left(
                \frac{(1\boxtimes e_\X)\cdot  (\Delta_{\I_p^n})_!(c)}{1\boxtimes e_\X}
            \right)\\
            &=(q_1)_!\circ (\Delta_{\I_p^n})_!(c)\\
            &= c
        \end{align*}
    \end{enumerate}
    Finally we use the four facts to prove the lemma.
    \begin{align*}
        \tau(c)&\overset{(d)}{=}
        \tau\left((1\otimes \sigma)\left((1\boxtimes c) K\right)\right)\\
        &=(\tau\otimes \sigma)\left(\left((1\boxtimes c) K\right)\right)\\
        &=\sigma\left((\tau\otimes 1)\left((c\boxtimes 1) K\right)\right)\\
        &\overset{(c)}{=}\sigma(c)
    \end{align*}
\end{proof}
\begin{proof}[Proof of \cref{thm:a Riemann--Roch}]
    Let $f\colon X\to \Y$ be any morphism. By projection formula we have
    for $a\in H^*(\X,E), b\in H^*(\Y,E)$
    \begin{equation}
        \label{eq:integration projection formula}
        t_\X^E(a f^* b)= t_\Y^E\circ f_! (a f^*b)
        =t_\Y^E \left((f_! a)\cdot b\right)
    \end{equation}
    Therefore, fix $a\in H^*(\X,E)$ and for any $b\in H^*(\Y,E)$,
    we compute
\begin{align*}
&t_{\I_p^n(\mathcal Y)}^L\left(
    \frac{\chi_{\mathcal Y}(f_!^E a)}{e_{\mathcal Y}}
    \chi_{\mathcal Y}(b)
\right) \\
&= \operatorname{ch}\left(
    t_{\mathcal Y}^E((f_!^E a)b)
\right)
&& \text{by \cref{lem:RR Y is a point case}} \\
&= \operatorname{ch}\left(
    t_{\mathcal X}^E(a f^*b)
\right)
&& \text{by \cref{eq:integration projection formula}} \\
&= t_{\I_p^n(\mathcal X)}^L\left(
    \frac{\chi_{\mathcal X}(a)}{e_{\mathcal X}}
    \chi_{\mathcal X}(f^*b)
\right)
&& \text{by \cref{lem:RR Y is a point case}} \\
&= t_{\I_p^n(\mathcal X)}^L\left(
    \frac{\chi_{\mathcal X}(a)}{e_{\mathcal X}}
    (\I_p^n f)^*\chi_{\mathcal Y}(b)
\right)
&& \text{by character naturality} \\
&= t_{\I_p^n(\mathcal Y)}^L\left(
    (\I_p^n f)_!^L\left(
        \frac{\chi_{\mathcal X}(a)}{e_{\mathcal X}}
    \right)
    \chi_{\mathcal Y}(b)
\right)
&& \text{by projection formula}
\end{align*}
In the end, the pairing $t_{\I_p^n \X}^L(a\cup b)$ is perfect\footnote{
Since we base-change to $L(E)$,
we have Kunneth formula (\cref{lem:Kunneth formula after base change}). Write $\delta_\X=\sum_i x_i\boxtimes y_i$.
Then by \eqref{eq:inverse kernel property}
\[
    \sum_i t_\X^E (a \cdot x_i)y_i=a
\]
So if $a\neq 0$ (after base change), there is
some $x_i$ such that $t_\X^E(a x_i)\neq 0$.
This means the pairing is perfect after base-change.}
and $\chi_\Y(b)$ spans $H^*\left((\I_p^n \X,L(E)\right)$. So we reach
\[
    \frac{\chi_{\mathcal Y}(f_!^E a)}{e_\Y}=(\I_p^n f)_!^L\left(
        \frac{\chi_{\mathcal X}(a)}{e_{\mathcal X}}
    \right)
\]
\end{proof}
\subsubsection{Torus-equivariant pushforward}
Let $\TT$ be a real torus and let $E$ be $K(n)$ or an HKR theory.
Suppose $[M/G]$ and $[M'/G]$ are stably almost complex global quotient
orbifolds with $\TT$-actions in the sense of
\cref{def:working definition of group action on orbifold}, using the same
extension $\widetilde G$.
Assume that $M$ and $M'$ have the $\widetilde G$-equivariant homotopy type
of finite $\widetilde G$-CW complexes.

Let $f\colon M\to M'$ be a smooth $\widetilde G$-equivariant map.
Assume either that $M,M'$ have no boundary and $f$ is proper, or that
they are compact with boundary and $f(\partial M)\subseteq\partial M'$.
We construct an equivariant pushforward
\[
    f_{!,\eq}\colon H^*\left(M_{h\widetilde{ G}},E\right)
    \to H^{*+\dim M'-\dim M}(M'_{h\widetilde{G}},E).
\]
\begin{definition}
    Let $k=\dim \TT$. 
    \begin{enumerate}[(1)]
        \item  For $N\in \N$, we define \[
   P_N=\underbrace{S^{2N+1}\times S^{2N+1} \times \cdots\times S^{2N+1}}_{k\text{ times}}=(S^{2N+1})^k.
\]
\[
    B_N=(\CP^N)^k.
\]
We equip $P_N$ with the obvious $\TT$-action such that $P_N\to B_N$
is a $\TT$ principal bundle. There are natural embeddings
$P_N\hookrightarrow P_{N+1}$, $B_N\hookrightarrow B_{N+1}$.
As $N\to \infty$, $B_N$ converges to $B\TT$ and $P_N$
converges to a contractible space.
        \item Define an $\widetilde{G}$-equivariant map
        \[
            f_N=f\times \id_{P_N}\colon M\times P_N\to M'\times P_N.
        \]
        The $\widetilde{G}$-action
on $M$($M'$) comes with definition, and the $\widetilde{G}$-action
on $P_N$ is via $\widetilde{G}\to \TT$ and $\TT$ acts on $P_N$.
    \end{enumerate}  

\end{definition}
There is a cartesian square
\[\begin{tikzcd}
	{[M/G]} & {[M\times P_N/\widetilde{G}]} \\
	& {B_N=(\CP^N)^k}
	\arrow[from=1-1, to=1-2]
	\arrow[from=1-2, to=2-2]
\end{tikzcd}\]
The stably almost complex structures on $[M/G]$ and $[M'/G]$ induce
stably almost complex structures on the vertical tangent bundles of the
finite-stage orbifolds over $B_N$. Together with the complex structure on
$B_N$, these give stably almost complex structures on
$[M\times P_N/\widetilde G]$ and $[M'\times P_N/\widetilde G]$.

The non-equivariant pushforwards $(f_N)_!$ are defined under the above
assumptions and commute with restriction by transverse base change:
\begin{equation}
\label{eq:compatibility of classifying space approximation}
\begin{tikzcd}[column sep=large, row sep=large]
    \cdots
    & {H^*\left((M\times P_N)_{h\widetilde{G}},E\right)}
    & {H^*\left((M\times P_{N+1})_{h\widetilde{G}},E\right)}
    & \cdots \\
    \cdots
    & {H^{*+\delta}\left((M'\times P_N)_{h\widetilde{G}},E\right)}
    & {H^{*+\delta}\left((M'\times P_{N+1})_{h\widetilde{G}},E\right)}
    & \cdots
    \arrow["{(f_N)_!}"', from=1-2, to=2-2]
    \arrow["{(f_{N+1})_!}", from=1-3, to=2-3]
    \arrow["{i_N^*}", from=1-3, to=1-2]
    \arrow["{(i'_N)^*}"', from=2-3, to=2-2]
\end{tikzcd}
\end{equation}
Here $\delta=\dim_\R M'-\dim_\R M$. 
Consequently, the maps $(f_N)_!$ induce a map on
inverse limits:
\[
f_{!,\eq}\colon
\varprojlim_N
H^*\left((M\times P_N)_{h\widetilde G},E\right)
\longrightarrow
\varprojlim_N
H^{*+\delta}\left((M'\times P_N)_{h\widetilde G},E\right).
\] 
$E^0$ is a complete Noetherian local ring, and 
by the finite-type assumption, $H^i((M\times P_N)_{h\widetilde G},E)$
is a finitely generated $E^0$-module for every $N$ and $i$,
and similarly for $M'$.
The following lemma is a special case of
\cite[Theorem 1.1]{BakerRichter2005}, and it
identifies these inverse limits with
$H^*(M_{h\widetilde{G}},E)$
\begin{lemma}[{\cite[Theorem 1.1]{BakerRichter2005}}]
    Let $(R,\m)$ be a complete Noetherian local ring and
    let $M_i$ be an inverse sequence of finitely generated
    $R$-modules. Then
    \[
        \limone_i M_i=0
    \]
\end{lemma}
\begin{definition}\label{def:equivariant pushforward}
Under the above assumptions, we define
\[
    f_{!,\eq}\colon H^*\left(M_{h\widetilde{ G}},E\right)
    \to H^{*+\dim M'-\dim M}(M'_{h\widetilde{G}},E)
\]
via the identifications
\[
    \varprojlim_N H^*\left((M\times P_N)_{h\widetilde G},E\right)
    \simeq H^*\left(M_{h\widetilde G},E\right),
\]
and similarly for $M'$.
\end{definition}
Next we consider $\TT$-equivariant HKR theory.
\cref{def:p-inertia stack construction} can be also generalized.
Let $[M/G]$ be $\TT$-equivariant in the sense of \cref{def:working definition of group action on orbifold}. Define (compare \cref{rmk:definition of p-inertia stack} (1))
    \[
        \Fix_{n,p}^G M=\{
        (\alpha,x)\in \Hom(\Lambda_m,G)\times M|
        \im(\alpha)\subset \Stab(x)
        \}
    \]
    $\Fix_{n,p}^G M$ admits a $\widetilde{G}$-action by
    $\widetilde{g}\cdot (\alpha,x)=(\widetilde{g}\alpha\widetilde{g}^{-1},\widetilde{g}x)$.
    This defines a $\TT$-action in the sense
    of \cref{def:working definition of group action on orbifold} on
    $[\Fix_{n,p}^G M/G]$.
\begin{definition}
    We call $[\Fix_{n,p}^G M/G]$ with the $\TT$-action the
    \emph{equivariant $n$-fold $p$-inertia stack}.
\end{definition}
There is a decomposition
\[
    [\Fix_{n,p}^G M/\widetilde{G}]=\bigsqcup_{[\alpha]}
    [M^{\im\alpha }/Z_{\widetilde{G}}(\alpha)]
\]
over $\widetilde{G}$-conjugacy classes of $\alpha\colon \Lambda_m\to G$, where $Z_{\widetilde{G}}(\alpha)=\{\widetilde{g}\in\widetilde{G}| \widetilde{g}\cdot \alpha(-)=\alpha(-)\cdot\widetilde{g}\}$.  
\begin{definition}
\label{def:equivariant universal map and equivariant normal bundle}
\begin{enumerate}[(1)]
    \item The \emph{equivariant universal map} $q\colon [\pt/\Lambda_m]\times
  [\Fix_{n,p}^G M/\widetilde{G}]\to [M/\widetilde{G}]$ on the component $[\pt/\Lambda_m]\times [M^{\im\alpha }/Z_{\widetilde{G}}(\alpha)]$ is defined to be
  \[
    M^{\im\alpha}\hookrightarrow M
  \]
  \[
  \Lambda_m\times Z_{\widetilde{G}}(\alpha)\to \widetilde{G},
  (u,\widetilde{g})\mapsto \alpha(u)\cdot \widetilde{g}
  \]
    \item We can define an $\TT$-equivariant vector orbi-bundle $\nu$ on $[\Fix_{n,p}^G M/G]$ as the normal
    bundle of $\Fix^G_{n,p}M\to M$. In fact, on each
    $M^{\im\alpha}$, $\nu|_{M^{\im\alpha}}=N_{M^{\im\alpha}/M}$ is
    the normal bundle of $M^{\im\alpha}\subset M$.
    Its Euler class in $L(E)$ cohomology is denoted by
    $e_{\mathrm{HKR},\eq}^L(\nu)\in H^*\left(
        (\Fix_{n,p}^G M)_{h\widetilde{G}}, L(E)
        \right)$
\end{enumerate}
\end{definition}

\begin{definition}
\label{equivariant HKR character}
    We define \emph{equivariant HKR character} to be the composition
    \[
        \chi_{n,p;\eq}\colon
        H^*(M_{h\widetilde{G}},E)
        \xrightarrow{q^*}H^*\left(B\Lambda_m\times (\Fix_{n,p}^G M)_{h\widetilde{G}},E\right)
        \to H^*\left(
        (\Fix_{n,p}^G M)_{h\widetilde{G}}, L(E)
        \right).
    \]
\end{definition}
The following proposition is the torus-equivariant form of \cref{thm:a Riemann--Roch}.
\begin{proposition}
\label{prop:Riemann Roch for equivariant HKR}
    Assume $[M/G]$ and $[M'/G]$ admit a $\TT$-action in the sense of
    \cref{def:working definition of group action on orbifold}, 
    and $f\colon M\to M'$ is $\widetilde{G}$-equivariant.
    There's a 
    commutative diagram
    \begin{equation}
    \label{eq:equivariant HKR Rieman Roch}
    \begin{tikzcd}[column sep={between origins,12em}, row sep=large]
        {H^*\left(M_{h\widetilde{G}},E\right)}
        &&
        {H^*\left(
            (\Fix_{n,p}^G M)_{h\widetilde{G}},L(E)
        \right)}
        \\
        {H^*\left(M'_{h\widetilde{G}},E\right)}
        &&
        {H^*\left(
            (\Fix_{n,p}^G M')_{h\widetilde{G}},L(E)
        \right)}
        \arrow[
            "{{e_{\mathrm{HKR},\eq}^L(\nu_M)}^{-1}
            \cdot\chi_{n,p;\eq}}",
            from=1-1,
            to=1-3
        ]
        \arrow[
            "{{e_{\mathrm{HKR},\eq}^L(\nu_{M'})}^{-1}
            \cdot\chi_{n,p;\eq}}",
            from=2-1,
            to=2-3
        ]
        \arrow["{f_{!,\eq}}", from=1-1, to=2-1]
        \arrow["{(f^{\Fix})_{!,\eq}}", from=1-3, to=2-3]
    \end{tikzcd}
\end{equation}
\end{proposition}
\begin{proof}
    Set
    \[
        \X_N=[(M\times P_N)/\widetilde{G}],\quad 
        \Y_N=[(M'\times P_N)/\widetilde{G}]
    \]
    and write $f_N\colon \X_N\to \Y_N$. This equivariant version
    of Riemann--Roch formula can be derived
    from \cref{thm:a Riemann--Roch} once we 
    have the following commutative cube diagram.
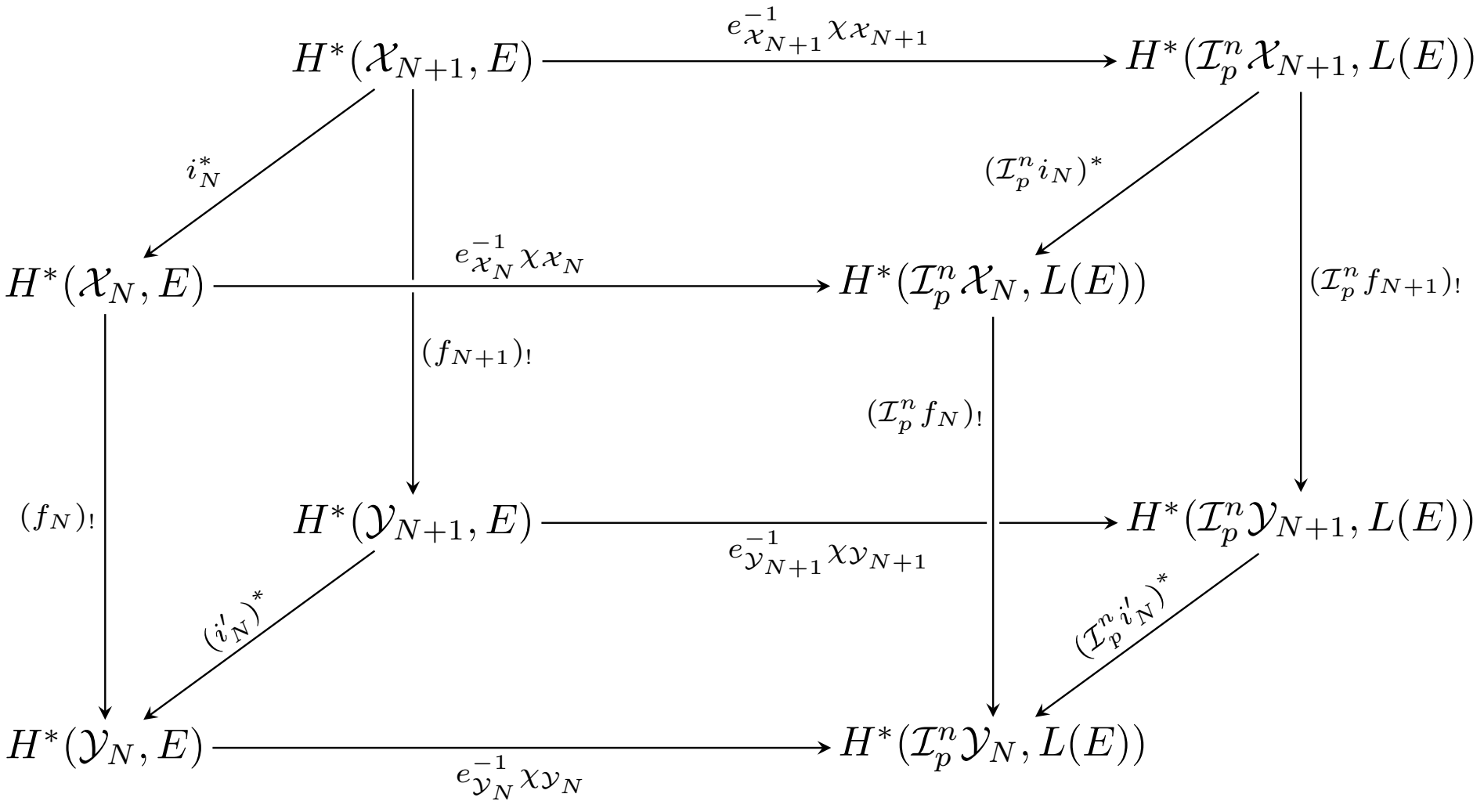
\begin{equation}
\label{eq:HKR-RR-cube}
\begin{tikzpicture}[
  >=stealth,
  every node/.style={inner sep=1.5pt},
  every path/.style={line width=.45pt},
  every edge quotes/.style={font=\scriptsize}
]
\node (XNp) at (2.6,6.5) {$H^*(\mathcal X_{N+1},E)$};
\node (INp) at (10.1,6.5) {$H^*(\I_p^n\mathcal X_{N+1},L(E))$};
\node (YNp) at (2.6,2.6) {$H^*(\mathcal Y_{N+1},E)$};
\node (JNp) at (10.1,2.6) {$H^*(\I_p^n\mathcal Y_{N+1},L(E))$};
\node (XN) at (0,4.6) {$H^*(\mathcal X_N,E)$};
\node (IN) at (7.5,4.6) {$H^*(\I_p^n\mathcal X_N,L(E))$};
\node (YN) at (0,0.7) {$H^*(\mathcal Y_N,E)$};
\node (JN) at (7.5,0.7) {$H^*(\I_p^n\mathcal Y_N,L(E))$};
\draw[->] (XNp) -- node[above,font=\scriptsize]
  {$e_{\mathcal X_{N+1}}^{-1}\chi_{\mathcal X_{N+1}}$} (INp);
\draw[->] (YNp) -- node[below,font=\scriptsize]
  {$e_{\mathcal Y_{N+1}}^{-1}\chi_{\mathcal Y_{N+1}}$} (JNp);
\draw[->] (XNp) -- node[right,font=\scriptsize,pos=.65]
  {$(f_{N+1})_!$} (YNp);
\draw[->] (INp) -- node[right,font=\scriptsize,pos=.48]
  {$(\I_p^n f_{N+1})_!$} (JNp);
\draw[->] (XNp) -- node[above left,font=\scriptsize,pos=.62] {$i_N^*$} (XN);
\draw[->] (INp) -- node[above left,font=\scriptsize,pos=.65] {$(\I_p^n i_N)^*$} (IN);
\draw[->] (YNp) -- node[above,sloped,font=\scriptsize,pos=.53] {$(i_N')^*$} (YN);
\draw[->] (JNp) -- node[above,sloped,font=\scriptsize,pos=.52] {$(\I_p^n i_N')^*$} (JN);
\draw[->,preaction={draw=white,line width=3pt}]
  (XN) -- node[above,font=\scriptsize]
  {$e_{\mathcal X_N}^{-1}\chi_{\mathcal X_N}$} (IN);
\draw[->] (YN) -- node[below,font=\scriptsize]
  {$e_{\mathcal Y_N}^{-1}\chi_{\mathcal Y_N}$} (JN);
\draw[->] (XN) -- node[left,font=\scriptsize,pos=.5]
  {$(f_N)_!$} (YN);
\draw[->,preaction={draw=white,line width=3pt}]
  (IN) -- node[left,font=\scriptsize,pos=.25]
  {$(\I_p^n f_N)_!$} (JN);
\end{tikzpicture}
\end{equation}
 \begin{enumerate}[(a)]
        \item Commutativity of the left and right faces
        comes from \cref{prop:properties of f_!}(2).
        \item Commutativity of the front and back
        faces comes from \cref{thm:a Riemann--Roch}.
        \item To prove the commutativity of the top and bottom
        faces, we claim that 
        \begin{equation}
            \label{eq:normal bundle identification}
             \nu_{\X_{N+1}}|_{\I_p^n\X_N}\simeq \nu_{\X_N}.
        \end{equation}
        Here $\nu_{\X_N},\nu_{\X_{N+1}}$ are defined in
        \cref{def:equivariant universal map and equivariant normal bundle}(2). Once \cref{eq:normal bundle identification} is true,
        for $a\in H^*(\X_{N+1},E)$ we have
        \[
        \begin{aligned}
            (\I_p^n i_N)^*\left(e_{\X_{N+1}}^{-1}\cdot \chi_{\X_{N+1}}\right)(a)
            &=(\I_p^n i_N)^*(e_{\X_{N+1}})^{-1} \cdot
            (\I_p^n i_N)^*\circ \chi_{\X_{N+1}}(a)\\
            &=e_{\X_{N}}^{-1} \cdot
             (\I_p^n i_N)^*\circ \chi_{\X_{N+1}}(a)
        \end{aligned}
        \]
        What remains is $(\I_p^n i_N)^*\circ \chi_{\X_{N+1}}=
        \chi_{\X_{N}}\circ i_N^*$, but this comes from naturality
        of HKR characters.
   
        To see \eqref{eq:normal bundle identification}, let $\alpha\colon \Lambda_r\to H$
        be an HKR sector, and $x\in (M\times P_N)^{\im\alpha}$
        represent a point in $\I_p^n\X_N$. We have
\[\begin{tikzcd}
	{(M\times P_N)^{\im\alpha}} & {M\times P_N} \\
	{(M\times P_{N+1})^{\im\alpha}} & {M\times P_{N+1}}
	\arrow[hook, from=1-1, to=1-2]
	\arrow[hook, from=1-1, to=2-1]
	\arrow[hook, from=1-2, to=2-2]
	\arrow[hook, from=2-1, to=2-2]
\end{tikzcd}\]
    Then $\nu_{\X_N}=\nu_{(M\times P_N)^{\im\alpha}/M\times P_N}$
    and $\nu_{\X_{N+1}}=\nu_{(M\times P_{N+1})^{\im\alpha}/M\times P_{N+1}}$. Their isomorphism is equivalent to another isomorphism
    \[
        \nu_{(M\times P_N)^{\im\alpha}/(M\times P_{N+1})^{\im\alpha}}
        \simeq \nu_{M\times P_N/M\times P_{N+1}}.
    \]
    This is true because $\nu_{M\times P_N/M\times P_{N+1}}\simeq \nu_{B_N/B_{N+1}}$ via $M\times P_N\to B_N$, and $\alpha\colon \Lambda_r\to G$ acts trivially on $B_N$ (because the action is via
    $\widetilde{G}\to \TT$). This means $\alpha$ fixes $\nu_{M\times P_N/M\times P_{N+1}}$, so $\nu_{M\times P_N/M\times P_{N+1}}$ must
    lie entirely in $\nu_{(M\times P_N)^{\im\alpha}/(M\times P_{N+1})^{\im\alpha}}$. 
    \end{enumerate}
    Taking inverse limit for $N\to \infty$, as $f_{!,\eq},(\I_p^n f)_{!,\eq}$ are all defined by inverse limits, and that
    $\nu_{\X_N}$ are the same for all $N$ and are restrictions
    of $\nu_M$ on the whole $[M/\widetilde{G}]$, we get
    the commutativity of \eqref{eq:equivariant HKR Rieman Roch}.
\end{proof}
\subsection{Equivariant localization and invertibility of normal Euler classes}
\label{section: equivariant localization}
Let $A$ be a compact abelian Lie group of the form
$A=\TT\times \Lambda_r$, where $\TT$ is a real torus and
$\Lambda_r=(C_{p^r})^n$. $E$ is an HKR theory at prime $p$ of height $n$.
When $A=\TT$, we also allow $E=K(n)$.
For nonzero characters $\chi\colon A\to S^1$,
the Euler class of the corresponding complex line bundle
over $BA$ is denoted $e_E(\chi)$.
\begin{definition}
    We denote $S_A\subset H^*(BA,E)$ the multiplicative set
    generated by $e_E(\chi)$ for all nonzero characters $\chi$.
\end{definition}
The HKR hypothesis is needed when the finite factor $\Lambda_r$ is nontrivial.
\begin{lemma}
    Under the above assumptions, $0\notin S_A$ and $S_A^{-1}H^*(BA,E)\neq 0$.
\end{lemma}
\begin{proof}
    Suppose that $E=K(n)$ and $A=\TT$.
    After a change of torus coordinates, the Euler class of a nontrivial character
    is $[g]_E t_1$ for some positive integer $g$, which is nonzero by the formula
    for the $p$-series. Since $K(n)^*\llbracket t_1,\ldots,t_d\rrbracket$
    is an integral domain, these Euler classes are non-zero-divisors.
    This proves the claim in this case.

    We now assume that $E$ is an HKR theory.

    When $A=\Lambda_r$, it is the same setting
    as in \cite{Generalized-Group-Characters-and-Complex-Oriented-Cohomology-Theories}. Write 
         $S_r:=S_{\Lambda_r}\subset H^*(B\Lambda_r,E)$.
         By \cite[Proposition 6.5]{Generalized-Group-Characters-and-Complex-Oriented-Cohomology-Theories}, 
         $S_r^{-1}H^*(B\Lambda_r,E)$ is faithfully flat
         over $p^{-1}E^*$. $p^{-1}E^*\neq 0$ by definition
         of HKR theories. Therefore, $S_r^{-1}H^*(B\Lambda_r,E)\neq 0$
         and $0\notin S_r$ in this case.

    Now consider general $A=\TT^d\times\Lambda_r$. There's an isomorphism
    $H^*(B\TT^d,E)\simeq E^*\llbracket t_1,t_2,\ldots,t_d\rrbracket$.
    Let $R=S_{\Lambda_r}^{-1}H^*(B\Lambda_r,E)$. There is
    a ring homomorphism
    \[
        q\colon H^*(B\TT^d\times B\Lambda_r,E)
        \to R\llbracket t_1,t_2,\ldots,t_d\rrbracket
    \]
    and the target is nonzero because $R\neq 0$ ($A=\Lambda_r$ case).

    Let $\chi\colon A\to S^1$ be a character of $A$, we consider
    the $q(e_E(\chi))\in R\llbracket t_1,t_2,\ldots,t_d\rrbracket$.
    We can distinguish between whether the restricted character $\chi'\colon \Lambda_r\to \TT^d\times \Lambda_r\xrightarrow{\chi} S^1$ is trivial.
    \begin{itemize}
        \item If $\chi'$ is nontrivial,
        then we claim $q(e_E(\chi))$ is a unit.
        In fact, the constant term of $q(e_E(\chi))\in R\llbracket t_1,t_2,\ldots,t_d\rrbracket$ is the image of
        $e_E(\chi')$. Since $e_E(\chi')$ is inverted in $R$,
        and in $R\llbracket t_1,t_2,\ldots,t_d\rrbracket$
        an element is invertible iff its constant term is,
        we get $q(e_E(\chi))$ is a unit.
        \item If $\chi'$ is trivial, then we claim
        $q(e_E(\chi))$ is not a zero-divisor.
        The character is the composition of 
        the projection $\TT^d\times \Lambda_r\to \TT^d$ with a
        homomorphism 
        $a\colon \TT^d\to S^1$. It's obvious that there is
        an automorphism $\phi$ of $\TT^d$ such that
        $a\circ \phi=(g,0,0,\ldots,0)\colon \TT^d\to S^1$ for
        some positive integer $g$. Therefore,
        \[
            (\phi\times\id_{\Lambda_r})^*(e_E(\chi))=
            [g]_Et_1.
        \]
        Since $g\neq 0$, the power series $[g]t_1$ is not
        a zero-divisor in $R\llbracket t_1,t_2,\ldots,t_d\rrbracket$.
    \end{itemize}
    For nonzero characters $\chi_1,\chi_2,\ldots,\chi_k$, 
    $q(e_E(\chi_i))$ is either a unit or a non-zero-divisor,
    so $\prod_i q(e_E(\chi_i))$ is a non-zero-divisor in $R\llbracket t_1,t_2,\ldots,t_d\rrbracket$.
    This implies $\prod_i e_E(\chi_i)\neq 0$.
\end{proof}
Consider a global quotient compact orbifold $[M/G]$ with $A$-action
in the sense of \cref{def:working definition of group action on orbifold}. Recall that there is an exact sequence
\[
    1\to G\to \widetilde{G}\to A\to 1
\]
and $\widetilde{G}$ extends the almost free action of $G$ on $M$.
In \cref{def:fixed points of orbifolds}, we defined
\[
    M^\fix=\{p\in M| \Stab_{\widetilde{G}}(p)\to A\text{ is surjective.}\}
\]
The following lemma is a direct analogue 
of Atiyah--Bott localization
theorem for torus actions 
\cite[Theorem~3.5]{AtiyahBottMomentMap}.
See \cite[Lemma~3.6]{bai2024equivariantformalitycomplexorientedtheories}
for the case when $A=\TT$ and $[M/G]$ is a manifold.
\begin{lemma}
\label{lem:combined form equivariant localzation}
    Let $i^* \colon H^*_{\widetilde{G}}(M,E)\to H^*_{\widetilde{G}}(M^\fix,E)$ be the restriction
    map on $E$-cohomology. Then after inverting
    $S_A$, we get an isomorphism
    \[
        S_A^{-1}H^*_{\widetilde{G}}(M,E)
        \xrightarrow{i^*}
        S_A^{-1} H^*_{\widetilde{G}}(M^\fix,E)
    \]
    In particular, if $M^\fix=\varnothing$, then $S_A^{-1}H^*_{\widetilde{G}}(M,E)=0$.
\end{lemma}
We only give a sketch of proof.
Assume a $\widetilde{G}$-CW structure on $M$
and argue by induction over cells. It suffices
to prove that
\[
    S_A^{-1}H^*(BH,E)=0
\]
for every stabilizer
$H=\Stab_{\widetilde{G}}(x)$
with $x\in M\setminus M^\fix$.
By definition of $M^\fix$,
the map $H\subset \widetilde{G}\to A$ is not surjective,
so there is a nontrivial character $\chi\colon A\to S^1$
such that $H\to A\xrightarrow{\chi}S^1$ is trivial.
This means the element $e_E(\chi)=0$ in $H^*(BH,E)$.
But $e_E(\chi)\in S_A$, so $S_A^{-1}H^*(BH,E)=0$.

For the next proposition, we write 
$\X=[M/G]$ and $\Fcal=[M^\fix/G]$, and call
$\Fcal$ the fixed locus of $\X$.
For $M^\fix\subset M$, let $\nu_\Fcal=\nu_{M^{\fix}/M}$ be
the normal bundle. Notice that $\nu$ is 
$\widetilde{G}$-equivariant, so $\nu$
can be also viewed as an $A$-equivariant vector orbi-bundle on
$\Fcal=[M^\fix/G]$ (at least in the working sense
of \cref{def:working definition of group action on orbifold}).

Let $i\colon \Fcal\to \X$ be the natural inclusion 
induced by $M^\fix\hookrightarrow M$. There are
pullback $i^*\colon H^*_{\widetilde{G}}(M,E)\to
H^*_{\widetilde{G}}(M^\fix,E)$ and equivariant pushforward $i_{!,\eq}\colon H^*_{\widetilde{G}}(M^\fix,E)
\to H^*_{\widetilde{G}}(M,E)$\footnote{
We don't need ambidexterity to define $i_{!,\eq}$ because
$\Fcal\to \X$ is representable, so the construction
relies only on \cref{lem:Pontryagin-Thom for orbifolds}.
}. It's classical that $i^*\circ i_{!,\eq}(a)=a\cdot e_{E,\eq}(\nu)\in
H^*_{\widetilde{G}}(M^\fix,E)$ for every $a\in H^*_{\widetilde{G}}(M^\fix,E)$, where $e_{E,\eq}(\nu)\in H^*_{\widetilde{G}}(M^\fix,E)$
is the $A$-equivariant Euler class of $\nu$. 
\begin{proposition}
\label{prop:invertibility of Euler class&isomorphism after S_A localization}
    After localizing at $S_A$, $i^*$ and $i_{!,\eq}$
        become isomorphisms. In other words, there are isomorphisms
        still denoted by $i^*,i_{!,\eq}$:
        \[
            i^*\colon S_A^{-1}H^*_{\widetilde{G}}(M,E)
            \to S_A^{-1}H^*_{\widetilde{G}}(M^\fix,E)
        \]
        \[
            i_{!,\eq}\colon S_A^{-1}H^*_{\widetilde{G}}(M^\fix,E)
            \to 
            S_A^{-1}H^*_{\widetilde{G}}(M,E)
        \]
        In particular, 
        the equivariant Euler class $e_{E,\eq}(\nu)\in
        H^*_{\widetilde{G}}(M^\fix,E)
        $ becomes invertible in
        \[
            S_A^{-1}H^*_{\widetilde{G}}(M^\fix,E).
        \]
\end{proposition}

\begin{proof}
     We can find an $\widetilde{G}$-equivariant
        open neighborhood $U$ of $M^\fix$ in $M$ such that
        $M^\fix\to U$ is homotopy equivariant. By Mayer--Vietoris
        we have a long exact sequence
        \[
            \begin{tikzcd}[column sep=small]
            \cdots
            &
            {H^*_{\widetilde G}(M,E)}
            &
            [3.5em]
            {H^*_{\widetilde G}(M-M^{\mathrm{fix}},E)
             \oplus H^*_{\widetilde G}(U;E)}
            &
            {H^*_{\widetilde G}(U-M^{\mathrm{fix}},E)}
            &
            \cdots
            \arrow[from=1-1, to=1-2]
            \arrow[
              "{\scriptstyle i_{M-M^{\mathrm{fix}}}^* \oplus i_U^*}",
              from=1-2,
              to=1-3
            ]
            \arrow[from=1-3, to=1-4]
            \arrow[from=1-4, to=1-5]
            \end{tikzcd}
        \]
        By \cref{lem:combined form equivariant localzation},
        since for $M-M^\fix$ and $U-M^\fix$, their ``fixed'' parts
        are empty, we have
        \[
        \begin{aligned}
        S_A^{-1}H^*_{\widetilde{G}}\left(M-M^\fix,E\right)&=0\\
        \text{and}\quad S_A^{-1}H^*_{\widetilde{G}}\left(U-M^\fix,E\right)&=0.
        \end{aligned}
        \]
        Equivalently, after $S_A$-localization
        the long exact sequence gives the isomorphism
        \[
            S_A^{-1}H^*_{\widetilde{G}}\left(M,E\right)
            \to S_A^{-1}H^*_{\widetilde{G}}\left(U,E\right)
            \simeq S_A^{-1}H^*_{\widetilde{G}}\left(M^\fix,E\right),
        \]
        which is exactly $i^*$.

        On the other hand, consider the cohomological long exact sequence
        for the space pair $(M,M-M^\fix)$:
        \[\begin{tikzcd}
	\cdots & {H^*_{\widetilde{G}}\left(M,M-M^\fix,E\right)} & {H^*_{\widetilde G}\left(M;E\right)} & {H^*_{\widetilde G}\left(M-M^\fix;E\right)} & \cdots
	\arrow[from=1-1, to=1-2]
	\arrow[from=1-2, to=1-3]
	\arrow[from=1-3, to=1-4]
	\arrow[from=1-4, to=1-5]
\end{tikzcd}\]
For the same reason, $S_A^{-1}H^*_{\widetilde G}\left(M-M^\fix;E\right)=0$, so we get an isomorphism
\[
    S_A^{-1}H^*_{\widetilde{G}}\left(M,M-M^\fix,E\right)
    \to S_A^{-1}H^*_{\widetilde G}\left(M;E\right).
\]
By Thom isomorphism, we have
$H^*_{\widetilde{G}}\left(M,M-M^\fix,E\right)\simeq 
H^*_{\widetilde{G}}\left(M^\fix,E\right)$.
Combining the last two isomorphisms
yields that $i_{!,\eq}$ is isomorphic after
$S_A$-localization.

$e_{E,\eq}(\nu)$ is invertible after $S_A$-localization because $e_{E,\eq}(\nu)\cdot (-)
=i^*\circ i_{!,\eq}(-)$.
\end{proof}
\begin{proposition}
\label{prop:Atiyah-Bott equivariant localization formula}
    Let $\pi_\X\colon M_{h\widetilde{G}}\to BA$
        and $\pi_\Fcal\colon (M^\fix)_{h\widetilde{G}}\to BA$
        be the natural projections as in \eqref{eq:homotopy quotient fiber sequence}, and $(\pi_\X)_{!,\eq}$ and $(\pi_\Fcal)_{!,\eq}$ be the equivariant pushforwards
        in $E$-cohomology.
        Then, for $a\in H^*_{\widetilde{G}}(M,E)$ we have
        \begin{equation}
            \label{eq:Atiyah-Bott formula}
            (\pi_\X)_{!,\eq}(a)=(\pi_\Fcal)_{!,\eq}\left(
            \frac{i^*a}{e_{E,\eq}(\nu)}\right)
            \in S_A^{-1}H^*(BA,E)
        \end{equation}
        Notice that LHS can be defined in $H^*(BA,E)$, while RHS
        is defined only in $S_A^{-1}H^*(BA,E)$ because $e_{E,\eq}(\nu)$ is invertible only after $S_A$-localization.
\end{proposition}
\begin{proof}
    We have a commutative diagram
\[\begin{tikzcd}
	{H^*_{\widetilde{G}}\left(M^\fix,E\right)} & {H^*_{\widetilde{G}}\left(M,E\right)} \\
	{H^*(BA,E)}
	\arrow["{i_{!,\eq}}", from=1-1, to=1-2]
	\arrow["{(\pi_\Fcal)_{!,\eq}}"', from=1-1, to=2-1]
	\arrow["{(\pi_\X)_{!,\eq}}", from=1-2, to=2-1]
\end{tikzcd}\]
For the rest of this proof, all the cohomologies are $S_A$-localized. The commutative diagram means that for $a\in H^*_{\widetilde{G}}(M,E)$, we have
\[
    (\pi_\X)_{!,\eq}(a)=(\pi_\Fcal)_{!,\eq}\left((i_{!,\eq})^{-1}(a)\right).
\]
Thus, it suffices to show 
\[
    (i_{!,\eq})^{-1}(a)=\frac{i^*(a)}{e_{E,\eq}(\nu)}.
\] 
Let $b=(i_{!,\eq})^{-1}(a)$ and use the invertibility
of $e_{E,\eq}(\nu)$, it is equivalent to showing that
\[
    b\cdot e_{E,\eq}(\nu)=i^*i_{!,\eq}(b).
\]
which is classical.
\end{proof}
\begin{remark}
    If $A=\TT$, the Euler class of every nontrivial character is a
    non-zero-divisor in $H^*(B\TT,E)$. Indeed, after a change of torus
    coordinates, it has the form $[g]_Et_1$ for some positive integer $g$.
    For an HKR theory, the finite-height condition makes this a Weierstrass
    series, and uniqueness in Weierstrass division implies that multiplication
    by it is injective
    \cite[Lemma~5.1]{Generalized-Group-Characters-and-Complex-Oriented-Cohomology-Theories}.
    Therefore, we can recover invariants in $E^*$ from computations that land in $L(E)^*$.
\end{remark}

\section{Global Kuranishi charts and equivariant localization}
Global Kuranishi charts for Gromov--Witten moduli spaces were constructed in \cite{abouzaid2021complexcobordismhamiltonianloops}
for $g=0$, and in \cite{hirschiGlobalKuranishiCharts2024b} and \cite{abouzaidGromovWittenInvariantsComplex2024a} for arbitrary genus. 
In this section we give an exposition of global Kuranishi charts, discussing their constructions, equivariant global Kuranishi charts in the presence of a Hamiltonian action on the symplectic manifold, and equivariant localization on equivariant global Kuranishi charts. 
Most of the content is from
\cite{abouzaid2021complexcobordismhamiltonianloops,abouzaidGromovWittenInvariantsComplex2024a} and \cite{Hirschi_Properties_of_Gromov_Witten}.

\subsection{General aspects of global Kuranishi charts}
\begin{definition}[{\cite[Definition 4.1]{abouzaid2021complexcobordismhamiltonianloops}}]
\label{def:global Kuranishi chart}
    A global Kuranishi chart is the following data: 
    $\Kur=\left(G,\T,\E,s\right)$
    where
    \begin{enumerate}[(1)]
        \item $G$ is a compact Lie group.
        \item $\T$ is a $G$-equivariant topological manifold, and
        we require that the $G$-action on $\T$ is almost free, i.e. stabilizers
        are finite groups.
        \item $\E\to \T$ is a $G$-equivariant (topological) vector bundle.
        \item $s\colon \T\to \E$ is a $G$-equivariant continuous section of $\E$.
    \end{enumerate}
    The virtual dimension of
    a global Kuranishi chart
    $\Kur=\left(G,\T,\E,s\right)$
    is $\vdim \Kur=\dim \T-\dim G-\rk \E$.
\end{definition}
In applications, $\T$ being a topological
 manifold is not enough, so we introduce more
 notions.
 \begin{definition}[{\cite[Definition 4.3]{abouzaid2021complexcobordismhamiltonianloops}}]
     Let $\Kur=(G,\T,\E,s)$ be a global
     Kuranishi chart.
     \begin{itemize}
         \item $\Kur$ is called smooth
         iff $\T$ is a smooth manifold
         with smooth $G$-action, and $\E$
         is a smooth vector bundle. But
         $s$ is NOT required to be smooth.
         \item $\Kur$ is called free
         iff at every $x\in s^{-1}(0)$,
         the stabilizer group $\Stab_{G}(x)$
         is trivial.
         \item $\Kur$ is called regular
         iff $s$ is smooth in a neighborhood of
         $s^{-1}(0)$ and transverse at every $x\in s^{-1}(0)$.
     \end{itemize}
 \end{definition}
 \begin{example}
    \begin{itemize}
        \item Let $M$ be a smooth manifold,
        and let $\Kur=(e,M,0,0)$ be a global
        Kuranishi chart where $G=e$ is the trivial
        Lie group, $\T=M$, $\E$ is the zero
        vector bundle with $s=0$. This
        global Kuranishi chart is smooth, free
        and regular.
        \item Let $\X=[M/G]$ be a global quotient
        orbifold, then similarly $\Kur=(G,M,0,0)$
        is a global Kuranishi chart. This global
        Kuranishi chart is smooth and regular,
        but not free unless $G$-action on $M$
        is free.
    \end{itemize}
\end{example}
Notice that when $\Kur$ is smooth, $[\E/G]\to [\T/G]$ is a vector orbi-bundle over 
the orbifold $[\T/G]$ and $s$ is a continuous section. We will call $[T(\T/G)]$ the
\emph{deformation bundle} and $[\E/G]$
the \emph{obstruction bundle}. 
\begin{definition}[{\cite[Remark 3.6]{abouzaidGromovWittenInvariantsComplex2024a}}]
    Let $\Kur=(G,\T,\E,s)$ be a smooth global
    Kuranishi chart. A stably almost complex 
    structure on $\Kur$ is a stably almost
    complex structure on the virtual vector orbi-bundle $[T(U/G)]-[\E|_{U}/G]$ 
    (\cref{def:stably almost complex structures on vector orbi-bundles}) for a $G$-equivariant
    open neighborhood $U\supset s^{-1}(0)$.
\end{definition}
\begin{convention}
\label{conv:global Kuranishi charts}
Unless otherwise stated, all global Kuranishi charts are assumed to be
smooth and equipped with a stably almost complex structure. Moreover,
their zero loci $s^{-1}(0)$ are assumed to be compact.
\end{convention}

For a given moduli space $\M$, its global Kuranishi charts are far from unique.
We should allow the following equivalence
of global Kuranishi charts.
\begin{definition}[{\cite[Definition 3.4]{abouzaidGromovWittenInvariantsComplex2024a}}]
\label{def:equivalence of global Kuranishi charts}
    Let $\Kur'=(G',\T',\E',s'),\Kur=(G,\T,\E,s)$
    be two global Kuranishi charts, and $\mathfrak{g}',\mathfrak{g}$ are the Lie algebras of $G',G$.
    \begin{enumerate}[(1)]
        \item A morphism of global Kuranishi
        charts $f\colon \Kur'\to \Kur$ is a group
        homomorphism $\rho\colon G'\to G$ and a $\rho$-equivariant smooth
        map $f\colon \T'\to \T$, $\widetilde{f}
        \colon \E'\to f^*\E$ such that $\widetilde{f}\circ s'=f^*s$ as
        sections of $f^*\E\to \T'$.
        \item A morphism $f\colon \Kur'\to\Kur$
        is an equivalence iff 
        \begin{itemize}
            \item $f$ induces an isomorphism
            of the underlying spaces and
            stabilizers, i.e. there is a
            homeomorphism of quotient spaces
            $s'^{-1}(0)/G'\xrightarrow{\simeq}
            s^{-1}(0)/G$ and $\Stab_{G'}(x')
            \xrightarrow{\simeq}\Stab_G(x)$
            for any $x'\in s'^{-1}(0),x=f(x')$.
            \item The following chain complex
        of vector bundles is exact in
        a neighborhood of $s'^{-1}(0)$:
        \[
            \mathfrak{g}'\to \mathfrak{g}\oplus
            T\T'\to T\T\oplus \E'
            \to \E
        \]
        The maps between these terms are
        as follows: $\mathfrak{g}'\to \mathfrak{g}$ is induced by $\rho$,
        $\mathfrak{g'}\to T\T'$ and $\mathfrak{g}\to T\T$ are induced by the $G'$-action on $\T'$
        and the $G$-action on $\T$, respectively, $T\T'\to
        T\T$ is induced by $f$, $T\T'\to \E'$
        and $T\T\to \E$ is induced by $Ds',Ds$.
        Note that the map $\mathfrak{g}\to \E'$
        is $0$.
        The stable identification
        \[
            T\T'-\mathfrak{g}'-\E'
            \simeq f^*(T\T-\mathfrak{g}-\E)
        \]
        induced by this exact sequence must identify the specified
        stably almost complex structures near $s'^{-1}(0)$.
        \end{itemize}
    \end{enumerate}
\end{definition}
In fact, any equivalence is a
zig-zag of simple moves on global
Kuranishi charts.
\begin{lemma}[{\cite[Lemma 3.5]{abouzaidGromovWittenInvariantsComplex2024a}}]
\label{lem:4 moves}
    Every equivalence factors as a composition of the following maps or their inverses:
    \begin{enumerate}[(1)]
        \item \emph{Free quotient}: $\Kur\to \Kur/H$
        for a normal
        subgroup $H\subset G$ such that $H$ acts freely
        on $\T$. This replaces $\Kur$
        by $\Kur/H=(G/H,\T/H,
        \E/H,s)$.
        \item \emph{Group enlargement}: $\Kur\to G'\times_G\Kur$ for an embedding
        $G\hookrightarrow G'$.  This replaces $\Kur$ by $G'\times_G\Kur=(G',G'\times_G\T,G'\times_G\E,s')$.
        \item \emph{Stabilization}: $\Kur\to \mathcal{N}$ associated to a $G$-equivariant vector bundle $p\colon \Fcal\to \T$.
        This replaces $\Kur$ by
        $\mathcal{N}=(G,\Tot(\Fcal),\E\oplus
        p^*\Fcal,p^*s\oplus\id_{\Fcal})$
        where $p\colon \Fcal\to \T$ is the
        projection map from $\Fcal$ to $\T$.
        \item \emph{Germ equivalence}: $\Kur'\hookrightarrow \Kur$ for a $G$-equivariant open neighborhood $U\supset s^{-1}(0)$
        inside $\T$. This replaces
        $\Kur$ by $\Kur'=(G,U,\E|_{U},s|_{U})$.
    \end{enumerate}
\end{lemma}
\begin{remark}
\label{rmk:split the stably almost complex condition}
    Let $\Kur$ be a smooth stably almost complex Kuranishi chart. After
    performing a stabilization and germ equivalence, we can assume that $T[\T/G]$ and $[\E/G]$ both admit stably almost complex structures.
    We use such representatives below.
\end{remark}
Besides equivalences, there is a
relation between global Kuranishi charts inducing the same virtual fundamental classes (see \cref{lem:invariance of vfc under cobordism}),
called ``cobordant''.
\begin{definition}[{\cite[Definition 3.8]{abouzaidGromovWittenInvariantsComplex2024a}}]
\label{def:cobordism of Kuranishi charts}
    Two global Kuranishi charts $\Kur_i=(G,\T_i,\E_i,s_i),i=1,2$
    are cobordant iff there is
    a $G$-equivariant cobordism $\mathcal{W}$ from $\T_1$ to $\T_2$ carrying
    a $G$-equivariant vector bundle $\E\to \W$ with
    a section $s$ whose restriction
    to $\T_i$ is isomorphic to $(\E_i,s_i)$.
    We require this cobordism to be \emph{proper}, meaning that the total
    zero locus $s^{-1}(0)\subset\W$ is compact.
\end{definition}
\begin{remark}
    The uniqueness of global Kuranishi charts
    for Gromov--Witten theory of a symplectic manifold $(X,\omega)$ is in the sense
    that any two charts (produced from
    certain constructions) are cobordant.
\end{remark}
 
Finally, let $E$ be a Morava $K$-theory or HKR theory. We recall the definition of virtual fundamental classes of a global Kuranishi chart
in $E$-homology. Let $\Kur=(G,\T,\E,s)$ be a smooth
global Kuranishi chart such that $s^{-1}(0)$ is compact and
that $T[\T/G]$ and $[\E/G]$ are both stably almost complex (\cref{rmk:split the stably almost complex condition}).
Let $X$ be a closed complex (or stably almost complex) manifold,
equipped with the trivial $G$-action. 
Assume there is a continuous $G$-equivariant map $f\colon
\T\to X$. We consider the following maps
\begin{itemize}
    \item Poincar\'e duality gives
    \[
        H_G^*(\T,E)\xrightarrow{\lambda_{G,\T}^{-1}}
        H_{\dim\T-\dim G-*}^{G,\mathrm{lf}}(\T;E).
    \]
    \item Let $Z=s^{-1}(0)$. Since $Z$ is compact, restriction followed
    by the map of pairs induced by $s$ gives
    \[
    \begin{aligned}
        H_{\dim\T-\dim G-*}^{G,\mathrm{lf}}(\T;E)
        &\xrightarrow{\mathrm{res}}
        H_{\dim\T-\dim G-*}^{G}(\T,\T\setminus Z;E)\\
        &\xrightarrow{s_*}
        H_{\dim\T-\dim G-*}^{G}(\E,\E\setminus\T;E).
    \end{aligned}
    \]
    If $\T$ has boundary, we require $Z\subset\operatorname{int}\T$.
    \item Thom isomorphism $H^G_{\dim\T-\dim G-*}(\E,\E-\T,E)
    \simeq H^G_{\dim\T-\dim G-\rk \E-*}(\T,E)$.
    \item $f_*\colon H^G_{\dim\T-\dim G-\rk \E-*}(\T,E)\to H_{\dim\T-\dim G-\rk \E-*}(X,E)$. It is the map induced by $ 
    EG\times_G\T\to \T/G\xrightarrow{f}X$.
    \item Another Poincare duality
    \[
        \lambda_X\colon
        H_{\dim\T-\dim G-\rk \E-*}(X,E)\to H^{\dim X+*-(\dim\T-\dim G-\rk \E)}(X,E).
    \]
\end{itemize}
Notice that $\dim\T-\dim G-\rk \E=\vdim\Kur$.
\begin{definition}[{\cite[Definition 5.23]{abouzaid2021complexcobordismhamiltonianloops}}]
\label{def:vfc as an pushforward}
    Define \[f_!^\Kur\colon H^*_{G}(\T,E)\to H^{\dim X-\vdim\Kur+*}(X,E)\] to be the composition of previous maps. 
    In particular, take $X=\pt$ and let $a\colon \T\to \pt$
    be the unique map. We call the map
    $a_!^\Kur\colon H^*_G(\T,E)\to H^{*-\vdim \Kur}(\pt,E)=E^{*-\vdim \Kur}$
    the virtual fundamental class of the global Kuranishi chart
    $\Kur$.
\end{definition}
We discuss the invariance of $f_!^\Kur$ and virtual fundamental
classes under moves described in \cref{lem:4 moves} and
the cobordism in \cref{def:cobordism of Kuranishi charts}.
\begin{lemma}
\label{lem:invariances of vfc under 4 moves}
    Let $\Kur=(G,\T,\E,s),\Kur'=(G',\T',\E',s')$ be smooth and satisfy \cref{rmk:split the stably almost complex condition}, $X$ stably almost complex,
    and $f\colon \T\to X$ be $G$-equivariant. If $\varphi\colon \Kur'\to \Kur$ is an equivalence of Kuranishi charts, then we have
    a commutative diagram
    \[\begin{tikzcd}
	{H^*_{G'}(\T',E)} & {H^*_{G}(\T,E)} \\
	& {H^{*+\dim X-\vdim \Kur}(X,E)}
	\arrow["{f'^{\Kur'}_!}"', from=1-1, to=2-2]
	\arrow["{\varphi^*}"', from=1-2, to=1-1]
	\arrow["{f^\Kur_!}", from=1-2, to=2-2]
\end{tikzcd}\]
\end{lemma}
\begin{proof}[Sketch of proof]
By \cref{lem:4 moves},
it suffices to consider the four simple moves therein.
The free-quotient and group-enlargement cases follow from the
change-of-group compatibility in 
\cite[Addendum 5.19]{abouzaid2021complexcobordismhamiltonianloops}. The stabilisation and open-embedding
cases are \cite[Lemmas 5.27 and 5.25]{abouzaid2021complexcobordismhamiltonianloops}, respectively.
See also
\cite[Lemma 3.18]{abouzaidGromovWittenInvariantsComplex2024a}
for the functoriality statement for Morava-local
counting theories.
\end{proof}

\begin{lemma}
\label{lem:invariance of vfc under cobordism}
Let $\Kur_i=(G,\T_i,\E_i,s_i)$, $i=1,2$, be smooth stably almost
complex Kuranishi charts, and let $f_i\colon \T_i\to X$ be $G$-equivariant smooth maps to a closed stably almost complex manifold $X$. Suppose that $\Kur_1$ and $\Kur_2$ are joined by a proper stably almost complex cobordism $(\W,\E)$ in the sense of \cref{def:cobordism of Kuranishi charts}, and that the $f_i$ extend to a $G$-equivariant smooth map $f\colon\W\to X$. Then, for every $a\in H_G^*(\W;E)$,
\[
(f_1)^{\Kur_1}_!(a|_{\T_1})=(f_2)^{\Kur_2}_!(a|_{\T_2}).
\]
\end{lemma}

\begin{proof}
Using collars of the two boundary components, we may modify the extension $f$, without changing $f_1$ or $f_2$, so that it is constant in the collar direction near $\T_1$ and $\T_2$.

Let $\widetilde{\W}=\W\cup_{\partial\W}\W$ be the double of $\W$. By the product structure near the boundary, the Kuranishi data and the stably almost complex structure double to a boundaryless stably almost complex Kuranishi chart $\widetilde{\Kur}$ on $\widetilde{\W}$. Its zero locus is compact, being the double of the compact total zero locus $s^{-1}(0)$. Let $r\colon\widetilde{\W}\to\W$ be the folding map. The map $f$ therefore doubles to $\widetilde f=f\circ r\colon\widetilde{\W}\to X$, and the class $a$ has a canonical extension $\widetilde a:=r^*a\in H_G^*(\widetilde{\W};E)$. Since $r|_{\T_i}=\id_{\T_i}$, we have $\widetilde a|_{\T_i}=a|_{\T_i}$.

Choose a smooth cobordism function $\pi\colon\W\to[0,1]$ which is the collar coordinate near the boundary and satisfies $\pi^{-1}(0)=\T_1$ and $\pi^{-1}(1)=\T_2$. Doubling $\pi$ gives a smooth map $\widetilde\pi\colon\widetilde{\W}\to S^1$ and points $q_1,q_2\in S^1$ such that $\widetilde\pi^{-1}(q_i)=\T_i$. Set $F=(\widetilde f,\widetilde\pi)\colon \widetilde{\W}\to X\times S^1$. Equip $S^1$ with its standard stable complex structure, so that $X\times S^1$ is a closed stably almost complex manifold.

By \cite[Lemma 4.5]{abouzaid2021complexcobordismhamiltonianloops}, after stabilizing $\widetilde{\Kur}$ we may assume that $F$ extends to a smooth $G$-equivariant submersion. The push-forward is unchanged by this stabilization by \cite[Lemma 5.27]{abouzaid2021complexcobordismhamiltonianloops}. Applying \cite[Lemma 5.34]{abouzaid2021complexcobordismhamiltonianloops} to the submanifolds $X\times\{q_i\}\subset X\times S^1$ therefore gives
\[
j_i^*F^{\widetilde{\Kur}}_!(\widetilde a)=(f_i)^{\Kur_i}_!(a|_{\T_i}), \quad i=1,2,
\]
where $j_i\colon X\cong X\times\{q_i\}\hookrightarrow X\times S^1$. Here the Kuranishi chart obtained by restricting the stabilized double to $X\times\{q_i\}$ is a stabilization of $\Kur_i$ near its zero locus, and \cite[Lemma 5.27]{abouzaid2021complexcobordismhamiltonianloops} identifies its push-forward with that of $\Kur_i$.

Finally, $j_1$ and $j_2$ are homotopic, hence $j_1^*=j_2^*\colon H^*(X\times S^1;E)\to H^*(X;E)$. It follows that
\[
(f_1)^{\Kur_1}_!(a|_{\T_1})
= j_1^*F^{\widetilde{\Kur}}_!(\widetilde a)
= j_2^*F^{\widetilde{\Kur}}_!(\widetilde a)
= (f_2)^{\Kur_2}_!(a|_{\T_2}),
\]
as required.
\end{proof}
Since the AMS construction of
global Kuranishi charts for Gromov--Witten moduli space
produces a ``unique up to equivalence and cobordism'' Kuranishi chart
with a stably almost complex structure, \cref{lem:invariances of vfc under 4 moves,lem:invariance of vfc under cobordism} together
prove that the $E$-cohomological virtual fundamental class for Gromov--Witten moduli space is well-defined regardless of choices
involved in constructions.

\subsection{The AMS construction of global Kuranishi charts}
\label{non-equivariant global Kuranishi chart construction}
Let $(M,\omega)$ be a closed symplectic manifold\footnote{In original constructions of global Kuranishi charts, $M$ is assumed to be closed, but actually we only require that
Gromov compactness holds}. Let $g,h\ge 0$,
and $\beta\in H_2(M,\Z)$ be a homology class.
Choosing a compatible almost complex structure
$J$ on $(M,\omega)$,
there is a topological space of 
``$J$-holomorphic stable maps from genus $g$ nodal (complex) curves with
$h$ marked points (to $M$) whose image is in homology class $\beta$'', denoted by $\Mb_{g,h}(M,J,\beta)$. It admits
the following structures.
\begin{itemize}
    \item As a topological space, $\Mb_{g,h}(M,J,\beta)$ is equipped with the Gromov topology. 
$\Mb_{g,h}(M,J,\beta)$ is compact by Gromov compactness. 
However, the problem of $\Mb_{g,h}(M,J,\beta)$ is that
it is generally not a manifold
because the Cauchy–Riemann equation is not transverse and because the domains vary and may have automorphisms.
    \item Besides the topological space structure, $\Mb_{g,h}(M,J,\beta)$ admits such a structure (compare this with orbifolds): for
a point $(\Sigma,u)\in \Mb_{g,h}(M,J,\beta)$, there's a group 
\[
    \Aut(\Sigma,u)=\{\phi\in\Aut(\Sigma)\colon u\circ \phi=u\}
\]
which is finite, as required by the stability condition, and it
is called the \emph{automorphism group} of the stable map $(\Sigma,u)$.
It is also called the \emph{stabilizer group}
of the point in $\Mb_{g,h}(M,J,\beta)$.
\end{itemize}

The AMS construction replaces this singular moduli space by a global Kuranishi chart respecting the two structures.
see \cite{abouzaid2021complexcobordismhamiltonianloops,abouzaidGromovWittenInvariantsComplex2024a}.
We write $\mathfrak a$ for a collection of admissible auxiliary choices in the AMS construction, satisfying the positivity and largeness conditions required therein. The individual components of $\mathfrak a$ will be introduced as they enter the construction below.
Note that in the following theorem, $J$ is still
present as part of the input.
\begin{theorem}
\label{thm:global Kuranishi charts existence and invariance}
    Let $(M,\omega,J,g,h,\beta)$ be the geometric input, and let $\mathfrak{a}$ be an admissible auxiliary choice.
    One can construct a smooth stably almost complex global Kuranishi chart
    $\Kur_\afrak=(G_\afrak,\T_\afrak,\E_\afrak,s_\afrak)$ together
    with an identification $s_\afrak^{-1}(0)/{G_\afrak}
    \simeq \Mb_{g,h}(M,J,\beta)$ which preserves
    stabilizer groups.

    For two admissible choices $\afrak_0$ and $\afrak_1$ with $J$ fixed, the resulting global Kuranishi charts are equivalent. Moreover, if $J_0$ and $J_1$ are joined by a path of $\omega$-compatible almost complex structures, the corresponding global Kuranishi charts are joined by a stably almost complex cobordism.
\end{theorem}
As a consequence, given only a symplectic manifold
$(M,\omega)$ (without specifying an almost complex structure $J$), there is a
unique equivalence/cobordism class of global
Kuranishi charts, which allows definition
of Gromov--Witten invariants.

\subsubsection{The AMS construction}
We now describe the construction in \cite{abouzaidGromovWittenInvariantsComplex2024a} of $\Kur_{\mathfrak a}$ in enough detail to explain the geometry of the thickening space $\mathcal T_{\mathfrak a}$, its smooth and stably almost complex structure, and the comparison of charts arising from different auxiliary choices. 

The construction has three main ingredients: a rigidification of the domains by projective embeddings, a finite-dimensional thickening of the Cauchy--Riemann equation, and a normalization of the projective framing. We begin with the first.

\paragraph{Projective parametrization of domains.}
\label{sec:projective parametrization of domains}
The first step is to place the varying domains of stable maps in a fixed finite-dimensional family. For this purpose, AMS use stable maps to projective space to give a projective presentation of the domain curves.

Fix for the moment an integer $d$, which will be determined later from the auxiliary line-bundle data. Let
$\Fcal=\Fcal_{g,h,d}\subset\Mb_{g,h}\left(\mathbb P^{d-g},d\right)$
be the locus of stable maps $\phi\colon \Sigma\to \mathbb P^{d-g}$ such that 
\begin{itemize}
    \item $H^1\left(\Sigma,\phi^*\O_{\P^{d-g}}(1)\right)=0$ and 
    \item $(\Sigma,\phi)$ has trivial automorphism group.
\end{itemize}
Let $\Ccal=\Ccal_{g,h,d}\to  \mathcal F$
be the universal curve.

These conditions imply that $\Fcal$ and $\Ccal$ are nonsingular quasi-projective varieties
\footnote{Terminology warning: the curve family $\Ccal\to \Fcal$
is not smooth because there are nodes}.
A point $\phi\in\Fcal$ determines a projectively presented marked curve $\Sigma_\phi:=\mathcal C|\phi$,
and the maps entering the thickened Gromov--Witten problem will be maps
$u\colon \Sigma_\phi\to M$.
In this way, the abstract varying domains are replaced by fibres of the fixed finite-dimensional family $\mathcal C\to\mathcal F$.

Changing homogeneous coordinates on $\P^{d-g}$ gives an action of
$G_{\C}=\GL(d-g+1,\C)$
on $\Fcal$ and $\Ccal$. Its maximal compact subgroup 
$G=U(d-g+1)$
will become the structure group of the global Kuranishi chart. 

Similarly, we define $\Fcal^{2}=\Fcal^{(2)}_{g,h,d}$ consisting of
 $u,v \colon\Sigma\to \P^{d-g}$ such that both $(\Sigma,u)$ and $(\Sigma,v)$ are in $\Fcal_{g,h,d}$. Let $\Delta\colon \Fcal\to \Fcal^{(2)}$
be the diagonal embedding. $\Fcal^{(2)}$
admits a natural $\G_\C$-action, too.
Define a vector bundle $H_{g,h,d}\to \Fcal_{g,h,d}$ with fiber
$(H_{g,h,d})_\phi=H^0\left(\Sigma_\phi,\phi^*T
\P^{d-g}\right)$, which can be identified
with the normal bundle $N_{\mathcal F_{g,h,d}/\mathcal F^{(2)}_{g,h,d}}$
Choose a $G$-equivariant tubular neighborhood $U_{g,h,d} \supset \Delta(\mathcal F_{g,h,d})$
and identify it with a neighborhood of the zero section in $H_{g,h,d}$. 

Later, a map $u\colon\Sigma_\phi\to M$, together with a choice of framing, will determine a second projective map $\iota_F\colon\Sigma_\phi\to\mathbb P^{d-g}$.
Whenever $(\phi,\iota_F)$ lies in $U_{g,h,d}$, the tubular-neighborhood identification determines an element of $H_{g,h,d}$, which vanishes precisely when the two projective presentations agree. This will give the $H_{g,h,d}$-component of the obstruction section.

\paragraph{Finite-dimensional thickening of the Cauchy--Riemann equation.}
The next step is to replace the Cauchy--Riemann equation by perturbed equations, simultaneously over the varying domains. Let $\Ffrak_{\Fcal}$ denote the space of pairs $(\phi,u)$, where $\phi\in\Fcal_{g,h,d}$ and
$u\colon\Sigma_\phi\to M$
is a smooth map representing $\beta\in H_2(M,\Z)$ and satisfying the conditions of
\cite[Definition~4.6]{abouzaidGromovWittenInvariantsComplex2024a}.
In particular, $u$ is required to be $J$-holomorphic in a neighborhood of the nodes; the remaining conditions impose positivity and a uniform lower bound on the energy of the relevant irreducible components, which we will not need explicitly.

To perturb the Cauchy--Riemann equation uniformly over this family, consider the $G$-equivariant vector bundle
\[
Y
=
\Omega^{0,1}_{\Ccal^\circ_{g,h,d}/\Fcal_{g,h,d}}
\otimes_{\C} TM
\to
\Ccal^\circ_{g,h,d}\times M .
\]
For $(\phi,u)\in\Ffrak_{\Fcal}$, let
\[
\Gamma_u\colon\Sigma_\phi^\circ
\to
\Ccal^\circ_{g,h,d}\times M,
\qquad
z\longmapsto (z,u(z)),
\]
denote the graph of $u$. Pulling a section of $Y$ back along $\Gamma_u$ gives a $(0,1)$-form on $\Sigma_\phi^\circ$ with values in $u^*TM$, and hence a perturbation of $\bar\partial_Ju$.

Following \cite[Definition~4.1]{abouzaidGromovWittenInvariantsComplex2024a}, choose a $G$-equivariant finite-dimensional approximation scheme
$\lambda_\mu\colon V_\mu\to C^\infty_c(Y),
\mu\in\mathbb N.$
Here the $V_\mu$ form an increasing sequence of finite-dimensional $G$-representations, the maps $\lambda_\mu$ are compatible with the inclusions, and
$\bigcup_\mu \lambda_\mu(V_\mu)$
is dense in $C^\infty(Y)$ with respect to the $C^\infty_{\loc}$-topology. For a fixed sufficiently large $\mu$, the associated pre-thickened moduli space is
\[
\T^{\mathrm{pre}}
=
\left\{
(\phi,u,e),
(\phi,u)\in\Ffrak_{\Fcal},
\ e\in V_\mu,
\ \pb_Ju\big|_{\Sigma_\phi^\circ}
+\lambda_\mu(e)\circ\Gamma_u=0
\right\}.
\]
Thus the auxiliary variable $e$ records a finite-dimensional perturbation of the Cauchy--Riemann equation. Although the original mapping space is infinite-dimensional, the failure of transversality of $\bar\partial_J$ is measured by a finite-dimensional cokernel. A point $(\phi,u,e)$ is called regular when the perturbation directions coming from $V_\mu$ surject onto the cokernel of the corresponding linearized Cauchy--Riemann operator. For every compact $K\subset\Fcal_{g,h,d}$, sufficiently large $\mu$ makes every point of the locus
$\left\{
(\phi,u,0)\in\T^{\mathrm{pre}}
\ \middle|\
\phi\in K
\right\}$
regular; the uniform choice of $\mu$ follows from Gromov compactness
\cite[Lemma~4.38]{abouzaidGromovWittenInvariantsComplex2024a}.

Moreover, on the regular locus the forgetful map
\[
\T^{\mathrm{pre,reg}}
\to
\Fcal_{g,h,d},
\qquad
(\phi,u,e)\longmapsto\phi,
\]
admits a fiberwise smooth $C^1_{\loc}$-structure in the sense of
\cite[Definition~4.27]{abouzaid2021complexcobordismhamiltonianloops};
see \cite[Lemma~4.39 and Appendix~A]{abouzaidGromovWittenInvariantsComplex2024a} and \cite[Appendix~B]{Pardon_2016}.
We will return to the role of this relative regularity when discussing the smooth structure on the global chart.

Finally, the locus $e=0$ is precisely the solution set of the unperturbed
Cauchy--Riemann equation.
Consequently, in the eventual global Kuranishi chart the projection to $e\in V_\mu$ becomes one component of the obstruction section: its vanishing removes the auxiliary perturbation and recovers the original moduli problem.
\paragraph{Framings and normalization.}
It remains to add a framing which produces a second projective presentation of the domain and, after normalization, reduces the change-of-framing group from $G_{\C}$ to its maximal compact subgroup $G$.

Recall that $\Ffrak_{\Fcal}$ parametrizes admissible maps $(\phi,u)$, with universal curve $\Cfrak_\Fcal\to \Ffrak_\Fcal$.
Choose a $G_{\C}$-invariant consistent domain metric $\mathfrak D$ in the sense of
\cite[Definition~4.9]{abouzaidGromovWittenInvariantsComplex2024a};
its existence is proved in
\cite[Lemma~4.17]{abouzaidGromovWittenInvariantsComplex2024a}.

Let $\Omega$ be an integral symplectic form on $M$ taming $J$, and let
$L\to M$ be a Hermitian line bundle with curvature $-2\pi\i\Omega$.
As in \cite[Definition~4.23 and Remark~4.24]{abouzaidGromovWittenInvariantsComplex2024a},
we assume that $L$ admits a root of order at least $3$ and fix $k\geq 3$.
Together with $\mathfrak D$, these choices form the line-bundle data
$\mathbb L=(L,k,\mathfrak D).$
For $(\phi,u,e)\in\T^{\mathrm{pre}}$, define the line bundle on $\Sigma_\phi$
\begin{equation}
\label{eq:polarizing line bundle in AMS}
       L_{u,\mathbb L}
=
\left(
\omega_{\Sigma_\phi}(p_1+\cdots+p_h)
\otimes_{\C}u^*L
\right)^{\otimes k}.
\end{equation}
The connection on $L$ induces a holomorphic structure on $L_{u,\mathbb L}$.
The integer $d$, which was left unspecified in the ``projective parametrization of domains'', is now fixed by
\[
    d=d_{\mathbb L}
    :=k\bigl(\langle[\Omega],\beta\rangle+2g-2+h\bigr)=\deg L_{u,\mathbb L}.
\]
The positivity assumptions above ensure that $L_{u,\mathbb L}$ is sufficiently positive on every irreducible component; in particular,
\begin{itemize}
    \item $H^1(\Sigma_\phi,L_{u,\mathbb L})=0$;
    \item $\dim_{\C}H^0(\Sigma_\phi,L_{u,\mathbb L})=d-g+1$.
\end{itemize}
and the associated complete linear system gives a projective presentation of the domain
\cite[Remark~4.24]{abouzaidGromovWittenInvariantsComplex2024a}.

The consistent domain metric $\mathfrak D$, together with the Hermitian metric on $L$, induces a Hermitian metric
$h_{L_{u,\mathbb L}}$ on $L_{u,\mathbb L}$ and a volume form
$\Omega_{\mathfrak D,\phi}$ on $\Sigma_\phi$. Hence there is an $L^2$-pairing
\begin{equation}
\label{eq:definition of metric on framings}
    \langle f_1,f_2\rangle_{L^2}
=
\int_{\Sigma_\phi}
h_{L_{u,\mathbb L}}(f_1,f_2),
\Omega_{\mathfrak D,\phi}
\end{equation}
on $H^0(\Sigma_\phi,L_{u,\mathbb L})$.

A holomorphic $\mathbb L$-framing of $(\phi,u,e)$ is a complex basis
$F=(F_0,\ldots,F_{d-g})$
of $H^0(\Sigma_\phi,L_{u,\mathbb L})$ as in
\cite[Definition~4.25]{abouzaidGromovWittenInvariantsComplex2024a}.
Associated to $F$ is the Hermitian matrix
$H_F
=
\bigl(
\langle F_i,F_j\rangle_{L^2}
\bigr)_{0\leq i,j\leq d-g}$.
A framing is called unitary when $H_F=I$.
The $G_{\C}$-invariance of the consistent domain metric makes this normalization compatible with changes of framing
\cite[Lemma~4.26]{abouzaidGromovWittenInvariantsComplex2024a}.
Thus arbitrary framings differ by the action of
$G_{\C}=\GL(d-g+1,\C)$,
whereas unitary framings differ by $G=U(d-g+1)$.
In this way the condition $H_F=I$ reduces the noncompact change-of-framing ambiguity to the compact group used in the global Kuranishi chart.
\paragraph{The global Kuranishi chart.}
\begingroup
\setlength{\parfillskip}{0pt plus 1fil}
The following choices are made during the preceding three steps:
\par
\endgroup
\begin{itemize}
    \item The line bundle and framing data: $\Omega,L,k,\mathfrak{D}$ 
    ($d$ is determined by $L,k$).
    \item A finite-dimensional approximation scheme $(V_\mu,\lambda_\mu)$ and $\mu$ large.
\end{itemize}
As in \cref{thm:global Kuranishi charts existence and invariance},
we use $\mathfrak{a}$ to denote a choice. A choice $\afrak$
is called admissible iff it satisfies the following conditions:
\begin{itemize}
    \item $L$ has a root of order$\ge 3$ and $k\ge 3$.
    \item $\mu$ is large enough such that $\T^{\pre,\mathrm{reg}}$
    contains all solutions to the unperturbed equation.
\end{itemize}
Let's introduce the AMS construction:
\begin{definition}
\label{def: AMS construction}
    With admissible auxiliary data $\afrak$, the global Kuranishi chart $\Kur=(G,\T,\E,s)$ is defined as follows
    \begin{enumerate}[(1)]
        \item $G_\afrak=U(d-g+1)$
        \item $\T_\afrak$ consists of $(\phi,u,e,F)$ where 
        $(\phi,u,e)\in \T^{\pre,\mathrm{reg}}$ and $F$ is a framing
        such that
        \begin{itemize}
            \item $H_F$ has positive eigenvalues and lies 
            in a neighborhood of the identity matrix,\newline
            so $\exp^{-1}(H_F)$ is well-defined.
            \item $(\phi,\iota_F)\in \Fcal_{g,h,d}^{(2)}$ lies
        in $U_{g,h,d}$, where $\phi$ is the original projective
        parametrization of domain and $\iota_F$ is the parametrization
        induced by $F$.
        \end{itemize}
        $\T_\afrak$ is called the \emph{thickened moduli space}.
        \item $\E_\afrak|_{(\phi,u,e,F)}=H_{g,h,d}\oplus V_\mu
        \oplus\H$, where $H_{g,h,d}$ is the normal bundle of $\Fcal_{g,h,d}$
        in $\Fcal_{g,h,d}^{(2)}$, $V_\mu$ is part of the finite-dimensional approximation scheme, and $\H$ is the vector space of Hermitian matrices
        \item The section $s_\afrak$ is defined as follows
        \begin{itemize}
            \item $H_{g,h,d}$ entry: since $(\phi,\iota_F)\in U_{g,h,d}$,
            we take the section to be its image under tubular neighborhood
            identification $U_{g,h,d}\simeq H_{g,h,d}$.
            \item $V_\mu$ entry: send $(\phi,u,e,F)$ to $e$.
            \item $\H$ entry: send $(\phi,u,e,F)$ to $\exp^{-1}(H_F)$
        \end{itemize}
        \end{enumerate}
\end{definition}
Since $\T^{\pre,\mathrm{reg}}\to \Fcal$ admits a fiberwise smooth $C^1_{\loc}$ structure, and $\T_\afrak$ is an open subset
of the framing principal bundle over $\T^{\pre,\mathrm{reg}}$,
we know $\T_\afrak\to \Fcal$ also admits a fiberwise
smooth $C^1_{\loc}$ structure.
\subsubsection{Smooth structure of {\texorpdfstring{$\T$}{T}}}
Recall that in \cref{conv:global Kuranishi charts}, we want
global Kuranishi charts to be smooth and stably
almost complex. By contrast, $\T$ is only a
$C^1_\loc$ fiberwise smooth bundle over $\Fcal$.
To get a smooth thickened moduli space,
one needs smoothing theory. 

Let $G$ be a compact Lie group. We use the language of
microbundles; see \cite[Secion 4.2]{abouzaid2021complexcobordismhamiltonianloops}.
For a $G$-equivariant topological manifold $M$, 
denote $T_\mu M$ to be its tangent $G$-microbundle.
If $M$ is smooth, then $T_\mu M$ admits a ($G$-equivariant)
vector lift $TM\to T_\mu M$. Lashof proved the following
equivariant smoothing result:
\begin{proposition}[{\cite{Lashof1979-Stable_G-smoothing}}]
\label{prop:Lashof G-smoothing}
    Let $M$ be a $G$-equivariant 
    topological manifold. Assume that
    $M$ has finitely many $G$-orbit types.
    Assume that there's a $G$-equivariant
    vector bundle reduction $\varphi\colon E\to T_\mu M$
    for some $G$-equivariant vector bundle $E\to M$,
    then there is a $G$-representation $V$ and
    a $G$-equivariant smooth structure $\alpha$ on
    $M\times V$ satisfying the following:
    \begin{itemize}
        \item There is an isomorphism of $G$-equivariant
        vector bundles $p_M^* E\oplus \underline{V}\simeq T(M\times V)_\alpha$, where $p_M\colon M\times V\to M$
        is the projection.
        \item The natural lift $T(M\times V)_\alpha
        \to T_\mu (M\times V)$, when restricted to
        $M\times 0$, is isotopic through G-equivariant lifts to the induced lift
        $\varphi\oplus\id_V\colon p_M^*E\oplus \underline{V}
        \to T_\mu(M)\oplus \underline{V}$.        
    \end{itemize}
    Furthermore, $\varphi\mapsto \alpha$ gives a bijection
    between
    homotopy classes of $G$-vector bundle lifts of $T_\mu M$
    and 
    stable isotopy classes\footnote{
    Look at \cite{Lashof1979-Stable_G-smoothing} or \cite[Definition 6.8]{bai2022arnoldconjectureintegers} for terminology involved.} of
    $G$-smoothings. 
\end{proposition}

Let $\Kur=(G,\T,\E,s)$ be a global
Kuranishi chart from AMS construction, where we omit
the auxiliary data $\afrak$ from notation.
$\T$ satisfies the premise of \cref{prop:Lashof G-smoothing}:
\begin{itemize}
    \item $\T$  has only finitely many $G$-orbit types
    possibly after a germ equivalence. It is a consequence of stable condition 
and Gromov-compactness that we can shrink
$\T$ so that $\T$ has only finitely
many $G$ orbit types.
\item $T_\mu \T$ admits a $G$-equivariant vector bundle reduction.
Since $\pi\colon \T\to \Fcal$ is a $C^1_\loc$ fiberwise smooth
bundle, the vertical tangent bundle $T^{vt}\T$ is well-defined
as a $G$-equivariant topological vector bundle. 
The vector bundle reduction is $\pi^*T\Fcal\oplus T^{vt} \T
\to T_\mu \T$. See details in \cite[Corollary~4.26 and Lemma~4.29]{abouzaid2021complexcobordismhamiltonianloops}.
\end{itemize}
\begin{remark}
    It is implicit in \cite{abouzaid2021complexcobordismhamiltonianloops,abouzaidGromovWittenInvariantsComplex2024a} that the $C^1_\loc$-structure from Pardon's gluing construction is compatible with the $G$-action.
    \cite{hirschiGlobalKuranishiCharts2024b}
    proved this $G$-equivariance in its section 5.4 (the exact statement is Proposition 5.1(c)).
\end{remark}
Therefore, one can conclude that:
\begin{proposition}[{\cite[Proposition 4.30]{abouzaid2021complexcobordismhamiltonianloops}}]
\label{prop:stable smoothing of global Kuranishi charts}
    Let $\Kur := (G, \T, \E, s)$ be a global Kuranishi chart for which $\T$ is the total space of a fiberwise smooth $C^1_\loc$ 
    $G$-bundle whose base $B$ is a smooth $G$-manifold.
    Also, suppose that the $G$-action on $\T$ has only finitely many orbit types.
    Then there is a stabilization $\Kur'$ of $\Kur$ which admits a smooth structure. The thickening of $\Kur'$ is $\T \times V$ for some $G$-representation $V$ and the restriction to $\T\times 0$ of the natural lift $T (\T\times V) \to T_\mu (\T\times V)$
    is isotopic through $G$-equivariant lifts to the induced lift 
    $l \oplus \id_V\colon  \E\oplus V \to  T_\mu\T \oplus V$ .
\end{proposition}
For two stabilization smoothings $\Kur_i=(G,\T_i,\E_i,s_i), i=1,2$ (where
$\T_i=\T\times V_i$ topologically) from \cref{prop:stable smoothing of global Kuranishi charts}, \cref{prop:Lashof G-smoothing}
argues that $\T_1,\T_2$ are $G$-equivariantly stably isotopic.
This means that after a further stabilization respectively,
we can assume that there is a homeomorphism of $G$-equivariant
\textit{topological} manifolds $\T_1\to \T_2$ and that
the $G$-smooth structures on $\T_1,\T_2$ are isotopic.
Hence, $\Kur_1,\Kur_2$ (after further stabilizations)
are related by a cobordism of global Kuranishi charts,
where $\W$ is topologically $\T_1\times [0,1]$ whose
smooth structure is given by the isotopy of smooth structures from
$\T_1$ to $\T_2$.

\subsubsection{Stably almost complex structure on
AMS construction}
Next, we argue that there is a canonical stably almost complex
structure on global Kuranishi charts from the AMS construction.
Recall that a stably almost complex structure on a \textit{smooth} global Kuranishi chart $\Kur=(G,\T,\E,s)$ is an equivalence
class of stable complex structures on the $G$-equivariant virtual
vector bundle $T\T-\mathfrak{g}-\E$. In the case of AMS construction,
we can define this before stabilization smoothing: we instead
consider $T^{vt}\T+\pi^*T\Fcal-\mathfrak{g}-\E$, where
$\pi\colon \T\to \Fcal$ is the natural projection, and $T^{vt}\T$
is the vertical tangent space. 

Write $N=d-g+1,G=U(N),\mathfrak{g}=\mathfrak{u}(N)$, and denote
$\H=\mathrm{Herm}_N$ the $N\times N$ Hermitian matrices as
in \cref{def: AMS construction}. Let $q\colon \T\to \T^\pre$ be
the map forgetting the framing. Then we have identities in $G$-equivariant
real $KO$-theory $\KO_G(\T)$:
\begin{equation}\label{eq:K-theory identity 1}
    T\T=q^*T(\T^\pre)+\underline{\mathfrak{gl}_N(\C)}=
    q^*T^{vt}(\T^\pre)+\pi^*T\Fcal+\underline{\mathfrak{gl}_N(\C)}
\end{equation}
Also, from \cref{def: AMS construction}
\begin{equation}\label{eq:K-theory identity 2}
    \E=\pi^* H_{g,h,d}+\underline{V_\mu}+\underline{\H}
\end{equation}
Combining \eqref{eq:K-theory identity 1}, \eqref{eq:K-theory identity 2},
and $\mathfrak{gl}_N(\C)=\mathfrak{g}\oplus \H$,
we get 
\begin{lemma}
\label{lem:computation of virtual deformation theory}
    In $\KO_G(\T)$ there is a distinguished isomorphism
    \[
        T\T-\mathfrak{g}-\E=q^*(T^{vt}\T^\pre-\underline{V_\mu})+\pi^*(T\Fcal-H_{g,h,d}).
    \]
\end{lemma}
The virtual bundle $T\Fcal-H_{g,h,d}$ is already complex.  It
therefore remains to construct a canonical $G$-equivariant stable
complex structure on $T^{vt}\T^\pre-\underline{V_\mu}$.

Fix a point $x=(\phi,u,e)\in\T^\pre$ with domain $\Sigma=\Ccal|_\phi$,
and let $\widetilde{\Sigma}$ be its normalization.
Since the projection
$\T^\pre\to\Fcal$ fixes the domain in each fibre, its vertical
deformation theory is governed by the linearisation of the perturbed
Cauchy--Riemann equation
\[
    \pb_Ju+\lambda_\mu(e)\circ\Gamma_u=0.
\]
The linearization of this operator is denoted by 
\[  
    \widehat{D}_x\colon \Omega^0(\Sigma,u^*TM)\oplus V_\mu\to 
    \Omega^{0,1}(\widetilde{\Sigma},u^*T M).
\]
We denote its restriction to $\Omega^0(\Sigma,u^*TM)$ by
\[
    D_x\colon \Omega^0(\Sigma,u^*TM)\to 
    \Omega^{0,1}(\widetilde{\Sigma},u^*T M).
\]
By assumptions of $V_\mu$ in \cref{def: AMS construction},
we may assume (by shrinking $\T^\pre$) $\widehat{D}_x$ is surjective $\forall x\in \T^\pre$. 

If we fix the domain type $\Sigma=\Sigma_\phi$, i.e. consider the fiber of $\T^\pre\to \Fcal$ over $\phi$, then by
standard Fredholm theory, solutions to perturbed Cauchy--Riemann equation
for a fixed $\Sigma$
form a smooth manifold whose tangent bundle can be identified with $\ker \widehat{D}_x$.
This is the reason that each fiber of $\T^\pre\to \Fcal$ is a smooth
manifold and $T^{vt}\T^\pre\simeq \ker(\widehat{D}_x)$ on
a single fiber. Therefore, $T^{vt}\T^\pre-\underline{V_\mu}=\ind_\R(D_x)$,
where $\ind_\R$ denotes $G$-equivariant family Fredholm
index, but $\ind_\R(D_x)$ only makes sense in $\KO_G(\T^\pre|_{G\cdot\phi})$ for $\phi\in \Fcal$.
 
The preceding description does not by itself give a global family
index across the nodal strata, since the Sobolev spaces vary when the
topological type of the domain changes.  The required compatibility 
is supplied by gluing results, like from \cite[Appendix B]{Pardon_2016}.

For a point $x_0=(\phi_0,u_0,e_0)$ and let $K_0=\ker(\widehat{D}_{x_0})$,
Pardon constructed a local parametrization of $\T^\pre$ near $x_0$ 
(\cite[Definition B.10.1]{Pardon_2016}):
\[
    U\times K_0 \to \T^\pre|_{U}
\]
where $U\subset \Fcal$ is some open neighborhood of $\phi_0$.
Working in the direction of \cite[Corollary 6.29]{abouzaid2021complexcobordismhamiltonianloops},
one can prove that for another trivialization
$U'\times K_0'\to \T^\pre|_V$, on their overlapping
part they give a map
$U\cap U' \to C^\infty(V\subset K_0,K_1)$ which becomes continuous if we use $C^1_{\loc}$-topology on RHS. This
is why $\T^\pre\to \Fcal$ admits a $C^1_\loc$ fiberwise smooth bundle
structure, and that $\ker(\widehat{D}_x)$
can be glued to a $G$-equivariant topological vector bundle
over $\T^\pre$.

The following route is the same as
\cite[Section 5.5]{hirschiGlobalKuranishiCharts2024b}, which has
its roots in \cite[]{Pardon_2016}.
For $D_x$
over a single $x$, we define its
complexification by 
\[  
    D^\C=\frac{D-J\circ D\circ J}{2}
    \colon \Omega^0(\Sigma,u^*TM)
    \to \Omega^{0,1}(\widetilde{\Sigma},u^*TM)
\]
and define a homotopy 
\[
    L_x(t)=t D_{x}+(1-t) D^\C_x
    =D^\C_x+t A_x,
\]
where $A_x= D_x-D^\C_x$ is a zeroth-order complex-antilinear differential operator.
By density of the approximation scheme, compactness of
$q(s^{-1}(0))\times[0,1]$, and openness of surjectivity under linear
gluing, we may choose $\nu\geq\mu$ so that adding $V_\nu$ makes
$L_x(t)$ surjective near this compact set. The linear gluing input is
essential here; see \cite[Section~5.5, footnote~14]{hirschiGlobalKuranishiCharts2024b}
and \cref{app:Proof of vector bundle structure of K}.
To ensure complex linearity of the augmented operator at $t=0$, set
\[
    E=V_\nu\otimes_\R\C,\qquad
    \lambda(v+iw)=\lambda_\nu(v)+J\lambda_\nu(w)
    \quad (v,w\in V_\nu).
\]
Thus $E$ is a complex $G$-representation and
$\lambda\colon E\to C^\infty_c(Y)$ is complex-linear and
$G$-equivariant. The real inclusion $V_\mu\subset V_\nu\subset E$
satisfies $\lambda|_{V_\mu}=\lambda_\mu$.
After shrinking $\T^\pre$ around $q(s^{-1}(0))$, the operators
\[
    \widehat{L}_x(t)\colon 
    \Omega^0(\Sigma,u^*TM)\oplus E\to 
    \Omega^{0,1}(\widetilde{\Sigma},u^*TM)
\]
\[
    \widehat{L}_x(t)(\xi,\dot e)=L_x(t)\xi+\lambda(\dot e)
    \circ \Gamma_u
\]
are surjective for every $(x,t)\in \T^\pre\times [0,1]$.
Indeed, complexification retains the surjective real augmentation.
At $t=1$ this operator restricts to $\widehat D_x$ on
$\Omega^0(\Sigma,u^*TM)\oplus V_\mu$.
\begin{proposition}[{\cite[Section 5.5]{hirschiGlobalKuranishiCharts2024b}}]
\label{prop:vector bundle structure of K}
    After shrinking $\T^\pre$ if necessary, the disjoint union
    \[
        \mathcal{K}:=\bigsqcup_{(x,t)\in \T^\pre\times [0,1]} K_{x,t},
        \quad K_{x,t}:=\ker \widehat{L}_{x}(t)
    \]
    admits a natural structure of a finite-rank $G$-equivariant topological real vector bundle over $\T^\pre\times[0,1]$.

    Moreover, the fiberwise complex structure on
    $\mathcal{K}_0=\mathcal{K}|_{\T^\pre\times\{0\}}$ induced by the
    complex-linear augmented operator $\widehat L_x(0)$ is continuous.
    Hence $\mathcal{K}_0-\underline E$ is a complex virtual $G$-bundle.
\end{proposition}
The proof strategy of \cref{prop:vector bundle structure of K}
is largely discussed throughout \cite[Section 5.5]{hirschiGlobalKuranishiCharts2024b}. See \cref{app:Proof of vector bundle structure of K} for discussion.

For the complex assertion, use evaluations at finitely many smooth points
of the domain together with the $E$-coordinate, and choose complex-linear
quotients giving an isomorphism on the kernel at the reference point.
With $J$-linear parallel transport these remain continuous, complex-linear
local trivializations near that point, as in
\cite[Lemma~5.11]{hirschiGlobalKuranishiCharts2024b}.

Write $\mathcal K_1=\mathcal K|_{\T^\pre\times\{1\}}$.
The compatibility $\lambda|_{V_\mu}=\lambda_\mu$ gives an exact sequence
of real $G$-vector bundles
\[
    0\longrightarrow T^{vt}\T^\pre
    \longrightarrow\mathcal K_1
    \xrightarrow{\ (\xi,a)\mapsto[a]\ }
    \underline{E/V_\mu}\longrightarrow0.
\]
The kernel is $\ker\widehat D=T^{vt}\T^\pre$.
For any $a\in E$, surjectivity of $\widehat D_x$ supplies
$(\xi,b)$ with $b\in V_\mu$ and
$\widehat D_x(\xi,b)=-\lambda(a)\circ\Gamma_u$;
then $(\xi,a+b)\in K_{x,1}$ maps to $[a]$.
The inclusion and projection are continuous in the gluing evaluation
coordinates and are $G$-equivariant. Consequently,
\[
    [\mathcal K_1]-[\underline E]
    =[T^{vt}\T^\pre]-[\underline{V_\mu}]
    \quad\text{in }\KO_G(\T^\pre).
\]
A $G$-equivariant splitting of the exact sequence and transport along
the bundle $\mathcal K$ over $[0,1]$ identify this real virtual bundle
stably with the underlying real bundle of
$\mathcal K_0-\underline E$. This supplies the required stable complex
structure, using the complex endpoint and its continuous transport,
not merely the equality of real $KO$-theory classes.
    
\subsubsection{Independence of choices}
We briefly recall why the preceding construction is independent of the
auxiliary data.  We refer to
\cite[Sections~4.7--4.9]{abouzaidGromovWittenInvariantsComplex2024a}
for the detailed comparison arguments.

Let $\afrak_0$ and $\afrak_1$ be two admissible choices with the almost
complex structure $J$ fixed.  The basic comparison device is to form a
common ``doubly framed'' thickening carrying simultaneously the
projective framings and finite-dimensional approximation data associated
to $\afrak_0$ and $\afrak_1$.  After group enlargements, stabilizations
and germ equivalences, both $\Kur_{\afrak_0}$ and $\Kur_{\afrak_1}$
are equivalent to this common chart.  Thus
\[
    \Kur_{\afrak_0}\simeq \Kur_{\afrak_1}.
\]
This is the content of
\cite[Section~4.9, in particular Proposition~4.66 and
Corollary~4.70]{abouzaidGromovWittenInvariantsComplex2024a}.
See also \cite[Section~6.1]{hirschiGlobalKuranishiCharts2024b}
for the analogous comparison by a doubly thickened chart.

If $J_t$, $t\in[0,1]$, is a path of compatible almost complex
structures, the same construction can be carried out parametrically.
After making the auxiliary choices uniformly over the compact
parametrized moduli space, one obtains a stably almost complex global
Kuranishi cobordism whose restrictions at $t=0,1$ are global charts for
$J_0$ and $J_1$; see
\cite[Section~4.7 and Corollary~4.70]
{abouzaidGromovWittenInvariantsComplex2024a}
and \cite[Section~6.2]{hirschiGlobalKuranishiCharts2024b}.

Consequently, by
\cref{lem:invariances of vfc under 4 moves,lem:invariance of vfc under cobordism},
the resulting virtual push-forward and virtual fundamental class are
independent of all the auxiliary choices and of $J$.

\subsection{Equivariant generalization}
When the symplectic manifold $M$ admits a Hamiltonian $K$-action,
we expect the $K$-action lifts to the moduli space $\Mb_{g,h}(M,\beta)$. The organization of this section is similar
to previous ones. First, we will introduce general theories of
equivariant global Kuranishi charts, especially
defining
virtual fundamental classes in this case. The adaptation
of equivariant localization from orbifolds to equivariant
global Kuranishi charts is discussed.
We then describe the $K$-action on the AMS chart of \cref{def: AMS construction}
(cf.~\cite[Lemma~5.14]{Hirschi_Properties_of_Gromov_Witten}).
\subsubsection{General theory of equivariant global Kuranishi charts and equivariant localization}
Let $K$ be a compact Lie group. Like \cref{def:working definition of group action on orbifold}, we can write a working definition
of $K$-equivariant global Kuranishi charts.
\begin{definition}[{\cite[Definition 5.1]{Hirschi_Properties_of_Gromov_Witten}}]
\label{def:equivariant global Kuranishi chart}
    A $ K$-equivariant global Kuranishi chart consists of the following data:
    \begin{enumerate}[(1)]
         \item An extension of compact Lie groups
        $1\to G\to \widetilde{G}\to K\to 1$.
        \item $\T$ is a $\widetilde{G}$-equivariant smooth manifold, and
        $\E\to \T$ is a $\widetilde{G}$-equivariant vector bundle.
        \item A $\widetilde{G}$-equivariant section $s\colon \T\to \E$.
    \end{enumerate}
    We require that, after forgetting the $K$-action, $(G,\mathcal T,\mathcal E,s)$ is an ordinary smooth global Kuranishi chart
    as in \cref{def:global Kuranishi chart}. In particular, this means the $G$-action on $\T$ has only finite stabilizers.
\end{definition}
\begin{remark}
    Recall that in \cref{conv:global Kuranishi charts}, we want to endow
    global Kuranishi charts with smooth and stably almost complex structures, and require that $s^{-1}(0)$ is compact.
    In equivariant case, 
    compactness
    and smoothness are easier to formulate. 
    For stably almost complex
    structure, we see $3$ vector bundles on $\T$:
    the tangent bundle $T\T$, the tangent bundle of the $G$-orbit $\underline{\mathfrak{g}}$,
    and the obstruction bundle $\E$. They are all $\widetilde{G}$-equivariant\footnote{The $\widetilde{G}$-action on $\mathfrak{g}$
    is the differential at identity of conjugation.}.  
    Thus, 
    we form a $\widetilde{G}$-equivariant virtual bundle
    $T\T-\mathfrak{g}-\E$ over $\T$. 
    Then,
    by ``a stably almost complex structure on the
    $K$-equivariant global Kuranishi chart'',
    we actually mean a stably almost complex structure 
    of the $\widetilde{G}$-equivariant virtual
    vector bundle $(T\T-\mathfrak{g})-\E$\footnote{
    In \cref{def:stably almost complex structures on vector orbi-bundles} we defined
    stably almost complex structures for virtual
    bundles over orbifolds. The definition
    works the same for virtual bundles over differentiable stacks. In our case,
    the differentiable stack is $[\T/\widetilde{G}]$
    }. More strongly, we
    want both $T\T-\mathfrak{g}$ and $\E$ to admit
    stably complex structures (e.g. after a stabilization).
\end{remark}

Let $\Kur=(G,\T,\E,s)$ be a $K$-equivariant global Kuranishi chart,
and $X$ be a closed complex manifold with a (holomorphic) $K$-action.
$\widetilde{G}$ acts on $X$ by $\widetilde{G}\to K$.
Let $f\colon \T\to X$ be a $\widetilde{G}$-equivariant smooth map.
For a torus $K=\TT$, we define the equivariant virtual pushforward as follows.
\begin{proposition}
\label{prop: equivariant vfc}
    Fix $K=\TT$. Let $\Kur=(G,\T,\E,s)$ and $X$ be as before, and $f\colon \T\to X$
    be a $\widetilde{G}$-equivariant smooth map. 
    Let $E$ be a Morava $K$-theory or an HKR theory.
    Then we can construct the equivariant
    virtual fundamental class as a map
    \[
        f_{!,eq}^\Kur\colon H_{\widetilde{G}}^*\left(\T,E\right)
        \to H_\TT^{*-\vdim\Kur+\dim X}(X,E)
    \]
\end{proposition}
\begin{proof}
    Similar to the construction of \cref{def:equivariant pushforward}, choose a sequence of finite-dimensional
    approximations of $B\TT$. In fact, if $\TT=T^{m}$,
    then we can take $B_N=(\CP^N)^m,P_N=(S^{2N+1})^m$.
    $P_N\to B_N$ is a $\TT$ principal bundle, and 
    there are inclusions $\iota_N\colon B_N\to B_{N+1}$.
    
    We define global Kuranishi charts
    \[\
            \Kur_N=\left(\widetilde{G}, P_N\times \T,\pi_{\T}^* \E=P_N\times \E, s_N=\pi_\T^* s\right)
    \]
    and stably almost complex manifolds
    \[
        X_N=P_N\times_\TT X
    \]
    The $\widetilde{G}$-action on $P_N\times \T$
    is diagonal, and we assume $X_N$ has trivial $\widetilde{G}$-action.
    There are smooth maps $f_N\colon P_N\times \T\to P_N\
    \times_\TT X$ which are $\widetilde{G}$-equivariant.

    Apply \cref{def:vfc as an pushforward}
    to $\Kur_N,X_N$ and $f_N\colon P_N\times \T\to X_N$,
    we have
    \[
        (f_N)_{!}^{\Kur_N}\colon H_{\widetilde{G}}^*\left(P_N\times \T,E\right)
        \to H^{*-\vdim\Kur_N+\dim X_N}(X_N,E).
    \]
    For $N$ and $N+1$, there is a commutative diagram
    \[\begin{tikzcd}
	{H^*_{\widetilde{G}}\left(P_{N}\times \T,E\right)} && {H^*_{\widetilde{G}}\left(P_{N+1}\times \T,E\right)} \\
	{H^{*-\vdim\Kur+\dim X}\left(X_N\right)} && {H^{*-\vdim\Kur+\dim X}\left(X_{N+1}\right)}
	\arrow["{(f_N)_!^{\Kur_N}}", from=1-1, to=2-1]
	\arrow["{\iota_N^*}"', from=1-3, to=1-1]
	\arrow["{(f_{N+1})_!^{\Kur_{N+1}}}", from=1-3, to=2-3]
	\arrow["{\iota_N^*}"', from=2-3, to=2-1]
\end{tikzcd}\]
This compatibility is from \cite[Lemma 5.34]{abouzaid2021complexcobordismhamiltonianloops}.

Finally, we define $f_{!,eq}^\Kur$ as follows: for
$a\in H^*_{\widetilde{G}}(\T,E)$, denote $a_N=\pi_\T^*(a)\in H^*_{\widetilde{G}}(P_N\times \T,E)$.
Then, define
\[
    f_{!,eq}^\Kur(a):=\left((f_N)_!^{\Kur_N}(a_N)\right)\in \varprojlim_N 
    H^{*-\vdim\Kur+\dim X}\left(P_N\times_\TT X\right)
\]
Each $X_N$ is a closed manifold, so $H^i(X_N,E)$ is finitely generated
over the complete Noetherian local ring $E^0$ for every $N$ and $i$.
As in the construction of \cref{def:equivariant pushforward}, Milnor's
$\lim^1$ term therefore vanishes, identifying the inverse limit with
$H_\TT^{*-\vdim\Kur+\dim X}(X,E)$.
\end{proof}
Similarly, we can define the ``fixed locus'' of
an equivariant global Kuranishi chart. Let $\Kur=(G,\T,\E,s)$ be a smooth $K$-equivariant global Kuranishi chart.  
We denote $\T^\fix$ as the fixed part defined in \cref{def:fixed points of orbifolds}. On $\T^\fix$
there are two vector bundles:
\begin{itemize}
    \item The normal bundle $\nu=\nu_{\T^\fix/\T}$. Notice that
    we have
    \[
        T\T|_{\T^\fix}=T\T^\fix\oplus \nu_{\T^\fix/\T}
    \]
    and that $\nu$ is the moving part of $T\T$ under
    $\Stab_{\widetilde{G}}(p)$.
    \item The restriction of $\E$ to $\T^\fix$. 
\end{itemize}
\begin{definition}
    For $p\in \T^\fix$,
    denote $\Stab_{\widetilde{G}}(p)$ to be its stabilizer group.
    Consider the action of the identity component
    $\Stab_{\widetilde{G}}(p)^0$, we can split it into
    the fixed and nontrivial part 
    \[
        \E_p=\E_p^\fix\oplus \E_p^\mov
    \]
    $\E^\fix$ and $\E^\mov$ are both vector bundles
    over $\T^\fix$.
\end{definition}
\begin{remark}
     Notice that they are all $\widetilde{G}$-equivariant. For example, for $g\in \Stab_{\widetilde{G}}(p)$
    and $h\in \Stab_{\widetilde{G}}(p)^0$, 
    $g^{-1}hg$ still lies in the identity component.
    Therefore, for $v\in \E^\fix_p$, we have $g^{-1}hg v=v$, or
    equivalently $h\cdot gv=gv$. This means $gv\in \E^\fix_p$,
    so $gv$ still lies in $\E^\mov_p$.
    However, our definition does not ensure that
    $\E^\fix$ is also fixed by elements from
    non-identity components in $\Stab_{\widetilde{G}}(p)$.
\end{remark}
It is easy to see that at $p\in \T^\fix$, the section $s(p)$ lies
entirely in $\E^\fix$. Therefore, we define a new global Kuranishi
chart.
\begin{definition}
\label{def:fixed locus of equivariant global Kuranishi chart}
    Assume $K$ is connected. For an $K$-equivariant global Kuranishi chart $\Kur=(G,\T,\E,s)$, define $\Kur^\fix$ to be
    \[
        \Kur^\fix=(G,\T^\fix, \E^\fix,s|_{\T^\fix}).
    \]
    There is a canonical map $i\colon \Kur^\fix\to \Kur$.
\end{definition}
Unlike orbifolds, the tangential differences between
$\Kur$ and $\Kur^\fix$ consist not only of the normal bundle
$\nu_{\T^\fix/\T}$, but also the moving part of the obstruction
$\E^\mov$. Recall \cref{section: equivariant localization} (and take
$A=\TT$), we have $S_\TT\subset H^*(B\TT,E)$. We need
the following invertibility of Euler classes of $\nu$ and $\E^\mov$.
\begin{lemma}
    The equivariant Euler classes $e_{E,eq}(\nu)$ and $e_{E,eq}(\E^\mov)$ are both invertible in
    \[
        S_\TT^{-1}H_{\widetilde{G}}\left(\T^\fix,E\right).
    \]
\end{lemma}
\begin{proof}
    The invertibility of $e_{E,eq}(\nu)$ is already proved
    in \cref{prop:invertibility of Euler class&isomorphism after S_A localization}.\footnote{
    There $M$ is a closed manifold, so we can induct on $\widetilde{G}$-cells. Here we need the assumption that
    $s^{-1}(0)$ is compact, so we can shrink $\T$ to make
    $\T$ also compact (with boundary but $s^{-1}(0)\cap \partial \T=\varnothing$).}
    
    For $\E^\mov$, consider the disk bundle of $\E^\mov$, denoted
    by $D(\E^\mov)$, which also admits a $\widetilde{G}$-action. The key observation is that
    $D(\E^\mov)^\fix=\T^\fix$. Let $p\in \T^\fix$ and $v\in \E^\mov_p$.
    Denote $H=\Stab_{\widetilde{G}}(p)$ and let $H^0$ be its
    identity connected component. By definition of $\T^\fix$, $H\to \TT$ is surjective. In fact, since $G\acts \T$ is almost free,
    $H\to \TT$ is a covering map.
    If $(p,v)$ lies in $D(\E^\mov)^\fix$,
    the map $\Stab_{\widetilde{G}}(p,v)\to \TT$ is also surjective.
    By computation of dimension, we have $\dim \Stab_{\widetilde{G}}(p,v)=\dim H=\dim \TT$. But $\Stab_{\widetilde{G}}(p,v)\subset H$,
    we have $H^0\subset \Stab_{\widetilde{G}}(p,v)$, so
    $v(p)\in \E^\mov_p$ is fixed by $H^0$. This contradicts
    our definition of $\E^\mov$.

    The normal bundle of $\T^\fix\subset D(\E^\mov)$ is exactly
    $\E^\mov$, so again by \cref{prop:invertibility of Euler class&isomorphism after S_A localization} we
    get the invertibility of $e_{E,eq}(\E^\mov)$
    in $S_\TT^{-1}H^*_{\widetilde{G}}\left(\T^\fix,E\right)$.
\end{proof}
We state the adapted version of \cref{prop:Atiyah-Bott equivariant localization formula}
as follows: 
\begin{proposition}
\label{prop:Atiyah-Bott localization for equivariant global Kuranishi charts}
    Let $\Kur$ be a $\TT$-equivariant global Kuranishi chart,
    and let $\Kur^\fix$ be its fixed chart as in \cref{def:fixed locus of equivariant global Kuranishi chart}. 
    Let $X$ be a $\TT$-equivariant stably almost complex
    manifold and $f\colon \T\to X$ be $\widetilde{G}$-equivariant.
    Then for $a\in H^*_{\widetilde{G}}(\T,E)$
    \begin{equation}
    \label{eq:Atiyah-Bott for equivariant gloal Kuranishi charts}
        f_{!,eq}^\Kur(a)=
        (f^\fix)_{!,eq}^{\Kur^\fix}\left(
            \frac{i^*(a)e_{E,eq}(\E^\mov)}{e_{E,eq}(\nu)}
        \right)
        \in S_{\TT}^{-1}H^*_\TT(X,E)
    \end{equation}

    Here $f^\fix\colon \T^\fix\to X$ is the restriction of $f$
    to $\T^\fix$.
\end{proposition}
Similar to \cref{prop:Atiyah-Bott equivariant localization formula},
the RHS 
of \eqref{eq:Atiyah-Bott for equivariant gloal Kuranishi charts}  makes
sense only after inverting $S_{\TT}$\footnote{
If $e\in S_\TT$, then $f_{!,eq}^\Kur(a\cup p_\Kur^* e)=f_{!,eq}^\Kur(a)\cup (p_X)^*e$ where $p_\Kur\colon \T \to \pt$
and $p_X\colon X\to \pt$ are the $\TT$-equivariant projection map to
a point. This allows us to extend definition of $f_{!,eq}^\Kur$
from $E$-cohomology to $S_\TT$-localized $E$-cohomology
in a compatible way. In particular, the RHS of
\eqref{eq:Atiyah-Bott for equivariant gloal Kuranishi charts} is defined this way.
}. 
    Choose an intermediate $\TT$-equivariant global Kuranishi chart
    $\Kur'=(G,\T^\fix,\E,s)$, and let $f'\colon \T^\fix\to X$
    be the restriction of $f$ to $\T^\fix$. Our proof splits into two
    goals in \cref{lem:two goals}.
    \begin{lemma}
    \label{lem:two goals}
        \begin{enumerate}[(1)]
        \item For $a\in H^*_{\widetilde{G}}(\T,E)$, we have
        \[
            f_{!,eq}^\Kur(a)=(f')_{!,eq}^{\Kur'}\left(
            \frac{i^*(a)}{e_{E,eq}(\nu)}
            \right)
        \]
        \item For $a\in H^*_{\widetilde{G}}(\T^\fix,E)$, we have
        \[
            (f')_{!,eq}^{\Kur'}(a)=
            (f^\fix)_{!,eq}^{\Kur^\fix}\left(
            a\cup e_{E,eq}(\E^\mov)
        \right)
        \]
    \end{enumerate}
    \end{lemma}
    \begin{proof}
        \begin{enumerate}[(1)]
            \item Recall that the equivariant pushforward is defined
            using an approximation \cref{prop: equivariant vfc}. We have the following
            commutative diagram,
            where the 
        vertical sequences come
        from \cref{def:vfc as an pushforward}
            \begin{center}
              \begin{adjustbox}{max width=\textwidth,center}
                    \begin{tikzcd}
	{H_{\widetilde{G}}^{*+\dim \T^\fix-\dim \T}\left(P_N\times\T^\fix,E\right)} & {H_{\widetilde{G}}^*\left(P_N\times\T,E\right)} \\
	{H^{\widetilde{G}}_{\dim (P_N\times\T)-\dim\widetilde{G}-*}\left((P_N\times\T^\fix)^c,E\right)} & {H^{\widetilde{G}}_{\dim (P_N\times\T)-\dim\widetilde{G}-*}\left((P_N\times\T^\fix)^c,E\right)} \\
	{H^{\widetilde{G}}_{\dim (P_N\times\T)-\dim\widetilde{G}-*}\left(P_N\times\E|_{\T^\fix},P_N\times(\E|_{\T^\fix}-\T^\fix),E\right)} & {H^{\widetilde{G}}_{\dim (P_N\times\T)-\dim\widetilde{G}-*}\left(P_N\times\E,P_N\times(\E-\T),E\right)} \\
	{H^{\widetilde{G}}_{\dim (P_N\times\T)-\dim\widetilde{G}-\rk \E-*}\left(P_N\times\T^\fix,E\right)} & {H^{\widetilde{G}}_{\dim (P_N\times\T)-\dim\widetilde{G}-\rk \E-*}\left(P_N\times\T,E\right)} \\
	{H_{\dim (P_N\times\T)-\dim\widetilde{G}-\rk \E-*}\left(X_N,E\right)} & {H_{\dim (P_N\times\T)-\dim\widetilde{G}-\rk \E-*}\left(X_N,E\right)} \\
	{H^{*-\vdim\Kur+\dim X}(X_N,E)} & {H^{*-\vdim\Kur+\dim X}(X_N,E)}
	\arrow["{i_{!}}", from=1-1, to=1-2]
	\arrow["{\lambda^{-1}_{\widetilde{G},P_N\times \T^\fix}}"', from=1-1, to=2-1]
	\arrow["{\lambda^{-1}_{\widetilde{G},P_N\times \T}}", from=1-2, to=2-2]
	\arrow["{i_*}", from=2-1, to=2-2]
	\arrow["{(s|_{\T^\fix})_*}", from=2-1, to=3-1]
	\arrow["{s_*}", from=2-2, to=3-2]
	\arrow["{i_*}", from=3-1, to=3-2]
	\arrow["{\operatorname{Thom}}", from=3-1, to=4-1]
	\arrow["{\operatorname{Thom}}", from=3-2, to=4-2]
	\arrow["{i_*}", from=4-1, to=4-2]
	\arrow["{(f')_*}", from=4-1, to=5-1]
	\arrow["{f_*}", from=4-2, to=5-2]
	\arrow["{=}", from=5-1, to=5-2]
	\arrow["{\lambda_{X_N}}", from=5-1, to=6-1]
	\arrow["{\lambda_{X_N}}", from=5-2, to=6-2]
	\arrow["{=}", from=6-1, to=6-2]
\end{tikzcd}
                \end{adjustbox}  
            \end{center}
               The commutativity of the top square is 
               by definition of $i_!$ in \cref{def:orbifold proper pushforward}. By this commutative diagram ,
               we have for $b\in H^*_{\widetilde{G}}(P_N\times \T^\fix,E)$
               \[
                    (f'_N)_{!}^{\Kur'_N}(b)
                    =(f_N)_{!}^{\Kur_N}(i_{!}(b))
               \]
               After localization at $S_\TT$,
               by \cref{prop:invertibility of Euler class&isomorphism after S_A localization} and $i^*\circ i_!=e_E(\nu)$
                we can prove the original formula by taking
               \[
                    b=\frac{i^*(a)}{e_{E}(\nu_{(P_N\times \T^\fix)/(P_N\times \T)})}
                \]
                because $i_!(b)=a$.
            \item The only difference
            between $\Kur'$ and $\Kur^\fix$
            is the obstruction bundle:
            there is an excess\footnote{Meaning: the
            section $s$ is zero on $\E^\mov$
            component} $\E^\mov$
            in $\Kur'$. We need two ingredients
            \begin{itemize}
                \item $\lambda^{-1}_{\widetilde{G},P_N\times \T^\fix}
                (a\cup e_{E}(\E^\mov))=e_E(\E^\mov)\cap \lambda^{-1}_{\widetilde{G},P_N\times \T^\fix}(a)$
                for any $ a\in H^*_{\widetilde{G}}(P_N\times \T^\fix,E)$.
                \item ${\operatorname{Thom}}_{\E^\fix}\left(s_*(e_E(\E^\mov)\cap b)\right)=\operatorname{Thom}_{\E^\fix\oplus \E^\mov}\left((s,0)_*(b)\right)$
                for $b\in H_*^{\widetilde{G}}\left((P_N\times \T^\fix)^c\right)$
            \end{itemize}
            Their proofs are as follows, let $Y=P_N\times \T^\fix$,
            $F=\E^\fix$, and $M=\E^\mov$ are $\widetilde{G}$-equivariant vector bundles over $Y$.
            \begin{itemize}
                \item Let $[Y^c]_E=\lambda^{-1}_{\widetilde{G},Y}(1)$
                be the equivariant $E$-fundamental class of $Y$, then
                $\lambda^{-1}_{\widetilde{G},Y}$ can be also defined 
                as $\lambda^{-1}_{\widetilde{G},Y}(a)=a\cap [Y^c]_E$.
                Therefore,
                \begin{align*}
                    \lambda^{-1}_{\widetilde{G},Y}(a\cup e_E(M))
                    &=(a\cup e_E(M))\cap [Y^c]_E\\
                    &=e_E(M)\cap (a\cap [Y^c]_E)\\
                    &=e_E(M)\cap \lambda^{-1}_{\widetilde{G},Y}(a).
                \end{align*}
                \item Let $q\colon F\to Y$ $p\colon F\oplus M\to Y$ be the bundle projections
                and let $z\colon F\to F\oplus M$ be the embedding. Let
                $u_F,u_M$ be the Thom classes. Then by \cref{lem:multiplicativity of the complex orientation}
                we have
                \[
                    u_{F\oplus M}=u_F\cup u_M
                \]
                Pulling back along $z$, we get
                \[
                    z^*u_{F\oplus M}=u_F\cup q^*e_E(M)
                \]
                Therefore,
                \begin{align*}
                    \operatorname{Thom}_{F\oplus M}((s,0)_*b)
                    &=\operatorname{Thom}_{F\oplus M}(z_*s_* b)\\
                    &=p_*\left(u_{F\oplus M}\cap z_*s_* b\right)\\
                    &=p_*z_*(z^*u_{F\oplus M}\cap s_* b)\\
                    &=q_*((u_F\cup q^*e_E(M))\cap s_*b)\\
                    &=q_*(u_F\cap (q^*e_E(M)\cap s_* b))
                \end{align*}
                Now use the projection formula for $s\colon Y\to F$,
                \[
                    q^*e_E(M)\cap s_*b=s_*(s^*q^*e_E(M)\cap b)
                \]
                The last term is $s_*(e_E(M)\cap b)$ because $q\circ s=\id_Y$. So finally we get
                \[
                    \operatorname{Thom}_{F\oplus M}((s,0)_*b)=
                    q_*(u_F\cap s_*(e_E(M)\cap b))=
                    \operatorname{Thom}_{F}(e_E(M)\cap b).
                \]
            \end{itemize}
            By the two ingredients, we can similarly prove
            the desired formula by comparing two definitional
            sequences, and one can notice that there are
            only two non-identical maps, covered by
            each ingredient.
        \end{enumerate}
    \end{proof}
\subsubsection{Equivariant lift of AMS construction} 
\label{sec:Equivariant lift of AMS construction}
Let $(M,\omega)$ be a symplectic manifold, and assume it admits a Hamiltonian
$K$-action for some connected compact Lie group $K$. 
Choose a compatible almost complex structure $J$ invariant under $K$-action.

Thus $K$ acts on the moduli space $\Mb_{g,h}(M,\beta)$ by postcomposition.
To lift this action to the AMS chart, we choose the auxiliary data as follows.
\begin{itemize}
\item We choose the auxiliary $J$-linear
connection $\nabla$ on $TM$ to be $K$-invariant, which can be obtained
by averaging any $J$-linear connection over $K$.
\item $\Ccal\to \Fcal$ and 
$\Cfrak_\Fcal\to \Ffrak_\Fcal$: note that $K$ does not
act on $\Ccal\to \Fcal$ but on
$\Cfrak_\Fcal
\to \Ffrak_\Fcal$. A $k\in K$ sends
$(\phi,u)\in \Ffrak_\Fcal$ to $(\phi,k\cdot u)$,
where $k\cdot u$ is defined using the $K$-action
on $M$.
    \item $\Omega,L$: 
    by \cite[Lemma~5.14]{Hirschi_Properties_of_Gromov_Witten},
using \cite[Corollary~1.4]{RIERA_2001}, we may choose the
Hermitian line bundle $L\to M$ such that the $K$-action on $M$
lifts to a linear unitary action on $L$. Choose a $K$-invariant
Hermitian connection on $L$ by averaging, and write its curvature
as $-2\pi\i\Omega$.
\item $k$: no difference.
\item $\mathfrak{D}$: we require the consistent domain metric
on $\Cfrak_\Fcal\to \Ffrak_\Fcal$ to be $\G_\C\times K$-invariant, with canonical $\G_C\times K$-action on $\Ffrak_\Fcal$.
\item $(V_\mu,\lambda_\mu)$: notice that 
\[\begin{tikzcd}
	{\Omega^{0,1}_{\Ccal^\circ/\Fcal}\otimes_\C TM} \\
	{\Ccal^\circ \times M}
	\arrow[from=1-1, to=2-1]
\end{tikzcd}\]
admits a canonical $G\times K$-action. We require
$V_\mu$ to be not only $G$-representations but also
$G\times K$-representations. This can be done by
applying \cite[Lemma 4.2]{abouzaidGromovWittenInvariantsComplex2024a}
where the group is taken to be $G\times K$.
Denseness is still required.
\end{itemize}
If we choose the auxiliary data $\mathfrak{a}$ in such a way, then
\cref{def: AMS construction} admits a $K$-equivariant lift.
We choose $\widetilde{G}=G\times K$, because the $G$ and $K$-action
on the data $(\phi,u,e,F)$ commute.
For $k\in K$, we can commute $k$ past $\pb_J$ and $\lambda_\mu$,
so the thickened moduli space
of framed curves $\T^\mathfrak{a}$ is preserved by
$\widetilde{G}$-action. On the obstruction bundle
$\E=H_{g,h,d}\oplus V_\mu\oplus \mathcal{H}$,
$K$ acts trivially on $H_{g,h,d}$ and $\mathcal{H}$, while $G$
acts on $H_{g,h,d}$ by the differential of its projective action
and on $\mathcal{H}$ by conjugation. Thus
\[
    (g,k)\cdot(\eta,v,H)
    =\bigl(g_*\eta,(g,k)\cdot v,gHg^{-1}\bigr),
\]
where the action on $V_\mu$ is part of $\mathfrak{a}$.
With a $G$-equivariant tubular-neighborhood identification,
these actions make $s_\afrak$ $G\times K$-equivariant.

For the smooth and stably almost complex structures, we retain the natural
$G\times K$-actions on the domain deformation and kernel parameters.
The domain gluing data are chosen equivariantly on local slices and
transported along group orbits. Metrics, cutoffs, and right inverses are
chosen compatibly using compact-group averaging. Pregluing then commutes
with the group action, and uniqueness of the Newton--Picard correction
gives equivariance of the gluing maps. The corresponding linear gluing
constructions respect the induced actions on the kernels. Therefore, all
constructions and choices involved in the gluing argument can be chosen
$G\times K$-invariant.

\newpage
\appendix
\section{Proof of \texorpdfstring{\cref{prop:an equivalence prop}}{Proposition \ref*{prop:an equivalence prop}}}
\label[appendix]{app:Proof of an equivalence}

  Before proving \cref{prop:an equivalence prop}, we introduce some notations.
    \begin{definition}
    \label{def:defining some new bases}
    \begin{enumerate}[(1)]
        \item Define the unitary frame bundle $B_1=\Fr_U(E_k)$, so
        $\Tot(K_1)\to B_1$ is a $U(r)$-equivariant
        vector bundle.
        We also define the frame bundle $B_\C=\Fr(E_k)$.
        \item  Define $B_2=B_\C\times_{\Fcal_{0,0,D}}U_{0,0,D}$. In other words, a point in $B_2$
        consists of $(\phi,f,F)$ such that
        \begin{itemize}
        \item $\phi\in \Fcal_{0,0,D}$.
            \item $f\colon \Sigma=\Sigma|_\phi\to \P^1$
            satisfy $\pb f=0$.
            \item $F$ is a (not necessarily
            unitary) framing of $H^0\left(
            \Sigma,\left(\omega_\Sigma\otimes
            f^*\O_{\P^1}(3)\right)^{\otimes k}
            \right)$\newline
            such that $(\phi,\phi_F)\in U_{0,0,D}$.
        \end{itemize}
        \item Define a map $j\colon B_1\to B_2$
        \[
            B_1\hookrightarrow B_\C
            =B_\C\times_{\Fcal_{0,0,D}}\Fcal_{0,0,D}
            \xhookrightarrow{\id_{B_\C}\times
            \Delta}B_\C\times_{\Fcal_{0,0,D}}
            U_{0,0,D}=B_2
        \]
        The first map is inclusion into complexification, and $\Delta\colon 
        \Fcal_{0,0,D}\to U_{0,0,D}$ is the diagonal map.
        $j$ sends $(f,F)\in B_1$ to 
        $(\phi_F,f,F)\in B_2$.
    \end{enumerate}
    \end{definition}
    \begin{lemma}
    \label{lem:submersiveness and normal bundle identification}
        \begin{enumerate}[(1)]
            \item $B_\C\to \Fcal_{0,0,D}$ and $U_{0,0,D}
            \to \Fcal_{0,0,D}$ are both submersive.
            Thus, $B_2\to \Fcal_{0,0,D}$ is also submersive.
            \item The map $j\colon B_1\to B_2$
            is an embedding with normal bundle
        \[
        N_{B_1/B_2}=\H\oplus H_{0,0,D}
        \]
         In fact, $\H$ is the normal bundle of $ B_1\hookrightarrow B_\C$, and $H_{0,0,D}$ is
    the normal bundle of $B_\C\hookrightarrow B_2$.
        \end{enumerate}
    \end{lemma}
    We omit the proof as it is just inspecting the 
    differential of each smooth map.
    
    Write $q\colon B_2\to \Fcal_{0,0,D}$. Note that
    for a point $(\phi,u,e,F)\in \T$, we can decompose
    the perturbed Cauchy--Riemann equation into $T\P^1$
    and $N$ direction, and $u=(f,\xi)$ where
    $f\colon \Sigma\to \P^1$ and $\xi\in \Omega^0(\Sigma,f^* N)$. Since $\lambda(e)$ is only in
    $N$ direction, on $T\P^1$ direction we actually
    have $\pb f=0$, i.e. $f\colon \Sigma\to \P^1$ is
    holomorphic.
    The chosen pullback Chern connection identifies the induced
    $\bar\partial$-operator on $u^*L=f^*\O_{\P^1}(3)$ with the usual
    Dolbeault operator, since $f$ is holomorphic.
    Thus the same $F$ is a holomorphic framing of the line bundle used in the
    definition of $B_2$, and determines the same projective map $\phi_F$.
    \begin{definition}
        Define a map $\rho\colon \T\to B_2$ to be
        \[
            (\phi,u=(f,\xi),e,F)\mapsto (\phi,f,F)        
        \]
    \end{definition}
    We have a commutative diagram
\[\begin{tikzcd}
	\T & {B_2} \\
	& {\Fcal_{0,0,D}}
	\arrow["\rho", from=1-1, to=1-2]
	\arrow["\pi"', from=1-1, to=2-2]
	\arrow["q", from=1-2, to=2-2]
\end{tikzcd}\]
Note that $\T\to \Fcal_{0,0,D}$ admits a $C^1_\loc$
structure, and $B_2\to \Fcal_{0,0,D}$ is a submersion.
\begin{lemma}
\label{lem:a C^1_loc compatibility lemma and surjective lemma}
The map $\rho$ is compatible with the $C^1_\loc$ structure
and is submersive in the sense that
\[
    D^{\operatorname{vert}} \rho\colon
    T^{\operatorname{vert}}(\T)\to 
    T^{\operatorname{vert}}(B_2)
\]
is surjective.
\end{lemma}
\begin{remark}
    $\rho$ is compatible with the $C^1_\loc$ structure means
    the following. Given any $C^1_\loc$ trivialization 
    $U\times V_1\to \T$
    over
    $U\subset \Fcal_{0,0,D}$ and a smooth trivialization
    $U\times V_2\to B_2$.
    Then on overlapping region, the coordinate transformation \[
        U\times V_1\to U\times V_2
    \]
    \[
        (b, v)\mapsto (b, f_b(v))
    \]
    \begin{itemize}
        \item gives a smooth map $v\mapsto f_b(v)$ for fixed $b$.
        \item $f_b\colon V_1\to V_2$ varies \emph{continuously} with $b$, in $C^1_\loc$
        topology of functions from $V_1$ to $V_2$.
    \end{itemize}
\end{remark}
\begin{remark}
\label{rmk:detailed discussion of C^1_loc trivialization}
    We recall how gluing arguments in \cite[Appendix B]{Pardon_2016} give the $C^1_\loc$
    structure over regular locus of pre-thickened moduli space. 
    
    Let $X$ be a symplectic manifold,
        which is $X_\loc$ in the case considered here.
        For a point $x_0=(\Sigma_0,u_0,e_0)\in \T^\pre$,
        let $U$ be a small neighborhood of $x_0\in \Fcal$,
        and let $b$ be local coordinates of $U$ near $x_0$.
        We have $\kappa\in \K$ where $K$ is the kernel of
        linearization of perturbed Cauchy--Riemann equation. 
        $K$ will become $T^{\operatorname{vert}}(\T^\pre/\Fcal)$. 
        
        For small $b$, one constructs
        pregluing maps $u_b^\pre\colon \Sigma_b\to X$,
        and uses the Newton--Picard iteration to define a map
        \begin{equation}
        \label{eq:gluing}
            \kappa\in K\mapsto \kappa_b\in \Gamma(\Sigma_b,u_b^\pre TX)
        \end{equation}
         to make $u_{b,\kappa}:=\exp_{u_b}(\kappa_b)$ really $J$-holomorphic. Such maps become the candidates for
         the structural $C^1_\loc$ trivializations of $U\times K\to \T^\pre$. See the details
         of gluing in \cite[B.6--B.11]{Pardon_2016}.

         Although it requires a complicated process to
         define $\kappa\mapsto u_{b,\kappa}$, there is
         a trick that makes the inverse map obvious.
         By appropriately choosing (\cite[Equation (B.7.3)]{Pardon_2016})
         \begin{itemize}
             \item Points $q_i\in \Sigma_0$ away from nodes and marked points.
             \item Subspaces $V_i\subset T_{u_0(q_i)}X$
             and $W''\subset W$.
         \end{itemize}
         such that the evaluation map
         \[
            L_0\colon K\subset \Omega^0\left(\Sigma_0, u_0^* TX\right)\oplus W\to \bigoplus_i (T_{u_0(q_i)}X/V_i)
            \oplus W/W''.
         \]
         is an equivalence, then the gluing
         construction can be performed so that
         the composition
         \[
            \kappa\mapsto (u_{b,\kappa}\colon \Sigma_b\to X)
    \xmapsto{\operatorname{ev}\footnotemark}\left(u_{b,\kappa}(q_i), {\operatorname{proj}}_{W\to W/W''}(e)\right).
        \]
        \footnotetext{
         To make sense of this evaluation map, we find coordinate charts near $u_0(q_i)\subset X$ such that
         the subspace $V_i$ becomes totally geodesic w.r.t. $\nabla$. Since all domains $\Sigma_b$ can be identified on regions away from nodes, we can extend
         the points $q_i\in \Sigma_0$ to a section of $\Ccal\to \Fcal$,
         also denoted by $q_i$. Then the evaluation is given by $u_{b,\kappa}\mapsto u_{b,\kappa}(q_i)$ in
         our chosen coordinate and to the quotient space
         $T_{u_0(q_i)}X/V_i$.
         }
        coincides with $L_0$. In other words, 
        the evaluation map followed by $L_0^{-1}$ is exactly
        the inverse of \eqref{eq:gluing}. 
\end{remark}
\begin{proof}[Proof of \cref{lem:a C^1_loc compatibility lemma and surjective lemma}]
        We write $\Fcal_{0,0,D}$ as $\Fcal$ for short.
        Since $\T$ is an open subset of the full $\GL(r,\C)$-framing bundle
        over $\T^{\pre,\mathrm{reg}}$, and gains its $C^1_\loc$ structure
        over $\Fcal$ from that of $\T^\pre/\Fcal$, it
        is convenient to define the corresponding
        prethickened object on the RHS. To be more precise,
        define\footnote{
        Another way to define $B_2^\pre$ is that it 
        is the prethickened moduli space in AMS
        construction if we don't choose perturbations
        (i.e. $V_\mu=0$), which is possible by $\P^1$ being convex.
        One might take another different approach from
        the current by comparing the AMS construction
        and gluing arguments of $X_\loc$ and $\P^1$.
        One chooses the data so the covariant derivative,
        the exponential map, and the parallel transport all respect
        the projection $p\colon X_\loc\to \P^1$:
        the $\nabla$ on $T\P^1$ and $TX_\loc$ (with $\nabla J=0$). The Cauchy--Riemann equation is compatible automatically. The only nontrivial part is
        choosing the evaluation map $L_0$ (\cite[Equation (B.7.3)]{Pardon_2016}) to
        be also compatible with $p$.
        Then one can claim that the map $\T^\pre\to B_2^\pre$
        becomes linear and even a projection to a vector bundle summand through trivializations from gluing (using these gluing data).
        }
        \[
            B_2^\pre=\{
                (\phi,f)\mid \phi\in \Fcal, f\colon
                \Sigma|_\phi\to \P^1 \text{ holomorphic}
            \}
        \]
        It is an algebro-geometric object, unobstructed, and does not 
        have nontrivial stabilizer groups since
        $\phi$ already fixes the domain. In fact,
        the projection $B_2^\pre\to \Fcal$ is smooth (in algebraic geometry) since $H^1\left(\Sigma,
        f^*T\P^1\right)=0$ for any genus $0$ $\Sigma$ and holomorphic map $f$. Finally,
        $B_2$ is locally an open subset of the full $\GL(r,\C)$-framing
        bundle over $B_2^\pre$. Since $\rho$ preserves the framing coordinate,
        restricting the local trivializations to these open subsets reduces
        the claim to proving that the map
        \[
            \rho^\pre\colon \T^\pre\to B_2^\pre
        \]
        \[
            (\phi,u,e)\mapsto (\phi,f)
        \]
        is compatible with the $C^1_\loc$ structure
        and submersive.

        \paragraph{$C^1_\loc$ compatibility.} We choose the 
        following
        $C^1_\loc$ or smooth trivializations. The details
        are recalled in \cref{rmk:detailed discussion of C^1_loc trivialization}.
        \begin{itemize}
            \item For $\T^\pre/\Fcal$: we choose
        the $C^1_\loc$ trivialization of $\T$
        \[
            U\times K\hookrightarrow \T
        \]
        \[
            (b,\kappa)\mapsto (u_{b,\kappa}\colon \Sigma_b\to X_\loc).
        \]
            \item For $B_2^\pre/\Fcal$, although its smooth structure comes from
         algebraic geometry directly, the evaluation still
         gives a smooth trivialization of fibers. 
         We can find $q_{i,2}\in \Sigma_b$ away from nodes and $V_i\subset T_{f_0(q_{i,2})}\P^1$ such that
         \[
            B_2^\pre\xrightarrow{\operatorname{ev}}\bigoplus
            _i T_{f_0(q_{i,2})}\P^1/ V_i
         \]
         becomes a smooth vertical coordinate.
        \end{itemize}  
        For all $\Sigma_b$, the region away from nodes can be
        identified; we denote it by $\Sigma'$.
        Thus, $C^1_\loc$ compatibility can be checked
        through 
        \[
            (b,\kappa)\in U\times K\mapsto u_{b,\kappa}\in W^{k,2}\left(\Sigma', X\right)\xmapsto{\operatorname{ev}} (u_{b,\kappa}(q_{i,2}))\in \bigoplus
            _i T_{f_0(q_{i,2})}\P^1/ V_i.
        \]
        Since the evaluation map is smooth for sufficiently
        large $k$, it suffices to prove the following two
        \footnote{This is also what is actually proved
        when checking different charts of $\T^\pre$ from gluing are $C^1_\loc$ compatible}
        \begin{itemize}
            \item $\kappa\mapsto u_{b,\kappa}\in W^{k,2}\left(\Sigma', X\right)$ is smooth for fixed $b$.
            \item As $b_j\to b_\infty$, the maps $\kappa\mapsto u_{b,\kappa}$ and $\kappa\mapsto D_{\kappa}u_{b,\kappa}$ converge uniformly
            in $C^0$ over any
            compact set.
        \end{itemize}
        These two facts come from (and are exactly the facts underlying the $C^1_\loc$ structure) the $C^1_\loc $ compatibility of
        two different charts for $\T^\pre $ from gluing,
        which was proved in \cite[Theorem 6.1, Corollary 6.29]{abouzaid2021complexcobordismhamiltonianloops}.

        \paragraph{Submersiveness.} The canonical map $\rho^\pre$ is $C^1_\loc$ compatible. At $e=0$,
        we can identify
        \begin{equation}
        \label{eq:vertical tangent of T}
             T^{\operatorname{vert}}(\T^\pre/\Fcal)|_{(\phi,u=(f,\xi),e=0)}
         = 
         \left\{
            \zeta\in \Gamma\left(\Sigma, u^* TX_\loc\right),
            e'\in W
            \mid 
            D_u \zeta+\lambda(e')\circ du=0
        \right\}
        \end{equation}
    \begin{equation}
\label{eq:vertical tangent of T equivalent description}
\begin{aligned}
T^{\operatorname{vert}}(\T^\pre/\Fcal)
|_{(\phi,u=(f,\xi),e=0)}
&=
\Bigl\{
\zeta_{\P^1}\in\Gamma(u^*T\P^1),\
\zeta_N\in\Gamma(u^*N),\
e'\in W
\ \Bigm|                                      \\
&\qquad
\pb_{f^*T\P^1}\zeta_{\P^1}=0,\quad
\pb_{f^*N}(\zeta_N)
+\lambda(e')\circ du=0
\Bigr\}.
\end{aligned}
\end{equation}
        Since the almost complex structure is just the complex structure
        (i.e. it is integrable), the linearization $D_u$ coincides with
        $\pb_{u^* TX_\loc}$. Thus we can split the equation into $T\P^1$
        and $N$ direction.
        On the other hand, we can write
        \[
            T^{\operatorname{vert}}(B_2^\pre/\Fcal)|_{(\phi,f)}
            =
            \left\{
                \zeta_{\P^1}\in \Gamma\left(\Sigma,f^*T\P^1\right)
                \mid \pb_{f^* T\P^1}\zeta_{\P^1}=0.
            \right\}
        \]
        It is obvious that the natural projection map
        \[
        D^{\operatorname{vert}}\rho^\pre\colon  T^{\operatorname{vert}}
        (\T^\pre/\Fcal)\to
        T^{\operatorname{vert}}(B_2^\pre/\Fcal)
        \]
        \[
            (\zeta_{\P^1},\zeta_N,e')\mapsto \zeta_{\P^1}
        \]
        is surjective, because we have already split
        $\zeta_{\P^1}$ and $(\zeta_{N},e')$. In fact, its kernel is $\ker(W\to H^1(\Sigma, f^* N))$. This completes the proof.
    \end{proof}
We notice that the kernel of $ D^{\operatorname{vert}}\rho^\pre$ is already defined on $B_2^\pre$. It lifts to a vector bundle $K_2\to B_2$
given by
\[
    K_2|_{(\phi,f,F)}=\ker\left(
    W\to H^1\left(\Sigma,f^* N\right).
    \right)
\]
Note that we have a pullback diagram of vector bundles
\begin{equation}
    \label{eq: pullback K_1 to K_2}
    \begin{tikzcd}
	{K_1} & {K_2} \\
	{B_1} & {B_2}
	\arrow["{\widetilde{j}}", from=1-1, to=1-2]
	\arrow[from=1-1, to=2-1]
	\arrow[from=1-2, to=2-2]
	\arrow["j", from=2-1, to=2-2]
\end{tikzcd}
\end{equation}
\begin{definition}
    Define $t\colon \T\to \Tot(K_2)$ by sending
    $(\phi,u,F,e)$ to $(\phi,f,F)\in B_2, e\in K_2|_{(\phi,f,F)}$
    \footnote{This is
because $u=(f,\xi),e\in W$ satisfy
the perturbed Cauchy--Riemann equation:
$\pb_{f^* N}\xi+\lambda(e)|_{(\phi,f)}=0$,
which implies $e\in K_2$.}.
\end{definition}
\begin{proposition}
\begin{enumerate}[(1)]
    \item $t$ is compatible with the $C^1_\loc$
    structure of $\T/\Fcal$ and
    $\Tot(K_2)/\Fcal$. The $C^1_\loc$
    structure on the latter comes from
    the vector bundle structure $K_2\to B_2$
    and the submersion $B_2\to \Fcal$.
    \item $t$ induces a homeomorphism of topological
    spaces (near the zero section) and an isomorphism on vertical tangent bundles. Thus,
    $t$ lifts the $C^1_\loc$ structure on $\T/\Fcal$ to a smooth submersion.
    \item Via the identification $t\colon \T\to \Tot(K_2)$,
    the inclusion $i\colon \Tot(K_1)\to \T$ can be viewed
    as the map induced by the vector bundle map
    $\widetilde{j}\colon K_1\to K_2$.
\end{enumerate}
\end{proposition}
\begin{proof}
    \begin{enumerate}[(1)]
        \item There is a smooth embedding $\Tot(K_2)\hookrightarrow
        B_2\times W$. The $C^1_\loc$ compatibility of $\rho\colon \T\to B_2$
        is already proved in \cref{lem:a C^1_loc compatibility lemma and surjective lemma}. The $C^1_\loc$ compatibility of $\T\to W$
        is because it is just forgetting other coordinates. Therefore,
        the map $t\colon \T\to \Tot(K_2)$ also respects the $C^1_\loc$
        structure.
        \item In the identification \eqref{eq:vertical tangent of T equivalent description}, $T^{\operatorname{vert}}(\T/B_2)
        \simeq T^{\operatorname{vert}}(\T^\pre/B_2^\pre) $
        can be identified with
        \[
            \left\{\zeta_N,e'\mid \pb_{f^* N}(\zeta_N)+\lambda(e')\circ du=0\right\}.
        \]
        The tangent map relative to $B_2$, $T^{\operatorname{vert}}(\T/B_2)\xrightarrow{dt}
        T^{\operatorname{vert}}(\Tot(K_2)/B_2)=K_2$ is identified with 
        \[
            (\zeta_N,e')\mapsto e'\in K_2
        \]
        This is a complex-linear isomorphism, which means $t$ induces an isomorphism $T^{\operatorname{vert}}(\T/B_2)\to T^{\operatorname{vert}}(\Tot(K_2)/B_2)$ and also
        an isomorphism
        $T^{\operatorname{vert}}(\T/\Fcal)\to T^{\operatorname{vert}}(\Tot(K_2)/\Fcal)$. After shrinking, $t\colon \T/\Fcal\to \Tot(K_2\to B_2)/\Fcal$ is an isomorphism of $C^1_\loc$
        fiber bundles. 
        \item We compose $\Tot(K_1)\xrightarrow{i} \T \xrightarrow{t}\Tot(K_2)$. Its effect on $(f,F,e)$ is
        \[
            (f,F)\in \Fr_U(E_k),e\in K_1\mapsto (\phi_F, u=(f,\xi),F,e\in W)\mapsto 
            (\phi_F,f,F,e\in K_2)
        \]
        Here $\xi$ is the unique solution of $\pb_{f^*N}\xi+\lambda(e)\circ du=0$.
        Omitting $e$, this gives exactly the base map $j\colon B_1\to B_2$ as defined in \cref{def:defining some new bases}.
        The total map preserves the fibers; in fact, it is the identity
        $e\mapsto e$. Therefore, the inclusion $i\colon \Tot(K_1)\to \T$
        is identified with $\Tot(K_1)\to \Tot(K_2)$ induced by
        \eqref{eq: pullback K_1 to K_2}.
    \end{enumerate}
\end{proof}
To summarize, we get a smoothing of $\T$ which is $\Tot(K_2)$.
\begin{proof}[Proof of \cref{prop:an equivalence prop}]
    Recall the morphism $\Kur_1\to \Kur^{\operatorname{AMS}}$ is
    $i\colon \Tot(K_1)\to \T\simeq \Tot(K_2)$ on the base,
    and $W\hookrightarrow H_{0,0,D}\oplus W\oplus \H$ on the obstruction
    bundle. Since the homotopy defining the stable complex structure is constant
    along the zero section, \cref{lem:computation of virtual deformation theory}
    and the complex-linear identifications above show that $i$ preserves the
    specified stable complex structures near it. Thus it remains to prove that
    the following complex is exact
    \[
        0\to T(\Tot(K_1))\xrightarrow{\operatorname{d}i} T \T\to H_{0,0,D}\oplus \H\to 0.
    \]
    Since \eqref{eq: pullback K_1 to K_2} is a pullback diagram,
    we have $\coker(\operatorname{d}i)\simeq N_{B_1/B_2}$,
    which is isomorphic to $H_{0,0,D}\oplus \H$ by \cref{lem:submersiveness and normal bundle identification}(2).
    In fact, one can check that the composition $T\T\to TB_2\to N_{B_1/B_2}\simeq H_{0,0,D}\oplus \H$ coincides with
    the map given by $s$, the section of $\Kur^{\operatorname{AMS}}$
\end{proof}

\section{Proof of \texorpdfstring{\cref{prop:vector bundle structure of K}}{Proposition \ref*{prop:vector bundle structure of K}}}
\label[appendix]{app:Proof of vector bundle structure of K}
We recall some context from \cite[Appendix B]{Pardon_2016}. Using notations there\footnote{In our notation, we use
$x\in\T^\pre$ to represent a general point corresponding to coordinates
$(\alpha,y,\kappa)$ and write
$\Sigma_x=C_\alpha$, $u_x=u_{\alpha,y,\kappa}$.},
    we denote $(\alpha,y)$ as the parameters
    of domain curves (near $\Sigma_0$), where $\alpha$ controls smoothing of nodal points
    and $y$ controls deformation of complex structures. Let $C_\alpha$ be the smooth
    nodal surface corresponding to $\alpha$,
    Pardon constructed a pregluing
    map $u_\alpha\colon C_\alpha\to M$
    \cite[Definition B.3.2]{Pardon_2016}.
    Then, for $\kappa\in K_0$ sufficiently small
    and $(\alpha,y)$ sufficiently close to 
    the origin,
    by Newton-Picard iteration, one
    can get a $\kappa_{\alpha,y}\in
    W^{k,2,\delta}(C_\alpha,u_\alpha^* TM)\oplus E$ such that 
    \[
        u_{\alpha,y,\kappa}:=\exp_{u_\alpha}(\kappa_{\alpha,y})
    \]
    \[
        e_{\alpha,y,\kappa}:=e_0+\mathrm{proj}_E
        \kappa_{\alpha,y}
    \]
    forms a solution of perturbed Cauchy--Riemann
    equation \cite[Proposition B.9.2]{Pardon_2016}, and this gives a local
    parametrization of $\T^\pre$ \cite[Proposition B.11.5]{Pardon_2016}.
    We take the $\nabla$ on $M$ to be
    \textit{$J$-linear}\footnote{
        This requirement is assumed in \cite[Appendix B.3]{Pardon_2016}.
    }, using which one defines parallel transport maps:
    \[
        \mathsf{PT}_{u_\alpha\to u_x}
        \colon W^{k,2,\delta}\left(C_\alpha,u_\alpha^*TM\right)\oplus E
        \to 
        W^{k,2,\delta}\left(C_\alpha, u_x^* TM\right)
        \oplus E
    \]
    \[
        \mathsf{PT}_{u_x\to u_\alpha}
        \colon W^{k-1,2,\delta}\left(\widetilde{C}_\alpha,
        \Omega^{0,1}_{\widetilde{C}_\alpha}
        \otimes_\C u_x^* TM\right)
        \to 
        W^{k-1,2,\delta}\left(\widetilde{C}_\alpha,
        \Omega^{0,1}_{\widetilde{C}_\alpha}
        \otimes_\C u_\alpha^* TM\right)
    \]

    Fix $(x_0,t_0)\in \T^\pre\times [0,1]$, 
    and let $K_0=K_{x_0,t_0}$.
    As in \cite[Equation (B.7.3)]{Pardon_2016}, we can choose
    finitely many points $q_i$ on smooth
    parts of $\Sigma_0$ with subspaces $V_i\subset T_{u_0(q_i)}M$ and a subspace
    $E''\subset E$ such that the projection
    map
    \[
        l_0\colon 
        W^{k,2,\delta}(\Sigma_0,u_0^*TM)\oplus
        E \to W:= \left(\bigoplus_{i}T_{u_0(q_i)M/V_i}\right)
        \oplus E/E''
    \]
    has a restriction to $K_0$ that is an
    isomorphism $l_0|_{K_0}\colon K_0\xrightarrow{\simeq} W$. Equivalently,
    the restriction 
    \[\widehat{L}_{x_0}(t_0)|_{\ker l_0}
    \colon \ker l_0\to \Omega^{0,1}(\widetilde{\Sigma_0}, u_0^* TM)\] is an isomorphism. Let
    $Q_0$ be its inverse.
    For $x$ sufficiently close to $x_0$ in $\T^\pre$, we can define\footnote{Note that
    $u_\alpha$ coincides with $u_0$ via identification of $C_\alpha$
    and $C_0$ away
    from nodes.} 
   $q_i(x)=\exp_{u_\alpha(q_i)}(\kappa_{\alpha,y})$
    and define $V_i(x)\subset T_{u_x(q_i)}M$
    by parallel transport.
    Define
    \[
        l_x\colon W^{k,2,\delta}(\Sigma_x, u_x^* TM)
        \oplus E\to W_x:=
        \left(\bigoplus_i 
        T_{u_x(q_i(x))}M/V_i(x)\right)
        \oplus 
        E/E''.
    \] 

    We extend diagram (B.7.7) in Pardon's
    paper as
    \begin{equation}
    \label{eq:Pardon B.7.7 extension}
        \begin{tikzcd}[
    column sep=7em,
    row sep=5.5em,
]
W^{k,2,\delta}\!\left(C_\alpha,u_x^*TX\right)\oplus E
    \arrow[r, "\widehat{L}_{x}(t_0)"]
&
W^{k-1,2,\delta}\!\left(
    \widetilde C_\alpha,
    \Omega^{0,1}_{\widetilde C_{\alpha,y}}
    \otimes_{\C} u_x^*TX
\right)
    \arrow[d, "\mathsf{PT}_{u_x\to u_\alpha}"]
\\
W^{k,2,\delta}\!\left(C_\alpha,u_\alpha^*TX\right)\oplus E
    \arrow[u, "\mathsf{PT}_{u_\alpha\to u_x}"']
    \arrow[r, "\widehat{L}_{\alpha,y}(t_0)"]
&
W^{k-1,2,\delta}\!\left(
    \widetilde C_\alpha,
    \Omega^{0,1}_{\widetilde C_{\alpha,y}}
    \otimes_{\C} u_\alpha^*TX
\right)
    \arrow[d, "\mathrm{break}"]
\\
W^{k,2,\delta}\!\left(C_0,u_{0|\alpha}^*TX\right)\oplus E
    \arrow[u, "\mathrm{glue}"']
    \arrow[r, "\widehat{L}_{0|\alpha,y}(t_0)"]
    \arrow[d, leftrightarrow, "\mathsf{PT}"']
&
W^{k-1,2,\delta}\!\left(
    \widetilde C_0,
    \Omega^{0,1}_{\widetilde C_{0,y}}
    \otimes_{\C} u_{0|\alpha}^*TX
\right)
    \arrow[d, leftrightarrow, "\mathsf{PT}"]
\\
W^{k,2,\delta}\!\left(C_0,u_0^*TX\right)\oplus E
    \arrow[r, "\widehat{L}_{0,y}(t_0)"]
    \arrow[d, equal, "\mathrm{id}"']
&
W^{k-1,2,\delta}\!\left(
    \widetilde C_0,
    \Omega^{0,1}_{\widetilde C_{0,y}}
    \otimes_{\C} u_0^*TX
\right)
    \arrow[d, leftrightarrow, "I_y"]
\\
W^{k,2,\delta}\!\left(C_0,u_0^*TX\right)\oplus E
    \arrow[r, shift left=0.6ex, "\widehat{L}_{0,0}(t_0)"]
&
W^{k-1,2,\delta}\!\left(
    \widetilde C_0,
    \Omega^{0,1}_{\widetilde C_{0,0}}
    \otimes_{\C} u_0^*TX
\right)
    \arrow[l, shift left=0.6ex, "Q_{0}"]
\end{tikzcd}
    \end{equation}
For each of the operators occurring in 
    \eqref{eq:Pardon B.7.7 extension},
    write
    \[
        \widehat L_\bullet(t)
=
\widehat D_\bullet+(1-t)\widehat A_\bullet,
\qquad
\widehat A_\bullet
:=
\widehat D_\bullet^{\C}-\widehat D_\bullet,
    \]
    where $\widehat A_\bullet$ is zeroth order and vanishes on the
finite-dimensional stabilization factor $E$.
Then, after shrinking the gluing neighborhood, the estimates of
\cite[Lemmas~B.7.2--B.7.4 and Equation~(B.8.10)]{Pardon_2016}
continue to hold with every $\widehat D_\bullet$
replaced by $\widehat L_\bullet(t)$, with constants independent of $t\in [0,1]$.
More precisely:
\begin{lemma}
\label{lem:generalizaing Pardon's estimate}
\begin{enumerate}[(1)]
    \item For the bottom square,
    \[
        \bigl|
\widehat L_{0,y}(t)
-
I_y\circ\widehat L_{0,0}(t)
\bigr|
\leq C|y|.
    \]
    \item For the middle square,
    \[
    \bigl|
\mathsf{PT}\circ\widehat L_{0,y}(t)
-
\widehat L_{0|\alpha,y}(t)\circ\mathsf{PT}
\bigr|
\leq
C\sum_{i=1}^d e^{-S_i}.
    \]
    \item For the glue--break square, if
    $\widehat L_{0|\alpha,y}(t)\xi
=
\operatorname{break}(\eta)$,
then
\[
    \left\|
\widehat L_{\alpha,y}(t)\operatorname{glue}(\xi)
-
\eta
\right\|_{k-1,2,\delta}
\leq
C\left\|\xi\right\|_{k,2,\delta}
\sum_{i=1}^d e^{-2\delta S_i}.
\]
\item Suppose
$u_x=\exp_{u_\alpha}\zeta$.
Using the $J$-linear parallel transports between $u_\alpha^*TM$
and $u_x^*TM$, one has
\[
    \begin{aligned}
&
\left\|
\widehat L_{\alpha,y}(t)
-
\mathsf{PT}{u_x\to u_\alpha}
\circ\widehat L_x(t)
\circ\mathsf{PT}_{u_\alpha\to u_x}
\right\|_{W^{k,2,\delta}\to W^{k-1,2,\delta}}
\\
&\hspace{4cm}
\leq C\left\|\zeta\right\|_{k,2,\delta}.
\end{aligned}
\]
\end{enumerate}
All four estimates hold uniformly for $t\in [0,1]$.
\end{lemma}

\begin{proof}
    \begin{enumerate}[(1)]
        \item Write
        \[
            \widehat{L}_{0,y}(t)-
            I_y \widehat{L}_{0,0}(t)
            =(\widehat{D}_{0,y}-I_y
            \widehat{D}_{0,0})
            +(1-t)(\widehat{A}_{0,y}
            -I_y
            \widehat{A}_{0,0}).
        \]
        Pardon proves
        \[
            \left\|
            \widehat{D}_{0,y}-I_y
            \widehat{D}_{0,0}
            \right\|
            \le C|y|.
        \]
        For the second term, $\widehat{A}_{0,y}$
        is a zeroth-order operator whose coefficients depend smoothly on
        $j_y,u_0,e_0$ and $I_y$ is complex
        linear which is smooth in $y$, then 
        \[
            \left\|
                \widehat{A}_{0,y}-I_y\widehat{A}_{0,0}
            \right\|
            \le C |y|.
        \]
        \item Expand
        \[
            \mathsf{PT}\circ 
            \widehat{L}_{0,y}(t)
            -\widehat{L}_{0|\alpha,y}(t)
            \circ \mathsf{PT}
            =
            \mathsf{PT}\circ 
            \widehat{D}_{0,y}(t)
            -\widehat{D}_{0|\alpha,y}(t)
            \circ \mathsf{PT}
            +(1-t)
            (\mathsf{PT}\circ 
            \widehat{A}_{0,y}(t)
            -\widehat{A}_{0|\alpha,y}(t)
            \circ \mathsf{PT}).
        \]
        Pardon proves 
        \[
            \mathsf{PT}\circ 
            \widehat{D}_{0,y}(t)
            -\widehat{D}_{0|\alpha,y}(t)
            \circ \mathsf{PT}
            \le C \sum_{i}e^{-S_i}.
        \]
        For the $A$-part, we write
        \[
        \begin{aligned}
        \widehat{A}_{0,y}
        &=-\frac{1}{2}\left(\widehat{D}_{0,y}+J_Y\circ \widehat{D}_{0,y}\circ J_X\right)\\
        \text{and}\quad \widehat{A}_{0|\alpha,y}
        &=-\frac{1}{2}\left(\widehat{D}_{0|\alpha,y}+J_{Y,\alpha}\circ \widehat{D}_{0|\alpha,y}\circ J_{X,\alpha}\right).
        \end{aligned}
        \]
        Since $\nabla J=0$,
        the parallel transport commutes with these $J_X,J_Y,J_{X,\alpha},J_{Y,\alpha}$.
        Then
        \begin{align*}
            \mathsf{PT} \circ \widehat{A}_{0,y}-
            \widehat{A}_{0|\alpha,y}\circ 
            \mathsf{PT}
            &= -\frac{1}{2}\Bigl(
            \mathsf{PT}\circ \widehat{D}_{0,y}
            +\mathsf{PT}\circ  J_Y \circ \widehat{D}_{0,y} J_X\\
            &\qquad
            -\widehat{D}_{0|\alpha,y}\circ 
            \mathsf{PT}- J_{Y,\alpha}
            \circ \widehat{D}_{0|\alpha,y}\circ J_{X,\alpha}\circ \mathsf{PT}
            \Bigr)\\
            &=-\frac{1}{2}\left(
            \Delta_{D}+ J_{Y,\alpha}
            \Delta_D J_{X,0}
            \right),
        \end{align*}
        where $\Delta_D= \mathsf{PT} D_{0,y}
        -D_{0|\alpha,y}\mathsf{PT}$. Suppose
        the Sobolev norms we used are induced by
        $J$-Hermitian metrics, then $\left\|J_{Y,\alpha}\right\|=1,
        \left\| J_{X}\right\|=1$. Therefore,
        \[
            \left\|\mathsf{PT} \circ \widehat{A}_{0,y}-
            \widehat{A}_{0|\alpha,y}\circ 
            \mathsf{PT}\right\|
            \le \left\|\Delta_D\right\|
        \]
        and RHS is controlled by Pardon's
        original estimate.
        \item 
        For the glue--break square, no new estimate is required.
Indeed, the defect in Pardon’s proof of Lemma B.7.4 is supported in $[2S_i, 2S_i+1]\cup 
[4S_i-1,4S_i]$ on each neck. On the larger region
the preglued map is constant and the perturbation term vanishes. Hence the corresponding linearized operator is simply the complex-linear Cauchy–Riemann operator
$\nabla^{0,1}$. In particular its complex-antilinear zeroth-order part $\widehat{A}_\bullet$ vanishes there. Thus
$\widehat{L}_\bullet(t)=\widehat{D}_\bullet$
on the support of the glue--break defect,
and Pardon's proof of Lemma B.7.4 applies verbatim, with the same estimate
\[
    \left\|\widehat{L}_{\alpha,y}(t)\operatorname{glue}
    (\xi)-\eta\right\|_{k-1,2,\delta}
    \le C\left\|\xi\right\|
    \sum_{i}e^{-2\delta S_i}.
\]
    \item The proof is similar to (2). 
\end{enumerate}
\end{proof}
Defining $T_{x,t_0}$ by
\[
    T_{x,t_0}=\mathsf{PT}_{u_\alpha\to u_x}\circ \operatorname{glue}\circ \mathsf{PT}\circ \id\circ Q_{0}\circ I_y \circ \mathsf{PT}\circ \operatorname{break}\circ\mathsf{PT}_{u_x\to u_\alpha}.
\]
\cref{lem:generalizaing Pardon's estimate} implies that for $t=t_0$ and $x\to x_0$, we have
\begin{equation}
    \label{eq:boundedness of T and asymptotic right inverse}
    \left\| T_{x,t_0}\right\|\le C,
    \quad
    \left\|\widehat{L}_x(t_0)\circ T_{x,t_0}-\operatorname{Id}\right\|
    \to 0
\end{equation}
For $t\neq t_0$, we have
\[
    \left\|\widehat{L}_{x}(t)T_{x,t_0}-\operatorname{Id}\right\|
    \le 
    \left\|\widehat{L}_x(t_0)\circ T_{x,t_0}-\operatorname{Id}\right\|+
    |t-t_0|\cdot
    \left\|\widehat{A}_xT_{x,t_0}\right\|.
\]
By taking $(x,t)$ in some open neighborhood of $(x_0,t_0)$, there is 
\begin{equation}
    \label{eq:pre right inverse estimate}
    \left\|\widehat{L}_{x}(t)T_{x,t_0}-\operatorname{Id}\right\|\le \frac{1}{2}
\end{equation}
Define 
\[
    Q_{x,t}=T_{x,t_0}\sum_{n\ge 0}
    \left(\operatorname{Id}- \widehat{L}_{x}(t)T_{x,t_0}\right)^n
\]
such that
\[
    \widehat{L}_{x}(t)Q_{x,t}=\operatorname{Id},
    \quad
    \left\|Q_{x,t}\right\|\le C
\]
By construction of gluing maps and
definition(by parallel transport) of $V_i(x)\subset T_{u_x(q_i(x))}M$,
we have $\im(Q_{x,t})\subset \ker l_x$. Moreover, by surjectivity of $\widehat{L}_x(t)$, we have
\[\codim \im Q_{x,t}=\dim \ker \widehat{L}_x(t)=\footnotemark
\dim \ker \widehat{L}_{x_0}(t_0)=\dim W=\dim W_x=\codim \ker 
l_x,\]\footnotetext{Here we used computation of index of Cauchy--Riemann operators to get $\dim \ker \widehat{L}_x(t)=
\dim \ker \widehat{L}_{x_0}(t_0)$. }
so 
\begin{equation}
\label{eq:im Q_x t=ker l_x}
    \im(Q_{x,t})=\ker l_x.
\end{equation}
\begin{proof}[Proof of \cref{prop:vector bundle structure of K}]
    For a fixed $(x_0,t_0)\in \T^\pre$ and $(x,t)$ 
    in some open neighborhood of $(x_0,t_0)$, consider
    \[
        l_x|{K_{x,t}}\colon K_{x,t}\to W_x.
    \]
    By \eqref{eq:im Q_x t=ker l_x}, this is
    an isomorphism of vector spaces. $W_x$ can be canonically identified with $W$ by backward
    parallel transport, so $l_x$ gives a
    local trivialization of $K_{x,t}$ in
    an open neighborhood of $(x_0,t_0)$.

    We prove that this defines a topological
    vector bundle structure on $\mathcal{K}$. Equivalently,
    for two points $a,b\in \T^\pre \times [0,1]$ 
    and open neighborhoods $U_a,U_b$
    where the aforementioned trivializations are defined,
    if $U_a\cap U_b\neq\varnothing$, the transport
    $\tau_x\colon W^a_x\xrightarrow{(l_x^{a})^{-1}} K_{x,t}
    \xrightarrow{(l^b_x)} W^b_x$ is continuous in $(x,t)$
    (using the canonical trivialization of $W_x\simeq W$).
    For simplicity of notation, we use $p$ to denote the
    point $(x,t)$.
Write
    $X_p=W^{k,2,\delta}\left(C_\alpha,u_x^* TM\right)\oplus 
    E$ and $Y_p=W^{k-1,2,\delta}\left(
\widetilde{C}_\alpha,\Omega^{0,1}\widetilde{C}_\alpha\otimes_\C u_{x}^* TM
    \right)$. Write $K_p=\ker D_p$ (which is $K_{x,t}$ before). We write $D_p:= \widehat{L}_{x}(t)$ and
    $Q_p^a,Q^b_p$ for $Q_{x,t}$ where $a$ or $b$ refer
    to the base point of construction (which is $(x_0,t_0)$ before).

    We can explicitly choose the inverse of $l_p^a\colon K_p\to W^a_p$
    as follows. 
    \begin{itemize}
        \item First choose a section $s_p^a\colon W^a_p
    \to X_p$ continuous in $p$. This is possible because
    $q_i$ are away from neck region, so we can choose $s_p^a$ to be bump functions near each $q_i$.
    \item Define $R_p^a:= \operatorname{Id} - Q_p^a D_p$. This
    is the projection $X_p\to K_p$ with regard to
    $X_p=K_p\oplus \im Q_p^a$.
    \item The inverse is defined to be $
    w\in W_p^a\mapsto R_p^a \circ s_p^a (w)\in K_p$.
    \end{itemize}
    Therefore, $\tau_p$ can be represented by
    \[
        \tau_p= l_p^b(\operatorname{Id}- Q_p^a D_p)s^a_p w
        =l_p^bs^a_p w-l^b_pQ_p^a D_p s_p^a w
    \]
    It suffices to prove: for fixed $w\in W^a$ (using
    parallel transport $W_p^a\simeq W^a$), 
    $l_p^bs^a_p w$ and $l^b_pQ_p^a D_p s_p^a w$
    are both continuous in $p$.
        
        For $l_p^bs^a_p w$, notice that $s^a_p w$
        as a function is supported away from neck
        regions/thin parts, and one can identify the thick parts of
        domain for different $p$,
        thus $s^a_p w$ can be viewed as a continuous family (in $p$) of functions over a single domain.
        Therefore, restricting to another set
        of points, which is what $l_p^b$ does, is
        continuous in $p$.
       
       For $l^b_pQ_p^a D_p s_p^a w$, we
        write
        \[
            l^b_pQ_p^a D_p s_p^a w=\sum_{n\ge 0}
            l_p^b T_p^a(\operatorname{Id}-D_p T_p^a)^nD_p s_p^n w
        \]
        It suffices to prove each term is continuous in $p$
        and this summation converges uniformly in $C^0$-norm. Let's examine every term in this expression:
        \begin{itemize}
            \item $s_p^a w$ varies continuously over thick parts.
            \item $D_p$: its coefficients vary continuously.
            \item break/glue/$\mathsf{PT}$/$I_y$:
            continuous with regard to gluing parameters.
            \item $Q_0^a$ fixed
            \item $l_p^b$: finite-dimensional evaluation.
        \end{itemize}
        So $p\mapsto l_p^b T_p^a(\operatorname{Id}-D_p T_p^a)^nD_p s_p^n w$ is continuous. Uniform
        convergence follows from \eqref{eq:pre right inverse estimate} and \eqref{eq:boundedness of T and asymptotic right inverse}. 

        Finally, we cite \cite[Lemma 5.11]{hirschiGlobalKuranishiCharts2024b}
        for lifting $K|_{\T^\pre\times 0}$ 
        to a topological complex vector bundle.
\end{proof}

\printbibliography

@misc{cheng2013poincaredualitymoravaktheory,
      title={Poincare duality in Morava K-theory for classifying spaces of orbifolds}, 
      author={Man Chuen Cheng},
      year={2013},
      eprint={1305.2731},
      archivePrefix={arXiv},
      primaryClass={math.AT},
      url={https://arxiv.org/abs/1305.2731}, 
}

@misc{abouzaid2021complexcobordismhamiltonianloops,
      title={Complex cobordism, Hamiltonian loops and global Kuranishi charts}, 
      author={Mohammed Abouzaid and Mark McLean and Ivan Smith},
      year={2021},
      eprint={2110.14320},
      archivePrefix={arXiv},
      primaryClass={math.SG},
      url={https://arxiv.org/abs/2110.14320}, 
}

@book {hori2003mirror,
    AUTHOR = {Hori, Kentaro and Katz, Sheldon and Klemm, Albrecht and
              Pandharipande, Rahul and Thomas, Richard and Vafa, Cumrun and
              Vakil, Ravi and Zaslow, Eric},
     TITLE = {Mirror symmetry},
    SERIES = {Clay Mathematics Monographs},
    VOLUME = {1},
      NOTE = {With a preface by Vafa},
 PUBLISHER = {American Mathematical Society, Providence, RI; Clay
              Mathematics Institute, Cambridge, MA},
      YEAR = {2003},
     PAGES = {xx+929},
      ISBN = {0-8218-2955-6},
   MRCLASS = {14J32 (14N35 32Q25 81T30 81T60)},
  MRNUMBER = {2003030},
MRREVIEWER = {Marcos\ Mari\~no},
}

@incollection {Enumeration_of_Rational_Curves_Via_Torus_Actions,
    AUTHOR = {Kontsevich, Maxim},
     TITLE = {Enumeration of rational curves via torus actions},
 BOOKTITLE = {The moduli space of curves ({T}exel {I}sland, 1994)},
    SERIES = {Progr. Math.},
    VOLUME = {129},
     PAGES = {335--368},
 PUBLISHER = {Birkh\"auser Boston, Boston, MA},
      YEAR = {1995},
      ISBN = {0-8176-3784-2},
   MRCLASS = {14N10 (14D22 14L30)},
  MRNUMBER = {1363062},
MRREVIEWER = {Anatoly\ Libgober},
       DOI = {10.1007/978-1-4612-4264-2\_12},
       URL = {https://doi.org/10.1007/978-1-4612-4264-2_12},
}

@article {strickland2000knlocaldualityfinitegroups,
    AUTHOR = {Strickland, N. P.},
     TITLE = {{$K(N)$}-local duality for finite groups and groupoids},
   JOURNAL = {Topology},
  FJOURNAL = {Topology. An International Journal of Mathematics},
    VOLUME = {39},
      YEAR = {2000},
    NUMBER = {4},
     PAGES = {733--772},
      ISSN = {0040-9383},
   MRCLASS = {55N15 (18D35 19L99 57R56)},
  MRNUMBER = {1760427},
MRREVIEWER = {Ulrike\ Tillmann},
       DOI = {10.1016/S0040-9383(99)00031-2},
       URL = {https://doi.org/10.1016/S0040-9383(99)00031-2},
}

@article {Greenlees1996TheTS,
    AUTHOR = {Greenlees, J. P. C. and Sadofsky, Hal},
     TITLE = {The {T}ate spectrum of {$v_n$}-periodic complex oriented
              theories},
   JOURNAL = {Math. Z.},
  FJOURNAL = {Mathematische Zeitschrift},
    VOLUME = {222},
      YEAR = {1996},
    NUMBER = {3},
     PAGES = {391--405},
      ISSN = {0025-5874,1432-1823},
   MRCLASS = {55N22},
  MRNUMBER = {1400199},
MRREVIEWER = {Mark\ Hovey},
       DOI = {10.1007/PL00004264},
       URL = {https://doi.org/10.1007/PL00004264},
}

@article {Generalized-Group-Characters-and-Complex-Oriented-Cohomology-Theories,
    AUTHOR = {Hopkins, Michael J. and Kuhn, Nicholas J. and Ravenel, Douglas
              C.},
     TITLE = {Generalized group characters and complex oriented cohomology
              theories},
   JOURNAL = {J. Amer. Math. Soc.},
  FJOURNAL = {Journal of the American Mathematical Society},
    VOLUME = {13},
      YEAR = {2000},
    NUMBER = {3},
     PAGES = {553--594},
      ISSN = {0894-0347,1088-6834},
   MRCLASS = {55N22 (55N20 55N91 55R40)},
  MRNUMBER = {1758754},
MRREVIEWER = {Andrew\ J.\ Baker},
       DOI = {10.1090/S0894-0347-00-00332-5},
       URL = {https://doi.org/10.1090/S0894-0347-00-00332-5},
}

@article {ojm/1200771835,
    AUTHOR = {Kawasaki, Tetsuro},
     TITLE = {The {R}iemann-{R}och theorem for complex {$V$}-manifolds},
   JOURNAL = {Osaka Math. J.},
  FJOURNAL = {Osaka Mathematical Journal},
    VOLUME = {16},
      YEAR = {1979},
    NUMBER = {1},
     PAGES = {151--159},
      ISSN = {0388-0699},
   MRCLASS = {58G10 (32C10)},
  MRNUMBER = {527023},
MRREVIEWER = {Harold\ Donnelly},
       URL = {http://projecteuclid.org/euclid.ojm/1200771835},
}

@article {abouzaidGromovWittenInvariantsComplex2024a,
    AUTHOR = {Abouzaid, Mohammed and McLean, Mark and Smith, Ivan},
     TITLE = {Gromov-{W}itten invariants in complex and {M}orava-local
              {$K$}-theories},
   JOURNAL = {Geom. Funct. Anal.},
  FJOURNAL = {Geometric and Functional Analysis},
    VOLUME = {34},
      YEAR = {2024},
    NUMBER = {6},
     PAGES = {1647--1733},
      ISSN = {1016-443X,1420-8970},
   MRCLASS = {53D45 (14N35 55N15)},
  MRNUMBER = {4823210},
       DOI = {10.1007/s00039-024-00697-4},
       URL = {https://doi.org/10.1007/s00039-024-00697-4},
}

@article {hirschiGlobalKuranishiCharts2024b,
    AUTHOR = {Hirschi, Amanda and Swaminathan, Mohan},
     TITLE = {Global {K}uranishi charts and a product formula in symplectic
              {G}romov-{W}itten theory},
   JOURNAL = {Selecta Math. (N.S.)},
  FJOURNAL = {Selecta Mathematica. New Series},
    VOLUME = {30},
      YEAR = {2024},
    NUMBER = {5},
     PAGES = {Paper No. 87, 74},
      ISSN = {1022-1824,1420-9020},
   MRCLASS = {53D45 (14N35)},
  MRNUMBER = {4807086},
       DOI = {10.1007/s00029-024-00982-y},
       URL = {https://doi.org/10.1007/s00029-024-00982-y},
}

@article {Hirschi_Properties_of_Gromov_Witten,
    AUTHOR = {Hirschi, Amanda},
     TITLE = {Properties of {G}romov-{W}itten invariants defined via global
              {K}uranishi charts},
   JOURNAL = {Ann. H. Lebesgue},
  FJOURNAL = {Annales Henri Lebesgue},
    VOLUME = {8},
      YEAR = {2025},
     PAGES = {811--871},
      ISSN = {2644-9463},
   MRCLASS = {53D45 (32Q65)},
  MRNUMBER = {5005386},
       DOI = {10.5802/ahl.249},
       URL = {https://doi.org/10.5802/ahl.249},
}

@article {RIERA_2001,
    AUTHOR = {Mundet i Riera, Ignasi},
     TITLE = {Lifts of smooth group actions to line bundles},
   JOURNAL = {Bull. London Math. Soc.},
  FJOURNAL = {The Bulletin of the London Mathematical Society},
    VOLUME = {33},
      YEAR = {2001},
    NUMBER = {3},
     PAGES = {351--361},
      ISSN = {0024-6093,1469-2120},
   MRCLASS = {57S15 (53D05)},
  MRNUMBER = {1817775},
MRREVIEWER = {Ian\ Hambleton},
       DOI = {10.1017/S0024609301007937},
       URL = {https://doi.org/10.1017/S0024609301007937},
}

@article{Ionel-Parker_Gopakumar-Vafa,
author = {Eleny-Nicoleta Ionel and Thomas Parker},
title = {{The Gopakumar-Vafa formula for symplectic manifolds}},
volume = {187},
journal = {Annals of Mathematics},
number = {1},
publisher = {Department of Mathematics of Princeton University},
pages = {1 -- 64},
year = {2018},
doi = {10.4007/annals.2018.187.1.1},
URL = {https://doi.org/10.4007/annals.2018.187.1.1}
}

@article{kontsevichGromovWittenClassesQuantum1994,
  title = {Gromov-{{Witten}} Classes, Quantum Cohomology, and Enumerative Geometry},
  author = {Kontsevich, M. and family=Manin, given=Yu., given-i={{Yu}}},
  date = {1994-08-01},
  journaltitle = {Communications in Mathematical Physics},
  shortjournal = {Communications in Mathematical Physics},
  volume = {164},
  number = {3},
  pages = {525--562},
  issn = {1432-0916},
  doi = {10.1007/BF02101490},
  url = {https://doi.org/10.1007/BF02101490}
}

@article {Li-Tian_algebraic,
    AUTHOR = {Li, Jun and Tian, Gang},
     TITLE = {Virtual moduli cycles and {G}romov-{W}itten invariants of
              algebraic varieties},
   JOURNAL = {J. Amer. Math. Soc.},
  FJOURNAL = {Journal of the American Mathematical Society},
    VOLUME = {11},
      YEAR = {1998},
    NUMBER = {1},
     PAGES = {119--174},
      ISSN = {0894-0347,1088-6834},
   MRCLASS = {14D20 (14D15 14N10)},
  MRNUMBER = {1467172},
MRREVIEWER = {Ralph\ Martin\ Kaufmann},
       DOI = {10.1090/S0894-0347-98-00250-1},
       URL = {https://doi.org/10.1090/S0894-0347-98-00250-1},
}

@article {intrisic-normal-cone,
    AUTHOR = {Behrend, K. and Fantechi, B.},
     TITLE = {The intrinsic normal cone},
   JOURNAL = {Invent. Math.},
  FJOURNAL = {Inventiones Mathematicae},
    VOLUME = {128},
      YEAR = {1997},
    NUMBER = {1},
     PAGES = {45--88},
      ISSN = {0020-9910,1432-1297},
   MRCLASS = {14F99 (14C15 14D20)},
  MRNUMBER = {1437495},
MRREVIEWER = {Tohru\ Nakashima},
       DOI = {10.1007/s002220050136},
       URL = {https://doi.org/10.1007/s002220050136},
}

@article {Ruan-Tian_1995,
    AUTHOR = {Ruan, Yongbin and Tian, Gang},
     TITLE = {A mathematical theory of quantum cohomology},
   JOURNAL = {J. Differential Geom.},
  FJOURNAL = {Journal of Differential Geometry},
    VOLUME = {42},
      YEAR = {1995},
    NUMBER = {2},
     PAGES = {259--367},
      ISSN = {0022-040X,1945-743X},
   MRCLASS = {58D29 (14C05 14N10 58D10)},
  MRNUMBER = {1366548},
MRREVIEWER = {Bernd\ Siebert},
       URL = {http://projecteuclid.org/euclid.jdg/1214457234},
}

@incollection {Li-Tian_1996,
    AUTHOR = {Li, Jun and Tian, Gang},
     TITLE = {Virtual moduli cycles and {G}romov-{W}itten invariants of
              general symplectic manifolds},
 BOOKTITLE = {Topics in symplectic {$4$}-manifolds ({I}rvine, {CA}, 1996)},
    SERIES = {First Int. Press Lect. Ser.},
    VOLUME = {I},
     PAGES = {47--83},
 PUBLISHER = {Int. Press, Cambridge, MA},
      YEAR = {1998},
      ISBN = {1-57146-019-5},
   MRCLASS = {53D45 (53D35 57R17 57R57)},
  MRNUMBER = {1635695},
}

@incollection {Fukaya-Ono_1999,
    AUTHOR = {Fukaya, Kenji and Ono, Kaoru},
     TITLE = {Arnold conjecture and {G}romov-{W}itten invariant for general
              symplectic manifolds},
 BOOKTITLE = {The {A}rnoldfest ({T}oronto, {ON}, 1997)},
    SERIES = {Fields Inst. Commun.},
    VOLUME = {24},
     PAGES = {173--190},
 PUBLISHER = {Amer. Math. Soc., Providence, RI},
      YEAR = {1999},
      ISBN = {0-8218-0945-8},
   MRCLASS = {53D40 (53D45 57R58)},
  MRNUMBER = {1733575},
MRREVIEWER = {David\ E.\ Hurtubise},
       DOI = {10.1090/fic/024/12},
       URL = {https://doi.org/10.1090/fic/024/12},
}

@book {McDuff-Salamon_2012,
    AUTHOR = {McDuff, Dusa and Salamon, Dietmar},
     TITLE = {{$J$}-holomorphic curves and symplectic topology},
    SERIES = {American Mathematical Society Colloquium Publications},
    VOLUME = {52},
   EDITION = {Second},
 PUBLISHER = {American Mathematical Society, Providence, RI},
      YEAR = {2012},
     PAGES = {xiv+726},
      ISBN = {978-0-8218-8746-2},
   MRCLASS = {53D45 (32Q65 53D35)},
  MRNUMBER = {2954391},
MRREVIEWER = {Mark\ Alan\ Branson},
}

@article {Pardon_2016,
    AUTHOR = {Pardon, John},
     TITLE = {An algebraic approach to virtual fundamental cycles on moduli
              spaces of pseudo-holomorphic curves},
   JOURNAL = {Geom. Topol.},
  FJOURNAL = {Geometry \& Topology},
    VOLUME = {20},
      YEAR = {2016},
    NUMBER = {2},
     PAGES = {779--1034},
      ISSN = {1465-3060,1364-0380},
   MRCLASS = {53D35 (37J10 53D37 53D40 53D42 53D45 54B40 57R17)},
  MRNUMBER = {3493097},
MRREVIEWER = {Sonja\ Hohloch},
       DOI = {10.2140/gt.2016.20.779},
       URL = {https://doi.org/10.2140/gt.2016.20.779},
}

@incollection {Givental_2000,
    AUTHOR = {Givental, Alexander},
     TITLE = {On the {WDVV} equation in quantum {$K$}-theory},
      NOTE = {Dedicated to William Fulton on the occasion of his 60th
              birthday},
   JOURNAL = {Michigan Math. J.},
  FJOURNAL = {Michigan Mathematical Journal},
    VOLUME = {48},
      YEAR = {2000},
     PAGES = {295--304},
      ISSN = {0026-2285,1945-2365},
   MRCLASS = {14N35 (19E08 53D45)},
  MRNUMBER = {1786492},
MRREVIEWER = {Gary\ P.\ Kennedy},
       DOI = {10.1307/mmj/1030132720},
       URL = {https://doi.org/10.1307/mmj/1030132720},
}

@article {Lee_2004,
    AUTHOR = {Lee, Y.-P.},
     TITLE = {Quantum {$K$}-theory. {I}. {F}oundations},
   JOURNAL = {Duke Math. J.},
  FJOURNAL = {Duke Mathematical Journal},
    VOLUME = {121},
      YEAR = {2004},
    NUMBER = {3},
     PAGES = {389--424},
      ISSN = {0012-7094,1547-7398},
   MRCLASS = {14N35 (19E08 53D45 55N15)},
  MRNUMBER = {2040281},
MRREVIEWER = {Andrew\ Kresch},
       DOI = {10.1215/S0012-7094-04-12131-1},
       URL = {https://doi.org/10.1215/S0012-7094-04-12131-1},
}

@book {Adamsbook,
    AUTHOR = {Adams, J. F.},
     TITLE = {Stable homotopy and generalised homology},
    SERIES = {Chicago Lectures in Mathematics},
 PUBLISHER = {University of Chicago Press, Chicago, Ill.-London},
      YEAR = {1974},
     PAGES = {x+373},
   MRCLASS = {55B20 (55E10)},
  MRNUMBER = {402720},
MRREVIEWER = {P.\ S.\ Landweber},
}

@book {Rugyak1998,
    AUTHOR = {Rudyak, Yuli B.},
     TITLE = {On {T}hom spectra, orientability, and cobordism},
    SERIES = {Springer Monographs in Mathematics},
      NOTE = {With a foreword by Haynes Miller},
 PUBLISHER = {Springer-Verlag, Berlin},
      YEAR = {1998},
     PAGES = {xii+587},
      ISBN = {3-540-62043-5},
   MRCLASS = {55-02 (55N22 55P42 57-02)},
  MRNUMBER = {1627486},
MRREVIEWER = {Donald\ M.\ Davis},
}

@unpublished{LurieChromatic252x,
  author = {Lurie, Jacob},
  title  = {Chromatic Homotopy Theory},
  note   = {Course notes for Mathematics 252x, Harvard University. Available at \url{https://people.math.harvard.edu/~lurie/252x.html}; accessed July 5, 2026}
}

@incollection {Wurgler-MoravaKtheory,
    AUTHOR = {W\"urgler, Urs},
     TITLE = {Morava {$K$}-theories: a survey},
 BOOKTITLE = {Algebraic topology {P}ozna\'n{} 1989},
    SERIES = {Lecture Notes in Math.},
    VOLUME = {1474},
     PAGES = {111--138},
 PUBLISHER = {Springer, Berlin},
      YEAR = {1991},
      ISBN = {3-540-54098-9},
   MRCLASS = {55N35 (19L99 55P42)},
  MRNUMBER = {1133896},
MRREVIEWER = {Mark\ Hovey},
       DOI = {10.1007/BFb0084741},
       URL = {https://doi.org/10.1007/BFb0084741},
}

@article {Moerdijk-Pronk1997,
    AUTHOR = {Moerdijk, I. and Pronk, D. A.},
     TITLE = {Orbifolds, sheaves and groupoids},
   JOURNAL = {$K$-Theory},
  FJOURNAL = {$K$-Theory. An Interdisciplinary Journal for the Development,
              Application, and Influence of $K$-Theory in the Mathematical
              Sciences},
    VOLUME = {12},
      YEAR = {1997},
    NUMBER = {1},
     PAGES = {3--21},
      ISSN = {0920-3036,1573-0514},
   MRCLASS = {22A22 (18B25 55N30 57S15)},
  MRNUMBER = {1466622},
MRREVIEWER = {A.\ Verschoren},
       DOI = {10.1023/A:1007767628271},
       URL = {https://doi.org/10.1023/A:1007767628271},
}

@article {Lerman-orbifolds_as_stacks,
    AUTHOR = {Lerman, Eugene},
     TITLE = {Orbifolds as stacks?},
   JOURNAL = {Enseign. Math. (2)},
  FJOURNAL = {L'Enseignement Math\'ematique. Revue Internationale. 2e
              S\'erie},
    VOLUME = {56},
      YEAR = {2010},
    NUMBER = {3-4},
     PAGES = {315--363},
      ISSN = {0013-8584},
   MRCLASS = {18D05 (22A22)},
  MRNUMBER = {2778793},
MRREVIEWER = {Chenchang\ Zhu},
       DOI = {10.4171/LEM/56-3-4},
       URL = {https://doi.org/10.4171/LEM/56-3-4},
}

@article {Pronk1996,
    AUTHOR = {Pronk, Dorette A.},
     TITLE = {Etendues and stacks as bicategories of fractions},
   JOURNAL = {Compositio Math.},
  FJOURNAL = {Compositio Mathematica},
    VOLUME = {102},
      YEAR = {1996},
    NUMBER = {3},
     PAGES = {243--303},
      ISSN = {0010-437X,1570-5846},
   MRCLASS = {18E35 (18D05)},
  MRNUMBER = {1401424},
MRREVIEWER = {A.\ Verschoren},
       URL = {http://www.numdam.org/item?id=CM_1996__102_3_243_0},
}

@incollection {Moerdijk-orbifolds_as_groupoids2002,
    AUTHOR = {Moerdijk, Ieke},
     TITLE = {Orbifolds as groupoids: an introduction},
 BOOKTITLE = {Orbifolds in mathematics and physics ({M}adison, {WI}, 2001)},
    SERIES = {Contemp. Math.},
    VOLUME = {310},
     PAGES = {205--222},
 PUBLISHER = {Amer. Math. Soc., Providence, RI},
      YEAR = {2002},
      ISBN = {0-8218-2990-4},
   MRCLASS = {22A22 (55N30 55P15 58H05)},
  MRNUMBER = {1950948},
MRREVIEWER = {Janez\ Mr\v cun},
       DOI = {10.1090/conm/310/05405},
       URL = {https://doi.org/10.1090/conm/310/05405},
}

@article {Behrend-Xu-differentiable_stacks,
    AUTHOR = {Behrend, Kai and Xu, Ping},
     TITLE = {Differentiable stacks and gerbes},
   JOURNAL = {J. Symplectic Geom.},
  FJOURNAL = {The Journal of Symplectic Geometry},
    VOLUME = {9},
      YEAR = {2011},
    NUMBER = {3},
     PAGES = {285--341},
      ISSN = {1527-5256,1540-2347},
   MRCLASS = {53C08 (19L50 22A22 53D50 58H05)},
  MRNUMBER = {2817778},
MRREVIEWER = {Chenchang\ Zhu},
       DOI = {10.4310/jsg.2011.v9.n3.a2},
       URL = {https://doi.org/10.4310/jsg.2011.v9.n3.a2},
}

@misc{bai2024equivariantformalitycomplexorientedtheories,
      title={Equivariant formality in complex-oriented theories}, 
      author={Shaoyun Bai and Daniel Pomerleano},
      year={2024},
      eprint={2405.05821},
      archivePrefix={arXiv},
      primaryClass={math.SG},
      url={https://arxiv.org/abs/2405.05821}, 
}

@article{virtual_localization,
    AUTHOR = {Aranha, Dhyan and Khan, Adeel A. and Latyntsev, Alexei and
              Park, Hyeonjun and Ravi, Charanya},
     TITLE = {Virtual localization revisited},
   JOURNAL = {Adv. Math.},
  FJOURNAL = {Advances in Mathematics},
    VOLUME = {479},
      YEAR = {2025},
     PAGES = {Paper No. 110434, 53},
      ISSN = {0001-8708,1090-2082},
   MRCLASS = {14D23 (14A20 14A30 14C15 14L30 14N35)},
  MRNUMBER = {4933128},
MRREVIEWER = {Victor\ Przyjalkowski},
       DOI = {10.1016/j.aim.2025.110434},
       URL = {https://doi.org/10.1016/j.aim.2025.110434},
}

@article {Ebert-Giansiracusa2011,
    AUTHOR = {Ebert, Johannes and Giansiracusa, Jeffrey},
     TITLE = {Pontrjagin-{T}hom maps and the homology of the moduli stack of
              stable curves},
   JOURNAL = {Math. Ann.},
  FJOURNAL = {Mathematische Annalen},
    VOLUME = {349},
      YEAR = {2011},
    NUMBER = {3},
     PAGES = {543--575},
      ISSN = {0025-5831,1432-1807},
   MRCLASS = {14D23 (14H10 55R40)},
  MRNUMBER = {2754996},
MRREVIEWER = {Andrea\ Bruno},
       DOI = {10.1007/s00208-010-0518-2},
       URL = {https://doi.org/10.1007/s00208-010-0518-2},
}

@article {Pardon2023-Orbifold_bordism,
    AUTHOR = {Pardon, John},
     TITLE = {Orbifold bordism and duality for finite orbispectra},
   JOURNAL = {Geom. Topol.},
  FJOURNAL = {Geometry \& Topology},
    VOLUME = {27},
      YEAR = {2023},
    NUMBER = {5},
     PAGES = {1747--1844},
      ISSN = {1465-3060,1364-0380},
   MRCLASS = {55P42 (55M05 55N91 55P25 55P91 55Q91 57R85 57R91)},
  MRNUMBER = {4621920},
MRREVIEWER = {Lisu\ Wu},
       DOI = {10.2140/gt.2023.27.1747},
       URL = {https://doi.org/10.2140/gt.2023.27.1747},
}

@book {Cheng-Deligne_Mumford_Stacks2011,
    AUTHOR = {Cheng, Man Chuen},
     TITLE = {A {D}uality {T}heorem for {D}eligne-{M}umford {S}tacks with
              {R}espect to {M}orava {K}-{T}heory},
      NOTE = {Thesis (Ph.D.)--Stanford University},
 PUBLISHER = {ProQuest LLC, Ann Arbor, MI},
      YEAR = {2011},
     PAGES = {79},
      ISBN = {979-8672-19498-1},
   MRCLASS = {99-05},
  MRNUMBER = {4172290},
       URL =
              {http://gateway.proquest.com/openurl?url_ver=Z39.88-2004&rft_val_fmt=info:ofi/fmt:kev:mtx:dissertation&res_dat=xri:pqm&rft_dat=xri:pqdiss:28168008},
}

@misc{BenMoshe2025Transfers,
  author       = {Ben-Moshe, Shay},
  title        = {Transfers in Transchromatic Homotopy Theory},
  year         = {2025},
  month        = oct,
  howpublished = {Lecture notes for the University of Toronto
                  Algebraic Geometry Seminar},
  note         = {Talk delivered October 22, 2025. Accessed July 15, 2026},
  url          = {https://shaybm.com/pdfs/notes/2025-10-22-Toronto-AG.pdf}
}

@article {AtiyahBottMomentMap,
    AUTHOR = {Atiyah, M. F. and Bott, R.},
     TITLE = {The moment map and equivariant cohomology},
   JOURNAL = {Topology},
  FJOURNAL = {Topology. An International Journal of Mathematics},
    VOLUME = {23},
      YEAR = {1984},
    NUMBER = {1},
     PAGES = {1--28},
      ISSN = {0040-9383},
   MRCLASS = {58F05 (55N91 57R20 57R22)},
  MRNUMBER = {721448},
MRREVIEWER = {G.\ J.\ Heckman},
       DOI = {10.1016/0040-9383(84)90021-1},
       URL = {https://doi.org/10.1016/0040-9383(84)90021-1},
}

@incollection {Lashof1979-Stable_G-smoothing,
    AUTHOR = {Lashof, Richard},
     TITLE = {Stable {$G$}-smoothing},
 BOOKTITLE = {Algebraic topology, {W}aterloo, 1978 ({P}roc. {C}onf., {U}niv.
              {W}aterloo, {W}aterloo, {O}nt., 1978)},
    SERIES = {Lecture Notes in Math.},
    VOLUME = {741},
     PAGES = {283--306},
 PUBLISHER = {Springer, Berlin},
      YEAR = {1979},
      ISBN = {3-540-09545-4},
   MRCLASS = {57R10 (57S10)},
  MRNUMBER = {557173},
MRREVIEWER = {Michael\ C.\ Bix},
}

@misc{bai2022arnoldconjectureintegers,
      title={Arnold conjecture over integers}, 
      author={Shaoyun Bai and Guangbo Xu},
      year={2022},
      eprint={2209.08599},
      archivePrefix={arXiv},
      primaryClass={math.SG},
      url={https://arxiv.org/abs/2209.08599}, 
}

@article {Graber-Pandharipande99,
    AUTHOR = {Graber, T. and Pandharipande, R.},
     TITLE = {Localization of virtual classes},
   JOURNAL = {Invent. Math.},
  FJOURNAL = {Inventiones Mathematicae},
    VOLUME = {135},
      YEAR = {1999},
    NUMBER = {2},
     PAGES = {487--518},
      ISSN = {0020-9910,1432-1297},
   MRCLASS = {14C17 (14D20 14N10 14N35)},
  MRNUMBER = {1666787},
MRREVIEWER = {Paolo\ Aluffi},
       DOI = {10.1007/s002220050293},
       URL = {https://doi-org.libproxy.mit.edu/10.1007/s002220050293},
}

@book {Rotman95,
    AUTHOR = {Rotman, Joseph J.},
     TITLE = {An introduction to the theory of groups},
    SERIES = {Graduate Texts in Mathematics},
    VOLUME = {148},
   EDITION = {Fourth},
 PUBLISHER = {Springer-Verlag, New York},
      YEAR = {1995},
     PAGES = {xvi+513},
      ISBN = {0-387-94285-8},
   MRCLASS = {20-01},
  MRNUMBER = {1307623},
MRREVIEWER = {Ronald\ Solomon},
       DOI = {10.1007/978-1-4612-4176-8},
       URL = {https://doi.org/10.1007/978-1-4612-4176-8},
}

@article {BakerRichter2005,
    AUTHOR = {Baker, Andrew and Richter, Birgit},
     TITLE = {On the {$\Gamma$}-cohomology of rings of numerical polynomials
              and {$E_\infty$} structures on {$K$}-theory},
   JOURNAL = {Comment. Math. Helv.},
  FJOURNAL = {Commentarii Mathematici Helvetici. A Journal of the Swiss
              Mathematical Society},
    VOLUME = {80},
      YEAR = {2005},
    NUMBER = {4},
     PAGES = {691--723},
      ISSN = {0010-2571,1420-8946},
   MRCLASS = {55P43 (13D03 55N15)},
  MRNUMBER = {2182697},
MRREVIEWER = {Sarah\ Whitehouse},
       DOI = {10.4171/CMH/31},
       URL = {https://doi.org/10.4171/CMH/31},
}

@article {Stapleton2015,
    AUTHOR = {Stapleton, Nathaniel},
     TITLE = {Subgroups of {$p$}-divisible groups and centralizers in
              symmetric groups},
   JOURNAL = {Trans. Amer. Math. Soc.},
  FJOURNAL = {Transactions of the American Mathematical Society},
    VOLUME = {367},
      YEAR = {2015},
    NUMBER = {5},
     PAGES = {3733--3757},
      ISSN = {0002-9947,1088-6850},
   MRCLASS = {55N20 (14L05 55N22)},
  MRNUMBER = {3314822},
MRREVIEWER = {Rui\ Miguel\ Saramago},
       DOI = {10.1090/S0002-9947-2014-06344-7},
       URL = {https://doi.org/10.1090/S0002-9947-2014-06344-7},
}

@misc{benmoshe2024highersemiadditivitytranschromatichomotopy,
      title={Higher Semiadditivity in Transchromatic Homotopy Theory}, 
      author={Shay Ben-Moshe},
      year={2024},
      eprint={2411.00968},
      archivePrefix={arXiv},
      primaryClass={math.AT},
      url={https://arxiv.org/abs/2411.00968}, 
}

@article{Ragnarsson2008RetractiveTransfers,
  author  = {Ragnarsson, Kari},
  title   = {Retractive transfers and {$p$}-local finite groups},
  journal = {Proceedings of the Edinburgh Mathematical Society},
  volume  = {51},
  year    = {2008},
  pages   = {465--487},
  doi     = {10.1017/S0013091505001719},
  url     = {https://doi.org/10.1017/S0013091505001719},
}

@misc{abouzaid2021arnoldconjecturemoravaktheory,
      title={Arnold Conjecture and Morava K-theory},
      author={Mohammed Abouzaid and Andrew J. Blumberg},
      year={2021},
      eprint={2103.01507},
      archivePrefix={arXiv},
      primaryClass={math.SG},
      url={https://arxiv.org/abs/2103.01507},
}

@misc{bai2025newtransversalityconditionorbifolds,
      title={A new transversality condition on orbifolds and integer-valued Gromov-Witten type invariants},
      author={Shaoyun Bai and Guangbo Xu},
      year={2025},
      eprint={2201.02688},
      archivePrefix={arXiv},
      primaryClass={math.SG},
      url={https://arxiv.org/abs/2201.02688},
}

@misc{abouzaid2025bordismresolutionsingularities,
      title={Bordism and resolution of singularities},
      author={Mohammed Abouzaid and Shaoyun Bai},
      year={2025},
      eprint={2412.04451},
      archivePrefix={arXiv},
      primaryClass={math.AT},
      url={https://arxiv.org/abs/2412.04451},
}

@misc{Gopakumar-Vafa_I,
      title={{M}-Theory and Topological Strings--{I}},
      author={Rajesh Gopakumar and Cumrun Vafa},
      year={1998},
      eprint={hep-th/9809187},
      archivePrefix={arXiv},
      primaryClass={hep-th},
      url={https://arxiv.org/abs/hep-th/9809187},
}

@misc{Gopakumar-Vafa_II,
      title={{M}-Theory and Topological Strings--{II}},
      author={Rajesh Gopakumar and Cumrun Vafa},
      year={1998},
      eprint={hep-th/9812127},
      archivePrefix={arXiv},
      primaryClass={hep-th},
      url={https://arxiv.org/abs/hep-th/9812127},
}

@article{ChouLee2025,
      author={Chou, You-Cheng and Lee, Y.-P.},
      title={Quantum {$K$}-invariants and {Gopakumar--Vafa} invariants {II}. {Calabi--Yau} threefolds at genus zero},
      journal={Geom. Topol.},
      volume={29},
      year={2025},
      pages={4665--4693},
      doi={10.2140/gt.2025.29.4665},
      eprint={2305.08480},
      archivePrefix={arXiv},
      primaryClass={math.AG},
      url={https://arxiv.org/abs/2305.08480},
}

@article{BravermanMaulikOkounkov2011,
      author={Braverman, Alexander and Maulik, Davesh and Okounkov, Andrei},
      title={Quantum cohomology of the {Springer} resolution},
      journal={Adv. Math.},
      volume={227},
      number={1},
      year={2011},
      pages={421--458},
      doi={10.1016/j.aim.2011.01.021},
      eprint={1001.0056},
      archivePrefix={arXiv},
      primaryClass={math.AG},
      url={https://arxiv.org/abs/1001.0056},
}

@misc{Large2022,
  author       = {Large, Timothy},
  title        = {Localization and {Morava K-theory}
                  {Gromov--Witten} Invariants},
  year         = {2022},
  howpublished = {Recorded lecture,
                  MSRI/SLMath Lecture Archive},
  url          = {https://vimeo.com/770979516},
  note         = {Video recording}
}
\end{document}